\input amstex
\documentstyle{amsppt}
\nopagenumbers
\nologo
\magnification\magstephalf
\catcode`@=11
\redefine\output@{%
  \def\break{\penalty-\@M}\let\par\endgraf
  \global\voffset=-20pt
  \ifodd\pageno\global\hoffset=5pt\else\global\hoffset=-34pt\fi  
  \shipout\vbox{%
    \ifplain@
      \let\makeheadline\relax \let\makefootline\relax
    \else
      \iffirstpage@ \global\firstpage@false
        \let\rightheadline\frheadline
        \let\leftheadline\flheadline
      \else
        \ifrunheads@ 
        \else \let\makeheadline\relax
        \fi
      \fi
    \fi
    \makeheadline \pagebody \makefootline
  }%
  \advancepageno \ifnum\outputpenalty>-\@MM\else\dosupereject\fi
}
\font\cpr=cmr7
\newcount\xnumber
\footline={\xnumber=\pageno
\divide\xnumber by 7
\multiply\xnumber by -7
\advance\xnumber by\pageno
\ifnum\xnumber>0\hfil\else\vtop{\vskip 0.5cm
\noindent\cpr CopyRight \copyright\ Sharipov R.A.,
1998, 2007.}\hfil\fi}
\def\setfirstpage{\global\firstpage@true}
\catcode`\@=\active
\fontdimen3\tenrm=3pt
\fontdimen4\tenrm=0.7pt

\def\leaderfill{\leaders\hbox to 0.3em{\hss.\hss}\hfill}
\font\tvbf=cmbx12
\font\tvrm=cmr12
\font\etbf=cmbx8
\font\tencyr=wncyr10
\def\negskp{\hskip -2pt}
\def\Cl{\operatorname{Cl}}
\def\Img{\operatorname{Im}}
\def\sign{\operatorname{sign}}
\def\Rbinary{\overset{\sssize R}\to\rightharpoondown}
\def\Pprec{\overset{\sssize P}\to\prec}
\def\Rsim{\overset{\sssize R}\to\sim}
\def\id{\operatorname{id}}
\def\compos{\,\raise 1pt\hbox{$\sssize\circ$} \,}
\accentedsymbol\updownarrows{\,\uparrow\hskip-1.6pt\downarrow\,}
\Monograph
\loadbold
\TagsOnRight
\newcount\chapternum
\def\blue#1{#1}
\def\red#1{#1}
\catcode`#=11\def\diez{#}\catcode`#=6
\catcode`_=11\catcode`_=8
\def\mycite#1{\cite{\blue{#1}}\immediate\special{ps:
     ShrHPSdict begin /ShrBORDERthickness 0 def}}

\def\mytag#1{%
    \tag#1}
\def\mythetag#1{\thetag{\blue{#1}}\immediate\special{ps:
     ShrHPSdict begin /ShrBORDERthickness 0 def}}

\def\myrefno#1{\no#1}
\def\myhref#1#2{\blue{#2}\immediate\special{ps:
     ShrHPSdict begin /ShrBORDERthickness 0 def}}
\def\myEarXivlink{\myhref{http://arXiv.org}{http:/\negskp/arXiv.org}}
\def\mytheorem#1{\csname proclaim\endcsname{Theorem #1}}
\def\mytheoremwithtitle#1#2{\csname proclaim\endcsname{Theorem #1#2}}
\def\mythetheorem#1{\blue{#1}\immediate\special{ps:
     ShrHPSdict begin /ShrBORDERthickness 0 def}}
\def\mythetheoremchapter#1#2{\blue{#1}\immediate\special{ps:
     ShrHPSdict begin /ShrBORDERthickness 0 def}}
\def\mylemma#1{\csname proclaim\endcsname{Lemma #1}}
\def\mylemmawithtitle#1#2{\csname proclaim\endcsname{Lemma #1#2}}
\def\mythelemma#1{\blue{#1}\immediate\special{ps:
     ShrHPSdict begin /ShrBORDERthickness 0 def}}
\def\mythelemmachapter#1#2{\blue{#1}\immediate\special{ps:
     ShrHPSdict begin /ShrBORDERthickness 0 def}}
\def\myproposition#1{\csname proclaim\endcsname{Proposition #1}}
\def\mypropositionwithtitle#1#2{\csname proclaim\endcsname{Proposition #1#2}}

\def\mycorollary#1{\csname proclaim\endcsname{Corollary #1}}

\def\myexercise#1{\csname proclaim\endcsname{Exercise #1}}

\def\mydefinition#1{\definition{Definition #1}}
\def\mythedefinition#1{\blue{#1}\immediate\special{ps:
     ShrHPSdict begin /ShrBORDERthickness 0 def}}
\def\mythedefinitionchapter#1#2{\blue{#1}\immediate\special{ps:
     ShrHPSdict begin /ShrBORDERthickness 0 def}}
\def\myaxiom#1{\csname proclaim\endcsname{\red{Axiom #1}}}
\def\myaxiomwithtitle#1#2{\csname proclaim\endcsname{\red{Axiom #1}#2}}
\def\mytheaxiom#1{\blue{#1}\immediate\special{ps:
     ShrHPSdict begin /ShrBORDERthickness 0 def}}
\def\SectionNum#1#2{\S
     \,#1.}
\pagewidth{10cm}
\pageheight{15.5cm}
\fontdimen3\tenrm=3pt
\fontdimen4\tenrm=0.7pt

\document

\chapternum=1
\vbox to\vsize{\centerline{\etbf MINISTRY OF GENERAL AND PROFESSIONAL}
\centerline{\etbf EDUCATION OF RUSSIAN FEDERATION}
\medskip 
\centerline{\etbf BASHKIR STATE UNIVERSITY}
\vskip 3cm
\centerline{SHARIPOV\ R.\,A.}
\vskip 1.5cm
\centerline{\tvbf FOUNDATIONS OF GEOMETRY}
\centerline{\tvbf FOR UNIVERSITY STUDENTS}
\centerline{\tvbf AND HIGH SCHOOL STUDENTS}
\vskip 1.3cm
\centerline{\tvrm The textbook}
\vskip 5.0cm
\centerline{Ufa 1998}

\vss}
\newpage
\vbox to 13.5cm{
UDC 541.1\par
Sharipov R. A. {\bf Foundations of geometry for university students 
and high-school students.} The textbook --- Ufa, 1998. --- 220 pages 
--- ISBN 5-7477-02494-1.
\bigskip
\bigskip
     This book is a textbook for the course of foundations of geometry. 
It is addressed to mathematics students in Universities and to High
School students for deeper learning the elementary geometry. It can
also be used in mathematics coteries and self-education groups.\par
     In preparing Russian edition of this book I used
computer typesetting on the base of \AmSTeX\ package and
I used Cyrillic fonts of Lh-family distributed by CyrTUG
association of Cyrillic \TeX\ users. English edition is also
typeset by \AmSTeX.\par
\medskip
\noindent
Referees:\ \ \ \
\vtop{\hsize 8.0cm\noindent Prof. R.~R.~Gadylshin, Bashkir State
Pedagogical University ({\tencyr BGPU});
\vskip 0.1cm
\noindent Prof\.~E.~M.~Bronshtein, Ufa State University for Aircraft 
and Technology ({\tencyr UGATU});
\vskip 0.1cm
\noindent Mrs. N.~V.~Medvedeva, Honored teacher of the Republic 
Bashkortostan, High School No. 42, Ufa.}
\vfil
\line{ISBN 5-7477-02494-1\hss\copyright\ Sharipov R.A., 1998}
\line{English Translation\hss\copyright\ Sharipov R.A., 2007}}
\newpage
\ \bigskip\medskip
\centerline{\bf CONTENTS.}
\bigskip
\line{CONTENTS.\ \leaderfill\ 3.}
\medskip
\line{PREFACE.\ \leaderfill\ \myhref{\diez pg6}{6}.}
\medskip
\line{CHAPTER~\uppercase\expandafter{\romannumeral 1}.
EUCLID'S GEOMETRY. ELEMENTS\hss}
\line{\qquad OF THE SET THEORY AND AXIOMATICS.\ \leaderfill\ 
\myhref{\diez s1pg7}{7}.}
\medskip
\line{\S\,1. Some initial concepts of the set theory.\ \leaderfill\ 
\myhref{\diez s1pg7}{7}.}
\line{\S\,2. Equivalence relations and breaking into\hss}
\line{\qquad equivalence classes.\ \leaderfill\ 
\myhref{\diez s2pg9}{9}.}
\line{\S\,3. Ordered sets.\ \leaderfill\ \myhref{\diez s3pg10}{10}.}
\line{\S\,4. Ternary relations.\ \leaderfill\ 
\myhref{\diez s4pg11}{11}.}
\line{\S\,5. Set theoretic terminology in geometry.\ \leaderfill\ 
\myhref{\diez s5pg11}{11}.}
\line{\S\,6. Euclid's axiomatics.\ \leaderfill\ \myhref{\diez s6pg12}{12}.}
\line{\S\,7. Sets and mappings.\ \leaderfill\ 
\myhref{\diez s7pg12}{12}.}
\line{\S\,8. Restriction and extension of mappings.\ \leaderfill\ 
\myhref{\diez s8pg16}{16}.}
\bigskip
\line{CHAPTER~\uppercase\expandafter{\romannumeral 2}.
AXIOMS OF INCIDENCE\hss}
\line{\qquad AND AXIOMS OF ORDER.\ \ \leaderfill\ \myhref{\diez s1pg18}{18}.}
\medskip
\line{\S\,1. Axioms of incidence.\ \leaderfill\ \myhref{\diez s1pg18}{18}.}
\line{\S\,2. Axioms of order.\ \leaderfill\ 
\myhref{\diez s2pg24}{24}.}
\line{\S\,3. Segments on a straight line. \ \leaderfill\ 
\myhref{\diez s3pg30}{30}.}
\line{\S\,4. Directions. Vectors on a straight line.\ \leaderfill\ 
\myhref{\diez s4pg34}{34}.}
\line{\S\,5. Partitioning a straight line and a plain.\ \leaderfill\ 
\myhref{\diez s5pg40}{40}.}
\line{\S\,6. Partitioning the space.\ \leaderfill\ 
\myhref{\diez s6pg45}{45}.}
\bigskip
\line{CHAPTER~\uppercase\expandafter{\romannumeral 3}.
AXIOMS OF CONGRUENCE.\ \ \leaderfill\ \myhref{\diez s1pg49}{49}.}
\medskip
\line{\S\,1. Binary relations of congruence.\ \leaderfill\ 
\myhref{\diez s1pg49}{49}.}
\line{\S\,2. Congruence of segments.\ \leaderfill\ 
\myhref{\diez s2pg49}{49}.}
\line{\S\,3. Congruent translation of straight lines. \ \leaderfill\ 
\myhref{\diez s3pg56}{56}.}\pagebreak
\line{\S\,4. Slipping vectors. Addition of vectors\hss}
\line{\qquad on a straight line.\ \leaderfill\ 
\myhref{\diez s4pg59}{59}.}
\line{\S\,5. Congruence of angles.\ \leaderfill\ 
\myhref{\diez s5pg64}{64}.}
\line{\S\,6. A right angle and orthogonality.\ \leaderfill\ 
\myhref{\diez s6pg73}{73}.}
\line{\S\,7. Bisection of segments and angles.\ \leaderfill\ 
\myhref{\diez s7pg78}{78}.}
\line{\S\,8. Intersection of two straight lines by a third one.
\ \leaderfill\ \myhref{\diez s8pg81}{81}.}
\bigskip
\line{CHAPTER~\uppercase\expandafter{\romannumeral 4}.
CONGRUENT TRANSLATIONS\hss}
\line{\qquad AND MOTIONS.\ \ \leaderfill\ \myhref{\diez s1pg84}{84}.}
\medskip
\line{\S\,1. Orthogonality of a straight line and a plane.\ \leaderfill\ 
\myhref{\diez s1pg84}{84}.}
\line{\S\,2. A perpendicular bisector of a segment and\hss}
\line{\qquad the plane of perpendicular bisectors.\ \leaderfill\ 
\myhref{\diez s2pg88}{88}.}
\line{\S\,3. Orthogonality of two planes. \ \leaderfill\ 
\myhref{\diez s3pg89}{89}.}
\line{\S\,4. A dihedral angle.\ \leaderfill\ 
\myhref{\diez s4pg93}{93}.}
\line{\S\,5. Congruent translations of a plane and the space.\ 
\leaderfill\ \myhref{\diez s5pg97}{97}.}
\line{\S\,6. Mirror reflection in a plane and in a straight line.\ 
\leaderfill\ \myhref{\diez s6pg105}{105}.}
\line{\S\,7. Rotation of a plane about a point.\ \leaderfill\ 
\myhref{\diez s7pg106}{106}.}
\line{\S\,8. The total rotation group and the group of pure\hss}
\line{\qquad rotations of a plane.\ \leaderfill\ 
\myhref{\diez s8pg111}{111}.}
\line{\S\,9. Rotation of the space about a straight line.\ 
\leaderfill\ \myhref{\diez s9pg113}{113}.}
\line{\S\,10. The theorem on the decomposition of rotations.\ 
\leaderfill\ \myhref{\diez s10pg116}{116}.}
\line{\S\,11. The total rotation group and the group of pure\hss}
\line{\qquad\ rotations of the space.\ \leaderfill\ 
\myhref{\diez s11pg121}{121}.}
\line{\S\,12. Orthogonal projection onto a straight line.\ 
\leaderfill\ \myhref{\diez s12pg122}{122}.}
\line{\S\,13. Orthogonal projection onto a plane.\ 
\leaderfill\ \myhref{\diez s13pg124}{124}.}
\line{\S\,14. Translation by a vector along a straight line.\ 
\leaderfill\ \myhref{\diez s14pg129}{129}.}
\line{\S\,15. Motions and congruence of complicated\hss}
\line{\qquad\ geometric forms.\ 
\leaderfill\ \myhref{\diez s15pg134}{134}.}
\bigskip
\line{CHAPTER~\uppercase\expandafter{\romannumeral 5}.
AXIOMS OF CONTINUITY.\ \ \leaderfill\ \myhref{\diez s1pg138}{138}.}
\medskip
\line{\S\,1. Comparison of straight line segments.\ 
\leaderfill\ \myhref{\diez s1pg138}{138}.}
\line{\S\,2. Comparison of angles.\ 
\leaderfill\ \myhref{\diez s2pg141}{141}.}
\line{\S\,3. Axioms of real numbers.\ 
\leaderfill\ \myhref{\diez s3pg146}{146}.}
\line{\S\,4. Binary rational approximations of real numbers.\ 
\leaderfill\ \myhref{\diez s4pg150}{150}.}
\line{\S\,5. The Archimedes axiom and Cantor's axiom\hss}
\line{\qquad in geometry.\ 
\leaderfill\ \myhref{\diez s5pg154}{154}.}
\line{\S\,6. The real axis.\ 
\leaderfill\ \myhref{\diez s6pg156}{156}.}
\line{\S\,7. Measuring straight line segments.\ 
\leaderfill\ \myhref{\diez s7pg161}{161}.}
\line{\S\,8. Similarity mappings for straight lines. Multiplication\hss}
\line{\qquad of vectors by a number.\ 
\leaderfill\ \myhref{\diez s8pg167}{167}.}
\line{\S\,9. Measuring angles.\ 
\leaderfill\ \myhref{\diez s9pg170}{170}.}
\bigskip
\line{CHAPTER~\uppercase\expandafter{\romannumeral 6}.
THE AXIOM OF PARALLELS.\ \leaderfill\ \myhref{\diez s1pg175}{175}.}
\medskip
\line{\S\,1. The axiom of parallels and the classical\hss}
\line{\qquad Euclidean geometry.\ 
\leaderfill\ \myhref{\diez s1pg175}{175}.}
\line{\S\,2. Parallelism of a straight line and a plane.\ 
\leaderfill\ \myhref{\diez s2pg178}{178}.}
\line{\S\,3. Parallelism of two planes.\ 
\leaderfill\ \myhref{\diez s3pg181}{181}.}
\line{\S\,4. The sum of angles of a triangle.\ 
\leaderfill\ \myhref{\diez s4pg184}{184}.}
\line{\S\,5. Midsegment of a triangle.\ 
\leaderfill\ \myhref{\diez s5pg185}{185}.}
\line{\S\,6. Midsegment of a trapezium.\ 
\leaderfill\ \myhref{\diez s6pg187}{187}.}
\line{\S\,7. Parallelogram.\ 
\leaderfill\ \myhref{\diez s7pg190}{190}.}
\line{\S\,8. Codirected and equal vectors in the space.\ 
\leaderfill\ \myhref{\diez s8pg192}{192}.}
\line{\S\,9. Vectors and parallel translations.\ 
\leaderfill\ \myhref{\diez s9pg197}{197}.}
\line{\S\,10. The group of parallel translations.\ 
\leaderfill\ \myhref{\diez s10pg199}{199}.}
\line{\S\,11. Homothety and similarity.\ 
\leaderfill\ \myhref{\diez s11pg202}{202}.}
\line{\S\,12. Multiplication of vectors by a number.\ 
\leaderfill\ \myhref{\diez s12pg211}{211}.}
\bigskip 
\line{REFERENCES.\ \ \leaderfill\ \myhref{\diez pg214}{214}.}
\bigskip 
\line{CONTACTS.\ \ \leaderfill\ \myhref{\diez pg215}{215}.}
\bigskip
\line{APPENDIX.\ \ \leaderfill\ \myhref{\diez pg216}{216}.}
\newpage
\newtoks\truehead
\truehead=\headline
\headline{\hfill}
\centerline{\bf PREFACE.}
\medskip
\bigskip
     The elementary geometry is a part of geometry which we usually 
first meet in school. It describes the structure of our everyday 
material environment. In very ancient ages people learned to 
discriminate some primitive constituents in the large variety of 
forms they observe in the world surrounding them. These are a 
{\it point}, a {\it straight line\/} and a {\it segment\/} of it, 
a {\it plane}, a {\it circle}, a {\it cylinder}, a {\it ball}, 
and some others. People began to study their properties. Geometers 
of the Ancient Greece succeeded in it better than others. They 
noted that the properties of the simplest geometric forms are not 
a collection of facts, but they are bound to each other by many 
logical bonds. Some of these properties can be deduced from some 
others.\par
     In 5-th century before Christ Euclid offered a list of simplest
basic properties of geometric forms, which now are called 
{\it postulates\/} or {\it axioms\/} of Euclid. The elementary geometry
or Euclidean geometry based on these axioms became the first axiomatic
theory in mathematics.\par
     The aim of this book is to explain the elementary geometry starting
from the Euclid's axioms in their contemporary edition. It is addressed 
to University students as textbook for the course of foundations of
geometry. It can also be recommended to High school students when they
wish to know better what is their subject of elementary geometry from 
the professional point of view of mathematicians. Especially for the
convenience of High school students in 
Chapter~\uppercase\expandafter{\romannumeral 1} of the book I give some
preliminary material from the set theory.\par
     In writing this book I used the books by I.~Ya.~Backelman \mycite{1}
and by  N.~V.~Efimov \mycite{2}. Proofs of some theorems are taken from 
the book by A.~V.~Pogorelov \mycite{3}.
\bigskip\bigskip
\line{\vbox{\hsize 7.5cm\noindent July, 1998;\newline January,
2007.}\hss R.~A.~Sharipov.}

\newpage
\setfirstpage
\topmatter
\title\chapter{1}
Euclid's geometry. Elements 
of the set theory and axiomatics.
\endtitle
\endtopmatter
\leftheadtext{CHAPTER~\uppercase\expandafter{\romannumeral 1}.
ELEMENTS OF THE SET THEORY \dots}
\document
\headline=\truehead
\head
\SectionNum{1}{7} Some basic concepts of the set theory.
\endhead
\rightheadtext{\S\,1. Some basic concepts of the set theory.}
     The set theory makes a ground for constructing the modern 
mathematics in whole. This theory itself is based on two very
simple concepts: the concept of a {\it set\/} and the concept of
an {\it element\/} of a set. Saying a set one usually understand
any collection of objects which for some reason should be grouped
together. Individual objects composing a set are called its 
{\it elements}. A set $A$ and its element $a$ are in the relation 
of {\it belonging}: $a\in A$. This writing says that the element
$a$ {\it belongs\/} to the set $A$ and the set $A$ 
{\it comprises\/} its element $a$. The transposed writing $A\owns a$
means the same.\par
     Let's mark a part of the elements in a set $A$. This marked part
of elements can be treated as another set $B$. The fact that $B$ is a
part of $A$ is denoted as $B\subset A$. If $B\subset A$, we say that
$B$ is a {\it subset\/} of $A$. One should clearly distinguish
two writings:
$$
\xalignat 2
&a\in A,&&B\subset A.
\endxalignat
$$
The inclusion sign $\ \ssize\subset\ $ relates two sets, while the
belonging sign $\ \ssize\in\ $ relates a set with its element.\par
     When composing the set $B$ above, we could mark all of the
elements of $A$. Then we would get $B=A$. But even in this special
case $B$ can be treated as a part of $A$. This means that the
writing $B\subset A$ does not exclude the possibility of coincidence
$B=A$. If we wish to show that $B$ is a subset of $A$ different from
$A$, we should write $B\varsubsetneq A$.\par
     Another special case of $B\subset A$ arises when $B$ contains no
elements at all. Such a set is called the {\it empty\/} set. It is
denoted by the special sign $\varnothing$. The empty set is treated
as a subset of an arbitrary set $A$, i\.\,e\. $\varnothing\subset A$.
\par
     Let $A$ and $B$ be two arbitrary sets. Some of their elements
could be common for them: $c\in A$ and $c\in B$. Such elements
constitute a set $C$ which is called the {\it intersection} of the
sets $A$ and $B$. This set is denoted $C=A\cap B$. If $A\cap B\neq
\varnothing$, then we say that the sets $A$ and $B$ {\it do intersect}.
Otherwise, if $A\cap B=\varnothing$, then we say that these sets 
{\it do not intersect}.\par
      Again let $A$ and $B$ be two arbitrary sets. Let's gather
into one set $C$ all of the elements taken from $A$ and $B$. The
resulting set $C$ in this case is called the {\it union\/} of the
sets $A$ and $B$. It is denoted $C=A\cup B$.\par
     Elements composing the set $A\cup B$ are divided into 
three groups (into three subsets). These are
\roster
\item elements that belong to the sets $A$ and $B$ simultaneously;
\item elements that belong to the set $A$, but do not belong to 
the set $B$;
\item elements that belong to the set $B$, but do not belong to 
the set $A$.
\endroster
The first group of elements constitutes the intersection $A\cap B$. 
The second group of elements constitutes the set which is called
the {\it difference\/} of the sets $A$ and $B$. It is denoted 
$A\setminus B$. Now it is clear that the third group of elements
constitutes the set being the difference $B\setminus A$. The sets
$A\cap B$, $A\setminus B$, and $B\setminus A$ do not intersect with 
each other. Their union coincides with the union of $A$ and $B$, 
i\.\,e\. we have the formula
$$
A\cup B=(A\cap B)\cup(A\setminus B)\cup(B\setminus A).
$$
\head
\SectionNum{2}{9} Equivalence relations and breaking 
into equivalence classes.
\endhead
\rightheadtext{\S\,2. Equivalence relations \dots}
     Let $M$ be some set. Let's consider ordered pairs of elements
$(a,b)$, where $a\in M$ and $b\in M$. Saying 
{\tencyr\char '074}ordered{\tencyr\char '076}, we mean that $a$ is 
the first element
in the pair and $b$ is the second element, the pair $(b,a)$ being
distinct from the pair $(a,b)$. The set of all ordered pairs of
the elements taken from $M$ is called the {\it Cartesian square\/} of
the set $M$. It is denoted $M\times M$.\par
    Assume that some pairs in the set $M\times M$ are somehow marked.
Then the marked pairs form a subset $R\subset M\times M$. If such a
subset $R$ is given, we say that a {\it binary relation\/} $R$ in $M$ 
is given. Indeed, each marked pair $(a,b)\in R$ can be understood as
a sign indicating that its elements $a$ and $b$ are related in some 
way so that such a relation is absent among the elements of a non-marked
pair. If a pair $(a,b)$ is marked, this fact can be denoted in a special
way, e\.\,g\. $a\Rbinary b$. The writing $a\Rbinary b$ is read as follows:
the element $a$ is in relation $R$ to the element $b$.\par
     The relation of equality and the relation of order among real 
numbers are well-known examples of binary relations. They are written
as $a=b$, $a<b$, or $b>a$.\par
\mydefinition{2.1} A binary relation $R$ in a set $M$ is called an
{\it equivalence relation} and is denoted by the sign $\Rsim$ if the
following conditions are fulfilled:
\roster
\item {\it reflectivity:} $a\Rsim a$ for any $a\in M$;
\item {\it symmetry:} $a\Rsim b$ implies $b\Rsim a$;
\item {\it transitivity:} $a\Rsim b$ and $b\Rsim c$ imply
$a\Rsim c$.
\endroster
\enddefinition\noindent
If a binary relation $R$ is implicitly known from the context, the
letter $R$ in the writing $a\Rsim b$ can be omitted and the relation 
$R$ among $a$ and $b$ is written as $a\sim b$.
\mydefinition{2.2} Assume that in a set $M$ an equivalence relation 
$R$ is given. An {\it equivalence class\/} of an element $a\in M$
is the set of all elements $x\in M$ equivalent to $a$, i\.\,e\.
$$
\Cl_{\sssize R}(a)=\{x\in M: x\Rsim a\}.
$$
\enddefinition
\mytheorem{2.1} If $a\Rsim b$, then $\Cl_{\sssize R}(a)=
\Cl_{\sssize R}(b)$. If the elements $a$ and $b$ are not equivalent,
then their classes do not intersect,\linebreak i\.\,e\. 
$\Cl_{\sssize R}(a)\cap\Cl_{\sssize R}(b)=\varnothing$.
\endproclaim
     Sometimes the classes determined by an equivalence relation 
$R$ in $M$ are considered as elements of another set. The set 
composed by all equivalence classes is called the {\it factorset\/}.
It is denoted $M/R$. The passage from $M$ to the factorset $M/R$
is called the {\it factorization}.\par
     The theorem~\mythetheorem{2.1} shows that if two equivalence
classes are distinct, they have no common elements, while each
element $a\in M$ belongs to at least one equivalence class. 
Therefore, each equivalence relation $R$ defines the division
of the set $M$ into the union of non-intersecting equivalence
classes:
$$
M=\bigcup_{\sssize Q\in M/R} Q.
$$
\myexercise{2.1} Using the properties of reflectivity, symmetry,
and transitivity of the equivalence relation $R$, prove the above
theorem~\mythetheorem{2.1}.
\endproclaim
\head
\SectionNum{3}{10} Ordered sets.
\endhead
\rightheadtext{\S\,3. Ordered sets.}
\mydefinition{3.1} A binary relation $P$ in a set $M$ is called
an {\it order relation\/} and is denoted by the symbol $\Pprec$
if the following conditions are fulfilled:
\roster
\item {\it non-reflectivity:} $a\Pprec b$ implies $a\neq b$;
\item {\it non-symmetry:} the condition $a\Pprec b$ excludes 
$b\Pprec a$;
\item {\it transitivity:} $a\Pprec b$ and $b\Pprec c$ imply
$a\Pprec c$.
\endroster
\enddefinition
     The writing $a\Pprec b$ is read as {\tencyr\char '074}{\it $a$
precedes $b$}{\tencyr\char '076} or as {\tencyr\char '074}$b$ follows
$a${\tencyr\char '076}. 
If a binary relation $P$ is implicitly known from the context, the
letter $P$ in the writing $a\Pprec b$ can be omitted and the relation 
$P$ among $a$ and $b$ is written as $a\prec b$.\par
     If one of the mutually exclusive conditions $a\prec b$ or 
$b\prec a$ is fulfilled, we say that the elements $a$ and $b$ are
{\it comparable}. A set $M$ equipped with an order relation $P$
is called a {\it linearly ordered\/} set if any two elements of 
$M$ are comparable. Otherwise, if there are non-comparable pairs 
of elements in $M$, it is called a {\it partially ordered\/} set.
\par
\head
\SectionNum{4}{11} Ternary relations.
\endhead
\rightheadtext{\S\,4. Ternary relations.}
    Along with binary relations, sometimes ternary (or triple) 
relations are considered. A simple example is given by the
addition operation for numbers. The equality $a+b=c$ means that
the ordered triple of numbers $(a,b,c)$ is distinguished as 
compared to other triples, for which such an equality is not
fulfilled. One can easily formalize this example.
\mydefinition{4.1} We say that a {\it ternary relation\/} $R$
in a set $M$ is defined if some subset $R$ in $M\times M\times M$
is fixed.
\enddefinition
\head
\SectionNum{5}{11} Set theoretic terminology in geometry.
\endhead
\rightheadtext{\S\,5. Set theoretic terminology \dots}
     The primary set, which is studied in Euclidean geometry, is 
the {\it space}. Its elements are called {\it points}. The geometric
space of Euclidean geometry is usually denoted by $\Bbb E$. Individual
points of this space by tradition are denoted by capital letters of 
the Roman alphabet. Apart from the whole space and individual points,
various other geometric forms are considered: {\it planes}, 
{\it straight lines}, {\it  segments of straight lines}, {\it rays},
{\it polygons}, {\it polyhedra} etc. All of these geometric forms are
subsets of the space, they consist of points.\par
    The relations of belonging and inclusion denoted by the signs
$\,\ssize\in\,$ and $\,\ssize\subset\,$ in geometry are expressed by
various words corresponding to their visual meaning. Thus, for example,
if a point $A$ belongs to a straight line $m$, then we say that $A$ 
lies on the line $m$, while the line $m$ passes through the point
$A$. Similarly, if a straight line $m$ is included into a plane
$\alpha$, we say that the line $m$ lies on the plane $\alpha$, while
the plane $\alpha$ passes through the line $m$. Usually such a
deliberate wording produces no difficulties for understanding and 
makes an explanation more vivid and visual.
\head
\SectionNum{6}{12} Euclid's axiomatics.
\endhead
\rightheadtext{\S\,6. Euclid's axiomatics.}
     The geometric space $\Bbb E$ consists of points. All points 
of this space are equipollent, none of them is distinguished. If 
a separate point is taken, it has no geometric properties by itself.
The properties of points reveal in their relation to other points.
For example, if we take three points, they can be lying on a straight
line and they can be not lying either. A triangle given by these 
points can be equilateral, or isosceles, or rectangular, or somewhat
else. When composing a geometric form the points of the space $\Bbb E$ 
come into some definite relations with each other. The basic properties
of such relations are formulated in Euclid's axioms. The total number
of Euclid's axioms in their contemporary edition is equal to twenty.
They are divided into five groups:
\roster
\item axioms of {\it incidence} (8 axioms A1--A8);
\item axioms of {\it order} (4 axioms A9--A12);
\item axioms of {\it congruence} (5 axioms A13--A17);
\item axioms of {\it continuity} (2 axioms A18 and A19);
\item axiom of {\it parallels} (1 axiom A20).
\endroster
In forthcoming chapters of this book we give a successive explanation
of the above axioms and the geometry based on them.
\head
\SectionNum{7}{12} Sets and mappings.
\endhead
\rightheadtext{\S\,7. Sets and mappings.}
     Let $X$ and $Y$ be two sets. A {\it mapping\/} of the set
$X$ to the set $Y$ is a rule that associates each element 
$x$ of the set $X$ with some definite element $y$ in the set $Y$. 
The mappings, as well as the sets, are denoted by various letters
(usually by small letters of the Roman alphabet). The writing 
$f\!:X\to Y$ means that $f$ is a mapping of the set $X$ to the 
set $Y$. If $x\in X$, then $f(x)$ denotes the result of applying 
the rule $f$ to the element $x$. The element $y=f(x)$ of the set 
$Y$ is called the {\it image\/} of the element $x$ from $X$. 
An element $x\in X$ such that $y=f(x)$ is called a {\it preimage} 
of the element $y$ from $Y$.\par
      For the rule $f$ to be treated as a mapping $f\!:X\to Y$
it should be unambiguous, i\.\,e\. each occasion of applying it to 
the same element $x\in X$ should yield the same result. In other
words, $x_1=x_2$ should imply $f(x_1)=f(x_2)$.\par
      The simplest example of a mapping is an {\it identical 
mapping\/} of a set $X$ to the same set $X$. It is denoted as
$\id_X\!:X\to X$. The identical mapping $\id_X$ associates each
element $x$ of the set $X$ with itself, i\.\,e\. $\id_X(x)=x$ 
for all $x\in X$.\par
     Let $f\!:X\to Y$ and $g\!:Y\to Z$ are two mappings. In this 
case we can construct the third mapping. Let's define a rule $h$
such that applying it to an element $x$ of $X$ consists in applying
$f$ to $x$ and then applying $g$ to $f(x)$. Ultimately this new rule
yields $g(f(x))$, i\.\,e\. $h(x)=g(f(x))$. The newly constructed
mapping $h\!:X\to Z$ is called the {\it composition\/} of the mappings
$g$ and $f$, it is denoted $h=g\compos f$. So we have
$$
\hskip -2em
g\compos f(x)=g(f(x))
\mytag{7.1}
$$
for all $x$ of $X$. Thus the relationship \mythetag{7.1} is a short
form for the definition of the composition $g\compos f$. The operation
of composition can also be understood as a multiplication, where the
multiplicands are two mappings.
\mytheorem{7.1} If three mappings $f\!:Z\to W$, $g\!:Y\to Z$ and
 $h\!:X\to Y$ are given, then the relationship
$$
\hskip -2em
(f\compos g)\compos h=f\compos (g\compos h)
\mytag{7.2}
$$
is valid. It expresses the {\it associativity\/} of the composition.
\endproclaim
\demo{Proof} We have mappings both in the left and in the right sides
of \mythetag{7.2}, i\.\,e\. \mythetag{7.2} is an equality of mappings.
Two mappings in our case are two rules that associate the elements of
$X$ with some elements in $W$. The statements of these rules can be
quite different, not similar to each other. However, these rules are
treated to be equal if the results of applying them to an element
$x$ do coincide for all $x\in X$. For this reason the proof 
\mythetag{7.2} reduces to verifying the equality
$$
(f\compos g)\compos h(x)=f\compos (g\compos h)(x)
\mytag{7.3}
$$
for all $x$ of $X$. Let's do it by means of direct calculations 
on the base of the formula \mythetag{7.1} defining the composition
of mappings:
$$
\aligned
&(f\compos g)\compos h(x)=f\compos g(h(x))=f(g(h(x))),\\
\vspace{1ex}
&f\compos (g\compos h)(x)=f(g\compos h(x))=f(g(h(x))).
\endaligned
$$
As a result of these rather simple calculations both left and 
right sides of \mythetag{7.3} are reduced to the same expression
$f(g(h(x)))$. The equality \mythetag{7.3} and, hence, the equality
\mythetag{7.2} are proved.
\qed\enddemo
\mydefinition{7.1} Let $f\!:X\to Y$ be a mapping of a set $X$ to 
a set $Y$ and let $A$ be some non-empty subset in $X$. Then the 
set $B\subset Y$ composed by the images of all elements of the 
set $A$ is called the {\it image\/} of the set $A$. It is denoted
$B=f(A)$.
\enddefinition
     According to this definition, the image of a non-empty set is
not empty. For the empty set we set $f(\varnothing)=\varnothing$. 
The image of the set $X$ under the mapping $f\!:X\to Y$ is sometimes
denoted by $\Img f$, i\.\,e\. $\Img f=f(X)$. The set $X$ is called
the {\it domain} of the mapping $f$, the set $Y$ is called its 
{\it domain of values}, and the set $\Img f$ is called the {\it image} 
of the mapping $f$. The domain of values and the image of a mapping $f$
often do not coincide.
\mydefinition{7.2} Let $f\!:X\to Y$ be a mapping of a set $X$ to a set 
$Y$ and let $y$ be some element of the set $Y$. The set composed by
all those elements of $X$ which are taken to the element $y$ by the
mapping $f$ is called the {\it total preimage\/} of the element $y$.
This set is denoted by $f^{-1}(y)$.
\enddefinition
\mydefinition{7.3} Let $f\!:X\to Y$ be a mapping of a set $X$ to 
a set $Y$ and let $B$ be a non-empty subset of the set $Y$. The 
set composed by all those elements $x$ of $X$ whose images $f(x)$
are in $B$ is called the {\it total preimage\/} of the set $B$. 
This set is denoted by $f^{-1}(B)$.
\enddefinition
     According to the definition~\mythedefinition{7.3} the total 
preimage of the set $Y$ coincides with $X$, i\.\,e\. $f^{-1}(Y)=X$. 
For the empty set we set $f^{-1}(\varnothing)=\varnothing$ by
definition. However, even for a non-empty set $B$ its preimage 
$f^{-1}(B)$ can be empty.
\mydefinition{7.4} A mapping $f\!:X\to Y$ is called {\it injective\/}
if for any $y\in Y$ the total preimage $f^{-1}(y)$ contains not more
than one element.
\enddefinition
\mydefinition{7.5} A mapping $f\!:X\to Y$ is called {\it surjective},
if for any $y\in Y$ the total preimage $f^{-1}(y)$ is not empty.
\enddefinition
\mydefinition{7.6} A mapping $f$ is called {\it bijective}, or a 
{\it one-to-one mapping\/} if it is injective and surjective
simultaneously.
\enddefinition
\mytheorem{7.2} A mapping $f\!:X\to Y$ is injective if and only if
$x_1\neq x_2$ implies $f(x_1)\neq f(x_2)$.
\endproclaim
\mytheorem{7.3} A mapping $f\!:X\to Y$ is surjective if and only if
$\,\Img f=Y$.
\endproclaim
\myexercise{7.1} Prove the theorems~\mythetheorem{7.2} and
\mythetheorem{7.3}, which are often used for checking injectivity
and surjectivity of mappings instead of the initial 
definitions~\mythedefinition{7.4} and \mythedefinition{7.5}.
\endproclaim
      Assume that a mapping $f\!:X\to Y$ is bijective. Then for
any element $y\in Y$ the total preimage $f^{-1}(y)$ is not empty,
but it contains not more then one element. Hence, it contains 
exactly one element. For this reason we can define the mapping
$h\!:Y\to X$ which associates each element $y$ of $Y$ with that
very unique element of the set $f^{-1}(y)$. Such mapping $h$
is called the {\it inverse mapping\/} for $f$. It is denoted
$h=f^{-1}$.
\mytheorem{7.4} The mapping $h\!:Y\to X$ inverse to a bijective
mapping $f\!:X\to Y$ is bijective and $h=f^{-1}$ implies $h^{-1}=f$.
\endproclaim
\mytheorem{7.5} A mapping $f\!:X\to Y$ and its inverse mapping
$f^{-1}\!:Y\to X$ are related to each other as follows:
$$
\xalignat 2
&f\compos f^{-1}=\id_Y,&&f^{-1}\compos f=\id_X.
\endxalignat
$$
\endproclaim
\mytheorem{7.6} The composition of two injective mappings is an
injective mapping.
\endproclaim
\mytheorem{7.7} The composition of two surjective mappings is a
surjective mapping.
\endproclaim
\mytheorem{7.8} The composition of two bijective mappings is a
bijective mapping.
\endproclaim
\myexercise{7.2} Prove the theorems~\mythetheorem{7.4}, 
\mythetheorem{7.5}, \mythetheorem{7.6}, \mythetheorem{7.7}, and
\mythetheorem{7.8}.
\endproclaim
\head
\SectionNum{8}{16} Restriction and extension of mappings.
\endhead
\rightheadtext{\S\,8. Restriction an extension of mappings.}
     Let $X'$ be a subset in a set $X$ and assume that two mappings
$f\!:X\to Y$ and $h\!:X'\to Y'$ are given. If $h(x)=f(x)$ for all
$x\in X'$, then we say that $h$ is a {\it restriction\/} of the
mapping $f$ to the subset $X'$. To the contrary, the mapping $f$
is called an {\it extension} or a {\it continuation} of the
mapping $h$ from the set $X'$ to the bigger set $X$.\par
     If a mapping $f\!:X\to Y$ is given, one can easily construct
its restriction to an arbitrary subset $X'\subset X$. It is 
sufficient to forbid applying $f$ to the elements not belonging to
$X'$. As a result we get the mapping $f:X'\to Y$ being a restriction
of the initial one.\par
     To extent a mapping $f:X'\to Y'$ from $X'$ to a bigger set
$X$ is usually more complicated. For this purpose one should
define the values $f(x)$ for those elements $x$ of $X$, which do not
belong to $X'$. This can be done in many ways. However, in a typical
case the mapping $f:X'\to Y'$ possesses some properties that should
be preserved in extending it to $X$. This makes the problem of extending
sapid, but substantially reduces the arbitrariness in choosing possible
extensions.\par
\newpage
\setfirstpage
\topmatter
\title\chapter{2}
Axioms of incidence and axioms of order.
\endtitle
\endtopmatter
\leftheadtext{CHAPTER \uppercase\expandafter{\romannumeral 2}.
AXIOMS OF INCIDENCE AND ORDER.}
\document
\chapternum=2
\head
\SectionNum{1}{18} Axioms of incidence.
\endhead
\rightheadtext{\S\,1. Axioms of incidence.}
\myaxiom{A1} Each straight line contains at least two points.
\endproclaim
\myaxiom{A2} For each two distinct points $A$ and $B$ there is 
a straight line passing through them and this line is unique.
\endproclaim
\myaxiom{A3} In the space there are at least three points 
which do not lie on one straight line.
\endproclaim
The axioms~\mytheaxiom{A1} and \mytheaxiom{A2} show that each
straight line can be fixed by fixing two points on it. This
fact is used for denoting straight lines: saying the line $AB$,
we understand the line passing through the points $A$ and $B$. 
Certainly, the line $AB$ can coincide with the line $CD$ for 
some other two points $C$ and $D$. This occasion is not excluded.
\par
    The axiom~\mytheaxiom{A3} shows that in the space there is at
least one triangle. However, it is not yet that very triangle in
usual sense because the axioms~\mytheaxiom{A1}, \mytheaxiom{A2},
and \mytheaxiom{A3}, taken separately, do not define a segment. On 
the base of these three axioms one cannot distinguish the interior
of a triangle from its exterior.
\myaxiom{A4} For any three points $A$, $B$, and $C$ not lying on
one straight line there is some plane passing through them. Such
a plane is unique.
\endproclaim
\myaxiom{A5} Each plane contains at least one point.
\endproclaim
\myaxiom{A6} If some two distinct points $A$ and $B$ of a straight
line $a$ lie on a plane $\alpha$, then the whole line $a$ lies on 
the plane $\alpha$.
\endproclaim
\myaxiom{A7} If two planes do intersect, their intersection contains
at least two points.
\endproclaim
\myaxiom{A8} In the space there are at least four points not lying 
on one plane.
\endproclaim
     The incidence axioms A1--A8 are yet too few in order to derive
complicated and sapid propositions from them. However, some simple
and visually evident facts can be proved on the base of these axioms.
\mytheorem{1.1} If two distinct straight lines do intersect, their
intersection consists exactly of one point.
\endproclaim
\demo{Proof} Let $a\neq b$ be two distinct straight lines with
non-empty intersection and let $A$ be a point of $a\cap b$. Provided 
the proposition of the theorem is not valid, one could find another 
point $B$ in the intersection $a\cap b$. Thus we would have two 
straight lines $a$ and $b$ passing through the points $A$ and $B$. 
This fact would contradict the axiom~\mytheaxiom{A2}.
\qed\enddemo
\mytheorem{1.2} If two distinct planes do intersect, then their
intersection is a straight line.
\endproclaim
\demo{Proof} Let $A$ be a common point of two distinct planes
$\alpha\neq\beta$. Let's apply the axiom~\mytheaxiom{A7}. 
According to this axiom, there is at least one more commom
point of these two planes $\alpha$ and $\beta$. We denote it
by $B$ and consider the straight line $AB$.\par
     The points $A$ and $B$ lie on the plane $\alpha$. Let's apply
the axiom~\mytheaxiom{A6} to them. This axiom says that the line
$AB$ in whole lies on the plain $\alpha$.\par
     Let's repeat these arguments for the plane $\beta$. As a result
we find that the line $AB$ in whole lies on the plain $\beta$. Thus,
the straight line $AB$ is a part of the intersection $\alpha\cap\beta$. 
It is a common line for these two planes.\par
     The rest is to prove that the intersection of planes $\alpha\cap
\beta$ contains no points other than those lying on the line $AB$. 
If such a point $C$ would exist, then we would have three points 
$A$, $B$, and $C$ not lying on one straight line, and we would have
two distinct planes $\alpha$ and $\beta$ passing through these three
points. But it contradicts to the axiom~\mytheaxiom{A4}. The 
contradiction obtained shows that the intersection $\alpha\cap\beta$
coincides with the line $AB$.
\qed\enddemo
\mytheorem{1.3} For a straight line and a point not lying on this 
straight line there is a plane passing through this line and through
this point. Such a plane is unique.
\endproclaim
\demo{Proof} Let $C$ be a point not lying on a line $a$. Let's apply
the axiom~\mytheaxiom{A1} to the line $a$. According to this axiom,
we can find two points $A$ and $B$ on the line $a$. Then the points
$A$, $B$, and $C$ appear to be three points not lying on one straight
line. Due to the axiom~\mytheaxiom{A4} there is exactly one plane
$\alpha$ passing through the points $A$, $B$, and $C$. Let's apply
the axiom~\mytheaxiom{A6} to the line $a$ and the plane $\alpha$.
From this axiom we derive that the line $a$ is contained within 
the plane $\alpha$. Thus, the plane $\alpha$ a required plane passing
through the point $C$ and the line $a$.\par
    Now let's show that the plane $\alpha$ passing through the point 
$C$ and the line $a$ is unique. Indeed, each plane of this sort should
pass through the above three points $A$, $B$, and $C$ not lying in one
straight line. Due to the axiom~\mytheaxiom{A4} it is unique. 
\qed\enddemo
\mytheorem{1.4} A straight line $a$ not lying on a plane $\alpha$
has not more than one common point with that plane.
\endproclaim
    If the intersection $a\cap\alpha$ is empty, then a straight line
$a$ is said to be {\it parallel} to a plane $\alpha$. Let's consider
the case, where this intersection is non-empty. Assume that $A$ is a
point from the intersection $a\cap\alpha$. If this intersection 
contains more than one point, then there is another point $B\in a\cap
\alpha$. The points $A$ and $B$ of the line $a$ both lie on the plane
$\alpha$. Applying the axiom~\mytheaxiom{A6}, we get that the whole
line $a$ should lie on the plane $\alpha$. However, this contradicts
the initial premise of the theorem: $a\not\subset\alpha$. The
contradiction obtained shows that the intersection $a\cap\alpha$ 
consists of exactly one point $A$. The theorem is proved.
\mytheorem{1.5} For a pair of intersecting, but not coinciding 
straight lines there is exactly one plane containing both of them.
\endproclaim
\demo{Proof} Let $a$ and $b$ be a pair of intersecting, but not 
coinciding straight lines. According to the theorem~\mythetheorem{1.1}
their intersection consists of one point, we denote this point 
by $A$. Then we apply the axiom~\mytheaxiom{A1} to the line $b$.
According to this axiom, there is another point $B$ on $b$ distinct
from $A$. The point $B$ does not lie on $a$ because it lies on the
other line $b$ and does not belong to the intersection $a\cap b$.
\par
    Now let's apply the theorem~\mythetheorem{1.3} to the aline $a$ 
and the point $B$. According to this theorem, there is exactly one
plane $\alpha$ passing through the line $a$ and the point $B$. The
points $A$ and $B$ of the line $b$ lie on the plane $\alpha$. 
Therefore we can apply the axiom~\mytheaxiom{A6}, which says that
that the line $b$ in whole should lie on the plane $\alpha$. Thus,
the plane $\alpha$ contains both lines $a$ and $b$.\par
     The rest is to prove that the plane $\alpha$ is unique. If
it is not unique and if $\beta$ is another plane containing both
lines $a$ and $b$, then $\beta$ passes through the line $a$ and
the point $B$ not lying on $a$. According to the 
theorem~\mythetheorem{1.3}, such a plane is unique. Therefore
the plane $\beta$ should coincide with the plane $\alpha$.
\qed\enddemo
\mylemma{1.1} For any plane $\alpha$ there is a point in the space
not lying on this plane.
\endproclaim
\demo{Proof} Let $\alpha$ be some arbitrary plane. Using the
axiom~\mytheaxiom{A8}, we find four points $A$, $B$, $C$, and $D$,
not lying on one plane. It is clear that at least one of these
four points does not lie on the plane $\alpha$. Otherwise they
would be lying on one plane $\alpha$ in spite of their choice.
\qed\enddemo
\mylemma{1.2} For any straight line $a$ there is a point in the
space not lying on this line.
\endproclaim
\demo{Proof} Let $a$ be some arbitrary straight line. Let's apply
the axiom~\mytheaxiom{A3} and find three points $A$, $B$, and $C$
not lying on one straight line. It is clear that at least one of 
these three points does not lie on the line $a$. Otherwise the 
points $A$, $B$, and $C$ would be lying on one line $a$ in spite 
of their choice.
\qed\enddemo
\mytheorem{1.6} On each plane there are at least three points 
not lying in one straight line.
\endproclaim
\demo{Proof} Assume that some arbitrary plane $\alpha$ is given. 
Applying the axiom~\mytheaxiom{A5}, we choose the first point $A$
on this plane. Then we use the lemma~\mythelemma{1.1}. Because of 
this lemma we can choose a point $X$ outside the plane $\alpha$. 
Therefore, the straight line $AX$ intersects with the plane $\alpha$ 
at the point $A$, but it does not lie on that plane. Now we
can apply the theorem~\mythetheorem{1.4}. This theorem says that
the point $A$ is the unique common point of the line $AX$ and the 
plane $\alpha$.\par
     Now let's apply the lemma~\mythelemma{1.2} to the line $AX$. 
According to this lemma, there is a point $Z$ not lying on the line
$AX$. The points $A$, $X$, $Z$ do not lie on one straight line. 
Therefore, according to the axiom~\mytheaxiom{A4}, they fix a unique
plane $\beta=AXZ$ passing through these three points.\par
     The planes $\alpha$ and $\beta$ do intersect and have the common 
point $A$. Let's apply the axiom~\mytheaxiom{A7} and conclude that, 
apart from the point $A$, there is at least one other common point 
of the planes $\alpha$ and $\beta$. We denote it $B$. As a result we 
have found that there are two distinct points $A$ and $B$ on the
plane $\alpha$.\par
     Now we apply the lemma~\mythelemma{1.1} to the plane $\beta$. 
Due to this lemma we can find a point $Y$ not lying on the plane
$\beta$. The line $AX$ lies in the plane $\beta$, while the point 
$Y$ is outside of this plane.  Therefore $Y\notin AX$. Hence the
three points $A$, $X$, and $Y$ do not lie on one straight line. 
According to the axiom~\mytheaxiom{A4} they determine a unique
plane $\gamma=AXY$ passing through them.\par
     The planes $\alpha$ and $\gamma$ do not coincide since there
is the point $X$ belonging to $\gamma$ and not belonging to $\alpha$. 
These planes do intersect since they have the common point $A$. 
Hence, we can apply the theorem~\mythetheorem{1.2}. It says that 
the intersection of the planes $\alpha$ and $\gamma$ is a straight 
line $a$ containing their common point $A$.\par
     By construction the straight line $a$ lies on the plane $\alpha$. 
It intersects the line $AX$ at the unique point $A$. This fact follows
from the theorem~\mythetheorem{1.1} and from $X\notin a$. Let's prove
that the straight line $a$ does not contain the point $B$. Remember
that $A$ is the unique point of intersection of the line $AX$ and the
plane $\alpha$. Therefore $B\notin AX$ and the three points $A$, $X$,
and $B$ do not lie on one straight line. If we admit that $B\in a$, 
then both planes $\beta$ and $\gamma$ pass through the three points 
$A$, $X$, and $B$. Die to the axiom~\mytheaxiom{A4} they should 
coincide: $\gamma=\beta$. However, by construction the plane $\gamma$
contains the point $Y$ not belonging to the plane $\beta$, i\.\,e\.
$\gamma\neq\beta$. This contradiction shows that $B\notin a$.\par
     As a result of the above considerations within the plane $\alpha$ 
we have constructed a straight line $a$ passing through the point
$A$, and we have constructed a point $B$ not lying on that line. 
Applying the axiom~\mytheaxiom{A1} to the line $a$, we find another
point $C\in a$ distinct from $A$. The three points $A$, $B$, and $C$
is a required triple of points of the plane $\alpha$ not lying on
one straight line.
\qed\enddemo
\myexercise{1.1} Draw figures illustrating the proofs of the 
above six theorems~\mythetheorem{1.1}--\mythetheorem{1.6} and 
two lemmas~\mythelemma{1.1} and \mythelemma{1.2}.
\endproclaim
    Let's some set consisting of four elements. For example, this
can be the set of four initial positive integers $\{1,2,3,4\}$. 
Let's call this set the {\it space}, while the numbers $1$, $2$, 
$3$, and $4$ are its points. The subsets 
$$
\{1,2\},\ \{1,3\},\ \{1,4\},\ \{2,3\},\ \{2,4\},\ \{3,4\}
$$
are called the {\it straight lines\/} in this space. For the planes
we choose the following four subsets:
$$
\{1,2,3\},\ \{1,2,4\},\ \{1,3,4\},\ \{2,3,4\}.
$$
The above sets constitute a {\it finite model\/} of geometry with 
the axioms of incidence.
\myexercise{1.2} Prove that the above finite model 
of geometry satisfies all eight axioms of 
incidence~\mytheaxiom{A1}--\mytheaxiom{A8}.
\endproclaim
\head
\SectionNum{2}{24} Axioms of order.
\endhead
\rightheadtext{\S\,2. Axioms of order.}
     Axioms of order constitute the second group of Euclid's axioms.
Mostly, they describe the inner structure of separate straight lines.
For any three distinct points $A$, $B$, and $C$ lying on one straight
line one of them lies between two others. If $B$ lies between $A$ and 
$C$, we write this fact as
$$
\hskip -2em
(A\blacktriangleright B\blacktriangleleft C).
\mytag{2.1}
$$
The axioms of order determine the properties of the ternary 
relation of the points of a fixed straight line written as 
\mythetag{2.1}.
\myaxiom{A9} If a point $B$ lies between the points $A$ and $C$, then
it lies between $C$ and $A$.
\endproclaim
\noindent Using the notation \mythetag{2.1}, this axiom can be written
as follows:
$$
\hskip -2em
(A\blacktriangleright B\blacktriangleleft C)\implies
(C\blacktriangleright B\blacktriangleleft A).
\mytag{2.2}
$$
The axiom~\mytheaxiom{A9} and the formula \mythetag{2.2} mean the
symmetry of the ternary relation, which is sometimes called the
{\tencyr\char '074}betweenness{\tencyr\char '076} relation,
under the exchange of its first and third arguments.\par
     Let $A$ and $B$ be two arbitrary distinct points. According 
to the axiom~\mytheaxiom{A2} they fix the straight line $AB$. 
An {\it open interval\/} (or simply an {\it interval}) is the
set of all points of the line $AB$ lying between the points $A$ 
and $B$:
$$
(AB)=\{X\in AB:\quad(A\blacktriangleright X\blacktriangleleft B)\}.
$$
The axiom~\mytheaxiom{A9} means that the interval $(AB)$ coincides
with the interval $(BA)$. The points $A$ and $B$ are called the
{\it ending points\/} of the interval $(AB)$. Joining the ending 
points to an open interval, we get a {\it closed interval\/} or
a {\it segment}:
$$
[AB]=\{A\}\cup\{B\}\cup(AB).
$$
According to the axiom~\mytheaxiom{A9} the segment $[AB]$ coincides
with the segment $[BA]$.\par
     The interval $(AB)$ is called the {\it interior\/} of the
segment $[AB]$, while $A$ and $B$ are its {\it ending points}. 
The points of the straight line $AB$ not belonging to the segment 
$[AB]$ constitute the {\it exterior\/} of the segment $[AB]$. 
Along with the open and closed intervals, sometimes one defines 
semi-open intervals:
$$
\xalignat 2
&[AB)=\{A\}\cup (AB),&&[BA)=\{B\}\cup (AB).
\endxalignat
$$
\myaxiom{A10} For any two points $A$ and $B$ on the straight line 
$AB$ there is a point $C$ such that $B$ lies between $A$ and $C$.
\endproclaim
\myaxiom{A11} For any three distinct points $A$, $B$, and $C$ lying
on one straight line only one of them can lie between two others.
\endproclaim
The axiom~\mytheaxiom{A11} means that not more than one of the 
following three conditions can be fulfilled:
$$
\xalignat 3
&\quad(A\blacktriangleright B\blacktriangleleft C),
&&(B\blacktriangleright C\blacktriangleleft A),
&&(C\blacktriangleright A\blacktriangleleft B).
\quad\qquad
\mytag{2.3}
\endxalignat
$$
Generally speaking, the axiom~\mytheaxiom{A11} does not exclude 
the case where none of the above conditions \mythetag{2.3} is 
fulfilled.\par
\mytheorem{2.1} The exterior of any segment $[AB]$ is not empty.
\endproclaim
\demo{Proof} Let's apply the axiom~\mytheaxiom{A10} to the points
$A$ and $B$. It yields the existence of a point $C$ lying on the
straight line $AB$ such that the condition $(A\blacktriangleright 
B\blacktriangleleft C)$ is fulfilled. Due to the 
axiom~\mytheaxiom{A11} this condition excludes the other condition
$(B\blacktriangleright C\blacktriangleleft A)$. Hence, the point
$C$ is not an inner point of the segment $[AB]$. It is not an
ending point either since it does coincide neither with $A$, nor
with $B$. Hence, $C$ is a point on the line $AB$ external for 
the segment $[AB]$. 
\qed\enddemo
Actually, the axioms~\mytheaxiom{A9}, \mytheaxiom{A10} and
\mytheaxiom{A11} can yield more. From them one can derive the 
existence of at least two points in the exterior of any segment 
$[AB]$.
\mytheorem{2.2} For any segment $[AB]$ there are two points $C_1$ 
and $C_2$ on the line $AB$ such that the conditions
$(A\blacktriangleright B\blacktriangleleft C_1)$ and
$(B\blacktriangleright A\blacktriangleleft C_2)$ are fulfilled.
\endproclaim
\demo{Proof} Note that the points $A$ and $B$ enter the statement
of the axiom~\mytheaxiom{A10} in an asymmetric way. For the
beginning we apply the axiom~\mytheaxiom{A10} in its standard
form. It yields the existence of a point $C_1$ on the line $AB$
such that the condition $(A\blacktriangleright B\blacktriangleleft
C_1)$ is fulfilled. Then we exchange $A$ and $B$ and apply the
axiom~\mytheaxiom{A10} once more. Now it yields the existence 
of a point $C_2$ on the line $AB$ such that $(B\blacktriangleright 
A\blacktriangleleft C_2)$.\par
     The points $C_1$ and $C_2$ both are in the exterior of the
segment $[AB]$. However, they cannot coincide. Indeed, if 
$C=C_1=C_2$, then, using the axiom~\mytheaxiom{A9}, we would 
derive that the conditions $(A\blacktriangleright B
\blacktriangleleft C)$ and $(C\blacktriangleright A
\blacktriangleleft B)$ are fulfilled simultaneously. But this
opportunity is prohibited by the axiom~\mytheaxiom{A11}.
\qed\enddemo
\myaxiom{A12} Let $A$, $B$ and $C$ be three points of a plane
$\alpha$ not lying on one straight line and let $a$ be a straight
line on the same plane $\alpha$ passing through none of these
three points. If the line $a$ intersects the segment $[AB]$ at its 
interior point, then it necessarily passes through an interior
point of at least one of the segments $[AC]$ or $[BC]$.
\endproclaim
The axiom~\mytheaxiom{A12} is known as Pasch's axiom. It is
important for all of the further constructions in this section.
\par
\mytheorem{2.3} The interior of any segment $[AB]$ is not empty.
\endproclaim
\demo{Proof} Let's apply the lemma~\mythelemma{1.2} and find a 
point $C$ not lying on the line $AB$ (see Fig\.~2.1 below). Then 
we apply the axiom~\mytheaxiom{A10} to the points $A$ and $C$. 
It yields the existence of a point $D$ lying on the line $AC$ 
and such that $C$ is in the interior of the segment $[AD]$. The 
next step is to draw the line $DB$ and apply the
theorem~\mythetheorem{2.1} to the segment $[DB]$. As a result 
we find a point $E$ on the line $DB$ lying outside the segment 
$[DB]$.\par
\parshape 18 0cm 10cm 0cm 10cm 0cm 10cm 0cm 10cm 
4.5cm 5.5cm 4.5cm 5.5cm 4.5cm 5.5cm 4.5cm 5.5cm
4.5cm 5.5cm 4.5cm 5.5cm 4.5cm 5.5cm 4.5cm 5.5cm
4.5cm 5.5cm 4.5cm 5.5cm 4.5cm 5.5cm 4.5cm 5.5cm
4.5cm 5.5cm 0cm 10cm 
     The line $CE$ crosses the line $AD$ at the unique point $C$, 
which is in the interior of the segment $[AD]$. It crosses the 
line $DB$ at the point $E$, which is distinct from $B$. Hence,
\vadjust{\vskip 5pt\hbox to 0pt{\kern 5pt
\includegraphics{Oris01.eps}\hss}\vskip -5pt}the 
line $CE$ contains none of the points $A$, $D$, and $B$. We
denote the line $CE$ by $a$ and apply Pasch's 
axiom~\mytheaxiom{A12} to the points $A$, $D$, and $B$. According
to this axiom, the line $CE$ should cross the segment $[AB]$ or the
segment $[DB]$ somewhere at an interior point. In our case it 
cannot cross the segment $[DB]$. Indeed, the line $CE$ crosses
the line $DB$ at the unique point $E$, which is in the exterior
of the segment $[DB]$ by construction. Therefore, the line $CE$
crosses the segment $[AB]$ at some its interior point $F$. This 
means that the interior of the segment $[AB]$ is not empty. The 
theorem~\mythetheorem{2.3} is proved.
\qed\enddemo
\mytheorem{2.4} For any three points $A$, $B$, and $C$ lying on
one straight line exactly one of them lies between two others.
\endproclaim
\demo{Proof}\parshape 14 0cm 10cm 0cm 10cm
4.5cm 5.5cm 4.5cm 5.5cm 4.5cm 5.5cm 4.5cm 5.5cm
4.5cm 5.5cm 4.5cm 5.5cm 4.5cm 5.5cm 4.5cm 5.5cm
4.5cm 5.5cm 4.5cm 5.5cm 4.5cm 5.5cm 0cm 10cm
Let $A$, $B$ and $C$ be three arbitrary points lying on one
\vadjust{\vskip 5pt\hbox to 0pt{\kern -5pt
\includegraphics{Oris02.eps}\hss}\vskip -5pt}straight 
line. Assume that the $A$ does not lie between $B$ and $C$
and assume also that $C$ does not lie between $A$ and $B$. Under
these assumptions let's prove that $B$ lies between $A$ and $C$. 
For the beginning we apply the lemma~\mythelemma{1.2} and find a
point $D$ not lying on the line $AC$ (see Fig\.~2.2). Then we apply 
the axiom~\mytheaxiom{A10} to the points $A$ and $D$ and find a 
point $E$ on the line $AD$ such that $D$ is an interior point of 
the segment $[AE]$. Let's draw the lines $DC$, $EC$, and $EB$. At 
the intersection of the lines $DC$ and $EB$ we get the point $F$. 
Through the points $A$ and $F$ we draw the line $AF$. At the
intersection of the line $AF$ with the line $EC$ we get the 
point $G$.\par
     Let's consider the triangle $ABE$. The line $DC$ intersects
its side $[AE]$ at the interior point $D$. According to 
Pasch's axiom~\mytheaxiom{A12} this line should intersect one 
of the segments $[AB]$ or $[EB]$ at some interior point. In our
case the segment $[AB]$ is excluded since the lines $DC$ and $AB$ 
intersect at the point $C$. This point, according to the assumption 
made in the beginning of the proof, does not lie on the segment 
$[AB]$ between the points $A$ and $B$. Hence, the point $F$ obtained
as the intersection of the lines $DC$ and $EB$ is an interior point 
of the segment $[EB]$.\par
     Now let's consider the triangle $EBC$. The line $AF$ crosses
its side $[EB]$ at the interior point $F$. According to 
Pasch's axiom~\mytheaxiom{A12} this line should intersect one 
of the segments $[BC]$ or $[EC]$ at some interior point. The
segment $[BC]$ is excluded since the lines $AF$ and $BC$ intersect
at the point $A$, while this point, according to the assumption
made in the beginning of the proof, does not lie between the points 
$B$ and $C$. The rest is the segment $[EC]$, which should intersect
the line $AF$ at its interior point $G$.\par
     The nest step is to consider the triangle $EDC$ and the line
$AG$ which intersect its side $EC$ at the interior point $G$. 
Applying Pasch's axiom~\mytheaxiom{A12} in this case, we get 
that the line $AG$ intersect one of the segments $[DC]$ or $[DE]$
at some interior point. The segment $[DE]$ is excluded. Indeed,
the line $AG$ crosses the line $DE$ at the point $A$, while the
point $E$ is chosen so that the condition $(A\blacktriangleright
D\blacktriangleleft E)$ is fulfilled. Due to the 
axiom~\mytheaxiom{A11} this condition excludes the condition
$(E\blacktriangleright A\blacktriangleleft D)$, i\.\,e\. 
the point $A$ cannot be an interior point of the segment $[DE]$.
Hence, $F$ is an interior point of the segment $[DC]$.\par
     In the last step we consider the triangle $ADC$ and the line
$EB$ which intersects its side $[DC]$ at the interior point $F$. 
Let's apply Pasch's axiom~\mytheaxiom{A12} in this case. From
this axiom we derive that the line $EB$ should intersect one of the
segments $[AD]$ or $[AC]$ at some interior point. The segment 
$[AD]$ is excluded. Indeed, the line $EB$ crosses the line $AD$ 
at the point $E$. The condition $(A\blacktriangleright
D\blacktriangleleft E)$ for the point $E$ excludes the condition
$(D\blacktriangleright E\blacktriangleleft A)$ and, hence, the
point $E$ does not lie in the interior of the segment $[AD]$. The
rest is the segment $[AC]$. The point $B$ at the intersection of
the lines $EB$ and $AC$ should be an interior point of the segment
$[AC]$. The theorem~\mythetheorem{2.4} is proved.
\qed\enddemo
     The theorem~\mythetheorem{2.4} proved just above strengthens
the axiom~\mytheaxiom{A11}. Now for any three points $A$, $B$, and
$C$ lying on one straight line one of the conditions \mythetag{2.3}
is necessarily fulfilled, thus excluding other two conditions
\mythetag{2.3}. Pasch's axiom~\mytheaxiom{A12} also can be 
strengthened.
\mytheorem{2.5} Let $A$, $B$, and $C$ be three points of a plane
$\alpha$ not lying on one straight line and let $a$ be a straight
line on the plane $\alpha$ passing through neither of these 
three points. Then if $a$ crosses the segment $[AB]$ at some
interior point, it passes through an interior point of exactly 
one of the segments $[AC]$ or $[BC]$.
\endproclaim
\demo{Proof} Assume that the proposition of the theorem is not
valid. Then the line $a$ crosses each of the three segments
$[AB]$, $[BC]$, and $[CA]$ at their interior points. Let's
denote these points by $P$, $Q$, and $R$. The points $P$, $Q$,
and $R$ lie on three distinct straight lines $AB$, $BC$, and $CA$
intersecting each other at three points $A$, $B$, and $C$. 
According to the statement of the theorem, none of the points
$A$, $B$, and $C$ lies on the line $a$, hence, none of the points
$P$, $Q$, and $R$ can coincide with another one.\par
    Let's prove that the point $R$ does not lie between $P$ and $Q$. 
If we admit that the point lies between $P$ and $Q$, then we can 
apply Pasch's axiom~\mytheaxiom{A12} to the line $AR$ and to
the triangle $PQB$. It says that the line $AR$ should intersect
one of the segments $[PB]$ or $[QB]$ at some interior point. 
However, we know that the line $AR$ crosses the line $PB$ at the
point $A$, and it crosses the line $QB$ at the point $C$. If $A$
is in the interior of the segment $[PB]$, this contradicts the 
fact that $P$ lies between $A$ and $B$. Similarly, if $C$ is 
in the interior of the segment $[QB]$, then $Q$ cannot lie
between $B$ and $C$.\par
     The above contradiction proves that the point $R$ cannot
lie between $P$ and $Q$. Similarly, one can prove that $Q$ does 
not lie between $R$ and $P$, while $P$ does not lie between $Q$
and $R$. Thus, none of the points $P$, $Q$, and $R$ on the line
$a$ lies between two others. This contradicts the previous 
theorem~\mythetheorem{2.4}. Therefore the initial assumption 
that the line $a$ crosses both segments $[AC]$ and $[BC]$ at
their interior points is invalid. The theorem~\mythetheorem{2.5} 
is proved.
\qed\enddemo
\head
\SectionNum{3}{30} Segments on a straight line.
\endhead
\rightheadtext{\S\,3. Segments on a straight line.}
\mylemma{3.1} Let $A$, $B$, $C$, and $D$ be a group of four points. 
Assume that the point $B$ lies between $A$ and $C$, while the point 
$C$ lies between $B$ and $D$. Then both points $B$ and $C$ lie 
between the points $A$ and $D$.
\endproclaim
\demo{Proof} \parshape 6 0cm 10.1cm 0cm 10.1cm
4.5cm 5.6cm 4.5cm 5.6cm 4.5cm 5.6cm 4.5cm 5.6cm
From $(A\blacktriangleright B\blacktriangleleft C)$ it follows
that the point $A$ lies on the line $BC$, while from
$(B\blacktriangleright C\blacktriangleleft D)$ it follows that
\vadjust{\vskip 5pt\hbox to 0pt{\kern -5pt
\includegraphics{Oris03.eps}\hss}\vskip -5pt}$D$ 
also lies on the line $BC$. Thus, under the assumptions of 
the lemma~\mythelemma{3.1} all of the four points $A$, $B$, $C$,
and $D$ lie on one straight line.\par
\parshape 7 4.5cm 5.5cm 4.5cm 5.5cm 4.5cm 5.5cm 4.5cm 5.5cm
4.5cm 5.5cm 4.5cm 5.5cm 0cm 10cm
Using the lemma~\mythelemma{1.2} we find a point $E$ not lying
on the line $AD$ (see Fig\.~3.1). Then we apply the 
axiom~\mytheaxiom{A10} to the points $C$ and $E$. As a result 
on the line $CE$ we find a point $F$ such that the point $E$ 
lies in the interior of the segment $[CF]$. Let's draw the lines 
$AE$ and $FB$, then consider the triangle $FBC$. The line $AE$ 
crosses the line $FC$ at the point $E$ which is an interior point 
for the segment $[FC]$. The intersection of the lines $AE$ and 
$BC$ coincides with the point $A$ which is outside the segment
$[BC]$. Therefore, according to Pasch's 
axiom~\mytheaxiom{A12}, the line $AE$ should cross the side
$[FB]$ of the triangle $FBC$ at some interior point $G$.\par
     Now let's consider the triangle $AEC$. The line $FB$ 
crosses the line $AC$ at the point $B$ lying in the interior 
of the segment $[AC]$. The same line $FB$ crosses the line $EC$ 
at the point $F$ outside the segment $[EC]$. Hence, according
to Pasch's axiom~\mytheaxiom{A12}, the point $G$ obtained as 
the intersection of the lines $FB$ and $AE$ should be an interior
point of the segment $[AE]$.\par
     In the next step we consider again the triangle $FBC$ and 
draw the line $GD$. This line crosses its side $[FB]$ at the
interior point $G$ and it has no common points with the side 
$[BC]$ since the point $D$ lies outside the segment $[BC]$. 
Hence, due to Pasch's axiom~\mytheaxiom{A12} we conclude 
that the line $GD$ crosses the segment $[FC]$ at some interior
point $H$.\par
     Now let's consider the triangle $GBD$ and the line $FC$. 
We use the fact that the point $C$ is in the interior of the 
segment $[BD]$ and that $F\notin [GB]$. Then from Pasch's
axiom~\mytheaxiom{A12} we derive that $H$ is an interior point 
of the segment $[GD]$.\par
     In the last step we consider the triangle $AGD$ and the line
$FC$. The line $FC$ crosses the side $[GD]$ of this triangle
at the interior point $H$ and it does not intersect the side $[AG]$
since the point $E$ lies outside the segment $[AG]$. Now from
Pasch's axiom~\mytheaxiom{A12} we find that $C$ is an interior point
for the segment $[AD]$. This is one of the propositions of the
lemma, which we had to prove.\par
     The second proposition of the lemma does not require a separate
proof. In order to prove that the point $B$ lies in the interior
of the segment $[AD]$ it is sufficient to exchange the notations
of the points $A$ with $D$ and $B$ with $C$. Thereafter the rest
is to use the first proposition, which is already proved, and 
then return to the initial notations. 
\qed\enddemo
\mylemma{3.2} Let $A$, $B$, $C$, and $D$ be a group of four points.
If the point $C$ lies in the interior of the segment $[AD]$ and
if the point $B$ lies in the interior of the segment $[AC]$, then
$B$ is in the interior of $[AD]$ and $C$ is in the interior of 
$[BD]$.
\endproclaim
\demo{Proof}\parshape 16 0cm 10cm 0cm 10cm 0cm 10cm
0cm 10cm 0cm 10cm 4.5cm 5.5cm 4.5cm 5.5cm 4.5cm 5.5cm
4.5cm 5.5cm 4.5cm 5.5cm 4.5cm 5.5cm
4.5cm 5.5cm 4.5cm 5.5cm 4.5cm 5.5cm
4.5cm 5.5cm 0cm 10cm
It is easy to see that if $C$ lies in the interior of the segment 
$[AD]$ and if $B$ lies in the interior of the segment $[AC]$, then
all of the four points $A$, $B$, $C$, and $D$ lie on one straight 
line $AD$. \vadjust{\vskip 5pt\hbox to 0pt{\kern -5pt
\includegraphics{Oris04.eps}\hss}\vskip -5pt}(see
Fig\.~3.2) Applying the lemma~\mythelemma{1.2}, we choose a point 
$E$ not lying on the line $AD$. Then we apply the axiom~\mytheaxiom{A10}
and find a point $F$ on the line $BE$ such that the point $E$ lies 
between $F$ and $B$. Let's draw the line $FC$ and consider the triangle 
$ABE$. The line $FC$ does not intersect the sides $[AB]$ and $[BE]$ of
this triangle. Indeed, the lines $AB$ and $FC$ intersect at the point
$C$ outside the segment $[AB]$. The point $F$ is the intersection of
the lines $BE$ and $FC$, it lies outside the segment $[BE]$. If the
line $FC$ would intersect the segment $[AE]$ at an interior point, 
then due to Pasch's axiom~\mytheaxiom{A12} it would intersect on of 
the segments $[AB]$ or $[BE]$. However, it is not so. Therefore the
line $FC$ has no common points with the segment $[AE]$.\par
     Now let's draw the line $DE$ that intersects the line $FC$ at
some point $G$. Then we consider the triangle $AED$. The line $FC$ 
intersects the side $[AD]$ of this triangle at the interior point
$C$, but it has no common points with the side $[AE]$. Applying
Pasch's axiom~\mytheaxiom{A12}, we find that $G$ is an interior point
of the segment $[ED]$.\par
     In the next step we consider the triangle $BED$. The line $FC$
crosses its side $[ED]$ at the interior point $G$. The intersection
of the line $FC$ with the line $BE$ is the point $F$. It is outside
the segment $[BE]$. Therefore, due to Pasch's axiom~\mytheaxiom{A12}
the point $C$, which is the intersection of the lines $FC$ and $BD$,
should be an interior point of the segment $[BD]$. Thus we have proved
the first proposition of the lemma~\mythelemma{3.2} saying that
$C$ lies between the points $B$ and $D$.\par
     Note that now we can apply the previous lemma~\mythelemma{3.1}.
Indeed, the point $B$ lies between $A$ and $C$, while the point $C$ 
lies between $B$ and $D$. From the lemma~\mythelemma{3.1} we derive 
that $B$ lies in the interior of the segment $[AD]$. Thus we have
proved the second proposition of the lemma~\mythelemma{3.2}.
\qed\enddemo
\mytheorem{3.1} If a point $B$ lies between two other points $A$ and
$C$, then the segments $[AB]$ and $[BC]$ are subsets of the segment
$[AC]$.
\endproclaim
\demo{Proof} Let's prove the proposition of the theorem concerning
the segment $[AB]$. Remember that a segment consists of two ending
points and of all points lying between these ending points. For the
point $A$ we have $A\in [AC]$. The point $B$ lies between the points
$A$ and $C$. Therefore $B$ is an interior point of the segment $[AC]$, 
i\.\,e\. $B\in [AC]$.\par
    Let $X$ be an arbitrary interior point of the segment $[AB]$. Then
$X$ lies between $A$ and $B$. While the point $B$ lies between $A$ and
$C$. In this case the lemma~\mythelemma{3.2} is applicable. It yields
$X\in [AC]$. Thus, we have proved that $[AB]\subset [AC]$.\par
     In order to prove the second proposition $[BC]\subset [AC]$ it is
sufficient to exchange the notations of the points $A$ and $C$ and then
apply the first proposition $[AB]\subset [AC]$, which is already
proved, and afterward return to the initial notations.
\qed\enddemo
\mytheorem{3.2} If a point $B$ lies between two other points $A$ 
and $C$, then the segment $[AC]$ is the union of the segments $[AB]$ 
and $[BC]$.
\endproclaim
\demo{Proof} According to the previous theorem, the segments$[AB]$ and
$[BC]$ are subsets of the segment $[AC]$. Therefore, we have
$$
\hskip -2em
[AB]\cup [BC]\subset [AC].
\mytag{3.1}
$$
Let's prove the opposite inclusion $[AC]\subset [AB]\cup [BC]$. The
ending points of the segment $[AC]$ and the point $B$ belong to the
union $[AB]\cup [BC]$. Therefore, we consider some arbitrary 
interior point of the segment $[AC]$, different from the point $B$. 
Let's denote it $X$.\par
     If $X\notin [AB]$, then due to the theorem~\mythetheorem{2.4}
exactly one of the following two conditions is fulfilled: $A\in [BX]$ 
or $B\in [AX]$. First of these condition combined with $X\in [AC]$
allows us to apply the lemma~\mythelemma{3.1}. From this lemma we
derive $A\in [BC]$ and $X\in [BC]$. But $A\in [BC]$ contradicts the
fact that $B$ is an interior point of the segment $[AC]$. Hence, we
should study the second condition $B\in [AX]$. When combined with
$X\in [AC]$, it allows us to apply the lemma~\mythelemma{3.2}. From
the lemma~\mythelemma{3.2} we derive $B\in [AC]$ and $X\in [BC]$.\par
     Thus, for an arbitrary interior point $X\neq B$ of the segment
$[AC]$ we have shown that $X\notin [AB]$ implies $X\in [BC]$. Hence,
the required inclusion $[AC]\subset [AB]\cup [BC]$ is proved. When
combined with \mythetag{3.1} it yields the equality $[AB]\cup [BC]
=[AC]$. The proof of the theorem~\mythetheorem{3.2} is complete.
\qed\enddemo
\mytheorem{3.3} If a point $B$ lies between two other points $A$ 
and $C$, then the intersection of the segments $[AB]$ and $[BC]$ 
consists of exactly one point $B$.
\endproclaim
\demo{Proof} The point $B$ is an ending point for both segments 
$[AB]$ and $[BC]$. Therefore, this point belongs to the intersection
$[AB]\cap [BC]$. The ending points $A$ and $C$ do not belong to the
intersection $[AB]\cap [BC]$ since $A\notin [BC]$ and $C\notin [AB]$. 
Hence, each point $X$ of the intersection $[AB]\cap [BC]$ distinct
from $B$ should be an interior point of the segments $[AB]$ and 
$[BC]$.\par
    Let $X$ be an interior point of the segment $[AB]$. Then from
the conditions $B\in [AC]$ and $X\in [AB]$, applying the 
lemma~\mythelemma{3.2}, we derive $X\in [AC]$ and $B\in [XC]$. 
Due to the theorem~\mythetheorem{2.4} the condition $B\in [XC]$ 
excludes the condition $X\in [BC]$. Thus, the segment $[AB]$ cannot
have common interior points with the segment $[BC]$.
\qed\enddemo
\head
\SectionNum{4}{34} Directions. Vectors on a straight line.
\endhead
\rightheadtext{\S\,4. Directions. Vectors on a straight line.}
     Let's consider a set of $n$ on some straight line $a$. We 
enumerate these points denoting them $A_1,\,\ldots,\,A_n$. Let's
call $A_1,\,\ldots,\,A_n$ a {\it monotonic sequence of points\/}
on a line if $n\ge 3$ and if each point $A_i$ lies between the
points $A_{i-1}$ and $A_{i+1}$ for all $i=2,\ldots,n-1$. The
points of a monotonic sequence $A_1,\,\ldots,\,A_n$ determine
the family of $n-1$ segments
$$
\hskip -2em
[A_1,A_2],\,\,[A_2,A_3],\,\,\ldots,\,\,[A_{n-1},A_n].
\mytag{4.1}
$$
Adjacent segments in \mythetag{4.1} have non-empty intersections 
consisting of one point: $[A_i,A_{i+1}]\cap [A_{i+1},A_{i+2}]
=\{A_{i+1}\}$. This fact follows from the theorem~\mythetheorem{3.3}.
Applying the theorem~\mythetheorem{3.2}, we get
$$
[A_i,A_{i+1}]\cup [A_{i+1},A_{i+2}]=[A_i,A_{i+2}].
$$
Applying this theorem several times, we find
$$
\hskip -2em
[A_i,A_{i+m}]=\bigcup^m_{q=1}[A_{i+q-1},A_{i+q}].
\mytag{4.2}
$$
From \mythetag{4.2} one can conclude that $A_i\in [A_{i-q},A_{i+k}]$. 
In other words, the point $A_i$ lies between the points $A_{i-q}$ and 
$A_{i+k}$. Moreover, the following relationship is valid:
$$
\hskip -2em
[A_i,A_{i+1}]\cap [A_j,A_{j+1}]=\varnothing\text{\quad for\ }
j\geqslant i+2.
\mytag{4.3}
$$
In order to prove \mythetag{4.3} we use the fact that for $j\geqslant 
i+2$ the point $A_j$ lies between the points $A_{i+1}$ and $A_{j+1}$. 
Therefore, from the theorem~\mythetheorem{3.1}, we derive $[A_j,A_{j+1}]
\subset [A_{i+1},A_{j+1}]$. And, applying the theorem~\mythetheorem{3.3},
we get
$$
[A_i,A_{i+1}]\cap [A_j,A_{j+1}]\subset [A_i,A_{i+1}]\cap
[A_{i+1},A_{j+1}]=\{A_{i+1}\}.
$$
But the point $A_{i+1}$ does not belong to the segment $[A_j,A_{j+1}]$. 
Therefore, the intersection of the segment $[A_i,A_{i+1}]$ and the
segment $[A_j,A_{j+1}]$ is empty, which is in concordance with the 
formula \mythetag{4.3}.
\mytheorem{4.1} Let $A_1,\,\ldots,\,A_n$ be a monotonic sequence of 
points on a straight line and let $B$ be some point of this line
coinciding with none of the points $A_1,\,\ldots,\,A_n$. Then one can
join the point $B$ to the points $A_1,\,\ldots,\,A_n$ and enumerate
the resulting set of points so that the monotonic sequence of points 
$A_1,\,\ldots,\,A_{n+1}$ will be formed.
\endproclaim
\demo{Proof} Let's consider the three points $A_1$, $A_n$, and $B$. 
According to the theorem~\mythetheorem{2.4}, exactly one of the
following three conditions is fulfilled:
$$
\xalignat 3
&(A_1\blacktriangleright A_n\blacktriangleleft B),
&&(A_1\blacktriangleright B\blacktriangleleft A_n),
&&(B\blacktriangleright A_1\blacktriangleleft A_n).
\qquad\quad
\mytag{4.4}
\endxalignat
$$
If the first condition is valid, we denote $B=A_{n+1}$ and immediately
get the required monotonic sequence $A_1,\,\ldots,\,A_{n+1}$.\par
    If the second condition \mythetag{4.4} is valid, then the point
$B$ lies in the interior of the segment $[A_1,A_n]$ and does not coincide
with $A_1,\,\ldots,\,A_n$. But from the relationship \mythetag{4.2} we
get
$$
[A_1,A_n]=\bigcup^{n-1}_{i=1}[A_i,A_{i+1}],
\mytag{4.5}
$$
the segments in the right hand side of this equality intersecting only
by their ending points. Hence, the point $B$ is an interior point for
exactly one of the segments in the right hand side of \mythetag{4.5}.
Assume that $B\in [A_q,A_{q+1}]$. We advance by one the numbers 
of the points $A_{q+1},\,\ldots,\,A_n$:
$$
A_{q+1}\to A_{q+2},\,\,\ldots,\,\,A_n\to A_{n+1}.
$$
Then assign $B=A_{q+1}$ and get the required monotonic sequence 
$A_1,\,\ldots,\,A_{n+1}$.\par
     In the case where the third condition \mythetag{4.4} is fulfilled
we need to advance the numbers in the whole sequence $A_1,\,\ldots,
\,A_n$:
$$
\pagebreak
A_1\to A_2,\,\,\ldots,\,\,A_n\to A_{n+1}.
$$
Then we assign $B=A_1$ and as a result we obtain the required 
monotonic sequence of points $A_1,\,\ldots,\,A_{n+1}$ including 
the point $B$ and all of the initial points $A_1,\,\ldots,\,A_n$.
\qed\enddemo
\mytheorem{4.2} Any set of $n$ points, where $n\geqslant 3$, lying 
on one straight line can be enumerated so that a monotonic sequence
of points $A_1,\,\ldots,\,A_n$ will be produced.
\endproclaim
\demo{Proof} Let's choose some three point from the given set of 
$n$ points on a straight line. According to the 
theorem~\mythetheorem{2.4}, exactly one of the chosen three points
lies between two others. We denote it $A_2$, while two other points
are denoted $A_1$ and $A_3$. As a result we get the monotonic sequence
of three points $A_1,\,\,A_2,\,\,A_3$. The rest is to add step by step
the other points of the given set, relying on the
theorem~\mythetheorem{4.1} in each step. 
\qed\enddemo
    For any set of $n\geqslant 3$ points on a straight line there 
are exactly two ways of numbering these points converting them into 
a monotonic sequence of points. If one of these two numberings
$A_1,\,\ldots,\,A_n$ is given, the other numbering $B_1,\,\ldots,
\,B_n$ is obtained from the first one as follows:
$$
\hskip -2em
B_1=A_n,\,\,B_2=A_{n-1},\,\,\ldots,\,\,B_n=A_1.
\mytag{4.6}
$$
\mydefinition{4.1} A segment $[AB]$ of a straight line is called
a {\it directed segment} or a {\it vector}, if one of its ending
points is somehow distinguished with respect to the other.
\enddefinition
\parshape 9 4cm 6cm 4cm 6cm 4cm 6cm 4cm 6cm 4cm 6cm 4cm 6cm 4cm 6cm
4cm 6cm 0cm 10cm
    The distinguished ending point of a vector is usually marked by an
arrow in drawing. This point is called the true {\it ending point}, 
while the other (not distinguished) ending point is called the 
{\it starting point\/} of a vector. An arrow is also used for denoting
vectors in writing: $\overrightarrow{AB\,}$\!. 
\vadjust{\vskip 5pt\hbox to 0pt{\kern -5pt
\includegraphics{Oris05.eps}\hss}\vskip -5pt}Note that $[AB]$ 
and $[BA]$ are two equivalent notations for the same segment. However, 
$\overrightarrow{AB\,}$ and $\overrightarrow{BA}$ are 
two different vectors.\par
    Each vector defines a direction on a straight line. It is 
impor\-tant to be able to compare two directions given by two 
vectors.
\mydefinition{4.2} Two vectors $\overrightarrow{AB\,}$ and 
$\overrightarrow{CD\,}$ lying on one straight line are called 
{\it codirected\/} if there is a monotonic sequence of points
$A_1,\,\ldots,\,A_n$ including the points $A=A_i$, $B=A_k$, 
$C=A_j$, and $D=A_q$ such that $\sign(k-i)=\sign(q-j)$.
\enddefinition
     Note that adding new points to a monotonic sequence of points 
$A_1,\,\ldots,\,A_n$ as described in the theorem~\mythetheorem{4.1} 
does not change the signs of $(k-i)$ and $(q-j)$. Renumbering the 
points $A_1,\,\ldots,\,A_n$ as described in \mythetag{4.6} changes 
these signs to opposite ones:
$$
\xalignat 2
&\sign(k-i)\to -\sign(k-i),
&&\sign(q-j)\to -\sign(q-j).
\endxalignat
$$
Thus, we see that the equality $\sign(k-i)=\sign(q-j)$ being valid 
or not does not depend on a particular choice of the monotonic 
sequence of points that includes the starting and ending points
of the vectors $\overrightarrow{AB\,}$ and $\overrightarrow{CD\,}$\!.
Therefore, in order to verify if two vectors $\overrightarrow{AB\,}$
and $\overrightarrow{CD\,}$ are codirected or not it is sufficient
to enumerate the set of starting and ending points of these vectors
as described in the theorem~\mythetheorem{4.2}.\par
    The codirectedness is a binary relation in the set of vectors 
lying on one line. This relation possesses the following properties:
\roster
\item  $\overrightarrow{AB\,}\upuparrows\overrightarrow{AB\,}$ for any
       vector $\overrightarrow{AB\,}$\!;
\item  $\overrightarrow{AB\,}\upuparrows\overrightarrow{CD\,}$
       implies $\overrightarrow{CD\,}\upuparrows\overrightarrow{AB\,}$\!;
\item  $\overrightarrow{AB\,}\upuparrows\overrightarrow{CD\,}$ and
       $\overrightarrow{CD\,}\upuparrows\overrightarrow{EF\,}$
       imply $\overrightarrow{AB\,}\upuparrows\overrightarrow{EF\,}$\!;
\item  if a vector $\overrightarrow{AB\,}$ is not codirected with
       $\overrightarrow{CD\,}$\!, while $\overrightarrow{CD\,}$ is not
       codirected with $\overrightarrow{EF\,}$\!, then
       $\overrightarrow{AB\,}\upuparrows\overrightarrow{EF\,}$.
\endroster
The properties \therosteritem{1}--\therosteritem{4} are easily proved
if one considers some monotonic sequence of points $A_1,\,\ldots,\,A_n$, 
including all of the points $A$, $B$, $C$, $D$, $E$, and $F$. The first
three of these properties mean that the codirectedness relation is
reflective, symmetric, and transitive. The fourth property shows that
if we factorize the vectors on a straight line with respect to this
relation, we get only two equivalence classes, each corresponding one
of two possible directions on this line.\par
     Assume that some vector $\overrightarrow{MN\,}$ on a straight
line $a$ is fixed. Let's agree to call {\it positive\/}  the direction
given by this vector. Then the opposite vector $\overrightarrow{NM\,}$
fixes the {\it negative\/} direction. With these prerequisites, for any
two points $X$ and $Y$ on the line $a$ we say that the point $X$ 
{\it precedes} the point $Y$ if if the vector $\overrightarrow{XY\,}$ 
is in positive direction, i\.\,e\. if $\overrightarrow{XY\,}\upuparrows
\overrightarrow{MN\,}$\!. The relation of {\it precedence\/} is denoted 
as $X\prec Y$. It possesses the following properties, which are easy to
verify:
\roster
\item $A\prec B$ implies $A\neq B$;
\item $A\prec B$ excludes $B\prec A$;
\item $A\prec B$ and $B\prec C$ imply $A\prec C$;
\item for any two points $A$ and $B$ exactly one of the
      two conditions $A\prec B$ or $B\prec A$ is fulfilled.
\endroster
The properties \therosteritem{1}--\therosteritem{4} show that
the precedence relation turns a line with a distinguished
vector $\overrightarrow{MN\,}$ into a linearly ordered set.
\par
\mytheorem{4.3} On a straight line withe a fixed direction on 
it a point $B$ lies between $A$ and $C$ if and only if one of 
the two conditions $A\prec B\prec C$ or $C\prec B\prec A$ is 
fulfilled.
\endproclaim
\myexercise{4.1} Prove that the only way of renumbering the
points of a monotonic sequence $A_1,\,\ldots,\,A_n$ preserving 
the property of being monotonic is given by the formula 
\mythetag{4.6}.
\endproclaim
\myexercise{4.2} Verify the properties 
\therosteritem{1}--\therosteritem{4} for the relation of
codirectedness of vectors.
\endproclaim
\myexercise{4.3} Verify the properties 
\therosteritem{1}--\therosteritem{4} for the relation
of precedence of points on a straight line with a fixed
direction.
\endproclaim
\myexercise{4.4} Prove the theorem~\mythetheorem{4.3}.
\endproclaim
\head
\SectionNum{5}{40} Partitioning a straight line and a plain.
\endhead
\rightheadtext{\S\,5. Partitioning a straight line and a plain.}
    Let's consider some point $O$ on a straight line $a$. According
to the axiom~\mytheaxiom{A1}, on the line $a$ there is at least one
point other than $O$. Let's denote it $E$. The vector 
$\overrightarrow{OE\,}$ fixes one of two possible directions 
on $a$ and defines the precedence relation for the points of $a$.
Let's consider two infinite intervals:
$$
\aligned
(O,+\infty)&=\{X\in a:\quad O\prec X\},\\
(-\infty,O)&=\{X\in a:\quad X\prec O\}.
\endaligned
$$
Using the properties \therosteritem{1}--\therosteritem{4} of the
binary relation of precedence, one can show that the intervals 
$(O,+\infty)$ and $(-\infty,O)$ do not intersect, while the whole
line $a$ is divided into three subsets:
$$
\hskip -2em
a=(-\infty,O)\cup\{O\}\cup(O,+\infty).\hskip -2em
\mytag{5.1}
$$
Joining the point $O$ to each of the infinite intervals $(-\infty,O)$
and $(O,+\infty)$, we get two sets which are called {\it half-lines}
or {\it rays}:
$$
\xalignat 2
&[O,-\infty)=(-\infty,O)\cup\{O\},
&&[O,+\infty)=\{O\}\cup(O,+\infty).
\endxalignat
$$
Thus, each point $O$ on a line $a$ determines the division
of this line into two rays with one common point $O$.\par
    Now let's consider a line $a$ lying on a plane $\alpha$. The
theorem~\mythetheorem{1.6} is applicable to the plane $\alpha$.
It says that on any plane there are at least three points not lying
on one line. Hence the set $\alpha\setminus a$ is not empty. Let's
define an equivalence relation on $\alpha\setminus a$ by setting
$A\sim B$ if $A=B$ or if the segment $[AB]$ has no common points
with the line $a$. The reflexivity ans symmetry of such binary
relation are obvious. The rest is to verify its transitivity.
\par
     Let $A\sim B$ and $B\sim C$. If $A=B$ or if $B=C$, then 
$A\sim C$ is a trivial consequence of one of the relations  
$A\sim B$ or $B\sim C$. The coincidence $A=C$ implies $A\sim C$
by itself. Therefore, we can assume that $A$, $B$, and $C$ are
three distinct points. Under this assumption let's consider two 
cases:
\roster
\item where the points $A$, $B$, and $C$ lie on one straight line;
\item where $A$, $B$, and $C$ do not lie on one straight line.
\endroster\par
     In the first case if we assume that the points $A$ and $C$ are
not equivalent, then the lines $AC$ and $a$ intersect at some 
point $O$ interior for the segment $[AC]$. Let's define a positive
direction on the line $AC$ by means of the vector $\overrightarrow{OA}$. 
Then $C\prec O\prec A$. The point $B$ does not lie on the line $a$, 
therefore, $B\neq O$. Hence, $B$ belongs to one of the intervals 
$(-\infty,O)$ or $(O,+\infty)$. If $B\in (-\infty,O)$, then 
$B\prec O\prec A$, which contradicts the condition $A\sim B$. If
$B\in (O,+\infty)$, then $C\prec O\prec B$, which contradicts the
condition $B\sim C$. In both cases the assumption of non-equivalence
of $A$ and $C$ leads to a contradiction. Therefore, the required
condition $A\sim C$ is fulfilled.\par
     In the second case, assuming that $A$ and $C$ are not equivalent,
we find that the line $a$ passing through none of the points $A$, $B$, 
and $C$ intersects the segment $[AC]$ at some interior point $O$. Then
due to Pasch's axiom~\mytheaxiom{A12} it should intersect one of the
segments $[AB]$ or $[BC]$ at an interior point. This contradicts to
the fact that both conditions $A\sim B$ and $B\sim C$ are fulfilled
simultaneously. The contradiction obtained proves that $A\sim C$.\par
     The above equivalence relation determines the division of the
set $\alpha\setminus a$ into classes. As appears, the number of such
classes is equal to two. Taking into account the axiom~\mytheaxiom{A1},
let's choose some point $O$ lying on the line $a$. Then we choose and 
fix some point $A$ lying on the plane $\alpha$, but not lying on the
line $a$. Let's draw the line $AO$ and apply the axiom~\mytheaxiom{A10}
to the points $A$ and $O$ on this line. As a result we find a point
$B$ on the line $AO$ such that the point $O$ lies in the interior of
the segment $[AB]$.\par
     The points $A$ and $B$ belong to the set $\alpha\setminus a$. They 
are not equivalent since the segment $[AB]$ intersects the line $a$ at 
the point $O$. Hence, the equivalence classes $\Cl(A)$ and $\Cl(B)$ are
distinct. Let's prove that an arbitrary point $X$ of the set $\alpha
\setminus a$ belongs to one of these classes. Let's study two cases:
\roster
\item where the point $X$ lies on the line $AO$;
\item where the point $X$ does not lie on the line $AO$.
\endroster
Naturally, we can assume that the point $X$ differs from $A$ and $B$.
In the first case the vector $\overrightarrow{OA}$ fixes one of the
two possible directions on the line $AO$ and determines a precedence
relation on this line. Due to the division \mythetag{5.1} the point 
$X$ is in one of the infinite intervals $(-\infty,O)$ or $(O,+\infty)$.
\par
     If $X\in(-\infty,O)$, then $X\prec O$ and $B\prec O$. Hence,
applying the theorem~\mythetheorem{4.3} to the points $X$, $B$, and 
$O$, we conclude that the point $O$ cannot lie in the interior of the
segment $[BX]$. Therefore, we have $X\in\Cl(B)$.\par
     If $X\in(O,+\infty)$, then $O\prec X$ and $O\prec A$. Hence, 
applying the theorem~\mythetheorem{4.3} again, we get $X\in\Cl(A)$.
\par
     In the case, where the point $X$ does not lie on the line
$OA$, we can consider the triangle $ABX$ lying on the plane $\alpha$.
The line $a$ lies on the same plane and does not passes through the
points $A$, $B$, and $X$. This line intersects the side $[AB]$ in the
interior point $O$. Let's apply the theorem~\mythetheorem{2.5}, which
strengthens Pasch's axiom. According to this theorem, the line $a$
intersects exactly one of two remaining sides of the triangle $ABX$
--- the side $[AX]$ or the side $[BX]$. If $a$ intersects $[AX]$,
then $a$ does not intersect $[BX]$ and $X\in\Cl(B)$. Otherwise, if
$a$ intersects $[BX]$, then $a$ does not intersect $[AX]$ and we
have $X\in\Cl(A)$.\par
     Let's denote $a_{+}=\Cl(A)$ and $a_{-}=\Cl(B)$. The above 
considerations show that the line $a$ lying on the plane $\alpha$
determines the division of this plane into three subsets:
$$
\hskip -2em
\alpha=a_{-}\cup a\cup a_{+}.
\mytag{5.2}
$$
The division \mythetag{5.2} is analogous to the division \mythetag{5.1}.
The subsets $a_{-}$ and $a_{+}$ are called {\it open half-planes}. 
Extending the analogy with \mythetag{5.1}, we define {\it closed
half-planes}:
$$
\xalignat 2
&\overline{a_{-}}=a_{-}\cup a,
&&\overline{a_{+}}=a_{+}\cup a.
\endxalignat
$$\par
\parshape 14 0cm 10cm 0cm 10cm 0cm 10cm 0cm 10cm
4.0cm 6cm 4cm 6cm 4cm 6cm 4cm 6cm 4cm 6cm 4cm 6cm
4cm 6cm 4cm 6cm 4cm 6cm 0cm 10cm
     Let's consider two non-coinciding straight lines $a$ and $b$
intersecting at a point $O$. According to the theorem~\mythetheorem{1.5},
such lines fix a unique plane $\alpha$ containing both of them.
Each of the lines $a$ and $b$ determines a division of the plane alpha
$\alpha$ into two half-planes. The intersection of two closed half-planes
\vadjust{\vskip 5pt\hbox to 0pt{\kern 5pt
\includegraphics{Oris06.eps}\hss}\vskip -5pt}is 
called an {\it angle}. Let's choose a point $A$ other than $O$ on the
line $a$. In a similar way, on the line $b$ we choose a point $B$, 
different from $O$. The point $B$ does not lie on the line $a$, hence,
it belongs only to one of the closed half-planes --- to $\overline{a_{+}}$
or to $\overline{a_{-}}$. For the sake of certainty assume that  
$\overline{a_{+}}$ is that very half-plane which contains the point 
$B$, while by $\overline{b_{+}}$ we denote the half-plane containing
the point $A$. The angle produced as the intersection of the closed
half-planes $\overline{a_{+}}$ and $\overline{b_{+}}$ is usually
denoted as follows:
$$
\hskip -2em
\angle AOB=\overline{a_{+}}\cap\overline{b_{+}}.
\mytag{5.3}
$$
Applying the axiom~\mytheaxiom{A10}, now we choose a point $C$ on the
line $a$ such that the point $O$ lies between $A$ and $C$. A point $D$
on the line $b$ is chosen in a similar way. The lines $a$ and $b$ 
define four angles at a time on the plane $\alpha$:
$$
\xalignat 2
&\angle AOB=\overline{a_{+}}\cap\overline{b_{+}},
&&\angle BOC=\overline{a_{+}}\cap\overline{b_{-}},\\
&\angle COD=\overline{a_{-}}\cap\overline{b_{-}},
&&\angle DOA=\overline{a_{-}}\cap\overline{b_{+}}.
\endxalignat
$$
The points $A$ and $B$ marking the half-planes $\overline{a_{+}}$ and
$\overline{b_{+}}$ play equal roles in defining the angle $\angle AOB$.
Therefore $\angle AOB$ and $\angle BOA$ are different notations for
the same angle.\par
     Let's consider the angle $\angle AOB$ from \mythetag{5.3}. The 
intersection of the open half-planes $a_{+}\cap b_{+}$ is called the
{\it interior\/} of the angle $\angle AOB$. The point $O$ determines 
the division of the lines $a$ and $b$ into four closed half-lines
(four rays). We denote them as follows:
$$
\xalignat 4
&[OA\rangle,&&[OB\rangle,&&[OC\rangle,&&[OD\rangle.
\endxalignat
$$
The previous notations \pagebreak $[O,+\infty)$ and $[O,-\infty)$ 
for rays are convenient only if we consider rays lying on one 
fixed line.
\mytheorem{5.1} Any angle $\angle AOB$ is the union of its interior
and two rays $[OA\rangle$ and $[OB\rangle$.
\endproclaim
     The rays $[OA\rangle$ and $[OB\rangle$ are called the 
{\it sides\/} of the angle $\angle AOB$, while the point $O$ is 
called its {\it vertex}. The angles $\angle AOB$ and $\angle COB$
on Fig\.~5.1 have the common vertex $O$ and the common side  
$[OB\rangle$, while the other sides of these angles $[OA\rangle$ 
and $[OC\rangle$ lie on one line $a$ and intersect at the unique
point $O$. Such angles are called {\it adjacent angles}.\par
     The union of two adjacent angles is a half-plane. Indeed, let's
consider the union of the angles $\angle AOB$ and $\angle COB$:
$$
\hskip -2em
\gathered
\angle AOB\cup\angle COB=(\overline{a_{+}}\cap\overline{b_{+}})
\cup(\overline{a_{+}}\cap\overline{b_{-}})=\\
=\overline{a_{+}}\cap(\overline{b_{+}}\cup\overline{b_{-}})=
\overline{a_{+}}\cap\alpha=\overline{a_{+}}.
\endgathered
\mytag{5.4}
$$
A closed half-plane $\overline{a_{+}}$ with a marked point $O$ 
on the line $a$ can be treated as an angle. Such an angle is 
called a {\it straight angle}. One should be careful when using 
the notation $\angle AOC$ for a straight angle since this notation
fits for $\overline{a_{+}}$ and for $\overline{a_{-}}$ either.\par
     Let's consider the angles $\angle AOB$ and $\angle COD$ on
Fig\.~5.1. The sides $[OC\rangle$ and $[OD\rangle$ of the second 
angle complement the sides $[OA\rangle$ and $[OB\rangle$ of the 
first one up to the lines $a$ and $b$. Such angles are called
{\it vertical angles}.
\mytheorem{5.2} Any three points $A$, $B$, and $O$ not lying on 
one straight line determine exactly one angle $\angle AOB$ with
the vertex at the point $O$.
\endproclaim
     Let's consider three points $A$, $B$, and $C$ not lying on
one straight line. We denote by $a$ the line $BC$, by $b$ the
line $AC$, and by $c$ the line $AB$. The lines $a$, $b$, and $c$ 
lie on the plane $\alpha$ determined by the points $A$, $B$, and
$C$ according to the axiom~\mytheaxiom{A4}. Let's denote by
$a_{+}$ the half-plane on the plane $\alpha$ determined by the
line $a$ and possessing the point $A$. In a similar way, let
$B\in b_{+}$ and $C\in c_{+}$. The {\it triangle $ABC$\/} is the
set of points of the plane $\alpha$ obtained as the intersection
of three closed half-planes $\overline{a_{+}}$, $\overline{b_{+}}$,
and $\overline{c_{+}}$. Before now, saying a triangle $ABC$, we
could understand the collection of three segments $[AB]$, $[BC]$,
and $[AC]$ connecting some three points $A$, $B$, and $C$ not lying 
on one straight line. Now a triangle $ABC$ is equipped with the
interior. The {\it interior\/} of a triangle $ABC$ is the 
intersection of three open half-planes $a_{+}\cap b_{+}\cap c_{+}$.
\mytheorem{5.3} A triangle $ABC$ is the union of its interior
and its three sides $[AB]$, $[BC]$, and $[AC]$.
\endproclaim
\myexercise{5.1} Prove the theorems~\mythetheorem{5.1}, 
\mythetheorem{5.2}, and \mythetheorem{5.3} by proving the
the following lemma before it .
\endproclaim
\mylemma{5.1} For any points $A$ and $B$ the intersection of
the rays $[AB\rangle$ and $[BA\rangle$ is the segment $[AB]$.
\endproclaim
\myexercise{5.2} Verify the calculations \mythetag{5.4} using some
set-theoretic considerations.
\endproclaim
\myexercise{5.3} Let $A$, $B$, and $C$ be three arbitrary points 
not lying on one straight line. Prove that the interior of the
triangle $ABC$ is not empty.
\endproclaim
\head
\SectionNum{6}{45} Partitioning the space.
\endhead
\rightheadtext{\S\,6. Partitioning the space.}
    Let $\alpha$ be a plane. According to the lemma~\mythelemma{1.1},
there is a point not lying on the plane $\alpha$. Therefore, the
set $\Bbb E\setminus\alpha$ is not empty. We define an equivalence
relation in $\Bbb E\setminus\alpha$ by setting $A\sim B$ if $A=B$
or if the segment $[AB]$ does not intersect the plane $\alpha$. 
The reflexivity and symmetry of this binary relation are obvious.\par
     Let's verify its transitivity. Assume that $A\sim B$ and $B\sim C$.
If $A=B$ or $B=C$, then $A\sim C$ is a trivial consequence of one of the
conditions $A\sim B$ or $B\sim C$. The coincidence $A=C$ implies
$A\sim C$ by itself. Therefore we should consider a general case, where
$A$, $B$, and $C$ are three distinct points. Under this assumption we
study two cases:
\roster
\item where the points $A$, $B$, and $C$ lie on one straight line;
\item where $A$, $B$, and $C$ do not lie on one straight line.
\endroster\par
     In the first case if we assume that $A$ and $C$ are not equivalent,
then the line $AC$ intersects the plane $\alpha$ at some point $O$ lying
in the interior of the segment $[AC]$. According to the 
theorem~\mythetheorem{1.4}, the point $O$ is the unique common point of
the line $AC$ and the plane $\alpha$. Let's define a positive direction 
on the line $AC$ by means of the vector $\overrightarrow{OA}$. Then
$C\prec O\prec A$. The point $B$ does not lie on the plane $\alpha$,
therefore, $B\neq O$. Hence, $B$ belongs to one of the intervals
$(-\infty,O)$ or $(O,+\infty)$ determined by the division \mythetag{5.1}. 
If $B\in (-\infty,O)$, then $B\prec O\prec A$, which contradicts the
condition $A\sim B$. If $B\in (O,+\infty)$, then $C\prec O\prec B$, which
contradicts the condition $B\sim C$. Thus, the assumption that $A$ and $C$
are not equivalent leads to a contradiction. Therefore the required 
condition $A\sim C$ is fulfilled.\par
    In the second case one can draw a plane $\beta$ through the points 
$A$, $B$, and $C$. This fact follows from the axiom~\mytheaxiom{A4}.
Here the assumption that $A$ and $C$ are not equivalent means that
the line $AC$ intersects the plane $\alpha$ at some interior point $O$
of the segment $[AC]$. But $[AC]\subset\beta$, therefore the 
non-coinciding planes $\alpha$ and $\beta$ have the common point $O$. 
Let's apply the theorem~\mythetheorem{1.2} and denote by $a$ the line
obtained as the intersection $\alpha\cap\beta$. The line $a$ lies
on the plane $\beta$ and passes through none of the points $A$, $B$,
and $C$. It intersects the segment $[AC]$ at the interior point $O$. 
Then due to Pasch's axiom~\mytheaxiom{A12} it should intersect one of
the segments $[AB]$ or $[BC]$ at some interior point. But this contradicts 
to the fact that the conditions $A\sim B$ and $B\sim C$ are fulfilled
simultaneously. The contradiction obtained proves that the points 
$A$ and $C$ are equivalent.\par
     The above equivalence relation divides the set
$\Bbb E\setminus\alpha$ into equivalence classes. Here, as in the
case of partitioning a plane, the number of equivalence classes is
equal to two. They are denoted $\alpha_{+}$ and $\alpha_{-}$ and 
are called {\it open half-spaces}. Thus, each plane $\alpha$ defines
the division of the space 
$$
\hskip -2em
\Bbb E=\alpha_{-}\cup\alpha\cup\alpha_{+}
\mytag{6.1}
$$
analogous to the divisions \mythetag{5.1} and \mythetag{5.2} for a 
line and for a plane respectively. Relying on \mythetag{6.1} we define
{\it closed half-spaces}
$$
\xalignat 2
&\overline{\alpha_{-}}=\alpha_{-}\cup\alpha,
&&\overline{\alpha_{+}}=\alpha_{+}\cup\alpha.
\endxalignat
$$
\mylemma{6.1} If $A$, $B$, $C$, and $D$ are four points not lying 
on one plane, then neither three of them can lie on one straight 
line.
\endproclaim
     Let $A$, $B$, $C$, and $D$ are some four points not lying on one
plane. The existence of at least one of such quadruples of points is
granted by the axiom~\mytheaxiom{A8}. Due to the lemma~\mythelemma{6.1}
and the axiom~\mytheaxiom{A4} each three of these four points determine
some plane. Let's denote these planes as follows:
$$
\xalignat 4
&\alpha=BCD,&&\beta=ACD,&&\gamma=ABD,&&\delta=ABC.
\endxalignat
$$
Each of the four planes $\alpha$, $\beta$, $\gamma$, and $\delta$
determine two half-spaces. Let's choose the notations for these
half-spaces so that the following conditions are fulfilled:
$$
\xalignat 4
&A\in\alpha_{+},&&B\in\beta_{+},&&C\in\gamma_{+},&&D\in\delta_{+}.
\endxalignat
$$\par
\parshape 3 0cm 10cm 0cm 10cm 4cm 6cm \noindent
The intersection of closed half-spaces $\overline{\alpha_{+}}$,
$\overline{\beta_{+}}$, $\overline{\gamma_{+}}$, and 
$\overline{\delta_{+}}$ is called a {\it tetrahedron}.
\vadjust{\vskip 5pt\hbox to 0pt{\kern -2pt
\includegraphics{Oris07.eps}\hss}\vskip -5pt}The 
points $A$, $B$, $C$, and $D$ are called {\it vertices\/} of
the tetrahedron $ABCD$, the segments $[AB]$, $[AC]$, $[AD]$,
$[BC]$, $[BD]$, and $[DC]$ are called its {\it edges}, while
the triangles $ABC$, $ABD$, $ACD$, and $BCD$ are called the
{\it faces} of the tetrahedron $ABCD$. The intersection of the
open half-spaces $\alpha_{+}$, $\beta_{+}$, $\gamma_{+}$, and
$\delta_{+}$ is called the {\it interior\/} of the tetrahedron
$ABCD$.\par
\parshape 3 4cm 6cm 4cm 6cm 0cm 10cm
     A tetrahedron is called a {\it three-di\-men\-sional simplex}. 
A triangle is a {\it two-dimensional simplex}, while a segment
is a {\it one-dimensional simplex}. A point is a {\it zero-dimensional
simplex}. Such a terminology is popular in algebraic topology
(see \mycite{4}).
\myexercise{6.1} Prove that the number of classes into 
which the set $\Bbb E\setminus\alpha$ is broken by the above 
equivalence relation is two.
\endproclaim
\myexercise{6.2} Prove the lemma~\mythelemma{6.1}. For this purpose
use the axioms of incidence and the results of \S\,1.
\endproclaim
\mylemma{6.2} Let $\angle AOB$ be the angle determined by three
points $A$, $O$, and $B$ not lying on one straight line. Then
a ray coming out from the point $O$ lies within the angle $\angle 
AOB$ if and only if it intersects the segment $[AB]$.
\endproclaim
\myexercise{6.3} Prove the lemma~\mythelemma{6.2}. For this purpose
complete the ray $[OA\rangle$ up to a whole straight line and upon 
choosing a point $C$ on the line $OA$ not belonging to the ray 
$[OA\rangle$ do draw the triangle $ABC$.
\endproclaim
\myexercise{6.4} Let $A$, $B$, $C$, and $D$ be four points not 
lying on one plane. Show that the interior of the tetrahedron 
$ABCD$ is not empty.
\endproclaim
\myexercise{6.5} Prove that each tetrahedron $ABCD$ is the uni\-on
of its interior and the triangles $ABC$, $ABD$, $ACD$,
and $BCD$.
\endproclaim
\newpage
\setfirstpage
\topmatter
\title\chapter{3}
Axioms of congruence.
\endtitle
\leftheadtext{CHAPTER \uppercase\expandafter{\romannumeral 3}.
AXIOMS OF CONGRUENCE.}
\endtopmatter
\document
\chapternum=3
\head
\SectionNum{1}{49} Binary relations of congruence.
\endhead
\rightheadtext{\S\,1. Binary relations of congruence.}
    The axioms of congruence form the third group of Euclid's axioms.
In formulating these axioms it is assumed that in the set of all 
straight line segments a binary relation is defined which is called 
the {\it congruence}. A similar binary relation is assumed to be given
in the set of all angles. It is also called the {\it congruence},
though the congruence of segments and the congruence of angles are
certainly two different binary relations. For denoting the congruence
of segments and the congruence of angles usually the same sign $\cong$ 
is used.\par
     A straight line segment is given by two points. An angle can 
be given by three points. Therefore, the congruence of segments 
can be treated as a tetrary (or quadruple) relation in the set of 
points, while the congruence of angles can be treated as a hexary 
(or sextuple) relation in the set of points. Such a treatment would 
be more consistent from the formal point of view. But it is less 
visual and, hence, is less convenient.
\head
\SectionNum{2}{49} Congruence of segments.
\endhead
\rightheadtext{\S\,2. Congruence of segments.}
\myaxiom{A13} Any straight line segment $[AB]$ is congruent to
itself and for any ray beginning at an arbitrary point $C$ 
there is a unique point $D$ on this ray such that $[AB]\cong 
[CD]$.
\endproclaim
\myaxiom{A14} The binary relation of congruence for segments is
transitive, i\.\,e\. $[AB]\cong [CD]$ and $[CD]\cong [EF]$ imply
$[AB]\cong [EF]$.
\endproclaim
    The reflexivity of the congruence of segments is stated 
explicitly in the axiom~\mytheaxiom{A13}, while the transitivity
of this relation forms the content of the axiom~\mytheaxiom{A14}. 
Let's prove its symmetry.
\mylemma{2.1} The binary relation of congruence for segments is
symmetric, i\.\,e\. $[AB]\cong [CD]$ implies $[CD]\cong [AB]$.
\endproclaim
\demo{Proof} Assume that $[AB]\cong [CD]$. Let's apply the
axiom~\mytheaxiom{A13} to the segment $[CD]$ and to the ray
$[AB\rangle$ beginning at the point $A$. From this axiom
we get that there is a point $E$ on the ray $[AB\rangle$ such
that $[CD]\cong [AE]$. Since $[AB]\cong [CD]$ and $[CD]\cong [AE]$,
due to the axiom~\mytheaxiom{A14} we derive $[AB]\cong [AE]$.\par
    Now we apply the axiom~\mytheaxiom{A13} to the segment $[AB]$ 
and to the ray $[AB\rangle$. It says that the point $E$ on the
ray $[AB\rangle$ such that $[AB]\cong [AE]$ is unique. But
$[AB]\cong [AB]$. Therefore, the point $E$ coincides with $B$. 
Hence, $[CD]\cong [AB]$. Lemma is proved.\qed\enddemo
     Thus, due to the axioms~\mytheaxiom{A13} and \mytheaxiom{A14}
and due to the lemma~\mythelemma{2.1} proved just above the
relation of congruence is an equivalence relation in the set
of all straight line segments.
\myaxiom{A15} Let $B$ be a point lying between $A$ and $C$ on a
straight line $AC$, while $L$ be a point lying between $K$ and $M$ 
on a straight line $KM$. Then the following propositions are valid:
\roster
\item $[AB]\cong [KL]$ and $[BC]\cong [LM]$ imply $[AC]\cong [KM]$;
\item $[AB]\cong [KL]$ and $[AC]\cong [KM]$ imply $[BC]\cong [LM]$.
\endroster
\endproclaim
     Note that the propositions \therosteritem{1} and \therosteritem{2} 
under the assumptions of the axiom~\mytheaxiom{A15} can be complemented
with one more proposition of the same sort:
\roster
\item[3] $[AC]\cong [KM]$ and $[BC]\cong [LM]$ imply $[AB]\cong [KL]$.
\endroster
The proposition~\therosteritem{3} is obtained from the 
proposition~\therosteritem{2} by reformulating it upon exchanging 
the notations of points: $A$ with $C$ and $K$ with $M$.\par 
\parshape 14 0cm 10cm 0cm 10cm 0cm 10cm 0cm 10cm
4.5cm 5.5cm 4.5cm 5.5cm 4.5cm 5.5cm 4.5cm 5.5cm
4.5cm 5.5cm 4.5cm 5.5cm 4.5cm 5.5cm 4.5cm 5.5cm
4.5cm 5.5cm 0cm 10cm
   Under the assumptions of the axiom~\mytheaxiom{A15} the point
$B$ breaks the segment $[AC]$ into two segments $[AB]$ and $[BC]$, 
whose intersection consists of the unique point $B$ and whose union
coincides with the whole segment $[AC]$ (see 
theorems~\mythetheoremchapter{3.2}{2}
and \mythetheoremchapter{3.3}{2} from 
Chapter~\uppercase\expandafter{\romannumeral 2}).
\vadjust{\vskip 5pt\hbox to 0pt{\kern 5pt
\includegraphics{Oris08.eps}\hss}\vskip -5pt}In 
this situation we say that the segment $[AC]$ is {\it composed\/}
of the segments $[AB]$ and $[BC]$ or, in other words, $[AC]$ is
the {\it sum\/} of $[AB]$ and $[BC]$. Therefore, the first 
proposition of the axiom~\mytheaxiom{A15} can be shorten to the
following one: a segment composed of segments congruent to 
$[AB]$ and $[BC]$ is congruent to their sum $[AC]$. If we call
$[BC]$ the {\it difference\/} of the segments $[AC]$ and $[AB]$, 
then the second proposition of the axiom~\mytheaxiom{A15} can be
stated as follows: the difference of segments congruent to $[AC]$ 
and $[AB]$ is congruent to their difference $[BC]$.\par
\mytheorem{2.1} Let the segment $[KM]$ be congruent to the segment
$[AC]$ and let $B$ be an arbitrary point of the line $AC$ distinct
from $A$ and $C$. Then on the line $KM$ there is a unique point
$L$ such that $[KL]\cong [AB]$ and $[LM]\cong [BC]$.
\endproclaim
\demo{Proof} \parshape 7 4.3cm 5.7cm 4.3cm 5.7cm 4.3cm 5.7cm
4.3cm 5.7cm 4.3cm 5.7cm 4.3cm 5.7cm 0cm 10.1cm
Let's consider three possible locations of the point 
\vadjust{\vskip 5pt\hbox to 0pt{\kern 5pt
\includegraphics{Oris09.eps}\hss}\vskip -5pt}$B$
relative to the points $A$ and $C$:
$$
\hskip -2em
\aligned
&(A\blacktriangleright B\blacktriangleleft C),\\
&(C\blacktriangleright A\blacktriangleleft B),\\
&(B\blacktriangleright C\blacktriangleleft A).
\endaligned
\mytag{2.1}
$$
According to theorem~\mythetheoremchapter{2.4}{2} from 
Chapter~\uppercase\expandafter{\romannumeral 2} exactly one of
the conditions \mythetag{2.1} is necessarily fulfilled.
If it is the first condition, we apply the 
axiom~\mytheaxiom{A13} to the segment $[AB]$ and to the ray
$[KM\rangle$. As a result we find a unique point $L$ on the
ray $[KM\rangle$ such that $[KL]\cong [AB]$ (see Fig\.~2.2).
Then we consider the ray coming out from the point $L$ in the
direction opposite to the ray $[LK\rangle$. On this ray we
find a point $\tilde M$ such that $[L\tilde M]\cong [BC]$. 
Then the following conditions are fulfilled:
$$
\xalignat 2
&\hskip -2em
(A\blacktriangleright B\blacktriangleleft C),
&&(K\blacktriangleright L\blacktriangleleft\tilde M).
\mytag{2.2}
\endxalignat
$$
Due to \mythetag{2.2} we can apply the item~\therosteritem{1}
of the axiom~\mytheaxiom{A15} to the points $A$, $B$, $C$, $K$, 
$L$, and $\tilde M$. It yields $[K\tilde M]\cong [AC]$. According
to the premise of the theorem, $[KM]\cong [AC]$. Moreover, by
construction both points $M$ and $\tilde M$ lie on the same ray
$[KL\rangle$ beginning at the point $K$. Hence, due to the
axiom~\mytheaxiom{A13} we derive the coincidence of points
$M=\tilde M$. Then
$$
\xalignat 2
&\hskip -2em
[KL]\cong [AB],&&[LM]\cong [BC].
\mytag{2.3}
\endxalignat
$$
\par\parshape 8 0cm 10cm 5.3cm 4.7cm 5.3cm 4.7cm 5.3cm 4.7cm
5.3cm 4.7cm 5.3cm 4.7cm 5.3cm 4.7cm 0cm 10cm\noindent
This means that $L$ is a required point on the line $KM$. 
\vadjust{\vskip 5pt\hbox to 0pt{\kern -5pt
\includegraphics{Oris10.eps}\hss}\vskip -5pt}Let's
prove its uniqueness. The second condition in \mythetag{2.3} 
admits the existence of
exactly two points $L$ and $\tilde L$ satisfying this condition.
The first of them lies on the ray $[MK\rangle$, it is the point
$L$. The second point $\tilde L$ lies on the opposite ray coming
out from the point $M$. If we assume that the point $L$ is not
fixed uniquely by the conditions \mythetag{2.3}, then the second
point $\tilde L$ should also satisfy both conditions \mythetag{2.3}
simultaneously. Under this assumption on the ray coming out from
the point $\tilde L$ in the direction opposite to the ray 
$[\tilde LK\rangle$ we choose a point $\tilde M$ such that 
$[\tilde L\tilde M]\cong [BC]$. From $[K\tilde L]\cong [AB]$ and
$[\tilde L\tilde M]\cong [BC]$ due to the item~\therosteritem{1}
of the axiom~\mytheaxiom{A15} we derive $[K\tilde M]\cong [AC]$.
But the point $M$, according to the premise if the theorem, 
satisfies the same condition $[KM]\cong [AC]$. This fact contradicts
the axiom~\mytheaxiom{A13}. The coincidence $M=\tilde M$ is
excluded since $M$ and $\tilde M$ by construction lie on different
sides with respect to the point $\tilde L$. The contradiction obtained 
proves the uniqueness of the point $L$ in the case where the first
condition \mythetag{2.1} is fulfilled.\par
\parshape 12 0cm 10.1cm 0cm 10.1cm 4.3cm 5.7cm 4.3cm 5.7cm
4.3cm 5.7cm 4.3cm 5.7cm 4.3cm 5.7cm 4.3cm 5.7cm 4.3cm 5.7cm
4.3cm 5.7cm 4.3cm 5.7cm 0cm 10.1cm
     Now let's study the second case of mutual disposition of the
points $A$, $B$, and $C$ in \mythetag{2.1}. Let's apply the
axiom~\mytheaxiom{A13} to the ray coming out from the point 
\vadjust{\vskip 5pt\hbox to 0pt{\kern 5pt
\includegraphics{Oris11.eps}\hss}\vskip -5pt}$K$ 
in the direction opposite to the ray $[KM\rangle$. As a result 
we get a point $L$ on the line $KM$ such that $[KL]\cong [AB]$. 
Upon combining $[KL]\cong [AB]$ and $[KM]\cong [AC]$ we apply the
item \therosteritem{1} of the axiom~\mytheaxiom{A15}. This yields
$[LM]\cong [BC]$. Thus, $L$ is a required point on the line $KM$. 
The rest is to prove the uniqueness of the point $L$.\par
\par\parshape 9 0cm 10cm 0cm 10cm 6.1cm 3.9cm 6.1cm 3.9cm 6.1cm 3.9cm
6.1cm 3.9cm 6.1cm 3.9cm 6.1cm 3.9cm 0cm 10cm
    The condition $[LM]\cong [BC]$ in \mythetag{2.3} admits the 
existence of exactly two points $L$ on the line $KM$ satisfying 
this condition. The first of them is the above point $L$. The second 
\vadjust{\vskip 5pt\hbox to 0pt{\kern -5pt
\includegraphics{Oris12.eps}\hss}\vskip -5pt}one 
is the point 
$\tilde L$ that lies on the ray coming out from the point $M$ in the 
direction opposite to the ray $[MK\rangle$. The assumption of 
non-uniqueness of the point $L$ means that both points $L$ and
$\tilde L$ satisfy both conditions \mythetag{2.3}. Under this assumption
we apply the axiom~\mytheaxiom{A13} to the ray $[KL\rangle$ and mark
a point $\tilde M$ on this ray such that $[K\tilde M]\cong [AC]$. Now
from $[K\tilde L]\cong [AB]$ and $[K\tilde M]\cong [AC]$, applying the
item \therosteritem{1} of the axiom~\mytheaxiom{A15}, we derive
$[\tilde L\tilde M]\cong [BC]$. Then for two points $M$ and $\tilde M$
on the ray $[\tilde LK\rangle$ we have $[\tilde L\tilde M]\cong
[BC]$ and $[\tilde LM]\cong [BC]$, which contradicts the 
axiom~\mytheaxiom{A13}. This contradiction proves the uniqueness
of $L$ for the second disposition of points in \mythetag{2.1}.\par
     The third disposition of the points $A$, $B$, and $C$ in
\mythetag{2.1} does not require a special treatment. This 
disposition is reduced to the second one by simultaneous exchanging 
the notations of the points: $A$ with $C$ and $K$ with $M$.
\qed\enddemo
     Assume that a segment $[AC]$ on a straight line $a$ is congruent 
to a segment $[KM]$ on another straight line $b$. Using the 
theorem~\mythetheorem{2.1} proved above, one can define a mapping
$f\!:a\to b$ by setting $f(A)=K$, $f(C)=M$, and determining $f(X)$ 
by means of the conditions $[AX]\cong [Kf(X)]$ and $[CX]\cong [Mf(X)]$ 
for all other points \par
     The mapping $h\!:b\to a$ is constructed in a similar way. For
it we set $h(K)=A$, $h(M)=C$, and we set $[KZ]\cong [Ah(Z)]$ and
$[MZ]\cong [Ch(Z)]$ for all other points $Z\in b$. Due to the 
uniqueness of the point $L$ in the theorem~\mythetheorem{2.1} the
mappings $f$ and $h$ appear to be inverse to each other. In particular,
this means that both of them are bijective.\par
     Note, that following the proof of the theorem~\mythetheorem{2.1},
one can establish the following fact characterizing $f$:
$$
\hskip -2em
\aligned
&(A\blacktriangleright X\blacktriangleleft C)
\text{\ \ implies \ }(K\blacktriangleright f(X)\blacktriangleleft M),\\
&(C\blacktriangleright A\blacktriangleleft X)
\text{\ \ implies \ }(M\blacktriangleright K\blacktriangleleft f(X)),
\\
&(X\blacktriangleright C\blacktriangleleft A)
\text{\ \ implies \ }(f(X)\blacktriangleright M\blacktriangleleft K).
\endaligned
\mytag{2.4}
$$
In other words, the mutual disposition of the points $K$, $M$, and 
$f(X)$ on the line $b$ mimics the disposition if the initial points 
$A$, $C$, and $X$ on the line $a$.
\mytheorem{2.2} Assume that a segment $[KM]$ of a straight line $b$ is
congruent to a segment $[AC]$ of a line $a$. Let's denote $f\!:a\to b$
the mapping given by the relationships $f(A)=K$, $f(C)=M$, and by the
conditions $[AX]\cong[Kf(X)]$ and $[CX]\cong [Mf(X)]$ for all $X\in a$
distinct from $A$ and $C$. Under these assumptions, if we introduce
distinguished directions on the lines $a$ and $b$ by virtue of the
vectors $\overrightarrow{AC}$ and $\overrightarrow{KM}$, then for any
two points $X$ and $Y$ on the line $a$ the following conditions are
fulfilled:
\roster
\item $X\prec Y$ implies $f(X)\prec f(Y)$;
\item the segment $[XY]$ is congruent to the segment $[f(X)f(Y)]$.
\endroster
\endproclaim
\demo{Proof} Let's study various cases of mutual disposition of the
points $X$, $Y$, $A$, and $C$ on the line $a$. If the point $X$ or
the point $Y$ coincides with one of the points $A$ or $C$, in this case
the proposition of the theorem easily follows from \mythetag{2.4} and 
from the way in which the mapping $f$ is defined. For this reason, 
without loss of generality we can assume that $X$, $Y$, $A$, and $C$ 
are four distinct points on the line $a$.\par
\vskip 1pt plus 1pt minus 1pt
     The points $A$ and $C$ break the line $a$ into three parts:
the ray $(-\infty,A]$, the segment $[AC]$, and the ray $[C,+\infty)$. 
If $X$ and $Y$ are in different parts of this partitioning, then the
interval $(XY)$ contains at least one of the points $A$ or $C$. In
this case the item \therosteritem{1} of the theorem can be derived from
\mythetag{2.4}. The item \therosteritem{2} of the 
theorem~\mythetheorem{2.2} after that is proved by applying the first
item from the axiom~\mytheaxiom{A15}.\par
\vskip 1pt plus 1pt minus 1pt
\parshape 12 0cm 10cm 0cm 10cm 4.5cm 5.5cm 4.5cm 5.5cm
4.5cm 5.5cm 4.5cm 5.5cm 4.5cm 5.5cm 4.5cm 5.5cm 4.5cm 5.5cm
4.5cm 5.5cm 4.5cm 5.5cm 0cm 10cm
    Now let's consider the case where the points $X$ and $Y$ lie on 
the segment $[AC]$. Then $X\prec Y$ implies $A\prec X\prec Y\prec C$.
\vadjust{\vskip 5pt\hbox to 0pt{\kern -5pt
\includegraphics{Oris13.eps}\hss}\vskip -5pt}Applying the
axiom~\mytheaxiom{A13}, we choose a point $\tilde Y$
on the ray $[f(X)M\rangle$ such that $[f(X)\tilde Y]\cong [XY]$. 
Then, using the same axiom~\mytheaxiom{A13}, we draw the segment
$[\tilde Y\tilde C]\cong [YC]$ on the ray coming out from the point
$\tilde Y$ in the direction opposite to the ray $[\tilde YK\rangle$.
Applying the item \therosteritem{1} of the axiom~\mytheaxiom{A15} to
the segments $[Kf(X)]$ and $[f(X)\tilde Y]$, we derive $[K\tilde Y]
\cong [AY]$. Then we apply the same item of the 
axiom~\mytheaxiom{A15} to the segments $[K\tilde Y]$ and 
$[\tilde Y\tilde C]$. It yields $[K\tilde C]\cong [AC]$.
Hence, $\tilde C=M$, which implies $[\tilde YM]\cong [YC]$. From
this relationship we derive the coincidence $\tilde Y=f(Y)$.\par
\vskip 1pt plus 1pt minus 1pt
     The coincidences $\tilde C=M$ and $\tilde Y=f(Y)$ due to the
above construction yield $K\prec f(X)\prec f(Y)\prec M$. This
proves the first proposition of the theorem $f(X)\prec f(Y)$. 
The second proposition $[f(X)f(Y)]\cong [XY]$ follows from
$[f(X)\tilde Y]\cong [XY]$ and from the coincidence $\tilde Y=f(Y)$.
\par
\vskip 1pt plus 1pt minus 1pt
\parshape 12 0cm 10cm 0cm 10cm 0cm 10cm 0cm 10cm 0cm 10cm
0cm 10cm 4.5cm 5.5cm 4.5cm 5.5cm 4.5cm 5.5cm 4.5cm 5.5cm 
4.5cm 5.5cm 4.5cm 5.5cm
     Now let's consider the case, where the points $X$ and $Y$ lie
on the ray $[C,+\infty)$. Here from $X\prec Y$ we derive $A\prec C
\prec X\prec Y$. Let's apply the axiom~\mytheaxiom{A13} to the ray
coming out from the point $f(X)$ in the direction opposite to
the ray $[f(X)M]$. We mark on this ray a point $\tilde Y$ such that
$[f(X)\tilde Y]\cong [XY]$. Then we apply the first item of the
axiom~\mytheaxiom{A15} to the segments $[Mf(X)]$ and $[f(X)\tilde Y]$. 
This yields $[M\tilde Y]\cong [CY]$, which in turn leads to the
coincidence of points $\tilde Y$ and $f(Y)$. Now the relationships 
$f(X)\prec f(Y)$ and $[f(X)f(Y)]\cong [XY]$ are 
\vadjust{\vskip 5pt\hbox to 0pt{\kern -5pt
\includegraphics{Oris14.eps}\hss}\vskip -5pt}fulfilled 
by construction of the point $\tilde Y$.\par
\parshape 4 4.5cm 5.5cm 4.5cm 5.5cm 4.5cm 5.5cm 0cm 10cm
     The last case, where the points $X$ and $Y$ lie on the ray
$(-\infty,A]$, does not require a special consideration. It is
reduced to the second case if we exchange the notations of points
$A$ with $C$ and $K$ with $M$ and if we change the distinguished 
directions on $a$ and $b$ for the opposite ones.
\qed\enddemo
\head
\SectionNum{3}{56} Congruent translation of straight lines.
\endhead
\rightheadtext{\S\,3. Congruent translation of straight lines.}
\mydefinition{3.1} A mapping $f\!:a\to b$ is called a 
{\it congruent translation} of a straight line $a$ to a straight 
line $b$ if for any two points $X$ and $Y$ on the line $a$ the
condition of congruence $[f(X)f(Y)]\cong [XY]$ is fulfilled.
\enddefinition
     Let $f$ and $h$ be two congruent translations of a straight
line $a$ to a straight line $b$. If at some two points $A$ and $B$ 
on te line $a$ these mappings coincide
$$
\xalignat 2
&f(A)=h(A),&&f(B)=h(B),
\endxalignat
$$
then they coincide at all points $X\in a$, i.\,e\. $f=h$. This fact
is easily derived from the theorem~\mythetheorem{2.1}. This theorem
together with the theorem~\mythetheorem{2.2} show that congruent
translations of lines do exist. Indeed, in order to define such a
mapping $f\!:a\to b$ it is sufficient to choose two points $A$ and
$B$ on the line $a$ and construct the segment $[KM]$ congruent to
$[AB]$ on the line $b$.\par
     Assume that on a straight line $a$ some point $O$ is marked
and one of two possible directions is distinguished. Then $a$ can
be broken into two rays $[O,-\infty)$ and $[O,+\infty)$. Assume that
oh another straight line $b$ some point $Q$ is marked and some
direction is distinguished. Then this line is also broken into two
rays $[Q,-\infty)$ and $[Q,+\infty)$. Now we choose a point $E_{+}$
on the ray $[O,+\infty)$ and, using the axiom~\mytheaxiom{A13}, on
the line $b$ we construct two segments $[QF_{+}]$ and $[QF_{-}]$ 
congruent to the segment $[OE_{+}]$. The segment $[QF_{+}]$ lies
on the ray $[Q,+\infty)$, while the segment $[QF_{-}]$ lies in the
opposite ray $[Q,-\infty)$. Because of the presence of two segments
congruent to $[OE_{+}]$ we can define two mappings of congruent 
translation of the line $a$ to the line $b$:
$$
\xalignat 2
&f^{+}(O)=Q,&&f^{+}(E_{+})=F_{+},\\
\vspace{-1.5ex}
&&&\mytag{3.1}\\
\vspace{-1.5ex}
&f^{-}(O)=Q,&&f^{-}(E_{+})=F_{-}.
\endxalignat
$$
\mytheorem{3.1} For any point $O$ on a straight line $a$ with a
distinguished direction and for any point $Q$ on another straight
line $b$ with a distinguished direction there are exactly two 
mappings $f\!:a\to b$ performing congruent translation of the
line $a$ to the line $b$. The first of them $f^{+}_{\sssize OQ}$ 
preserves the precedence of points, i\.\,e\. $X\prec Y$ implies
$f^{+}_{\sssize OQ}(X)\prec f^{+}_{\sssize OQ}(Y)$. The second 
mapping inverts the precedence of points, i\.\,e\. $X\prec Y$ 
implies $f^{-}_{\sssize OQ}(Y)\prec f^{-}_{\sssize OQ}(X)$.
\endproclaim
\myexercise{3.1} Prove the theorem~\mythetheorem{3.1} by showing
that the mappings $f^{+}_{\sssize OQ}$ and $f^{-}_{\sssize OQ}$ 
do not depend on a particular choice of the point $E_{+}\in [O,+
\infty)$ in formulas \mythetag{3.1} defining them. Moreover, 
show that these mappings remain unchanged if one changes 
simultaneously the distinguished directions on the lines 
$a$ and $b$ for opposite ones.
\endproclaim
     Assume that the line $b$ coincides with the line $a$. We choose
two points $O$ and $Q$ and fix one of two possible directions on this
line. In this case the mapping $f^{+}_{\sssize OQ}$ is called the {\it
congruent translation by the vector $\overrightarrow{OQ}$}. It is 
denoted $f^{+}_{\sssize OQ}=p_{\sssize OQ}$. The coincidence $O=Q$ 
makes a special case. In this special case the points $O$ and $Q$ do
not define a vector (understood as an arrowhead segment), while the
the mapping of congruent translation $p_{\sssize OO}$ appears to be
the identical mapping: $p_{\sssize OO}=\id$. The mapping 
$f^{-}_{\sssize OQ}$ is not identical even in the case of coinciding
points $O$ and $Q$. For $O=Q$ the mapping $f^{-}_{\sssize OO}$ is called
the {\it inversion with respect to the point $O$}. It is denoted as 
$f^{-}_{\sssize OO}=i_{\sssize O}.$
\mytheorem{3.2} The mappings of congruent translation by vectors and
the mappings of inversion satisfy the relationships
$$
\xalignat 2
&p_{\sssize BC}\compos p_{\sssize AB}=p_{\sssize AC},
&&p_{\sssize AB}\compos i_{\sssize C}=i_{\sssize C}\compos p_{\sssize BA},\\
&i_{\sssize C}\compos i_{\sssize C}=\id,
&&i_{\sssize A}\compos i_{\sssize B}=p_{\sssize BC}\text{, \ where \ }
C=i_{\sssize A}(B).
\endxalignat
$$
\endproclaim
\mytheorem{3.3} Assume that on each of two straight lines $a$ and $b$ 
with distinguished directions two points are fixed: $O,\tilde O\in a$ 
and $Q, \tilde Q\in b$. Then we have
$$
\xalignat 2
&f^{+}_{\sssize\tilde O\tilde Q}\compos=p_{\sssize Q\tilde Q}\compos
f^{+}_{\sssize OQ}\compos p_{\sssize \tilde OO},
&&f^{-}_{\sssize\tilde O\tilde Q}=p_{\sssize Q\tilde Q}\compos
f^{-}_{\sssize OQ}\compos p_{\sssize \tilde OO}.
\endxalignat
$$
\endproclaim
\noindent{\bf A remark.} By means of the small circle in
theorems~\mythetheorem{3.2} and \mythetheorem{3.3} we denote 
the operation of composing two mappings: $f\compos h(x)=f(h(x))$
(see \myhref{\diez s7pg12}{\S\,7} in
Chapter~\uppercase\expandafter{\romannumeral 1}).
\myexercise{3.2} Relying on the theorems~2.1, 2.2, and 3.1, prove
the properties of mappings of congruent translations of lines
stated in theorems~\mythetheorem{3.2} and \mythetheorem{3.3}.
\endproclaim
     The domain of any of the above mappings of congruent translation 
is some straight line. At this moment we have no tools for extending 
this domain. The only exception is the inversion $i_{\sssize O}$. Assume 
that some point $O$ in the space is fixed. For any point $X$ different
from $O$ there is a unique line $a=OX$ passing through $O$ and $X$.
On this line the inversion mapping $i_{\sssize O}$ is defined. Let's
set $i(X)=i_{\sssize O}(X)$. For the point $O$ itself we set $i(O)=O$. 
As a result we get the mapping $i\!:\Bbb E\to\Bbb E$ which is called the 
{\it inversion\/} or the {\it central symmetry\/} with the center at the
point $O$.
\head
\SectionNum{4}{59} Slipping vectors. Addition of vectors\\
on a straight line.
\endhead
\rightheadtext{\S\,4. Slipping vectors. Addition of vectors \dots}
\mydefinition{4.1} Two vectors $\overrightarrow{AB}$ and
$\overrightarrow{CD}$ lying on one straight line are called {\it
equal\/} if they are codirected and if the segment $[AB]$ is 
congruent to the segment $[CD]$.
\enddefinition
     The equality of points, straight lines, planes, segments, and many
other geometric forms is understood as pure coincidence. The equality of
vectors, according to the definition~\mythedefinition{4.1}, is of different 
nature.
\myexercise{4.1} Verify that the relation of equality of vectors is 
a binary relation of equivalence.
\endproclaim
\noindent Vectors understood as arrowhead segments are sometimes called 
{\it geometric vectors}. They have strictly fixed positions in the space.
In contrast to geometric vectors, a {\it slipping vector\/} on a straight
line is a class of mutually equivalent vectors in the sense of the
definition~\mythedefinition{4.1}. A slipping vector has many representatives
lying on a given line. They are called {\it geometric realizations\/} of
this slipping vector.
\mytheorem{4.1} For any four points $A$, $B$, $C$, and $D$ lying on one
straight line $\overrightarrow{AB}=\overrightarrow{CD}$ implies
$\overrightarrow{AC}=\overrightarrow{BD}$ and, conversely,
$\overrightarrow{AC}=\overrightarrow{BD}$ implies $\overrightarrow{AB}
=\overrightarrow{CD}$.
\endproclaim
\demo{Proof} Let's consider the first proposition of the theorem. 
Assume that $\overrightarrow{AB}=\overrightarrow{CD}$. We choose
the direction of the vector $\overrightarrow{AB}$ for the positive
direction on the line where both vectors $\overrightarrow{AB}$ and
$\overrightarrow{CD}$ lie. Then the following relationships are
fulfilled:
$$
\xalignat 2
&\hskip -2em
A\prec B,&&C\prec D.
\mytag{4.1}
\endxalignat
$$
From \mythetag{4.1} we derive the complete list of possible mutual
dispositions of the points $A$, $B$, $C$, and $D$:
$$
\xalignat 2
&\hskip -2em
A\prec B\prec C\prec D,&&C\prec D\prec A\prec B,\qquad
\mytag{4.2}\\
&\hskip -2em
A\prec C\prec B\prec D,&&C\prec A\prec D\prec B,\qquad
\mytag{4.3}\\
&\hskip -2em
A\prec C\prec D\prec B,&&C\prec A\prec B\prec D.\qquad
\mytag{4.4}
\endxalignat
$$
All of the above cases in \mythetag{4.2}, \mythetag{4.3}, and 
\mythetag{4.4} are grouped into pairs. We can study only one
case in each pair since the other case is obtained by transposition 
of vectors $\overrightarrow{AB}$ and $\overrightarrow{CD}$, which 
does not change the proposition of the theorem in whole.\par
\parshape 9 0cm 10cm 0cm 10cm 5.1cm 4.9cm 5.1cm 4.9cm 5.1cm 4.9cm
5.1cm 4.9cm 5.1cm 4.9cm 5.1cm 4.9cm 0cm 10cm
     Let's show that the case \mythetag{4.4} is impossible. It is not
compatible with the condition $[AB]\cong [CD]$ which follows
\vadjust{\vskip 5pt\hbox to 0pt{\kern -5pt
\includegraphics{Oris15.eps}\hss}\vskip -5pt}from the 
equality $\overrightarrow{AB}=\overrightarrow{CD}$. Applying the
axiom~\mytheaxiom{A13}, on the right of $B$ we choose 
a point $\tilde D$ such that $[B\tilde D]\cong [AC]$. Combining
this condition with the condition $[CD]\cong [AB]$, from the 
axiom~\mytheaxiom{A15} we derive the congruence of the segments
$[AD]\cong [A\tilde D]$. This relationship contradicts the
axiom~\mytheaxiom{A13} since both points $D$ and $\tilde D$ lie
on the same ray coming out from the point $A$. The contradiction 
obtained excludes the case \mythetag{4.4} from our further
consideration.\par
\parshape 9 0cm 10.1cm 0cm 10.1cm 4.3cm 5.7cm 4.3cm 5.7cm 4.3cm 5.7cm
4.3cm 5.7cm 4.3cm 5.7cm 4.3cm 5.7cm 0cm 10.1cm
    Note that the first proposition of the theorem for the case
\mythetag{4.3} is equivalent to the second proposition for the 
case \mythetag{4.2}. Therefore, it is sufficient to consider 
\vadjust{\vskip 5pt\hbox to 0pt{\kern -5pt
\includegraphics{Oris16.eps}\hss}\vskip -5pt}only the 
case $A\prec B\prec C\prec D$ and prove both propositions of the 
theorem for this case. The codirectedness conditions
$\overrightarrow{AB}\upuparrows\overrightarrow{CD}$ and
$\overrightarrow{AC}\upuparrows\overrightarrow{BD}$ follow from the
disposition of points $A\prec B\prec C\prec D$. From $[AB]\cong [CD]$ 
and from the obvious relationship $[BC]\cong [BC]$, applying the
first item of the axiom~\mytheaxiom{A15}, we derive $[AC]\cong [BD]$. 
Conversely, from $[AC]\cong [BD]$ and $[BC]\cong [BC]$, upon applying
the item~\therosteritem{2} of the axiom~\mytheaxiom{A15}, we get
$[AB]\cong [CD]$. The theorem is proved.
\qed\enddemo
\mytheorem{4.2} The equality $p_{\sssize AB}=p_{\sssize CD}$ is valid
if and only if $\overrightarrow{AB}=\overrightarrow{CD}$ in the sense
of the definition~\mythedefinition{4.1}.
\endproclaim
\demo{Proof} Assume that $p_{\sssize AB}$ and $p_{\sssize CD}$ are
the mappings of congruent translation on a straight line $a$ and
assume that $p_{\sssize AB}=p_{\sssize CD}=p$. Let's define the
positive direction on the line $a$ by means of the vector
$\overrightarrow{AC}$ and thus define a relation of precedence for
the points of this line $a$. Then $A\prec C$. Let's apply the mapping
$p$ to the points $A$ and $C$ and use the theorem~\mythetheorem{2.2}:
$$
\xalignat 2
&p(A)\prec p(C),
&&[p(A)p(C)]\cong [AC].
\endxalignat
$$
But $p(A)=p_{\sssize AB}(A)=B$ and $p(C)=p_{\sssize CD}(C)=D$. Hence,
$A\prec C$ and $B\prec D$, which means that the vectors
$\overrightarrow{AC}$ and $\overrightarrow{BD}$ are codirected. 
Moreover, $[BD]\cong [AC]$, therefore, $\overrightarrow{AC}=
\overrightarrow{BD}$. Applying the theorem~\mythetheorem{4.1}, we
derive the required equality of vectors $\overrightarrow{AB}
=\overrightarrow{CD}$.\par
     Now, conversely, assume that $\overrightarrow{AB}
=\overrightarrow{CD}$. From this equality, applying the
theorem~\mythetheorem{4.1}, we derive $\overrightarrow{AC}
=\overrightarrow{BD}$. Hence, we have $[AC]\cong [BD]$, while from 
$A\prec C$ it follows that $B\prec D$. Let's apply the mapping
$p_{\sssize AB}$ to the point $C$ and denote $\tilde D
=p_{\sssize AB}(C)$. Then, according to the theorem~\mythetheorem{2.2},
we get $[AC]\cong [B\tilde D]$, while $A\prec C$ implies $B\prec\tilde D$.
From $[AC]\cong [B\tilde D]$ and $[AC]\cong [BD]$ we derive 
$[BD]\cong [B\tilde D]$ and from $B\prec D$ and $B\prec\tilde D$ 
we conclude that the points $D$ and $\tilde D$ lie on the same ray
coming out from the point $B$. Hence, $D=\tilde D$, which follows
from the axiom~\mytheaxiom{A13}. Now $D=p_{\sssize AB}(C)$. This
fact yields the coincidence $p_{\sssize AB}=p_{\sssize CD}$.
\qed\enddemo
\mytheorem{4.3} Any two mappings of congruent translation on the same
straight line do commute: $p_{\sssize AB}\compos p_{\sssize CD}=
p_{\sssize CD}\compos p_{\sssize AB}$.
\endproclaim
\demo{Proof} Let's choose some arbitrary point $E$ on the line $a$
on which the vectors $\overrightarrow{AB}$ and $\overrightarrow{CD}$
lie. Then denote $F=p_{\sssize AB}(E)$, $G=p_{\sssize CD}(F)$, and
$H=p_{\sssize CD}(E)$. As a result we get
$$
\xalignat 2
&\hskip -2em
p_{\sssize AB}=p_{\sssize EF},&&p_{\sssize CD}=p_{\sssize FG}
=p_{\sssize EH}.
\mytag{4.5}
\endxalignat
$$
To the last equality $p_{\sssize FG}=p_{\sssize EH}$ in 
\mythetag{4.5} the previous theorem~\mythetheorem{4.2} is 
applicable. It yields $\overrightarrow{FG}=\overrightarrow{EH}$. 
We apply the theorem~\mythetheorem{4.1} to this equality.
It yields $\overrightarrow{EF}=\overrightarrow{HG}$. Hence,
$p_{\sssize EF}=p_{\sssize HG}$. By means of direct calculations
we derive
$$
\hskip -2em
\aligned
&p_{\sssize CD}\compos p_{\sssize AB}=p_{\sssize FG}\compos
p_{\sssize EF}=p_{\sssize EG},\\
&p_{\sssize AB}\compos p_{\sssize CD}=p_{\sssize HG}\compos
p_{\sssize EH}=p_{\sssize EG}.
\endaligned
\mytag{4.6}
$$
Here in \mythetag{4.6} we used the theorem~\mythetheorem{3.2}. The 
rest to compare the right hand sides of the formulas \mythetag{4.6}, 
which immediately yields the required result $p_{\sssize AB}\compos
p_{\sssize CD}=p_{\sssize CD}\compos p_{\sssize AB}$.
\qed\enddemo
    Some fixed mapping of congruent translation on a straight 
line can be given by various pairs of points. However, the 
theorem~\mythetheorem{4.2} show that all such pairs of correspond
to geometric vectors equal to each other in the sense of the
definition~\mythedefinition{4.1}. Therefore, passing from geometric
vectors to slipping vectors, we get the one-to-one correspondence
of the set of congruent translations and the set of slipping vectors
on a line: $p=p_a$.\par
     The identical mapping $\id$ is also a mapping of congruent
translation: $\id=p_{\sssize AA}$. But a single point $A$ does not
define an arrowhead segment. Especially for to describe this 
situation the concept of {\it zero vector\/} is introduced. The
zero vector $\bold 0$ is a formal object complementing the set 
of slipping vectors on a line so that $p_{\bold 0}=\id$. Any 
one-point set on a line treated as a 
{\tencyr\char '074}degenerate{\tencyr\char '076} arrowhead segment
$\overrightarrow{AA}$ can be taken for a geometric realization of
the zero vector.\par
      The set of mappings of congruent translation is naturally 
equipped with the operation of composition. According to the
theorem~\mythetheorem{3.2}, the composition of two congruent
translations is a congruent translation. Let's set by definition
$$
\hskip -2em
p_{\bold a}\compos p_{\bold b}=p_{\bold a+\bold b}.
\mytag{4.7}
$$
The formula \mythetag{4.7} is the definition of the addition operation
for slipping vectors on a straight line.
\mytheorem{4.4} The addition operation of slipping vectors on a straight 
line possesses the following properties:
\roster
\item it is commutative, i\.\,e\. $\bold a+\bold b=\bold b+\bold a$;
\item it is associative, i\.\,e\. $(\bold a+\bold b)+\bold c=
      \bold a+(\bold b+\bold c)$;
\item there is a vector $\bold 0$, such that $\bold 0+\bold a=
      \bold a+\bold 0=\bold a$ for an arbitrary vector $\bold a$;
\item for any vector $\bold a$ there is an opposite vector $\bold a'$
      such that $\bold a+\bold a'=\bold a'+\bold a=\bold 0$.
\endroster
\endproclaim
     The first property follows from the theorem~\mythetheorem{4.3}.
The associativity follows from the formula \mythetag{4.7} since the
property of associativity is peculiar to the composition of any 
mappings. The third property follows from the formula \mythetag{4.7}
and from the definition of the zero vector $p_{\bold 0}=\id$. The 
rest is to prove the fourth property. Let $\overrightarrow{AB}$ be
a geometric realization for a slipping vector $\bold a$. Let's denote
by $\bold a'$ the slipping vector whose geometric realization is
the vector $\overrightarrow{BA}$. Then $p_{\bold a+\bold a'}
=p_{\bold a}\compos p_{\bold a'}=p_{\sssize AB}\compos 
p_{\sssize BA}=p_{\sssize BB}=\id=p_{\bold 0}$. This means that, 
$\bold a+\bold a'=\bold 0$.\par
     The addition is an algebraic operation on the set of slipping 
vectors. Sets equipped with various algebraic operations are 
studied in course of general algebra (see, for instance, \mycite{5}).
Let's recall the definition of a group --- it is a set with one
algebraic operation, which is usually called the {\it group 
multiplication}.
\mydefinition{4.2} A set $G$ is called a {\it group\/} if for
any two elements $a$ and $b$ of this set a third element of this 
set $a\cdot b$, which is called the {\it product\/} of $a$ and
$b$, is assigned so that the following three conditions are 
fulfilled:
\roster
\item $(a\cdot b)\cdot c=a\cdot (b\cdot c)$ i\.\,e\. 
      the group multiplication is associative;
\item there is an element $e\in G$ such that $e\cdot a=a\cdot e=a$
      for any element $a\in G$;
\item for any element $a\in G$ there is an element $a'$ such that
      $a\cdot a'=a'\cdot a=e$.
\endroster
The element $e$ is called the {\it unity\/} of the group $G$, while 
the element $a'$ is called the {\it inverse element\/} for an 
element $a\in G$.
\enddefinition
\mydefinition{4.3} A group $G$ is called {\it commutative\/} or
{\it Abelian} if the group multiplication in it is commutative, 
i\.\,e\. $a\cdot b=b\cdot a$.
\enddefinition
     Comparing the definitions~\mythedefinition{4.2} and
\mythedefinition{4.3} with the properties of the addition of vectors
in the theorem~\mythetheorem{4.4} shows that the set of slipping 
vectors is an Abelian group with respect to the addition. The
multiplication sign in it is replaced by the plus sign, while the
zero vector plays the role of the unity.
\head
\SectionNum{5}{64} Congruence of angles.
\endhead
\rightheadtext{\S\,5. Congruence of angles.}
    Let $h$ and $k$ be two rays coming out from one point and not lying
on one straight line. Let's choose a point $A$ on the ray $h$ and a point
$B$ on the ray $k$. Assume that the points $A$ and $B$ are distinct from 
the point $O$. Let's apply the theorem~\mythetheoremchapter{5.2}{2} from
Chapter~\uppercase\expandafter{\romannumeral 2} to the points $A$, $B$,
and $O$ and construct the angle $\angle AOB$. It is clear that this angle
does not depend on a particular choice of the points $A$ and $B$ on the
rays $h$ and $k$. It is determined by the rays $h$ and $k$ themselves. 
For this reason we shall denote such an angle as $\angle hOk$ or even as
$\angle hk$.
\myaxiom{A16} Any angle $\angle hk$ is congruent to itself and for any 
half-plane $a_{+}$ with a ray $m$ lying on the boundary line $a$ there 
is a unique ray $n$ within the half-plane $a_{+}$ such that 
$\angle hk\cong\angle mn$.
\endproclaim
\myaxiom{A17} Let $A$, $B$, and $C$ be three points not lying on one 
straight line and let $\tilde A$, $\tilde B$, and  $\tilde C$ be other
three points also not lying on one straight line. If the conditions
$$
\xalignat 3
&[AB]\cong [\tilde A\tilde B],
&&[AC]\cong [\tilde A\tilde C],
&&\angle BAC\cong\angle \tilde B\tilde A\tilde C
\endxalignat
$$
are fulfilled, then the other two conditions $\angle ABC\cong
\angle \tilde A\tilde B\tilde C$ and $\angle ACB\cong\angle\tilde A
\tilde C\tilde B$ are also fulfilled.
\endproclaim
    The axiom~\mytheaxiom{A16} for angles is analogous to the
axiom~\mytheaxiom{A13} for segments. The axiom~\mytheaxiom{A17}
can be formulated as follows: if an angle and two sides forming 
this angle in some triangle are congruent to an angle and to
the corresponding sides of some other triangle, then the remaining
two angles of the first triangle are congruent to the corresponding
angles of the second triangle. Let's define the concept of congruence
for triangles.
\mydefinition{5.1} Two triangles are called {\it congruent\/} if their
vertices are in one-to-one correspondence so that the angles at the
vertices and the sides of one triangle are congruent to the corresponding
angles and sides of the other triangle.
\enddefinition
For example, the triangle $ABC$ is congruent to the other triangle
$FGH$ if the following six conditions are fulfilled:
$$
\xalignat 3
[AB]&\cong [FG],
&[BC]&\cong [GH],
&[CA]&\cong [HF],\\
\angle ABC&\cong\angle FGH,
&\angle BCA&\cong\angle GHF,
&\angle CAB&\cong\angle HFG.
\endxalignat
$$
\mytheorem{5.1} If for triangles $ABC$ and $\tilde A\tilde B\tilde C$ 
the conditions $[AB]\cong [\tilde A\tilde B]$, $[AC]\cong [\tilde A
\tilde C]$, and $\angle BAC\cong\angle \tilde B\tilde A\tilde C$ are
fulfilled, then the triangle $ABC$ is congruent to the triangle 
$\tilde A\tilde B\tilde C$.
\endproclaim
\demo{Proof}\parshape 12 0cm 10cm 0cm 10cm 4.5cm 5.5cm 4.5cm 5.5cm 
4.5cm 5.5cm 4.5cm 5.5cm 4.5cm 5.5cm 4.5cm 5.5cm 4.5cm 5.5cm 
4.5cm 5.5cm 4.5cm 5.5cm 0cm 10cm
The congruence of three corresponding angles and the congruence of two
corresponding sides in the triangles $ABC$ and $\tilde A\tilde B\tilde C$
easily follow from the axiom~\mytheaxiom{A17}. The rest is to prove the
congruence $[BC]\cong[\tilde B\tilde C]$. Using the
axiom~\mytheaxiom{A13}, on the ray $[\tilde B\tilde C\rangle$ we find 
\vadjust{\vskip 5pt\hbox to 0pt{\kern -5pt
\includegraphics{Oris17.eps}\hss}\vskip -5pt}a point $C'$ 
such that $[BC]\cong [\tilde BC']$ and consider the
triangle $\tilde A\tilde BC'$. For this triangle the conditions
$$
\align
&[AB]\cong [\tilde A\tilde B],\\
&[BC]\cong [\tilde BC'],\\
&\angle ABC\cong\angle\tilde A\tilde BC',
\endalign
$$
are fulfilled. Due to these conditions one can apply the 
axiom~\mytheaxiom{A17} to the angle $\angle\tilde A\tilde BC'$ 
and to the sides $[\tilde A\tilde B]$ and $[\tilde BC']$ of the
triangle $\tilde A\tilde BC'$. It yields $\angle BAC\cong
\angle\tilde B\tilde AC'$. Moreover, in the statement of the theorem 
we have $\angle BAC\cong\angle\tilde B\tilde A\tilde C$. Therefore, 
if we assume that $\tilde C\neq C'$, we would have two rays
$[\tilde A\tilde C\rangle$ and $[\tilde AC'\rangle$ on one 
half-plane bounded by the line $\tilde A\tilde B$ which form
with the ray $[\tilde A\tilde B\rangle$ two angles congruent to
the angle $\angle BAC$. But this would contradict the
axiom~\mytheaxiom{A16}, hence, $\tilde C=C'$. As a result we get
the required congruence of segments $[BC]\cong [\tilde B\tilde C]$. 
Thus, the theorem is proved.
\qed\enddemo
The theorem~\mythetheorem{5.1} is known as the congruence criterion
for triangles by two sides and the angle between them. The next theorem
is called the congruence criterion for triangles by a side and two 
angles adjoint to this side.
\mytheorem{5.2} If for triangles $ABC$ and $\tilde A\tilde B\tilde C$
the conditions $[AB]\cong [\tilde A\tilde B]$, $\angle ABC\cong\angle
\tilde A\tilde B\tilde C$, and $\angle BAC\cong\angle\tilde B\tilde A
\tilde C$ are fulfilled, then the triangle $ABC$ is congruent to the
triangle $\tilde A\tilde B\tilde C$.
\endproclaim
\demo{Proof} Applying the axiom~\mytheaxiom{A13}, on the ray 
$[\tilde B\tilde C\rangle$ we choose a point $C'$ such that
$[BC]\cong [\tilde BC']$. Then to the triangles $ABC$ and $\tilde A
\tilde BC'$ the previous theorem~\mythetheorem{5.1} is applicable. 
It means that the triangle $ABC$ is congruent to the triangle
$\tilde A\tilde BC'$, hence, $\angle BAC\cong\angle\tilde B
\tilde AC'$. From the statement of the theorem we get 
$\angle BAC\cong\angle\tilde B\tilde A\tilde C$. Now, if we assume that 
$C'\neq\tilde C$, then on the half-plane bounded by the line 
$\tilde A\tilde B$ we would have two rays $[\tilde A\tilde C\rangle$ 
and $[\tilde AC'\rangle$ which form with the ray $[\tilde A\tilde B
\rangle$ two angles congruent to the angle $\angle BAC$. This would 
be a contradiction to the axiom~\mytheaxiom{A16}. Hence, $C'=\tilde C$ 
and the triangle $ABC$ is congruent to the triangle $\tilde A
\tilde B\tilde C$.
\qed\enddemo
\mytheorem{5.3} Let $h$, $k$, and $l$ be three distinct rays coming 
out from one point $O$ and lying on one plane. Let $h'$, $k'$, and 
$l'$ be other three distinct rays \pagebreak coming out from one point 
$O'$ and lying on one plane. If the ray $l$ is inside the angle 
$\angle hk$ and if the ray $l'$ is inside the angle $\angle h'k'$, 
then
\roster
\item $\angle hl\cong\angle h'l'$ and $\angle lk\cong\angle l'k'$
      imply $\angle hk\cong\angle h'k'$;
\item $\angle hk\cong\angle h'k'$ and $\angle hl\cong\angle h'l'$
      imply $\angle lk\cong\angle l'k'$;
\item $\angle hk\cong\angle h'k'$ and $\angle lk\cong\angle l'k'$
      imply $\angle hl\cong\angle h'l'$.
\endroster
\endproclaim
\demo{Proof}\parshape 7 0cm 10cm 0cm 10cm 0cm 10cm 0cm 10cm
5cm 5cm 5cm 5cm 5cm 5cm
Let's choose a point $A$ on the ray $h$ and a point $C$ on the ray
$k$. We connect them by means of the segment $[AC]$ and then we apply 
the lemma~\mythelemmachapter{6.2}{2} from 
Chapter~\uppercase\expandafter{\romannumeral 2} to the ray $l$ passing
within the angle $\angle AOC$. Let's denote by $B$ the point at which
the ray $l$ crosses the segment $[AC]$. The point $B$ is an interior
point of the segment $[AC]$
\vadjust{\vskip 5pt\hbox to 0pt{\kern -5pt
\includegraphics{Oris18.eps}\hss}\vskip -5pt}since the ray 
$l$ does coincide neither with the ray $h$ nor with the ray $k$.\par
\parshape 6 5cm 5cm 5cm 5cm 5cm 5cm 5cm 5cm 5cm 5cm 
0cm 10cm
    The condition $\angle hl\cong\angle h'l'$ is common for the first
and the second items of the theorem. Let's use it in the following 
constructions. Applying the axiom~\mytheaxiom{A13}, on the ray $h'$ 
we choose a point $A'$ such that $[OA]\cong [O'A']$. Then in the
same way on the ray $l'$ we choose a point $B'$ such that $[OB]\cong
[O'B']$. Due to the condition $\angle hl\cong\angle h'l'$ we can
apply the theorem~\mythetheorem{5.1}. As a result we find that the
triangle $AOB$ is congruent to the triangle $A'O'B'$. Hence, we have
$$
\xalignat 2
&\hskip -2em
[AB]\cong [A'B'],&&\angle OAB\cong\angle O'A'B'.
\mytag{5.1}
\endxalignat
$$
The point $B'$ divides the line $A'B'$ into two rays, the ray 
$[B'A'\rangle$ being one of them. On the ray opposite to $[B'A'\rangle$
we choose a point $C'$ such that $[BC]\cong [B'C']$. Applying the
first item of the axiom~\mytheaxiom{A15} to the segments $[A'B']$
and $[B'C']$, we get $[AC]\cong [A'C']$. From this relationship and
from \mythetag{5.1} due to the theorem~\mythetheorem{5.1} we conclude
that the triangle $AOC$ is congruent to the triangle $A'O'C'$.
Hence, the following relationships are valid:
$$
\hskip -2em
\aligned
&[OC]\cong [O'C'],\\
&\angle ACO\cong\angle A'C'O',\\
&\angle AOC\cong\angle A'O'C'.
\endaligned
\mytag{5.2}
$$
Let's combine the first two relationships \mythetag{5.2} with
$[BC]\cong [B'C']$. Applying the theorem~\mythetheorem{5.1} to
this combination, we derive that the triangle $BOC$ is congruent
to the triangle $B'O'C'$. From this congruence, in particular,
we get
$$
\hskip -2em
\angle BOC\cong\angle B'O'C'.
\mytag{5.3}
$$\par
     On order to prove the first item of the theorem~\mythetheorem{5.3}
we consider the angles $\angle l'k'$ and $\angle B'O'C'$. They lie on
one half-plane bounded by the line $O'B'$ and have the common side
$l'=[O'B'\rangle$. For this angles from the statement of the theorem 
and from the formula \mythetag{5.3} we extract the relationships 
$$
\xalignat 2
&\angle lk\cong\angle l'k',
&&\angle lk\cong\angle B'O'C'.
\endxalignat
$$
Applying the axiom~\mytheaxiom{A16} now yields the coincidence of 
rays $[O'C'\rangle=k'$. Being combined with the last relationship
\mythetag{5.2}, this coincidence immediately yields the required 
result $\angle hk\cong h'k'$.\par
     In order to prove the second item of the theorem~\mythetheorem{5.3}
we derive the coincidence of rays $[O'C'\rangle=k'$ by considering the
angles $\angle h'k'$ and $\angle A'O'C'$, which lie on one half-plane
bounded by the line $O'A'$ and have the ray $h'$ as their common side.
From the statement of the theorem and from \mythetag{5.2} for these
angles we get
$$
\xalignat 2
&\angle hk\cong\angle h'k',
&&\angle hk\cong\angle A'O'C'.
\endxalignat
$$
Applying the axiom~\mytheaxiom{A16} to the above relationships, we derive
$[O'C'\rangle=k'$. Now from $[O'C'\rangle=k'$ and \mythetag{5.3} we get
the required result \pagebreak $\angle lk\cong l'k'$.\par
     The third item of the theorem~\mythetheorem{5.3} does not require 
a separate proof. It is reduced to the second item upon exchanging the 
notations of the rays: $h$ with $k$ and $h'$ with $k'$.
\qed\enddemo
\mydefinition{5.2} A triangle $ABC$ is called {\it isosceles\/} if some
two sides of it are congruent. For example, $[AB]\cong [AC]$. The side
$[BC]$ in this case is called the {\it base\/} of the isosceles triangle
$ABC$, while the congruent sides $[AB]$ and $[AC]$ are called the 
{\it lateral sides\/} of this isosceles triangle.
\enddefinition
\mytheorem{5.4} The angles at the base of an isosceles triangle are
congruent to each other.
\endproclaim
\demo{Proof} Let $ABC$ be an isosceles triangle with lateral sides 
$[AB]$ and $[AC]$. Let's introduce the duplicate notations for the
vertices of this triangle: $\tilde A=A$, $\tilde B=C$, and $\tilde C=B$.
Then from $[AB]\cong [AC]$ due to the axiom~\mytheaxiom{A16} we get
$$
\xalignat 3
&[AB]\cong [\tilde A\tilde B],
&&[AC]\cong [\tilde A\tilde C],
&&\angle BAC\cong\angle\tilde B\tilde A\tilde C.
\endxalignat
$$
In this situation the axiom~\mytheaxiom{A17} is applicable. Applying
this axiom we get $\angle ABC\cong\angle\tilde A\tilde B\tilde C$ and
$\angle ACB\cong\angle \tilde A\tilde C\tilde B$. When taking into 
account the above duplicate notations for the vertices it means that 
$\angle ABC\cong\angle ACB$ and $\angle ACB\cong\angle ABC$. The
required congruence of angles at the base of the isosceles triangle 
$ABC$ is proved.
\qed\enddemo
\mylemma{5.1} Let $ABC$ and $ADC$ be two triangles with the common side
$[AC]$ lying on two different half-planes of the same plane separated by 
the line $AC$. In this case if $[AB]\cong [AD]$ and if $[CB]\cong [CD]$,
then $\angle ABC\cong\angle ADC$ and $\angle ADC\cong\angle ABC$.
\endproclaim
\demo{Proof} The points $B$ and $D$ lie on different sides of the
line $AC$, hence, the segment $[BD]$ intersects this line at some 
interior point $S$. There are the following five cases for the 
disposition of the point $S$ relative to the points $A$ and $C$:
\roster
\item the point $A$ lies between $S$ and $C$;
\item the point $S$ coincides with the point $A$;
\item the point $S$ lies between $A$ and $C$;
\item the point $S$ coincides with the point $C$;
\item the point $C$ lies between $A$ and $S$.
\endroster
The first three cases of mutual disposition of points are shown on
Fig\.~5.3, Fig\.~5.4 and Fig\.~5.5. The fourth case is reduced to
the second one and the fifth case is reduced to the first one when
exchanging the notations of the points $A$ and $C$. These two cases
do not require a separate consideration.\par
\vskip 1cm plus 1pt minus 1pt
\hbox to 0pt{\kern -5pt\includegraphics{Oris19.eps}\hss}
\vskip 4.7cm plus 1pt minus 1pt
    Let's consider the first case. From $[AB]\cong [AD]$ and from
$[CB]\cong [CD]$ we conclude that the triangles $BCD$ and $BAD$ on
Fig\.~5.3 are isosceles. Therefore, we have $\angle SBC\cong\angle 
SDC$ and $\angle SBA\cong\angle SDA$. The ray $[BA\rangle$ lies 
inside the angle $\angle SBC$, while the ray $[DA\rangle$ is inside 
the angle $\angle SDC$. Hence, we can apply the second item of the
theorem~\mythetheorem{5.3} and derive the required relationship
$\angle ABC\cong\angle ADC$. The second relationship $\angle ADC\cong
\angle ABC$ is derived similarly.\par
    In the second case the required relationships are derived 
immediately since due to the congruence $[CB]\cong [CD]$ the triangle 
$BCD$ on Fig\.~5.4 is isosceles.\par
    And finally, let's consider the third case. From $[AB]\cong [AD]$ 
and $[CB]\cong [CD]$ we conclude that the triangles $BCD$ and $BAD$ 
on Fig\.~5.5 are isosceles, which means that $\angle SBC\cong\angle SDC$
and $\angle SBA\cong\angle SDA$. Due to these relationships we can apply
the first item of the theorem~\mythetheorem{5.3}. It yields $\angle ABC
\cong\angle ADC$. The second relationship $\angle ADC\cong\angle ABC$
is derived similarly.
\qed\enddemo
\mylemma{5.2} For any triangle $ABC$ in the plane of this triangle
there is exactly one point $D$ different from $B$ and not lying on
the line $AC$ such that $[AB]\cong [AD]$ and $[CB]\cong [CD]$.
\endproclaim
\demo{Proof} For the beginning we prove the existence of the point
$D$. Let's denote by $h$ the ray $[AC\rangle$ lying on the line
$AC$. This line divides the plane $ABC$ into two half-planes. The 
triangle $ABC$ lies on one of them. Applying the axiom~\mytheaxiom{A16},
in the other half-plane we draw a ray $k$ coming out from the point $A$
and such that $\angle CAB\cong hk$. Applying the axiom~\mytheaxiom{A13},
on the ray $k$ we choose a point $D$ satisfying the condition
$[AB]\cong [AD]$. Then $\angle CAB\cong\angle CAD$ and due to the 
theorem~\mythetheorem{5.1} the triangle $ABC$ appears to be congruent 
to the triangle $ADC$. Hence, $[CB]\cong [CD]$. Thus, a required 
point $D$ is constructed. The points $B$ and $D$ lie on different 
half-planes outside the boundary line $AC$. Therefore $D\neq B$.\par
    Let's prove the uniqueness of the point $D$. Assume that $\tilde D$
is another point satisfying the conditions of the lemma. According to
the statement of the lemma, the point $\tilde D$ does not lie on the
line $AC$, therefore it lies on one of the open half-planes determined
my this line.\par
    Let's consider that half-plane which contains the point $D$. Applying
the lemma~\mythelemma{5.1}, we get $\angle ABC\cong\angle A\tilde DC$.
Combining this congruence with $[AB]\cong [A\tilde D]$ and $[CB]\cong
[C\tilde D]$, we conclude that the triangle $ABC$ is congruent to the 
triangle $A\tilde DC$. Hence, we have the following congruence of angles:
$$
\xalignat 2
&\hskip -2em
\angle CAB\cong\angle CA\tilde D,
&&\angle ACB\cong\angle AC\tilde D.
\mytag{5.4}
\endxalignat
$$
The similar relationships are fulfilled for the point $D$ too:
$$
\xalignat 2
&\hskip -2em
\angle CAB\cong\angle CAD,
&&\angle ACB\cong\angle ACD.
\mytag{5.5}
\endxalignat
$$
They are derived from $[AB]\cong [AD]$ and $[CB]\cong [CD]$ by means of
the lemma~\mythelemma{5.1}. Comparing \mythetag{5.4} with \mythetag{5.5} 
and applying the axiom~\mytheaxiom{A13}, we prove the coincidence of the
rays $[AD\rangle=[A\tilde D\rangle$ and $[CD\rangle=[C\tilde D\rangle$. 
Two non-coinciding straight lines $AD$ and $CD$ can have at most one
common point. Therefore $\tilde D=D$.\par
    If we assume that the point $\tilde D$ lies on the same half-plane
as the point $B$, we get $\tilde D=B$. This fact is derived with the use
of the lemma~\mythelemma{5.1} on the base of considerations quite similar
to the above ones. But the coincidence $\tilde D=B$ is excluded by the
provisions of the lemma. Therefore, the point $D$ constructed above 
is unique. 
\qed\enddemo
\mytheorem{5.5} If for triangles $ABC$ and $\tilde A\tilde B\tilde C$
the conditions $[AB]\cong [\tilde A\tilde B]$, $[AC]\cong [\tilde A
\tilde C]$, and $[BC]\cong [\tilde B\tilde C]$ are fulfilled, then the
triangle $ABC$ is congruent to the triangle $\tilde A\tilde B\tilde C$.
\endproclaim
\demo{Proof}\parshape 15 0cm 10cm 0cm 10cm 3.7cm 6.4cm
3.7cm 6.4cm 3.7cm 6.4cm 3.7cm 6.4cm 3.7cm 6.4cm 3.7cm 6.4cm
3.7cm 6.4cm 3.7cm 6.4cm 3.7cm 6.4cm 3.7cm 6.4cm 3.7cm 6.4cm 
3.7cm 6.4cm 0cm 10cm
Let's denote by $\tilde h$ the ray $[\tilde A\tilde C\rangle$. 
The line $\tilde A\tilde C$ divides the plane of the second 
triangle $\tilde A\tilde B\tilde C$ into two half-planes. The 
triangle $\tilde A\tilde B\tilde C$ itself lies on one of them. 
Applying the axiom~\mytheaxiom{A16} to both half-planes, we draw
two rays $\tilde k$ and $\tilde k'$ coming out from the point
$\tilde A$ and such that $\angle CAB\cong\angle\tilde h\tilde k$ 
and $\angle CAB\cong\angle\tilde h\tilde k'$.
\vadjust{\vskip 5pt\hbox to 0pt{\kern 0pt
\includegraphics{Oris20.eps}\hss}\vskip -5pt}The ray 
$\tilde k'$ is chosen to be lying in the same half-plane
as the triangle $\tilde A\tilde B\tilde C$. The ray $\tilde k$ lies
in the other half-plane. Let's apply the axiom~\mytheaxiom{A13} to
the rays $\tilde k$ and $\tilde k'$ and choose two points $\tilde D$ 
and $\tilde B'$ on them such that $[AB]\cong [\tilde A\tilde D]$ and
$[AB]\cong [\tilde A\tilde B']$. By construction the triangle $ABC$
appears to be congruent to the triangles $\tilde A\tilde B'\tilde C$ 
and $\tilde A\tilde D\tilde C$. This fact follows from the 
theorem~\mythetheorem{5.1} if we take into account the congruence
$[AC]\cong [\tilde A\tilde C]$. Hence, for the points $\tilde D$ 
and $\tilde B'$ we get
$$
\xalignat 2
&[AB]\cong[\tilde A\tilde B'],&&[CB]\cong[\tilde C\tilde B'],\\
&[AB]\cong[\tilde A\tilde D],&&[CB]\cong[\tilde C\tilde D].
\endxalignat
$$
According to the statement of the theorem, exactly the same conditions 
are fulfilled for the point $\tilde B$ lying on the same half-plane
as the point $\tilde B'$:
$$
\xalignat 2
&[AB]\cong[\tilde A\tilde B],&&[CB]\cong[\tilde C\tilde B].
\endxalignat
$$
Hence, the lemma~\mythelemma{5.2} yields $\tilde B=\tilde B'$. This 
means that the  triangle $ABC$ is congruent to the triangle $\tilde A
\tilde B\tilde C$.
\qed\enddemo
    The theorem~\mythetheorem{5.5} proved just above is known as the
congruence criterion for triangles by three sides.
\mytheorem{5.6} The congruence of angles is a reflexive, symmetric,
and transitive binary relation. Therefore, it is the an equivalence
relation.
\endproclaim
\myexercise{5.1} The reflexivity of the congruence of angles is 
explicitly stated in the axiom~\mytheaxiom{A16}. Prove the symmetry
and transitivity of this binary relation with the use of the 
theorem~\mythetheorem{5.5}.
\endproclaim
\head
\SectionNum{6}{73} A right angle and orthogonality.
\endhead
\rightheadtext{\S\,6. A right angle and orthogonality.} 
\mytheorem{6.1} The congruence of two angles imply the congruence
of their adjacent angles.
\endproclaim
\demo{Proof}\parshape 20 0cm 10cm 0cm 10cm 0cm 10cm 0cm 10cm 
0cm 10cm 0cm 10cm 5cm 5cm 5cm 5cm 5cm 5cm 5cm 5cm 5cm 5cm 
5cm 5cm 5cm 5cm 5cm 5cm 5cm 5cm 5cm 5cm 5cm 5cm 5cm 5cm 
5cm 5cm 0cm 10cm
Assume that an angle $\angle hk$ with the vertex at a point $O$ is
congruent to an angle $\angle\tilde h\tilde k$ with the vertex at 
a point $\tilde O$. Let's complement the ray $h$ with the ray $l$ 
up to a whole straight line. Similarly, we complement the ray 
$\tilde h$ with the ray $\tilde l$ up to a whole straight line.
As a result we get two angles $\angle kl$ and $\angle\tilde k
\tilde l$. Their congruence should be proved. For this purpose 
we choose some points $A$, $B$, and $C$ on the rays $h$, $k$, and 
$l$ respectively.
\vadjust{\vskip 5pt\hbox to 0pt{\kern -7pt
\includegraphics{Oris21.eps}\hss}\vskip -5pt}Then, 
applying the axiom~\mytheaxiom{A13}, we mark three points
$\tilde A$, $\tilde B$, and $\tilde C$ on the rays $\tilde h$, 
$\tilde k$ and $\tilde l$ so that the conditions 
$$
\hskip -2em
\aligned
&[OA]\cong[\tilde O\tilde A],\\
&[OB]\cong[\tilde O\tilde B],\\
&[OC]\cong[\tilde O\tilde C]
\endaligned
\mytag{6.1}
$$
are fulfilled. From the first two conditions \mythetag{6.1} 
complemented\linebreak with $\angle hk\cong\angle\tilde h\tilde k$ 
we derive that the triangles $AOB$ and $\tilde A\tilde O\tilde B$ 
are congruent. This yields $[AB]\cong [\tilde A\tilde B]$ and 
$\angle BAO\cong\angle\tilde B\tilde A\tilde O$. Moreover, from
the first and the last conditions \mythetag{6.1}, applying the
axiom~\mytheaxiom{A15}, we derive $[AC]\cong [\tilde A\tilde C]$. 
Combining the three conditions obtained and applying the 
theorem~\mythetheorem{5.1}, we derive that the triangles 
$ABC$ and $\tilde A\tilde B\tilde C$ are congruent. This congruence
yields $[BC]\cong [\tilde B\tilde C]$. Now we can apply the
theorem~\mythetheorem{5.5} to the triangles $BOC$ and $\tilde B
\tilde O\tilde C$. This theorem yields the congruence of these
triangles, which, in turn, yields the required relationship
$\angle BOC\cong\angle\tilde B\tilde O\tilde C$. The theorem is
proved.
\qed\enddemo
\mytheorem{6.2} Vertical angles are congruent to each other.
\endproclaim
    The proof of this theorem is obvious. Two vertical angles always 
have a common adjacent angle (see Fig\.~5.1 in 
Chapter~\uppercase\expandafter{\romannumeral 2}). Therefore, due to 
the previous theorem~\mythetheorem{6.1} and the 
axiom~\mytheaxiom{A16} these angles are congruent.
\mydefinition{6.1} An angle is called a {\it right angle}, if it is
congruent to its adjacent angle.
\enddefinition
\mylemma{6.1} Right angles do exist.
\endproclaim
\demo{Proof} \parshape 13 0cm 10cm 0cm 10cm 3.7cm 6.3cm 
3.7cm 6.3cm 3.7cm 6.3cm 3.7cm 6.3cm 3.7cm 6.3cm 3.7cm 6.3cm 
3.7cm 6.3cm 3.7cm 6.3cm 3.7cm 6.3cm 3.7cm 6.3cm 
0cm 10cm
Let $\angle hk$ be an arbitrary angle formed by two rays coming out from
a point $A$. \vadjust{\vskip 5pt\hbox to 0pt{\kern 5pt
\includegraphics{Oris22.eps}\hss}\vskip -5pt}If this angle
appears to be a right angle, then the proof is over. Assume that this
angle is not a right angle. Let's complete the ray $h$ up to a whole 
straight line. This line divides the plane of the angle $\angle hk$
into two half-planes. The angle $\angle hk$ belongs to one of these
half-planes. Applying the axiom~\mytheaxiom{A16}, on the other 
half-plane we draw a ray $\tilde k$ coming out from the point $A$ 
and such that $\angle hk\cong\angle h\tilde k$. Then we choose a point
$B$ on the ray $k$ and, applying the axiom~\mytheaxiom{A13}, on the ray
$\tilde k$ we find a point $\tilde B$ satisfying the condition
$[AB]\cong [A\tilde B]$.\par
\parshape 7 0cm 10cm 0cm 10cm 0cm 10cm 0cm 10cm 0cm 10cm 0cm 10cm 
3.7cm 6.3cm 
    The points $B$ and $\tilde B$ lie on different sides of the line
containing the ray $h$. Therefore, the segment $[B\tilde B]$ intersects
this line at some its interior point $O\neq A$ (the equality $O=A$
leads to the case, where $\angle hk$ and $\angle h\tilde k$ both are
right angles). The angles $\angle AOB$ and $\angle AO\tilde B$ are 
adjacent angles. They are congruent since the triangles $AOB$ and 
$AO\tilde B$ are congruent, this fact follows from the
theorem~\mythetheorem{5.1} due to $[AB]\cong[A\tilde B]$ and
\vadjust{\vskip 5pt\hbox to 0pt{\kern -5pt
\includegraphics{Oris23.eps}\hss}\vskip -5pt}$\angle OAB
\cong\angle OA\tilde B$. Hence, the angles $\angle AOB$ and
$\angle AO\tilde B$ are right angles.
\qed\enddemo
\mylemma{6.2}\parshape 1 3.7cm 6.3cm 
All right angles are congruent to each other.
\endproclaim
\demo{Proof}\parshape 6 3.7cm 6.3cm 3.7cm 6.3cm 3.7cm 6.3cm 3.7cm 6.3cm 
3.7cm 6.3cm 0cm 10cm 
In order to prove the lemma it is sufficient to show that all right 
angles are congruent to one of them. As such a reference model we 
choose the right angle $\angle AOB$ constructed in the proof of the
lemma~\mythelemma{6.1}. Let's complement Fig\.~6.2 with one more point
$\tilde A$ lying on the line $OA$. Let's determine this point $\tilde A$ 
by means of the condition $[O\tilde A]\cong [OA]$. Then we draw the
segments $[\tilde AB]$ and $[\tilde A\tilde B]$. As a result we get four 
right angles with the common vertex at the point $O$.\par
     Let $\angle h'q'$ be some arbitrary right angle. The line 
$A\tilde A$ divides the plane of Fig\.~6.3 into two half-planes. 
In that half-plane which contains the point $B$ we draw a ray $q$
coming out from the point $O$ so that the angle $\angle hq$ is congruent
to the angle $\angle h'q'$. Then the angle $\angle hq$ is also a right
angle.\par
     Let's prove that the ray $q$ coincides with the ray $[OB\rangle$. 
If it is not so, the ray $q$ lies within one of the angles $\angle AOB$ 
or $\angle\tilde AOB$. For the sake of certainty assume that it lies
within the angle $\angle AOB$ (the second case is reduced to this one 
by exchanging the notations of the points $A$ and $\tilde A$). Applying
the lemma~\mythelemmachapter{6.2}{2} from 
Chapter~\uppercase\expandafter{\romannumeral 2} to the ray $q$, we find
that it intersects the segment $[AB]$ at some interior point $C$. Let's
complement $q$ up to the whole line with the ray $\tilde q$ lying inside
the angle $\tilde AO\tilde B$. Then we perform the congruent translation
$f$ of the line $AB$ to the line $A\tilde B$ such that $f(A)=A$ and 
$f(B)=\tilde B$. Let's denote $\tilde C=f(C)$. Then $\tilde C$ is an
interior point of the segment $[A\tilde B]$ and $[A\tilde C]\cong [AC]$.
This relationship complemented with $\angle OAB\cong\angle OA\tilde B$ 
yields the congruence of the triangles $AOC$ and $AO\tilde C$. Hence,
the angle $\angle AO\tilde C$ is congruent to the angle $\angle AOC
=\angle hk$. The angle $\angle h\tilde q$ is an adjacent angle for
$\angle hq$. Therefore, $\angle h\tilde q\cong\angle hq$ since $\angle hq$
is a right angle. From $\angle AO\tilde C\cong\angle hq$ and
$\angle h\tilde q\cong\angle hq$ we derive the coincidence of the rays 
$[O\tilde C\rangle$ and $\tilde q$. But such a coincidence is forbidden 
since $[O\tilde C\rangle$ lies within the angle $\angle AO\tilde B$, while
$\tilde q$ is within its adjacent angle $\angle\tilde AO\tilde B$.\par
    The contradiction obtained just above proves that $C=B$ and 
$q=[OB\rangle$. Hence, an arbitrary right angle $\angle h'q'$ is congruent
to the reference right angle $\angle AOB$. The lemma is proved.
\qed\enddemo
\mydefinition{6.2} Two intersecting straight lines are called 
{\it perpendicular\/}  to each other if all four angles formed by them
at the intersection point are right angles.
\enddefinition
     The perpendicularity of lines  $a$ and $b$ is denoted as $a\perp b$.
Actually, for two lines to be perpendicular it is sufficient that one of
the four angles formed by them at the intersection point is a right angle.
Then other three angles are also right angles due to the
theorems~\mythetheorem{6.1} and \mythetheorem{6.1}.
\mytheorem{6.3} Let $a$ be some straight line lying on a plane $\alpha$.
Then for any point $A\in a$ there is exactly one straight line $b$ passing
trough this point, lying on the plane $\alpha$, and being perpendicular to
the line $a$.
\endproclaim
\mytheorem{6.4} Through any point of a plane one can draw at most two
straight lines lying on this plane and being perpendicular to each other.
\endproclaim
    The theorems~\mythetheorem{6.3} and \mythetheorem{6.4} are easily
derived from the lemmas~\mythelemma{6.1} and \mythelemma{6.2}. The
theorem~\mythetheorem{6.4} means the two-dimensionality of a plane.
\par
    Let $a$ be some straight line and assume that $B$ is some point
outside this line. A segment $[BA]$ connecting the point $B$ with
a point $A\in a$ is called a {\it perpendicular\/} dropped from $B$ 
onto the line $a$ if the line $AB$ is perpendicular to the line$a$.
\mytheorem{6.5} For any line $a$ and for any point $B\notin a$ there 
is exactly one perpendicular dropped from $B$ onto the line $a$.
\endproclaim
\demo{Proof}\parshape 14 0cm 10cm 0cm 10cm 4cm 6cm 4cm 6cm
4cm 6cm 4cm 6cm 4cm 6cm 4cm 6cm 4cm 6cm 4cm 6cm
4cm 6cm 4cm 6cm 4cm 6cm 0cm 10cm
For the beginning let's prove the existence of a perpendicular dropped
from the point $B$ onto the line $a$. For this purpose we draw the plane
$\alpha$ passing through this line and through the point $B$.
\vadjust{\vskip 5pt\hbox to 0pt{\kern -5pt
\includegraphics{Oris24.eps}\hss}\vskip -5pt}Let's choose
some point $O$ on the line $a$. If $OB\perp a$, the required 
perpendicular is found. Otherwise, we consider two rays produced by
the point $O$ on the line $a$. Let's denote one of them through $h$ 
and let $k=[OB\rangle$. The angle $\angle hk$ lies on one of the 
half-planes produced by the line $a$ on the plane $\alpha$. On the
other half-plane we draw a ray coming out from the point $O$ and such
that $\angle h\tilde k\cong\angle hk$. On this ray we choose a point 
$\tilde B$ satisfying the condition $[O\tilde B]\cong [OB]$. Then the
segment $[B\tilde B]$ intersects the line $a$ at some interior point
$A$, while the triangles $OAB$ and $OA\tilde B$ are congruent. Hence,
the angle $\angle OAB$, which is congruent to its adjacent angle 
$\angle OA\tilde B$, is a right angle. Thus, $[BA]$ is a required
perpendicular.\par
    Now let's prove the uniqueness of the perpendicular dropped from 
the point $B$ onto the line $a$. Assume that it is not unique and
consider two perpendiculars $[BA]$ and $[B\tilde A]$. As above, we 
draw the plane $\alpha$ passing through the point $B$ and the line
$a$. Let's choose a point $O$ on the line $a$ outside the segment
$[A\tilde A]$ and denote by $h$ the ray $[OA\rangle=[O\tilde A\rangle$.
On the lines $AB$ and $\tilde AB$ we choose the points $\tilde B$ and
$\tilde B'$ distinct from $B$ and satisfying the conditions
$[AB]\cong [A\tilde B]$ and $[\tilde AB]\cong[\tilde A\tilde B']$.
Now note that the right angle $\angle OAB$ is congruent to its 
adjacent angle $\angle OA\tilde B$, while the right angle $\angle O
\tilde AB$ is congruent to the angle $\angle O\tilde A\tilde B'$. 
Hence, we derive that $\triangle OAB\cong\triangle OA\tilde B$ and 
that $\triangle O\tilde AB\cong\triangle O\tilde A\tilde B'$. 
As a result we have
$$
\xalignat 2
\angle AO\tilde B&\cong\angle AOB,&\angle AO\tilde B'&\cong\angle AOB,\\
[O\tilde B]&\cong [OB],&[O\tilde B']&\cong [OB].
\endxalignat
$$
The first pair of the above relationships yields the coincidence of
the rays $[O\tilde B\rangle=[O\tilde B'\rangle$. Then from the second
pair we derive $\tilde B=\tilde B'$. Hence, $A=\tilde A$ and 
$[BA]=[B\tilde A]$.
\qed\enddemo
\mytheorem{6.6} A triangle cannot have two right angles.
\endproclaim
\myexercise{6.1} Prove the theorems~\mythetheorem{6.3} and
\mythetheorem{6.4} on the base of the lemmas~\mythelemma{6.1}
and \mythelemma{6.2}.
\endproclaim
\myexercise{6.2} Derive the theorem~\mythetheorem{6.6} from the
theorem~\mythetheorem{6.5}.
\endproclaim
\head
\SectionNum{7}{78} Bisection of segments and angles.
\endhead
\rightheadtext{\S\,7. Bisection of segments and angles.}
\mydefinition{7.1} A point $O$ is called a {\it center\/} of
a segment $[AB]$ if it lies in the interior of this segment and
if $[AO]\cong [OB]$.
\enddefinition
\mydefinition{7.2} Let $\angle hk$ be an angle with the vertex at 
a poit $O$. A ray $l$ coming out from the point $O$ is called a 
{\it bisector\/} of the angle $\angle hk$ if it lies within this
angle and if $\angle hl\cong\angle lk$.
\enddefinition
\mydefinition{7.3} A segment $[AO]$ connecting a vertex $A$ of 
a triangle $ABC$ with a point $O$ on the line $BC$ is called 
a {\it height\/} of this triangle if the lines $AO$ and $BC$
are perpendicular.
\enddefinition
\mydefinition{7.4} A segment $[AO]$ connecting a vertex $A$ of
a triangle $ABC$ with a center of the side $[BC]$ is called a
{\it median\/} of this triangle.
\enddefinition
\mytheorem{7.1} In an isosceles triangle $ABC$ a median $AO$ drawn
from a vertex $A$ to a center of the base $[BC]$ is a height and 
a bisector of the angle $\angle BAC$ simultaneously.
\endproclaim
     The proof of the theorem~\mythetheorem{7.1} is sufficiently
simple provided a of the triangle $ABC$ is already drawn. The
problem of existing a median is not considered in this theorem
at all.
\mytheorem{7.2} For any segment $[AB]$ there is a point $O$ being
its center.
\endproclaim
\demo{Proof}\parshape 14 0cm 10cm 0cm 10cm 4cm 6cm 4cm 6cm 4cm 6cm 
4cm 6cm 4cm 6cm 4cm 6cm 4cm 6cm 4cm 6cm 4cm 6cm 4cm 6cm 4cm 6cm 
0cm 10cm
Assume that a straight line segment $[AB]$ is given. Let's consider
some plane $\alpha$ containing the line $AB$. 
\vadjust{\vskip 5pt\hbox to 0pt{\kern 5pt
\includegraphics{Oris25.eps}\hss}\vskip -5pt}Applying 
the theorem~\mythetheorem{6.3}, on the plane $\alpha$ we 
draw two lines being perpendicular to the line $AB$ and passing through 
the points $A$ and $B$ respectively. These two lines do not 
intersect each other since if they intersect at some point $M$, the
triangle $ABM$ would have two right angles. Let's choose some point
$C$ on the line passing through the point $A$. Then on the line
passing through the point $B$ we choose a point $D$ such that
$C$ and $D$ are on different half-planes separated by the line $AB$
and such that $[BD]\cong [AC]$. Then the segment $[CD]$ intersects
the line $AB$ at some interior point $O$.\par
     Let's consider the triangle $BOD$ and the line $AC$. The line
$AC$ cannot intersect the side $[BD]$ of this triangle since the
lines $AC$ and $BD$ have no common points at all. The line $AC$ 
intersects the line $OD$ at the point $C$ outside the segment $[OD]$.
Hence the line $AC$ cannot intersect the third side $[OB]$ of the
triangle $BOD$ at an interior point (see Pasch's 
axiom~\mytheaxiom{A12}). This fact excludes the disposition
$(O\blacktriangleright A\blacktriangleleft B)$ of the points
$O$, $A$, and $B$ on the line $AB$. The second disposition 
$(A\blacktriangleright B\blacktriangleleft O)$ is excluded by
considering the line $BD$ and the triangle $AOC$. The rest is the
third disposition $(A\blacktriangleright O\blacktriangleleft B)$.
In this disposition $O$ is an interior point of the segment $[AB]$.
\par
     Now let's consider the triangles $ABC$ and $BAD$. They are
congruent due to the theorem~\mythetheorem{5.1} since $\angle BAC
\cong\angle ABD$ and $[AC]\cong [BD]$. From this congruence of 
the triangles $ABC$ and $BAD$ we derive $[AD]\cong [BC]$ and 
$\angle ABC\cong\angle BAD$. Applying the item \therosteritem{1}
of the theorem~\mythetheorem{5.3}, we get $\angle CBD\cong\angle DAC$, 
which means that $\triangle CBD\cong\triangle DAC$. Hence, we have
$$
\xalignat 2
&\angle DCB\cong\angle CDA,&&\angle CDB\cong\angle DCA.
\endxalignat
$$
From these two relationships, applying the theorem~\mythetheorem{5.2},
we derive that the triangle $AOC$ is congruent to the triangle $BOD$, 
while the triangle $COB$ is congruent to the triangle $DOA$. These
congruences of triangles yield the relationships
$$
\xalignat 2
&[AO]\cong [OB],&&[CO]\cong [OD].
\endxalignat
$$
They mean that the point $O$ is a required center of the segment
$[AB]$, and simultaneously, it is a center of the segment $[CD]$.
\qed\enddemo
    The theorem~\mythetheorem{7.2} complemented with the
theorem~\mythetheorem{7.1} yields an algorithm for bisecting
angles. Indeed, assume that an angle $\angle hk$ with the vertex
at a point $O$ is given. On the sides of this angle we choose 
two points $A$ and $B$ such that $[OA]\cong [OB]$. Then the
triangle $AOB$ is isosceles. Bisecting its base $[AB]$ by a
point $C$, we construct its median being a bisector of the
angle $\angle AOB$ at the same time. In other words we have
the following theorem.
\mytheorem{7.3} For any angle $\angle hk$ there is a ray $l$ 
being a bisector of this angle.
\endproclaim
\mytheorem{7.4} The center of any segment is unique.
\endproclaim
\mytheorem{7.5} Each angle has exactly one bisector.
\endproclaim
\myexercise{7.1} Using the theorem~\mythetheorem{2.1}, prove
the theorem~\mythetheorem{7.4} for a segment $[AB]$. In other words, 
show that on the line $AB$ there is exactly one point $O$ satisfying
the condition $[AO]\cong [OB]$. Then from the 
theorem~\mythetheorem{7.4} derive the theorem~\mythetheorem{7.5}.
\endproclaim
\head
\SectionNum{8}{81} Intersection of two straight lines 
by a third one.
\endhead
\rightheadtext{\S\,8. Intersection of two lines by a third one.}
\parshape 14 0cm 10cm 0cm 10cm 4cm 6cm 4cm 6cm 4cm 6cm 4cm 6cm 4cm 
6cm 4cm 6cm 4cm 6cm 4cm 6cm 4cm 6cm 4cm 6cm 4cm 6cm 0cm 10cm
     Let's consider three straight lines $a$, $b$, and $c$ lying on 
one plane. \vadjust{\vskip 5pt\hbox to 0pt{\kern -5pt
\includegraphics{Oris26.eps}\hss}\vskip -5pt}Assume that 
the line $c$ intersects the lines $a$ and $b$ at the points $A$ and 
$B$ respectively. At the intersection point $A$ we have four angles. 
Other four angles arise at the point $B$. The angles $\angle CAB$ and
$\angle ABD$ have the sides $[AB\rangle$ and $[BA\rangle$ being two
oppositely-directed rays lying on one line $c$ and intersecting along 
the segment $[AB]$. The angles $\angle CAB$ and $\angle ABD$ lie on
different half-planes separated by the line $c$. Such angles are called
{\it inner crosswise lying\/} angles. Apart from $\angle CAB$ and 
$\angle ABD$, by the intersection of the lines $a$ and $b$ with the third
line $c$ another pair of inner crosswise lying angles arise. These are the
angles $\angle EAB$ and $\angle ABF$. The angle $\angle EAB$ is an adjacent
angle for the angle $\angle CAB$, while $\angle ABF$ is an adjacent angle
for $\angle ABD$.
\mydefinition{8.1} Two straight lines $a$ and $b$ are called 
{\it parallel\/} if they coincide $a=b$ or if they lie on one 
plane and do not intersect each other.
\enddefinition
     The relation of {\it parallelism\/} of two straight lines
is reflexive and symmetric by definition. It is denoted as
$a\parallel b$. In order to prove the transitivity of this 
relation one should use the axiom~\mytheaxiom{A20}, which is 
not yet considered.
\mytheorem{8.1} Assume that $a$ and $b$ are two straight lines
lying on one plane and intersecting with a third straight line
$c$ at the points $A$ and $B$. If the inner crosswise lying 
angles at the points $A$ and $B$ are congruent, then the lines
$a$ and $b$ are parallel.
\endproclaim
\demo{Proof} Note that it is sufficient to require the congruence
of one pair of inner crosswise lying angles, e\.\,g\. $\angle EAB
\cong\angle ABF$. Then the relationship $\angle CAB\cong\angle ABD$ 
follows from the above relationship due to the
theorem~\mythetheorem{6.1}.\par
\parshape 13 0cm 10cm 0cm 10cm 4.2cm 5.8cm 4.2cm 5.8cm 4.2cm 5.8cm
4.2cm 5.8cm 4.2cm 5.8cm 4.2cm 5.8cm 4.2cm 5.8cm 4.2cm 5.8cm 
4.2cm 5.8cm 4.2cm 5.8cm 0cm 10cm
Let's prove the theorem by contradiction. Assume that the angles
$\angle EAB$ and $\angle ABF$ are congruent, but the lines $a$ and
$b$ are not parallel.
\vadjust{\vskip 5pt\hbox to 0pt{\kern -5pt
\includegraphics{Oris27.eps}\hss}\vskip -5pt}Then they
intersect. For the sake of certainty we can take the points
$C$ and $F$ to be coinciding with the intersection point of the
lines $a$ and $b$. We choose the point $E$ to be satisfying the 
condition $[AE]\cong [BF]$. Such a choice is possible due to the
axiom~\mytheaxiom{A13}. Let's connect the point $E$ with the point
$B$ by means of the segment $[EB]$ and consider the triangles $EAB$ 
and $ABF$. From the statement of the theorem and from the above
additional drawings we derive
$$
\xalignat 3
&[AE]\cong [BF],&&[AB]\cong [BA],&&\angle EAB\cong\angle ABF.
\endxalignat
$$
Applying the theorem~\mythetheorem{5.1} to these relationships, we 
get the congruence $\triangle EAB\cong\triangle ABF$. This congruence 
immediately yields $\angle EBA\cong\angle BAC$. But $\angle BAC\cong
\angle ABD$, which, as we already mentioned above, follows from
$\angle EAB\cong\angle ABF$. Therefore, $\angle EBA\cong\angle DBA$. 
Hence, due to the axiom~\mytheaxiom{A16} we get the coincidence of
the rays $[BE\rangle$ and $[BD\rangle$. This means that the point
$E$ should lie on the line $b$, which is impossible since the lines
$a$ and $b$ do not coincide. The contradiction obtained shows that 
the lines $a$ and $b$ cannot intersect, i\.\,e\. they are parallel.
\qed\enddemo
\mytheorem{8.2} Any two perpendiculars to one straight line lying
on one plane are parallel.
\endproclaim
\myexercise{8.1} Derive the theorem~\mythetheorem{8.2} from the
theorem~\mythetheorem{8.1}.
\endproclaim
\newpage
\setfirstpage
\topmatter
\title\chapter{4}
Congruent translations and motions.
\endtitle
\leftheadtext{CHAPTER \uppercase\expandafter{\romannumeral 4}.
TRANSLATIONS AND MOTIONS.}
\endtopmatter
\document
\chapternum=4
\head
\SectionNum{1}{84} Orthogonality of a straight line and a plane.
\endhead
\rightheadtext{\S\,1. Orthogonality of a line and a plane.}
\mydefinition{1.1} Assume that a straight line $a$ intersects 
a plane $\alpha$ at a point $O$. The line $a$ is said to be 
{\it perpendicular\/} to the plane $\alpha$ if it is perpendicular
to all straight lines lying on the plane $\alpha$ and passing through
the point $O$.
\enddefinition
\mytheorem{1.1} A line $a$ intersecting a plane $\alpha$ at a point
$O$ is perpendicular to this plane if and only if it is perpendicular 
to some two distinct straight lines lying on the plane $\alpha$ and
passing through the point $O$.
\endproclaim
\demo{Proof}\parshape 3 0cm 10cm 0cm 10cm 5.3cm 4.7cm 
The necessity of the condition formulated in the theorem is obvious:
\vadjust{\vskip 5pt\hbox to 0pt{\kern -15pt
\includegraphics{Oris28.eps}\hss}\vskip -5pt}if the line 
$a$ is 
perpendicular to the plane $\alpha$, then it is perpendicular to all
lines on this plane passing through the point $O$, including those
two of them mentioned in the theorem.\par
\parshape 8 5.3cm 4.7cm 5.3cm 4.7cm 5.3cm 4.7cm 5.3cm 4.7cm 5.3cm 4.7cm 
5.3cm 4.7cm 5.3cm 4.7cm 0cm 10.1cm
    Let's prove the sufficiency. Assume that the line $a$ is 
perpendicular to the lines $b$ and $c$ lying on the plane $\alpha$ 
and intersecting at the point $O$. Let's consider some arbitrary 
straight line $x$ lying on the plane $\alpha$ and passing through 
the point $O$. The point $O$ divides the line $x$ into two rays. 
Let's consider one of these rays $[O,+\infty)$. It lies inside 
one the four angles formed by the lines $a$ and $b$ at the 
intersection point $O$. On the sides of this angle we mark two 
points $B$ and $C$, then we draw the segment $[BC]$. According to
the lemma~\mythelemmachapter{6.2}{2} from 
Chapter~\uppercase\expandafter{\romannumeral 2}, the ray $[O,+\infty)$
intersects the segment $[BC]$ at some interior point $X$.\par
     Let's choose some arbitrary point $A$ different from $O$ on the
line $a$ and then, applying the axiom~\mytheaxiom{A13}, on the ray
opposite to $[OA\rangle$ we mark a point $\tilde A$ such that
$[O\tilde A]\cong [OA]$. Since $a\perp b$, we conclude that the angles
$\angle AOB$ and $\angle\tilde AOB$ are right angles. They are congruent
to each other $\angle AOB\cong\angle\tilde AOB$. Moreover, $[OB]\cong [OB]$. 
From these three relationships 
$$
\xalignat 3
&[OA]\cong [O\tilde A],&&[OB]\cong [OB],&&\angle AOB\cong
\angle\tilde AOB
\endxalignat
$$
we derive that the triangles $AOB$ and $\tilde AOB$ are congruent. 
Hence, $[AB]\cong [\tilde AB]$. In a similar way from $a\perp c$ 
we derive the congruence $[AC]\cong [\tilde AC]$. Let's complement
these two relationships with one more:
$$
\xalignat 3
&[AB]\cong [\tilde AB],&&[AC]\cong [\tilde AC],&&[BC]\cong [BC].
\endxalignat
$$
Now from these three relationships we derive that the triangles
$ABC$ and $\tilde ABC$ are congruent, which yields 
$\angle ABC\cong\angle\tilde ABC$.\par
     Let's consider the triangles $ABX$ and $\tilde ABX$. For their
sides and angles the following relationships are fulfilled:
$$
\xalignat 3
&[AB]\cong [\tilde AB],&&[BX]\cong [BX],&&\angle ABX\cong\angle
\tilde ABX.
\endxalignat
$$
These relationships yield the congruence of the triangles $ABX$ and
$\tilde ABX$. Hence, for the segments $[AX]$ and $[\tilde AX]$ we have
$[AX]\cong [\tilde AX]$. In other words, the triangle $AX\tilde A$
is isosceles, while the segment $[OX]$ is a median of it since $O$ 
is the center of the segment $[A\tilde A]$ by our choice of the
points $A$ and $\tilde A$. The rest is to apply the
theorem~\mythetheoremchapter{7.1}{3} from 
Chapter~\uppercase\expandafter{\romannumeral 3}. According to this 
theorem the median $[OX]$ in the isosceles triangle $\tilde AXA$ 
is its height at the same time. Therefore, we get the required
relationship $a\perp x$.
\qed\enddemo
\mytheorem{1.2} For any line $a$ and a point $O$ on this line
there is exactly one plane $\alpha$ passing through the point
$O$ and being perpendicular to the line $a$.
\endproclaim
\demo{Proof} For the beginning we prove the existence of the
required plane $\alpha$, then we prove its uniqueness. Let's
choose a point $B$ outside the line $a$ and draw the plane
$\beta$ passing through the line $a$ and the point $B$. On the
plane $\beta$ we apply the theorem~\mythetheoremchapter{6.3}{3}
from Chapter~\uppercase\expandafter{\romannumeral 3} to the
line $a$ and the point $O\in a$. This yields a line $b\in\beta$
passing through the point $O$ and being perpendicular to the 
line $a$.\par
      In the next step we choose a point $C$ not lying on the
plane $\beta$. Through the line $a$ and the point $C$ we draw
a plane $\gamma$. It is clear that the planes $\beta$ and $\gamma$
are distinct,the line $a$ being the intersection of these planes.
Applying the theorem~\mythetheoremchapter{6.3}{3} from
Chapter~\uppercase\expandafter{\romannumeral 3} to the line 
$a\in\gamma$, and the point $O$, we get a line $c$ on the plane
$\gamma$ passing through the point $O$ and being perpendicular 
to the line $a$.\par
     The lines $b$ and $c$ belong to the different planes $\beta$ 
and $\gamma$ and pass through the point $O$. They are two distinct
straight lines intersecting at the point $O$. There is a plane
$\alpha$ passing through such two lines (see 
theorem~\mythetheoremchapter{1.5}{2} in 
Chapter~\uppercase\expandafter{\romannumeral 2}). But $a\perp b$
and $a\perp c$ by construction. Therefore, according to the above
theorem~\mythetheorem{1.1}, we have $a\perp\alpha$. A required
plane $\alpha$ is constructed.\par
     Now let's prove the uniqueness of the constructed plane $\alpha$.
Assume that there is another plane $\tilde\alpha\perp a$ passing through
the point $O$. The planes $\alpha$ and $\tilde\alpha$ are distinct but 
have the common point $O$. Therefore, the intersect along some line
$b$ which passes trough the point $O$. Since $\alpha\perp a$, we have
$b\perp a$.\par
     The lines $a$ and $b$ are perpendicular to each other. They 
intersect at the point $O$. There is a plane $\beta$ passing through 
such two lines $a$ and $b$. Let's choose a point $C$ not lying on the 
plane $\beta$ and let's draw a plane $\gamma$ through the point $C$ 
and the line $a$. Note that such a plane $\gamma$ is different from
$\alpha$ and from $\tilde\alpha$. It follows from the fact that the 
plane $\gamma$ intersects with the line $b$ at the single point $O$, 
while the planes $\alpha$ and $\tilde\alpha$ contain the line $b$ 
in whole. The plane $\gamma$ intersecting with the planes $\alpha$ 
and $\tilde\alpha$ yields two lines $c$ and $\tilde c$. The lines
$c$ and $\tilde c$ lie on the plane $\gamma$ and pass through the
point $O$. From $\alpha\perp a$ and $\tilde\alpha\perp a$ we derive
$c\perp a$ and $\tilde c\perp a$ for them. If $c\neq\tilde c$, then
on the plane $\gamma$ we would have two perpendiculars to the
line $a$ passing through the point $O\in a$. This would contradict
the theorem~\mythetheoremchapter{6.3}{3} from
Chapter~\uppercase\expandafter{\romannumeral 3}. Hence, $c=\tilde c$.
From $c=\tilde c$ we easily derive $\alpha=\tilde\alpha$. The uniqueness
of the plane $\alpha$ is proved.
\qed\enddemo
\mytheorem{1.3} For any plane $\alpha$ and a point $O\in\alpha$ there 
is exactly one straight line $a$ passing through the point $O$ and 
being perpendicular to the plane $\alpha$.
\endproclaim
\demo{Proof} Let's begin with proving the existence of a required line.
On the plane $\alpha$ we choose some arbitrary point $B$ different from
$O$ and draw the line $OB$. Let's denote this line through $b$. Then
we apply the above theorem~\mythetheorem{1.2} to the line $b$ and the 
point $O$. As a result we get a plane $\gamma$ passing through the point
$O$ and being perpendicular to the line $b$. Intersecting with $\alpha$, 
the plane $\gamma$ produces a line $c$ passing through the point $O$. To
$\gamma$, $c$, and $O$ we apply the theorem~\mythetheoremchapter{6.3}{3} 
from Chapter~\uppercase\expandafter{\romannumeral 3}. As a result we get
a line $a\in\gamma$ passing through the point $O$ and being perpendicular
to the line $c$. From $a\subset\gamma$ and $\gamma\perp b$ we derive
$a\perp b$. Thus, for the constructed line $a$ we have
$$
\xalignat 2
&a\perp b,&&a\perp c.
\endxalignat
$$
This means that the line $a$ passes through the point $O$ and 
perpendicular to the lines $b$ and $c$ lying on the plane $\alpha$ 
and intersecting at the point $O$. According to the
theorem~\mythetheorem{1.1}, the line $a$ is perpendicular to the 
plane $\alpha$. The existence of a required line $a$ is proved.
\par
     The rest is to prove the uniqueness of the constructed line $a$.
Assume that there is another line $\tilde a$ passing through the point
$O$ and being perpendicular to the plane $\alpha$. If $a\neq\tilde a$, 
then the pair of lines $a$ and $\tilde a$ intersecting at the point $O$ 
defines a plane $\beta$. Intersecting with $\alpha$, the plane $\beta$ 
produces a line $b$ passing through the point $O$. From $b\subset\alpha$,
from $a\perp\alpha$, and from $\tilde a\perp\alpha$ we conclude that on
the plane $\beta$ there are two perpendiculars $a$ and $\tilde a$ to
the line $b$ passing through the point $O\in b$. This fact contradicts
the theorem~\mythetheoremchapter{6.3}{3} from
Chapter~\uppercase\expandafter{\romannumeral 3}. The contradiction
obtained proves the coincidence $a=\tilde a$ and, thus, it proves 
the uniqueness of the required line $a$.
\qed\enddemo
\mytheorem{1.4} Assume that a line $b$ intersects a line $a$ at a point
$O$. The line $b$ is perpendicular to the line $a$ if and only if it lies
in a plane $\alpha$ passing through the point $O$ and being perpendicular 
to the line $a$.
\endproclaim
\mytheorem{1.5} A plane $\alpha$ passing through a point $O$ on a line
$a$ and being perpendicular to this line is the union of all straight 
lines passing through the point $O$ and being perpendicular to $a$.
\endproclaim
\myexercise{1.1} Draw figures illustrating the proofs of the 
theorems~\mythetheorem{1.2} and \mythetheorem{1.3}.
\endproclaim
\myexercise{1.2} Prove the theorems~\mythetheorem{1.4} and
\mythetheorem{1.5}.
\endproclaim
\head
\SectionNum{2}{88} A perpendicular bisector of a segment and
the plane of perpendicular bisectors.
\endhead
\rightheadtext{\S\,2. A perpendicular bisector of a segment \dots}
\mydefinition{2.1} A straight line $a$ passing through the center
of a segment $[AB]$ and being perpendicular to it is called a 
{\it perpendicular bisector\/} of this segment.
\enddefinition
     Let $O$ be the center of a segment $[AB]$ and let $M$ be some 
point on its perpendicular bisector $a$ distinct from the point $O$. 
Then from $[AO]\cong [BO]$ and $[OM]\cong [OM]$ it follows that
the triangles $AOM$ and $BOM$ are congruent to each other. Hence, 
we get $[AM]\cong [BM]$.\par
      Conversely, assume that $[AM]\cong [BM]$. Then the triangle
$AMB$ is isosceles. Its median $[MO]$ is its height at the same time
(see theorem~\mythetheoremchapter{7.1}{3} in
Chapter~\uppercase\expandafter{\romannumeral 3}). Hence, the line
$AM$ is a perpendicular bisector for the segment $[AB]$. The conclusion
is that a point $M$ satisfies the relationship $[AM]\cong [BM]$ if and 
only if it lies on some perpendicular bisector of the segment $[AB]$.
\par
     Let's consider the set of all perpendicular bisectors of the given
segment $[AB]$. This is the set of all straight lines being perpendicular
to the line $AB$ and passing through the point $O$. According to the
theorem~\mythetheorem{1.5}, such a set is a plane passing through the 
point $O$ and being perpendicular to the segment $[AB]$. This plane is
called the {\it plane of perpendicular bisectors\/} of the segment $[AB]$. 
Now we can formulate the following theorem.
\mytheorem{2.1} For any two points $A$ and $B$ a point $M$ satisfies the
condition $[AM]\cong [BM]$ if and only if it lies on the plane of 
perpendicular bisectors of the segment $[AB]$.
\endproclaim
\head
\SectionNum{3}{89} Orthogonality of two planes.
\endhead
\rightheadtext{\S\,3. Orthogonality of two planes.}
\mylemma{3.1} Assume that two planes $\alpha$ and $\beta$ have a 
common point $O$. Under this assumption if the plane $\beta$ 
contains the perpendicular to the plane $\alpha$ passing through 
the point $O$, then the plane $\alpha$ contains the perpendicular 
to the plane $\beta$ passing through the same point $O$.
\endproclaim
\demo{Proof} Let $a$ be the perpendicular to the plane $\alpha$
passing through the point $O$ and assume that $b$ is the perpendicular 
to the plane $\beta$ also passing through the point $O$. Let's denote
by $c$ the intersection of the planes $\alpha$ and $\beta$. According 
to the statement of the lemma, $a\subset\beta$. Hence, $b\perp a$ and
$b\perp c$. Moreover, $a\perp c$, since the line $c$ lies on the plane
$\alpha$, while $a$ is perpendicular to the plane $\alpha$.\par
     Two lines $b$  and $c$ intersecting at the point $O$ define some
plane $\gamma$. From $a\perp b$ and $a\perp c$ due to the
theorem~\mythetheorem{1.1} we get $a\perp\gamma$. But the plane 
passing through the point $O$ and being perpendicular to the line $a$ 
is unique. Therefore, the plane $\gamma$ should coincide with the plane 
$\alpha$. Hence, we get the required inclusion $b\subset\alpha$.
\qed\enddemo
\mydefinition{3.1} Assume that two planes $\alpha$ and $\beta$ have
a common point $O$. The plane $\alpha$ is said to be {\it perpendicular
to the plane $\beta$ at the point $O$\/} if it contains the perpendicular
to the plane $\beta$ passing through the point $O$.
\enddefinition
The above lemma~\mythelemma{3.1} shows that the relation of orthogonality
of planes at a point is symmetric, i\.\,e\. $\alpha\perp\beta$ implies 
$\beta\perp\alpha$.
\mytheorem{3.1} Assume that two planes $\alpha$ and $\beta$ intersect
along a line $c$. In this case if the plane $\beta$ is perpendicular 
to the plane $\alpha$ at some point $A\in c$, then the plane $\beta$ 
is perpendicular to the plane $\alpha$ at any other point $B$ of the 
line $c$.
\endproclaim
\demo{Proof}\parshape 15 0cm 10cm 0cm 10cm 4.5cm 5.5cm
4.5cm 5.5cm 4.5cm 5.5cm 4.5cm 5.5cm 4.5cm 5.5cm 4.5cm 5.5cm
4.5cm 5.5cm 4.5cm 5.5cm 4.5cm 5.5cm 4.5cm 5.5cm 4.5cm 5.5cm 
4.5cm 5.5cm 
0cm 10cm
Let's denote by $a$ and $b$ the perpendiculars to the planes 
$\alpha$ and $\beta$ respectively passing through the point $A$. 
\vadjust{\vskip 5pt\hbox to 0pt{\kern -5pt
\includegraphics{Oris29.eps}\hss}\vskip -5pt}Since $\alpha$ 
is perpendicular to $\beta$ at the point $A$, we have $a\subset\beta$ 
and $b\subset\alpha$. Let's apply the
theorem~\mythetheoremchapter{6.3}{3} from
Chapter~\uppercase\expandafter{\romannumeral 3} to the line $AB$ 
and the point $B$ on the plane $\beta$. It yields a straight line 
lying on the plane $\beta$, passing through the point $B$, and
perpendicular to the line $AB$. Let's show that this line is a
perpendicular to the plane $\alpha$. For this purpose we mark some
point $C$ on it, then we choose a segment $[AD]$ on the line $b$ 
such that $[AD]\cong [BC]$. For the rectangular triangles $ABC$ 
and $BAD$ we have the relationships
$$
\xalignat 3
&[AB]\cong [BA],&&[BC]\cong [AD],&&\angle ABC\cong\angle BAD.
\endxalignat
$$
From these relationships we derive the congruence of these 
triangles $ABC$ and $BAD$. Hence, $[AC]\cong [BD]$.\par
     Now let's consider the triangles $CAD$ and $DBC$. These
triangles appear to be congruent due to the following congruence
relationships for their sides:
$$
\xalignat 3
&[AD]\cong [BC],&&[AC]\cong [BD],&&[CD]\cong [DC].
\endxalignat
$$
Since the triangles $CAD$ and $DBC$ are congruent, we have the 
congruence of the angles $\angle CAD\cong\angle DBC$. But $b$ 
is a perpendicular to the plane $\beta$, whic means that the 
angle $\angle CAD$ is a right angle. Hence, the angle $\angle DBC$ 
is also aright angle. This yields $BC\perp BD$ and $BC\perp AB$. 
In other words, the line $BC$ is perpendicular to the pair of lines 
$AB$ and $BD$ lying on the plane $\alpha$. Therefore, it is a
perpendicular to the plane $\alpha$ passing through the point $B$. 
The line $BC$ belongs to the plane $\beta$, hence, we obtain the
required orthogonality of the plane $\alpha$ and the plane $\beta$ 
at the point $B$.
\qed\enddemo
     The theorem~\mythetheorem{3.1} shows that the orthogonality 
of planes is their global property: if it takes place at some point,
it is present at all other points of the intersection of two planes.
\mytheorem{3.2} For any two perpendiculars to a given plane $\alpha$ 
there is a plane $\beta$ containing both of them. This plane $\beta$
is perpendicular to the plane $\alpha$.
\endproclaim
     The theorem~\mythetheorem{3.2} does not require a separate 
proof. The required plane $\beta$ was constructed in proving the
theorem~\mythetheorem{3.1}.
\mytheorem{3.3} Any two perpendiculars to a given plane $\alpha$
are parallel to each other.
\endproclaim
     The theorem~\mythetheorem{3.3} is easily derived from the
theorem~\mythetheorem{3.2} and from the 
theorem~\mythetheoremchapter{8.2}{3} from 
Chapter~\uppercase\expandafter{\romannumeral 3}.
\mydefinition{3.2} Assume that $B$ is some point not lying on a
plane $\alpha$. A segment $[BA]$ connecting the point $B$ with
some point $A\in\alpha$ is called a {\it perpendicular\/} dropped
from the point $B$ onto the plane $\alpha$ if the line $AB$ 
is perpendicular to the plane $\alpha$. The point $A\in\alpha$ 
is called the {\it foot\/} of the perpendicular or the
{\it orthogonal projection\/} of the point $B$ onto the plane
$\alpha$.
\enddefinition
\mytheorem{3.4} From any point $B\notin\alpha$ one can drop exactly
one perpendicular onto the plane $\alpha$.
\endproclaim
\demo{Proof} For the beginning let's prove the existence of a
perpendicular dropped from the point $B$ onto the plane $\alpha$. 
Let's choose some point $O\in\alpha$. If the line $OB$ is
perpendicular to the plane $\alpha$, then the segment $[BO]$ 
is a required perpendicular.\par
     Let's conside the case where $[BO]$ is not a perpendicular
to the plane $\alpha$. Using the theorem~\mythetheorem{1.3}, we 
draw the line $a$ being perpendicular to the plane $\alpha$ and 
passing through the point $O$. The segment $[BO]$ does not lie on the line
$a$, hence, $B\notin a$. Mow we draw a plane $\beta$ passing through 
the line $a$ and the point $B$. Let's denote by $c$ the line being
the intersection of the planes $\alpha$ and $\beta$ and then apply
the theorem~\mythetheoremchapter{6.5}{3} from
Chapter~\uppercase\expandafter{\romannumeral 3} to the point $B$ 
and the line $c$. As a result we get a point $A$ on the line $c$
such that the line $AB$ is perpendicular to the line $c$.\par
     Let's show that the segment $[BA]$ is a required perpendicular
dropped from the point $B$ onto the plane $\alpha$. For this purpose
note that the plane $\beta$ contains the perpendicular to the plane 
$\alpha$ drawn at the point $O$. According to the 
theorem~\mythetheorem{3.1} it contains a perpendicular to $\alpha$
drawn  at the point $A$ either. Let's denote this perpendicular by
$\tilde a$. From $\tilde a\subset\beta$ and $\tilde a\perp c$, applying
the theorem~\mythetheoremchapter{6.5}{3} from 
Chapter~\uppercase\expandafter{\romannumeral 3}, we derive the 
coincidence $\tilde a=AB$.\par
     Let's show that the perpendicular $[BA]$ constructed above
is unique. If we assume that another perpendicular $[BA']$ does
exist, then in the triangle $ABA'$ we would have two right angles,
which contradicts the theorem~\mythetheoremchapter{6.6}{3} from
Chapter~\uppercase\expandafter{\romannumeral 3}. The uniqueness
of the perpendicular $[BA]$ can be derived from the
theorem~\mythetheorem{3.3} too.
\qed\enddemo
\myexercise{3.1} Prove the theorems~\mythetheorem{3.2} and
\mythetheorem{3.3}.
\endproclaim
\head
\SectionNum{4}{93} A dihedral angle.
\endhead
\rightheadtext{\S\,4. A dihedral angle.}
    Let $\alpha$ and $\beta$ be two planes intersecting along some
straight line $a$. Each of these two planes $\alpha$ and $\beta$ 
divides the space $\Bbb E$ into two half-spaces. We describe this 
fact as follows:
$$
\xalignat 2
&\Bbb E=\alpha_{-}\cup\alpha\cup\alpha_{+},
&&\Bbb E=\beta_{-}\cup\beta\cup\beta_{+}.
\endxalignat
$$
\mydefinition{4.1} A {\it dihedral angle\/} is the intersection
of two closed half-spaces determined by two intersecting planes. 
The intersection of the corresponding open half-spaces is called
the {\it interior\/} of a dihedral angle. The straight line produced 
as the intersection of two planes confining a dihedral angle is 
called an {\it edge\/} of this dihedral angle.
\enddefinition
\parshape 1 4.5cm 5.5cm
By the intersection of two planes $\alpha$ and $\beta$ four dihedral
angles arise at a time. \vadjust{\vskip 5pt\hbox to 0pt{\kern -10pt
\includegraphics{Oris30.eps}\hss}\vskip -5pt}These are 
$$
\xalignat 2
&\overline{\alpha_{+}}\cap\overline{\beta_{+}},
&&\overline{\alpha_{+}}\cap\overline{\beta_{-}},\\
&\overline{\alpha_{-}}\cap\overline{\beta_{-}},
&&\overline{\alpha_{-}}\cap\overline{\beta_{-}}.
\endxalignat
$$
Note that each dihedral angle is the union of its interior and two
closed half-planes cut by the edge $a$ on the planes $\alpha$ and 
$\beta$. These two half-planes ate called the {\it sides\/} of a 
dihedral angle.\par
\parshape 11 0cm 10cm 0cm 10cm 0cm 10cm 0cm 10cm 0cm 10cm 0cm 10cm
0cm 10cm 0cm 10cm 0cm 10cm 0cm 10cm 4.5cm 5.5cm
     Let's consider some dihedral angle with the sides on two 
planes $\alpha$ and $\beta$ (see Fig\.~4.2 below). Let's choose 
some arbitrary point $O$ on its edge $a$ and draw the plane 
$\gamma$ passing through the point $O$ and being perpendicular 
to the line $a$. The plane $\gamma$ intersecting with the sides 
of the dihedral angle yields two rays $h$ and $k$ lying on the 
planes $\alpha$ and $\beta$ and being perpendicular to the edge 
of the dihedral angle. They form the angle $\angle hk$ which 
is called the {\it plane angle\/} of the dihedral angle at the 
point $O$. The plane angle of a dihedral angle depends on a 
point $O$ on its edge. \vadjust{\vskip 5pt\hbox to 0pt{\kern -10pt
\includegraphics{Oris31.eps}\hss}\vskip -5pt}However, 
the following theorem shows that all plane angles of a dihedral 
angle are equipollent.
\mytheorem{4.1}\parshape 1 4.5cm 5.5cm All plane angles of a given 
dihedral angle are congruent to each other.
\endproclaim
\parshape 19 4.5cm 5.5cm 4.5cm 5.5cm 4.5cm 5.5cm 4.5cm 5.5cm
4.5cm 5.5cm 4.5cm 5.5cm 0cm 10cm 0cm 10cm 
3.5cm 6.5cm 3.5cm 6.5cm 3.5cm 6.5cm 3.5cm 6.5cm 3.5cm 6.5cm
3.5cm 6.5cm 3.5cm 6.5cm 3.5cm 6.5cm 3.5cm 6.5cm 3.5cm 6.5cm
0cm 10cm 
Before proving this theorem we consider some additional constructions 
on a plane. Let $c$ be some straight line lying on a plane $\alpha$
and assume that $A$ and $B$ are two points of this line. 
\vadjust{\vskip 5pt\hbox to 0pt{\kern 0pt
\includegraphics{Oris32.eps}\hss}\vskip -5pt}On the plane
$\alpha$ we draw two lines passing through these points and being
perpendicular to the line $c$. According the
theorem~\mythetheoremchapter{8.2}{3} from
Chapter~\uppercase\expandafter{\romannumeral 3}, such two lines are
parallel. On these lines we choose the rays lying in one half-plane 
with respect to the line $c$. We denote these rays through $a$ and
$b$. On the ray $a$ coming out from the point $A$ we mark some point
$C$. Then on the ray $b$ we choose a point $D$ so that the condition
$[BD]\cong [AC]$ is fulfilled. Such a choice is enabled by the 
axiom~\mytheaxiom{A13}. Moreover, due to this axiom such a point $D$
is unique. Let's connect the points $C$ and $D$ with the points $A$ 
and $B$ by means of the segments $[AD]$ and $[BC]$.\par
     Since the lines $AC$ and $BD$ are parallel, the segment
$[BD]$ does not intersect the line $AC$. Hence, the points $B$ and
$D$ lie on one half-plane with respect to the line $AC$. Moreover,
the points $C$ and $D$ by construction lie on one half-plane with
respect to the line $AB$. Hence, the point $D$ and the ray
$[AD\rangle$ lie inside the angle $\angle BAC$. Applying the
lemma~\mythelemmachapter{6.2}{2} from 
Chapter~\uppercase\expandafter{\romannumeral 2}, we find that
the ray $[AD\rangle$ intersects the segment $[BC]$ at some its
interior point $M$.\par
     The similar considerations can be applied to the ray $[BC\rangle$
and the angle $\angle ABD$. They prove that the point $M$ is an interior
point of the segment $[AD]$.\par
     Let's consider the rectangular triangles $ABD$ and $BAC$. For these
triangles the following conditions are fulfilled:
$$
\xalignat 3
&[AB]\cong [BA],&&[AC]\cong [BD],&&\angle BAC\cong\angle ABD.
\endxalignat
$$
From these relationships we derive that the triangles $ABD$ and 
$BAC$ are congruent. Hence, $\angle DAB\cong\angle CBA$. We apply
this congruence to the triangle $AMB$. According to the 
theorem~\mythetheoremchapter{5.2}{3} from
Chapter~\uppercase\expandafter{\romannumeral 3}, it is congruent 
to itself under exchanging its vertices $A$ and $B$. This yields
$[AM]\cong [BM]$, i\.\,e\. the triangle $AMB$ is isosceles.\par
     Apart from $\angle DAB\cong\angle CBA$, the congruence of
the triangles $ABD$ and $BAC$ yields $[AD]\cong [BC]$. Combining
this relationship with the relationship $[AM]\cong [BM]$, we get
$[CM]\cong [DM]$. This result follows from the 
axiom~\mytheaxiom{A15}. It means that the triangle $CMD$ is also
isosceles. In Euclidean geometry the isosceles triangles $AMB$ 
and $CMD$ are congruent. However, in order to prove this congruence
one should use the axiom~\mytheaxiom{A20}, which is not yet
considered.\par
\demo{Proof of the theorem~\mythetheorem{4.1}} Let's consider
a dihedral angle and mark two arbitrary points $A$ and $B$ on
its edge (see Fig\.~4.4 below). Then we draw two planes passing
through these points and being perpendicular to the edge of our
dihedral angle. These planes cut out four rays on the sides of
the dihedral angle, which determine two plane angles with the 
vertices at the points $A$ and $B$. On the sides of one of these
plane angles  we choose the points $C$ and $E$. On the sides of
the other plane angle we mark two points $D$ and $F$ such that
the conditions
$$
\xalignat 2
&[AC]\cong [BD],
&&[AE]\cong [BF]
\endxalignat
$$
are fulfilled. \pagebreak Now let's draw the segments $[AD]$, $[BC]$,
$[AF]$,
and $[BE]$. Then on each side of the dihedral angle we get the pattern
shown on Fig\.~4.3. This pattern is already studied in details. We are
going to apply the above results.\par
\parshape 16 0cm 10cm 0cm 10cm 5.7cm 4.3cm 5.7cm 4.3cm 5.7cm 4.3cm
5.7cm 4.3cm 5.7cm 4.3cm 5.7cm 4.3cm 5.7cm 4.3cm 5.7cm 4.3cm 
5.7cm 4.3cm 5.7cm 4.3cm 5.7cm 4.3cm 5.7cm 4.3cm 5.7cm 4.3cm
0cm 10cm
     At the intersection of the segments $[AD]$ and $[BC]$ we have the
point $M$ and, similarly, the point $N$ is at the intersection of the
segments $[AF]$ and $[BE]$. \vadjust{\vskip 5pt\hbox to 0pt{\kern -5pt
\includegraphics{Oris33.eps}\hss}\vskip -5pt}As we know, 
the triangles $AMB$ and $ANB$ are isosceles. This fact yields the
relationships
$$
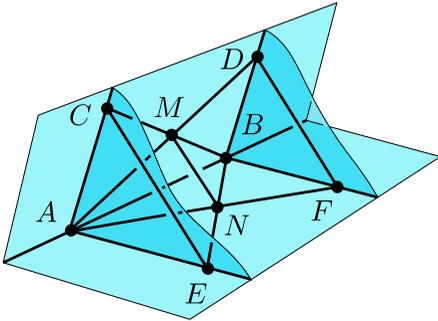
\align
&[AM]\cong [BM],\\
&[AN]\cong [BN].
\endalign
$$
Let's complement these relationships with 
$$
[MN]\cong [MN].
$$
Then, applying the theorem~\mythetheoremchapter{5.5}{3} from
Chapter~\uppercase\expandafter{\romannumeral 3}, we derive the 
congruence of the triangles $AMN$ and $BMN$. Hence, we have 
$\angle MAN\cong\angle MBN$ or, equivalently, $\angle DAF\cong
\angle CBE$.\par
     Now let's use the relationships $[AD]\cong [BC]$ and $[AF]\cong [BE]$,
which arise in considering the rectangular triangles $CAB$, $ABD$, $EAB$,
and $ABF$. Complementing them with $\angle DAF\cong\angle CBE$ and applying
the theorem~\mythetheoremchapter{5.1}{3} from
Chapter~\uppercase\expandafter{\romannumeral 3}, we derive the congruence
of the triangles $DAF$ and $CBE$. Now we have
$$
\xalignat 3
&[CE]\cong [DF],&&[AC]\cong [BD],&&[AE]\cong [BF],
\endxalignat
$$
which implies the congruence of the triangles $CAE$ and $DBF$. Hence,
$\angle CAE\cong\angle DBF$, which means that two plane angles of our
dihedral angle are congruent.
\qed\enddemo
\head
\SectionNum{5}{97} Congruent translations of a plane and the space.
\endhead
\rightheadtext{\S\,5. Congruent translations \dots}
    The concept of congruent translation for straight lines was 
introduced in \S\,3 of Chapter~\uppercase\expandafter{\romannumeral 3} 
(see definition~\mythedefinitionchapter{3.1}{3}). It is easily 
generalized for the case of planes and for the whole space.
\mydefinition{5.1} A mapping $f\!:\alpha\to\beta$ is called
a {\it congruent translation} of a plane $\alpha$ to a plane
$\beta$ if for any two points $X$ and $Y$ on the plane $\alpha$ 
the condition $[f(X)f(Y)]\cong [XY]$ is fulfilled.
\enddefinition
\mydefinition{5.2} A mapping $f\!:\Bbb E\to\Bbb E$ is called
a {\it congruent translation} of the space if for any two points 
$X$ and $Y$ the condition congruence $[f(X)f(Y)]\cong [XY]$ is 
fulfilled.
\enddefinition
\mytheorem{5.1} Let $f$ be a congruent translation of a plane to
another plane or a congruent translation of the space in whole.
In both cases the following propositions are valid:
\roster
\rosteritemwd=0pt
\item for any three points $X$, $Y$, and $Z$ from the domain 
      of the mapping $f$ if they lie on one straight line, then 
      their images $f(X)$, $f(Y)$, and $f(Z)$ also lie on one 
      straight line so that $(X\blacktriangleright Y
      \blacktriangleleft Z)$ implies $(f(X)\blacktriangleright 
      f(Y)\blacktriangleleft f(Z))$;
\item for any three points $X$, $Y$, and $Z$ from the domain 
      of $f$ if they do not lie on one straight line, 
      then their images $f(X)$, $f(Y)$, and $f(Z)$ also do not 
      lie on one straight line and the triangle $XYZ$ is
      congruent to the triangle $f(X)f(Y)f(Z)$.
\endroster
\endproclaim
     Before proving this theorem we shall formulate and prove the
following auxiliary lemma.
\mylemma{5.1} Assume that $A$, $B$, and $C$ are three points 
not lying on one straight line. Then for any point $\tilde B$ 
on the line $AC$ the conditions $[AB]\cong [A\tilde B]$ and 
$[BC]\cong [\tilde BC]$ cannot be fulfilled simultaneously.
\endproclaim
\demo{Proof}\parshape 17 0cm 10cm 0cm 10cm 0cm 10cm 0cm 10cm
0cm 10cm 0cm 10cm 3.6cm 6.4cm 3.6cm 6.4cm 3.6cm 6.4cm 3.6cm 
6.4cm 3.6cm 6.4cm 3.6cm 6.4cm 3.6cm 6.4cm 3.6cm 6.4cm 3.6cm 
6.4cm 3.6cm 6.4cm 0cm 10cm
The proof is by contradiction. Assume that $\tilde B\neq B$ 
is a point for which both conditions $[AB]\cong [A\tilde B]$ 
and $[BC]\cong [\tilde BC]$ are fulfilled. Applying the 
theorem~\mythetheorem{2.1}, we find that the points $A$ and
$C$ lie on the plane of perpendicular bisectors of the
segment $[B\tilde B]$. The intersection of this plane with 
the plane of the triangle $ABC$ is some particular 
perpendicular bisector of the segment $[B\tilde B]$ containing
both points $A$ and $C$. \vadjust{\vskip 5pt\hbox to 0pt{\kern 0pt
\includegraphics{Oris34.eps}\hss}\vskip -5pt}If $O$ 
is the center of the segment, it belongs to any perpendicular 
bisector, in particular, we have $O\in AC$. If we assume that 
$\tilde B\in AC$, then from $\tilde B\in AC$ and from $O\in AC$ 
we would conclude that the lines $\tilde BO$ and $AC$ do coincide. 
Hence, $B\in AC$, which contradicts the premise of the lemma: 
$A$, $B$, and $C$ are three points not lying on one straight line. 
The contradiction obtained shows that $\tilde B\notin AC$. The
proof of the lemma is complete.\qed\enddemo
\vskip 0pt plus 1pt minus 1pt
\demo{Proof of the theorem~\mythetheorem{5.1}}\parshape 12
0cm 10cm 0cm 10cm 4cm 6cm 4cm 6cm 4cm 6cm 4cm 6cm 4cm 6cm 4cm 6cm 
4cm 6cm 4cm 6cm 4cm 6cm 0cm 10cm 
Let's begin with proving the item \therosteritem{1} of the theorem. 
Assume that the points $X$, $Y$, $Z$ taken from the domain of the 
mapping $f$ lie on one straight line. If the main proposition in
\vadjust{\vskip 5pt\hbox to 0pt{\kern -5pt
\includegraphics{Oris35.eps}\hss}\vskip -5pt}the item
\therosteritem{1} is not valid, then $f(X)$, $f(Y)$, and
$f(Z)$ does not lie on one straight line at least for one particular
choice of the points $X$, $Y$, $Z$ lying on one straight line. Then
they define a triangle $f(X)f(Y)f(Z)$ for whose sides the relationships
$[XY]\cong [f(X)f(Y)]$, $[YZ]\cong [f(Y)f(Z)]$, and $[XZ]\cong [f(X)f(Z)]$
are fulfilled. Relying upon the last relationship $[XZ]\cong [f(X)f(Z)]$,
we apply the theorem~\mythetheoremchapter{2.1}{3} from
Chapter~\uppercase\expandafter{\romannumeral 3} to the points $X$, $Y$,
and $Z$. Due to this theorem we can find a point $\tilde Y$ on the line
$[f(X)f(Z)]$ such that $[XY]\cong [f(X)\tilde Y]$ and $[YZ]\cong 
[\tilde Yf(Z)]$. Comparing these relationships with $[XY]\cong [f(X)f(Y)]$
and $[YZ]\cong [f(Y)f(Z)]$, we note that we got exactly in a situation
forbidden by the above lemma~\mythelemma{5.1}. Hence, our assumption
that the points $f(X)$, $f(Y)$, and $f(Z)$ do not lie on one straight
line is wrong. Thus, we have proved the main proposition 
of the item \therosteritem{1} in the theorem~\mythetheorem{5.1}. Now 
the relationship $(f(X)\blacktriangleright f(Y)\blacktriangleleft f(Z))$ 
is derived from $(X\blacktriangleright Y\blacktriangleleft Z)$ by means 
of the theorem~\mythetheoremchapter{2.2}{3} from 
Chapter~\uppercase\expandafter{\romannumeral 3}.\par
\parshape 12 0cm 10cm 0cm 10cm 4cm 6cm 4cm 6cm 4cm 6cm 4cm 6cm 4cm 6cm 
4cm 6cm 4cm 6cm 4cm 6cm 4cm 6cm 
0cm 10cm
     Let's prove the item \therosteritem{2} of the theorem. Now the 
points $X$, $Y$, and $Z$ do not lie on one straight line. Assume that 
their images $f(X)$, $f(Y)$, and $f(Z)$ appear to be lying on one 
straight line. \vadjust{\vskip 5pt\hbox to 0pt{\kern -5pt
\includegraphics{Oris36.eps}\hss}\vskip -5pt}Let's use the
relationship $[XZ]\cong [f(X)f(Z)]$ and apply the
theorem~\mythetheoremchapter{2.1}{3} from
Chapter~\uppercase\expandafter{\romannumeral 3}. Due to this theorem
we can find a point $\tilde Y$ lying on the segment $[XZ]$ and such 
that $[X\tilde Y]\cong [f(X)f(Y)]$ and $[\tilde YZ]\cong [f(Y)f(Z)]$.
Now, if we recall the relationship $[XY]\cong [f(X)f(Y)]$ and
$[YZ]\cong [f(Y)f(Z)]$, we see that again we are in a situation
forbidden by the lemma~\mythelemma{5.1}. Hence, our preliminary 
assumption that the points $f(X)$, $f(Y)$, and $f(Z)$ lie on one 
straight line is not valid. The main proposition of the item
\therosteritem{2} in the theorem~\mythetheorem{5.1} is proved.
The rest is to prove the congruence of the triangles $XYZ$ and
$f(X)f(Y)f(Z)$. It follows from $[XY]\cong [f(X)f(Y)]$, $[YZ]\cong
[f(Y)f(Z)]$, and $[XZ]\cong [f(X)f(Z)]$.
\qed\enddemo
    As a corollary of the theorem~\mythetheorem{5.1} and the
theorem~\mythetheoremchapter{2.2}{3} from 
Chapter~\uppercase\expandafter{\romannumeral 3} we find that 
each congruent translations maps a straight line onto a straight 
line and a ray onto a ray. Under such a mapping each angle is
mapped onto a congruent angle. In particular, this means that
congruent translations preserve orthogonality of lines.\par
     Let $f\!:\Bbb E\to\Bbb E$ be a congruent translation of the
space. Using the theorem~\mythetheorem{1.5}, we conclude that 
such a mapping takes a plane onto a plane preserving the 
orthogonality of planes and preserving the orthogonality of a 
plane and a straight line. Hence, the restriction of a congruent 
translation of the space to some plane appears to be a congruent 
translation of planes, while the restriction of a congruent 
translation of a plane to some straight line appears to be a 
congruent translation of straight lines.\par
     Let's consider a congruent translation $f$ mapping a plane 
$\alpha$ onto a plane $\beta$. Let $a$ be some straight line on
the plane $\alpha$. Then, as we noted above, the points of the 
line $a$ are mapped to the points of some straight line $b$ lying
on the plane $\beta$. Let $X$ and $Y$ be two points of the plane
$\alpha$ lying in different half-planes with respect to the line
$a$. Then the segment $[XY]$ crosses the line $a$ at some its
interior point $O$. Due to the theorem~\mythetheorem{5.1} the
point $f(O)$ is the intersection point of the lines $f(X)f(Y)$ 
and $b$. From $(X\blacktriangleright O\blacktriangleleft Y)$ we 
derive $(f(X)\blacktriangleright f(O)\blacktriangleleft f(Y))$. 
This means that the points $f(X)$ and $f(Y)$ lie on different 
sides of the line $b$ on the plane $\beta$.\par
    Now let's consider two points $X$ and $Y$ lying on the same
half-plane with respect to the line $a$ on the plane $\alpha$.
Let $Z$ be a point lying on the other half-plane. Then due to
above considerations the points $f(X)$ and $f(Z)$ lie on different
half-planes with respect to the line $b$ on the plane $\beta$. 
The points $f(Y)$ and $f(Z)$ are also on different half planes.
Therefore, the points $f(X)$ and $f(Y)$ lie on the same side of 
the line $b$ on the plane $\beta$. Thus, we have the following
result.
\mytheorem{5.2} Any congruent translation of planes and any
congruent translation of the whole space are 
{\tencyr\char '074}half-planes preserving maps{\tencyr\char '076}, 
i\.\,e\. they take a half-plane onto a half-plane. 
\endproclaim
\mytheorem{5.3} Each congruent translation of the whole space
takes a half-space onto a half-space.
\endproclaim
     The proof of the theorem~\mythetheorem{5.3} is analogous
to the proof of the theorem~\mythetheorem{5.2}. We do not give 
it here.\par
     Let $f\!:\alpha\to\beta$ be a congruent translation of planes. 
Let's consider some arbitrary straight line $a$ on the plane $\alpha$. 
The points of this line are mapped into some definite line $b$ on the
plane $\beta$. The line $a$ divides the plane $\alpha$ into two
half-planes $a_{+}$ and $a_{-}$. Let's denote by $b_{+}$ that 
half-plane on the plane $\beta$ to which the points of $a_{+}$ are
mapped. Let $X$ be some point of the half-plane $a_{+}$. Within the
plane $\alpha$ we can drop the perpendicular onto the line 
$a$. Let's denote by $X_0$ the foot of this perpendicular. Any
congruent translation preserves the orthogonality, therefore, the
point $f(X_0)$ is the foot of the perpendicular dropped from 
$f(X)$ onto the line $b$ within the plane $\beta$. Being more
precise, the following three conditions are fulfilled for the
point $f(X)$:
$$
\xalignat 3
&[f(X)f(X_0)]\perp b,&&[f(X)f(X_0)]\cong [XX_0],
&&f(X)\in b_{+}.
\endxalignat
$$
Note that these conditions fix uniquely the point $f(X)$ on the 
plane $\beta$ provided the point $f(X_0)$ on the line $b$ is given. 
\vadjust{\vskip 5pt\hbox to 0pt{\kern -10pt
\includegraphics{Oris37.eps}\hss}\vskip 230pt}In a 
similar way, for an arbitrary point $Y$ from the half-plane $a_{-}$ 
and for the foot $Y_0$ of the perpendicular dropped from $Y$ onto 
the line $a$ the following conditions are fulfilled:
$$
\xalignat 3
&[f(Y)f(Y_0)]\perp b,&&[f(Y)f(Y_0)]\cong [YY_0],
&&f(Y)\in b_{-}.
\endxalignat
$$
These conditions fix uniquely the point $f(Y)$ on the plane $\beta$ 
upon fixing the point $f(Y_0)$ on the line $b$. These observations 
lead to the following theorem.
\mytheorem{5.4} Let $a$ be some straight line on a plane $\alpha$
dividing this plane into two half-planes $a_{+}$ and $a_{-}$. Let 
$b$ be some straight line on a plane $\beta$ dividing this plane 
into two half-planes $b_{+}$ and $b_{-}$. Then any congruent 
translation of lines $f\!: a\to b$ has a unique extension up to
a congruent translation of planes $f\!:\alpha\to\beta$ such that
the half-plane $a_{+}$ is mapped onto $b_{+}$ and the half-plane 
$a_{-}$ is mapped onto $b_{-}$.
\endproclaim
\demo{Proof} Let $X\in a_{+}$. We determine $f(X)$ by means of the
following construction. From the point $X$ we drop the perpendicular
onto the line $a$ and denote by $X_0$ the foot of this perpendicular.
Applying the congruent translation of lines $f\!: a\to b$ to the
point $X_0$, we get the point $f(X_0)$ on the line $b$. Then on the
plane $\beta$ we draw the line passing through the point $f(X_0)$ 
and being perpendicular to the line $b$. On this line we choose the
ray coming out from the point $f(X_0)$ and lying on the closed 
half-plane $\overline{b_{+}}$. Then the point $f(X)$ on this ray is 
determined by the condition $[XX_0]\cong [f(X)f(X_0)]$. For a point
$Y\in a_{-}$ the procedure of constructing the point $f(Y)$ differs 
only by the choice of the ray $[f(Y_0)f(Y)\rangle$ lying not on
$\overline{b_{+}}$, but on $\overline{b_{-}}$ (see Fig\.~5.4 above).
\par
     The above construction yields a mapping $f\!:\alpha\to\beta$.
For the points $X\in a$ this mapping coincide with the initial
mapping $f\!: a\to b$. Let's prove that this mapping is a congruent
translation of planes. For this purpose we need to prove that
$[XY]\cong [f(X)f(Y)]$ for any two points $X$ and $Y$ on the plane
$\alpha$. Let's consider the following four cases:
\roster
\rosteritemwd=5pt
\item both point $X$ and $Y$ lie on the line $a$;
\item only one of the points $X$ or $Y$ lies on the line $a$;
\item the points $X$ and $Y$ do not belong to the line $a$ and 
      are on the same side with respect to this line;
\item the points $X$ and $Y$ are on different sides with respect
      to $a$.
\endroster\par
\parshape 3 0cm 10cm 0cm 10cm 6.2cm 3.8cm
     In the first case the relationship  $[XY]\cong [f(X)f(Y)]$ 
follows from the fact that the initial mapping $f\!: a\to b$ is
a congruent translation of lines.\par
\parshape 11 6.2cm 3.8cm 6.2cm 3.8cm 6.2cm 3.8cm 6.2cm 3.8cm 6.2cm 3.8cm
6.2cm 3.8cm 6.2cm 3.8cm 6.2cm 3.8cm 6.2cm 3.8cm 6.2cm 3.8cm 0cm 10cm
     In the second case we assume that $Y\in a$ and $X\notin a$ for the
sake of certainty. Then on the planes $\alpha$ and $\beta$ we have two
rectangular triangles (see Fig\.~5.5). 
\vadjust{\vskip 5pt\hbox to 0pt{\kern -10pt
\includegraphics{Oris38.eps}\hss}\vskip -5pt}In these
triangles for their sides $[X_0Y]$ and $[f(X_0)f(Y)]$ we have the
relationship $[X_0Y]\cong [f(X_0)f(Y)]$ following from the fact
that $f\!: a\to b$ is a congruent translation of lines. Moreover, 
by construction of the point $f(X)$ we have $[XX_0]\cong [f(X)f(X_0)]$.
Combining this relationship with the congruence of the right angles
$\angle XX_0Y\cong\angle f(X)f(X_0)f(Y)$, we get the congruence of
the triangles $XX_0Y$ and $f(X)f(X_0)f(Y)$. This congruence yields
the required relationship $[XY]\cong [f(X)f(Y)]$.\par
\parshape 15 0cm 10cm 0cm 10cm 5.8cm 4.2cm 5.8cm 4.2cm 5.8cm 4.2cm
5.8cm 4.2cm 5.8cm 4.2cm 5.8cm 4.2cm 5.8cm 4.2cm 5.8cm 4.2cm 5.8cm 4.2cm
5.8cm 4.2cm 5.8cm 4.2cm 5.8cm 4.2cm 0cm 10cm
    Let's consider the third case. The congruence of the rectangular
triangles $X_0Y_0Y$ and $f(X_0)f(Y_0)f(Y)$ on this case is proved 
just \vadjust{\vskip 5pt\hbox to 0pt{\kern 5pt
\includegraphics{Oris39.eps}\hss}\vskip -5pt}like above 
in the previous case. From this congruence we derive 
that the angle $\angle YX_0Y_0$ and the angle $\angle f(Y)f(X_0)f(Y_0)$ 
are congruent. We also derive the congruence of the segments $[X_0Y]$ 
\ and \ $[f(X_0)f(Y)]$. Then we take into account the congruence of 
the right angles \ $\angle f(X)f(X_0)f(Y_0)$ and $\angle XX_0Y_0$ \ and 
\ apply the theorem~\mythetheoremchapter{5.3}{3} from
Chapter~\uppercase\expandafter{\romannumeral 3}. As a result we
obtain the congruence of the angle $\angle f(X)f(X_0)f(Y)$ and
the angle $\angle XX_0Y$. From the relationship $[XX_0]\cong [f(X)f(X_0)]$ 
and from the relationship $[X_0Y]\cong [f(X_0)f(Y)]$ now we derive that
the triangles $XX_0Y$ and $f(X)f(X_0)f(Y)$ are congruent. This congruence
yields the required relationship $[XY]\cong [f(X)f(Y)]$ for $X$ and $Y$.
\par
\parshape 12 5.8cm 4.2cm 5.8cm 4.2cm 5.8cm 4.2cm 5.8cm 4.2cm 5.8cm 4.2cm
5.8cm 4.2cm 5.8cm 4.2cm 5.8cm 4.2cm 5.8cm 4.2cm 5.8cm 4.2cm 5.8cm 4.2cm
0cm 10cm
      The proof of the congruence $[XY]\cong [f(X)f(Y)]$
in the fourth case almost literally the same as the proof of this fact in
the third case. 
\vadjust{\vskip 5pt\hbox to 0pt{\kern 5pt
\includegraphics{Oris40.eps}\hss}\vskip -5pt}However, we have
quite \ different \ picture \ in this case since the points $X$ and $Y$ 
are on different sides of the line $a$. By construction their images $f(X)$
and $f(Y)$ are also on different sides of the corresponding line $b$.\par
      In the cases \therosteritem{2}, \therosteritem{3}, and 
\therosteritem{4} there are three degenerate subcases, where $X_0=Y_0$.
We leave the proof of the relationship $[XY]\cong [f(X)f(Y)]$ in 
these subcases to the reader as an exercise.\par
      Thus we have constructed a mapping $f\!:\alpha\to\beta$ extending
the initial mapping $f\!: a\to b$ and have proved that it is a congruent
translation. The uniqueness of such an extension follows from the 
considerations preceding the statement of the theorem~\mythetheorem{5.4}.
\qed\enddemo
\mytheorem{5.5} Let $\alpha$ be some plane dividing the space into two
half-spaces $\alpha_{+}$ and $\alpha_{-}$. Let $\beta$ be another plane
dividing the space into two half-spaces $\beta_{+}$ and $\beta_{-}$.
Then any congruent translation of planes $f\!:\alpha\to\beta$ has a 
unique extension up to a congruent translation of the whole space 
$f\!:\Bbb E\to\Bbb E$ such that $\alpha_{+}$ is mapped onto $\beta_{+}$ 
and $\alpha_{-}$ is mapped onto $\beta_{-}$.
\endproclaim
     The construction of the required mapping $f\!:\Bbb E\to\Bbb E$
is analogous to that we used above in proving the 
theorem~\mythetheorem{5.4}. Other details of the proof for the 
theorem~\mythetheorem{5.5} are also very similar to those for the
theorem~\mythetheorem{5.4}.
\myexercise{5.1} Prove the theorems~\mythetheorem{5.3} and
\mythetheorem{5.5}.
\endproclaim
\myexercise{5.2} Prove that the inversion $i_O\!:\Bbb E\to\Bbb E$
with respect to a point $O$ is a congruent translation of the whole
space.
\endproclaim
\head
\SectionNum{6}{105} Mirror reflection in a plane and 
in a straight line.
\endhead
\rightheadtext{\S\,6. Mirror reflection \dots}
\parshape 13 0cm 10cm 0cm 10cm 0cm 10cm 4.3cm 5.7cm 4.3cm 5.7cm
4.3cm 5.7cm 4.3cm 5.7cm 4.3cm 5.7cm 4.3cm 5.7cm 4.3cm 5.7cm 4.3cm
5.7cm 4.3cm 5.7cm 0cm 10cm
     Let's apply the theorems~\mythetheorem{5.4} and \mythetheorem{5.5}
in order to construct some particular mappings of congruent translation.
Let's consider a plane $\alpha$ and denote by $\alpha_{+}$ and 
$\alpha_{-}$ the half-planes produced in the space $\Bbb E$ by this 
plane. \vadjust{\vskip 5pt\hbox to 0pt{\kern -5pt
\includegraphics{Oris41.eps}\hss}\vskip -5pt}Let's set 
$\beta=\alpha$, $\beta_{+}=\alpha_{-}$, and $\beta_{-}=\alpha_{+}$. 
For the initial mapping we take the identical mapping of the plane 
$\alpha$ onto itself. Since $\beta=\alpha$, we can treat it as
$f\!:\alpha\to\beta$. Applying the theorem~\mythetheorem{5.5} to
$f\!:\alpha\to\beta$, we get the mapping of congruent translation
$z_\alpha\!:\Bbb E\to\Bbb E$. This mapping exchanges the half-spaces 
$\alpha_{+}$ and $\alpha_{-}$, leaving stable the points of the plane 
$\alpha$. Such a mapping $z_\alpha$ is called a {\it mirror reflection 
in a plane $\alpha$}.\par
\parshape 5 0cm 10cm 0cm 10cm 0cm 10cm 0cm 10cm 3.5cm 6.5cm 
     Let $a$ be some straight line lying on the plane $\alpha$. It divides
this plane into two half-planes $a_{+}$ and $a_{-}$. Let's denote $b=a$,
$b_{+}=a_{-}$, and $b_{-}=a_{+}$. For the initial mapping $f\!:a\to b$ we
take the identical mapping of the line $a$ onto itself. Then, applying
the theorem~\mythetheorem{5.4}, we get a mapping of congruent translation 
of plane $z_a\!:\alpha\to\alpha$ that exchanges the half-planes $a_{+}$ 
\vadjust{\vskip 5pt\hbox to 0pt{\kern -15pt
\includegraphics{Oris42.eps}\hss}\vskip -5pt}and $a_{-}$. 
It is called the {\it mirror reflection of a plane $\alpha$ in a line
$a$}.\par
\parshape 6 3.5cm 6.5cm 3.5cm 6.5cm 3.5cm 6.5cm 3.5cm 6.5cm 3.5cm 6.5cm
0cm 10cm
    The mapping $z_a$ defined just above can be extended up to 
a mapping of the whole space. For this purpose we denote again 
$\beta=\alpha$, $\beta_{+}=\alpha_{-}$, and $\beta_{-}=\alpha_{+}$. 
For the initial mapping $f\!:\alpha\to\beta$ mapping now we take 
the mirror reflection of the plane $\alpha$ in the line $a$. Applying 
the theorem~\mythetheorem{5.5}, we define a mapping $z_a\!:\Bbb E\to
\Bbb E$, which is called the {\it mirror reflection of the space 
$\Bbb E$ in a line $a$}. The plane $\alpha$ plays an auxiliary role 
in defining this mapping $z_a\!:\Bbb E\to\Bbb E$. There is the 
following theorem that yields a different way for constructing the 
mapping $z_a\!:\Bbb E\to\Bbb E$.
\mytheorem{6.1} For any point $X\notin a$ the segment connecting
$X$ with its mirror symmetric point $z_a(X)$ intersects the line 
$a$ at the point $X_0$ being its center and this segment is 
perpendicular to the line $a$.
\endproclaim
\myexercise{6.1} Verify that the theorem~\mythetheorem{6.1} fixes
uniquely the point $z_a(X)$ provided the point $X$ is given.
\endproclaim
\myexercise{6.2} Prove the theorem~\mythetheorem{6.1} and show
that the mirror reflections in a line $z_a$ and in a plane $z_\alpha$ 
satisfy the identities $z_\alpha\compos z_\alpha=\id$ and $z_a\compos
z_a=\id$.
\endproclaim
\myexercise{6.3} Let a straight line $a$ be a perpendicular to a plane
$\alpha$ passing through some point $O\in\alpha$. In this case prove that
$z_a\compos z_\alpha=z_\alpha\compos z_a=i_O$.
\endproclaim
\head
\SectionNum{7}{106} Rotation of a plane about a point.
\endhead
\rightheadtext{\S\,7. Rotation of a plane about a point.}
\parshape 13 0cm 10cm 0cm 10cm 0cm 10cm 4.3cm 5.7cm 4.3cm 5.7cm
4.3cm 5.7cm 4.3cm 5.7cm 4.3cm 5.7cm 4.3cm 5.7cm 4.3cm 5.7cm 
4.3cm 5.7cm 4.3cm 5.7cm 0cm 10cm
    Assume that in a plane $\alpha$ an angle $\angle hk$ with the vertex
at a point $O$ is given. Let's extend the ray $h$ up to the whole line
$a$, and extend the ray $k$ up to the whole line $b$. The angle 
$\angle hk$ is the intersection of two half-planes determined by 
the lines $a$ and $b$. \vadjust{\vskip 5pt\hbox to 0pt{\kern 5pt
\includegraphics{Oris43.eps}\hss}\vskip -5pt}For the sake 
of certainty let's set $\angle hk=\overline{a_{+}}\cap\overline{b_{-}}$. 
According to the
result of \S\,3 in Chapter~\uppercase\expandafter{\romannumeral 3},
there is a unique mapping of congruent translation of lines
$f^{+}_{OO}\!:a\to b$ taking the point $O$ to itself and mapping the
ray to the ray $k$. Applying the theorem~\mythetheorem{5.4}, we extend
this mapping up to a congruent translation of planes $\theta_{hk}\!:
\alpha\to\alpha$ mapping $a_{+}$ to $b_{+}$ and $a_{-}$ to $b_{-}$.
Such a mapping is called the {\it rotation of a plane $\alpha$ about
a point $O$ by the angle $\angle hk$ from the ray $h$ toward the 
ray $k$}.\par
     Exchanging the rays $h$ and $k$ we do not change the angle 
$\angle hk$. However, in this case we get the rotation in the 
opposite direction from the ray $k$ toward the ray $h$. The 
mapping $\theta_{kh}$ is inverse for $\theta_{hk}$, i\.\,e\. we
have $\theta_{hk}\compos\theta_{kh}=\theta_{kh}\compos\theta_{hk}
=\id_\alpha$.\par
     Assume that in a plane $\alpha$ some rotation angle $\angle hk$
is given. Let's study the procedure of constructing the point 
$Y=\theta_{hk}(X)$ for some arbitrary point $X\in\alpha$. For the
sake of certainty e assume that $X\in a_{+}$. On the rays $h$ and 
$k$ we choose two points $A$ and $B$ so that the relationships
$[OA]\cong [OX]$ and $[OB]\cong [OX]$ are fulfilled. Then  
$B=\theta_{hk}(A)$. Let's consider the isosceles triangle $AOX$. 
It lies on the closed half-plane $\overline{a_{+}}$, which is
taken to $\overline{b_{+}}$ under the rotation $\theta_{hk}$.
Let's construct the triangle $BOY$ congruent to $AOX$ in the
half-plane $\overline{b_{+}}$. For this purpose in $\overline{b_{+}}$
we choose a ray forming with the ray $k$ an angle congruent to the 
angle $\angle XOA$. Afterwards, on this ray we mark a point $Y$ 
such that $[OY]\cong [OX]$. Now from $O=\theta_{hk}(O)$, from
$B=\theta_{hk}(A)$, and from the congruence of the triangles $AOX$ 
and $BOY$ we derive $Y=\theta_{hk}(X)$. For the case, where 
$X\in a_{-}$, the procedure of constructing the point $Y=\theta_{hk}(X)$
is analogous to the above one. The only difference is that the 
triangle $BOY$ is chosen in the other half-plane $\overline{b_{-}}$.
\par
\mytheorem{7.1}\parshape 1 4cm 6cm The \ rotation \ mapping
$\theta_{hk}\!:\alpha\to\alpha$ has exactly one stable point. 
\vadjust{\vskip 5pt\hbox to 0pt{\kern 0pt
\includegraphics{Oris44.eps}\hss}\vskip -5pt}This is 
the point $O=\theta_{hk}(O)$, about which the rotation is performed.
\endproclaim
\demo{Proof}\parshape 5 4cm 6cm 4cm 6cm 4cm 6cm 4cm 6cm 0cm 10cm
The mapping $\theta_{hk}$ is an extension of the mapping
$f^{+}_{OO}\!:a\to b$ due to the theorem~\mythetheorem{5.4}. The point
$O$ is a stable point for $f^{+}_{OO}\!:a\to b$, hence it is a stable
point for $\theta_{hk}$ either. Let's show that the mapping  
$\theta_{hk}\!:\alpha\to\alpha$ has no other stable points.\par
     Assume that it is not so. If $X$ is another stable point, then the
line $OX$ consists of the stable points of the mapping $\theta_{hk}$
(this fact follows from the theorem~\mythetheoremchapter{2.1}{3} in
Chapter~\uppercase\expandafter{\romannumeral 3}). The line $OX$ is 
distinct from $a$ (since $a$ is mapped to $b$), therefore, we can 
assume that the stable point $X\neq O$ is initially chosen on the
half-plane $a_{+}$. Let's apply to $X$ the above procedure of 
constructing the point $Y=\theta_{hk}(X)$. The condition $X=Y$ leads 
to the situation shown on Fig\.~7.2. From $Y\in b_{+}$ and $A\in b_{-}$
we conclude that the points $X$ and $A$ are on different sides of the
line $b$. Therefore the segment $[AX]$ crosses the line $OB$ at some
its interior point $C$ and $A$ is an external point with respect to 
the segment $[XC]$. Hence, from $X\in a_{+}$ it follows that $C\in a_{+}$.
But $B$ lies on the half-plane $a_{+}$. Hence, the points $B$ and
$C$ lie on one side with respect to the line $OA$ and the point $O$ 
is outside the segment $[BC]$. From this fact it follows that the 
points $B$ and $C$ lie on one side with respect to the line $OX$.
Now, applying the congruence of angles $\angle OYB\cong\angle OXA$,
we derive that the angles $[XA\rangle$ and $[XB\rangle$ do coincide.
Then from $[AX]\cong [BY]$ we conclude that $A=B$. But this 
contradicts the fact that the rays $h\neq k$ form the angle 
$\angle hk$. The contradiction obtained shows that the mapping
$\theta_{hk}$ has no stable points other than the point $O$.
\qed\enddemo
\mytheorem{7.2} Let $h$ and $k$ be two rays coming out from the point
$O$ and lying on a plane $\alpha$. There are exactly two mappings of
congruent translation $f\!:\alpha\to\alpha$ with the stable point $O$
that take the ray $h$ to the ray $k$. The first of them $f=\theta_{hk}$ 
is the rotation by the angle $\angle hk$ about the point $O$ and the 
second one $f=z_m$ is the mirror reflection of the plane $\alpha$ in
the line $m$ containing the bisector of the angle $\angle hk$.
\endproclaim
\demo{Proof} It is easy to see that both mappings $\theta_{hk}$ and $z_m$
take the ray $h$ to the ray $k$. Let $a$ be the line containing the ray
$h$ and let $b$ be the line containing the ray $k$. The restrictions of
$\theta_{hk}$ and $z_m$ to the line $a$ coincide with the mapping of 
congruent translation of lines $f^{+}_{OO}\!:a\to b$. But due to the
theorem~\mythetheorem{5.4} we have exactly two extensions of this 
mapping up to a congruent translation of planes $f\!:\alpha\to\alpha$. 
The first extension takes $a_{+}$ to $b_{+}$ and $a_{-}$ to $b_{-}$, 
it is the rotation $f=\theta_{hk}$. The second one takes $a_{+}$ to 
$b_{-}$ and $a_{-}$ to $b_{+}$, it is the mirror reflection $f=z_m$.
\qed\enddemo
{\bf A remark}. In theorem~\mythetheorem{7.2} there are two special
dispositions of the rays $h$ and $k$ where they lie on one straight 
line $n$. If $h=k$, then we set $\theta_{hk}=\theta_{hh}=\id$ 
by definition. In this case the line $m$ coincides with $n$. If the rays
$h$ and $k$ are opposite to each other, they define a straight angle,
the bisector of this angle is the perpendicular to the line $n$ passing
through the point $O$. The rotation $\theta_{hk}$ in this case is set to
be coinciding with the inversion of the plane $\alpha$ with respect to 
the point $O$. With these additional provisions the 
theorem~\mythetheorem{7.2}
remains valid both special dispositions of the rays $h$ and $k$.
\mytheorem{7.3} For any ray $l$ coming out from a point $O$ and 
lying on a plane $\alpha$ the equality $\theta_{hk}(l)=q$ implies
$\theta_{hk}=\theta_{lq}$.
\endproclaim
\demo{Proof} Indeed, $\theta_{hk}(l)=q$ means that the mapping of 
congruent translation $\theta_{hk}$ takes the ray $l$ to the ray $q$.
It has the unique stable point $O$. Therefore $\theta_{hk}$ coincides
with $\theta_{lq}$.
\qed\enddemo
\mytheorem{7.4} The mappings of rotation of a plane $\alpha$ about
a point $O$ and the mappings of mirror reflections in straight lines
passing through this point possess the following properties:
$$
\xalignat 2
&\theta_{kl}\compos\theta_{hk}=\theta_{hl};
&&\theta_{hk}\compos z_m=z_m\compos \theta_{kh};\\
&z_m\compos z_m=\id;&&z_m\compos z_n=\theta_{hk},
\text{ \ where \ }h\subset n,\,k=z_m(h).
\endxalignat
$$
\endproclaim
\demo{Proof} Let $h$, $k$, and $l$ be three rays coming out from 
a point $O$ and lying on a plane $\alpha$. Assume that none two
of them lie on one straight line. Then they define three angles 
$\angle hk$, $\angle kl$, $\angle hl$ and three rotations  
$\theta_{hk}$, $\theta_{kl}$, $\theta_{hl}$. Let's denote by $f$ 
the composition $f=\theta_{kl}\compos\theta_{hk}$. It is easy to 
see that the mapping $f$ takes the ray $h$ to the ray $l$. It has 
the stable point $O$. Due to the theorem~\mythetheorem{7.2} we 
have two options: $f=\theta_{hl}$ or $f=z_m$.\par
     Let's prove that $f$ has the unique stable point $O$. If we 
assume that there is another stable point $X$, then the ray
$q=[OX\rangle$ consists of stable points for the mapping $f$.
In this case we would have $\theta_{kl}\compos\theta_{hk}(q)=q$ 
or $\theta_{hk}(q)=\theta_{lk}(q)$. We denote $\theta_{hk}(q)=p$ 
and apply the theorem~\mythetheorem{7.3}. As a result we get 
$\theta_{hk}=\theta_{qp}$ and $\theta_{lk}=\theta_{qp}$. 
Hence, $\theta_{kh}=\theta_{kl}$, which leads to the coincidence
of the rays $h=l$. But this contradicts to the initial assumption 
that none of the rays $h$, $k$, and $l$ lies on one straight line.
This contradiction excludes the option $f=z_m$ and proves the
first relationship $\theta_{kl}\compos\theta_{hk}=\theta_{hl}$.\par
     Let $m$ be some straight line lying on the plane $\alpha$ and
passing through the point $O$. The point $O$ divides the line $m$ 
into two rays. Let's denote one of them by $l$ and denote
$q=\theta_{hk}(l)$. Then $\theta_{hk}=\theta_{lq}$, which follows
from the theorem~\mythetheorem{7.3}. Let's denote $p=\theta_{ql}(l)$. 
By construction the rays $p$ and $q$ lie on different sides of the
line $m$. From $p=\theta_{ql}(l)$ and $l=\theta_{ql}(q)$ we get
these rays form two congruent angles $\angle lq$ and $\angle lp$
with the ray $l$ lying on the line $m$. Hence, $p=z_m(q)$ and
$q=z_m(p)$. This yields
$$
\gathered
\theta_{hk}\compos z_m(l)=\theta_{hk}(l)=\theta_{lq}(l)=q,\\
z_m\compos\theta_{kh}(l)=z_m\compos\theta_{ql}(l)=z_m(p)=q.
\endgathered
$$
Two mappings $f=\theta_{hk}\compos z_m$ and $g=z_m\compos\theta_{kh}$
take the ray $l$ to the ray $q$. The equality $f=\theta_{lq}$ is
excluded since $f=\theta_{hk}\compos z_m=\theta_{lq}$ implies $z_m=\id$.
The equality $g=\theta_{lq}$ is also excluded since 
$g=z_m\compos\theta_{kh}=\theta_{lq}$ would lead to $z_m=\theta_{lq}
\compos\theta_{pl}=\theta_{pq}$. Hence, due to the
theorem~\mythetheorem{7.2} the mappings $f$ and $g$ coincide with 
the mirror reflection in the line containing the bisector of the
angle $\angle lq$. Therefore, $f=g$. This fact proves the second 
relationship $\theta_{hk}\compos z_m=z_m\compos \theta_{kh}$
in the theorem~\mythetheorem{7.4}.\par
     The third relationship $z_m\compos z_m=\id$ follows immediately
from the definition of the mirror reflection $z_m\!:\alpha\to\alpha$
of a plane in a line (see \S\,6 above).\par
     Let $m\neq n$ be two straight lines lying on the plane $\alpha$ 
and intersecting each other at the point $O$. The point $O$ divides
the line $n$ into two rays. We denote one of them by $h$ and set
$k=z_m(h)$. Then for the mapping $f=z_m\compos z_n$ we have $f(h)=k$. 
Then the theorem~\mythetheorem{7.2} provides two options: 
$f=\theta_{hk}$ or $f=z_u$, where $u$ is the line containing the 
bisector of the angle $\angle hk$.\par
     Let's show that the mapping $f$ has the unique stable point 
$O$. If we assume that there is another stable point $X$, then from  
$z_m\compos z_n(X)=X$ we derive $z_m(X)=z_n(X)$. Let's denote
$Y=z_m(X)=z_n(X)$. The coincidence $X=Y$ is excluded since the point
$O$ is the unique common stable point for the mappings $z_m$ and
$z_n$. Then from $Y=z_m(X)$ we get that $m$ with the perpendicular
bisector of the segment $[XY]$ lying on the plane $\alpha$. Due to
$Y=z_n(X)$ the line $n$ also coincides with this perpendicular
bisector, which contradicts the initial assumption $m\neq n$. Hence, 
$f$ nas no stable points other than $O$. Therefore, $f\neq z_u$. 
In this case we have $f=\theta_{hk}$, which completes the proof of 
the fourth relationship and completes the proof of the theorem in
whole.\qed\enddemo
\myexercise{7.1} Compare the theorem~\mythetheorem{7.4} 
with the theorem~\mythetheoremchapter{3.2}{3} in
Chapter~\uppercase\expandafter{\romannumeral 3}.
\endproclaim
\myexercise{7.2} Prove that the theorem~\mythetheorem{7.3} and 
the theorem~\mythetheorem{7.4} remain valid in special cases 
where some two straight lines do coincide or some rays appear
to be lying on one straight line.
\endproclaim
\head
\SectionNum{8}{111} The total rotation group and the group of pure
rotations of a plane.
\endhead
\rightheadtext{\S\,8. The total rotation group and \dots}
     Let $f\!:\alpha\to\alpha$ be some mapping of the congruent 
with a stable point $O$. Such a mapping is sometimes called 
a {\it generalized rotation of a plane $\alpha$ about a point
$O$}. Let's choose some arbitrary ray $h\subset\alpha$ coming 
out from the point $O$. The mapping $f$ takes it to the ray
$k=f(h)$ coming out from the same point. Now, applying the 
theorem~\mythetheorem{7.2}, we conclude that any generalized 
rotation $f$ of a plane $\alpha$ about a point $O$ is some
rotation $\theta_{hk}$ about this point or a reflection $z_m$
in a line passing through the point $O$. As a simple consequence
of this fact we get that such a mapping $f$ is bijective and has
the inverse mapping $f^{-1}$, which is also a generalized 
rotation about the point $O$. The set of all generalized 
rotations of a plane $\alpha$ about a point $O$ is a group with
respect to the composition (see 
definition~\mythedefinitionchapter{4.2}{3} in 
Chapter~\uppercase\expandafter{\romannumeral 3}). The unity of this 
{\it total rotation group\/} is the identical mapping $\id$, which
can be interpreted as a special case of the rotation: $\id=\theta_{hh}$
(see the remark to the theorem~\mythetheorem{7.2}).\par
     The set of pure (non-generalized) rotations of a plane $\alpha$ 
about a point $O$ is also a group. This fact follows from the
theorem~\mythetheorem{7.3} and from the first relationship in the
theorem~\mythetheorem{7.4}.
\mytheorem{8.1} Let $h$, $k$, $l$, and $q$ be four rays coming out 
from one point and lying on one plane $\alpha$. The equality 
$\theta_{hk}=\theta_{lq}$ takes place if and only if the bisectors of 
the angles $\angle lk$ and $\angle hq$ lie on one straight line.
\endproclaim
\demo{Proof} Let $\theta_{hk}=\theta_{lq}$. We denote by $m$ the
bisector of the angle formed by the rays $l$ and $k$. Then $z_m(l)=k$
and $z_m(k)=l$. Let's consider the mapping $f=z_m\compos\theta_{kh}
\compos z_m$ and calculate $f(l)$:
$$
f(l)=z_m(\theta_{kh}(z_m(l)))=z_m(\theta_{kh}(k))=z_m(h).
$$
On the other hand, applying the second and the third relationships
from the theorem~\mythetheorem{7.4}, for the mapping $f$ we get
$$
f=(z_m\compos\theta_{kh})\compos z_m=(\theta_{hk}\compos z_m)
\compos z_m=\theta_{hk}\compos (z_m\compos z_m)=\theta_{hk}.
$$
Therefore, $z_m(h)=f(l)=\theta_{hk}(l)=\theta_{lq}(l)=q$. This means
that the ray $q$ is produced  from the ray $l$ with the use of the 
mirror reflection in the line $m$. Hence, the bisector of the angle 
$\angle hq$ lies on the same line $m$ as the bisector of the angle 
$\angle lk$.
The necessity of the proposition stated in the theorem is proved.\par
     Let's prove its sufficiency. Assume that the bisectors of the 
angles $\angle lk$ and $\angle hq$ lie on one straight line. We denote
this line by $m$. Then $z_m(h)=q$ and $z_m(l)=k$. Let's consider the
mapping $f=z_m\compos\theta_{kh}\compos z_m=\theta_{hk}$ again and 
calculate $\theta_{hk}(l)$:
$$
\pagebreak
\theta_{hk}(l)=f(l)=z_m(\theta_{kh}(z_m(l)))=z_m(\theta_{kh}(k))
=z_m(h)=q.
$$
The relationship obtained $\theta_{hk}(l)=q$ and the 
theorem~\mythetheorem{7.3} yield the required result $\theta_{hk}
=\theta_{lq}$.\qed\enddemo
\mytheorem{8.2} The equality $\theta_{hk}=\theta_{lq}$ implies
$\theta_{hl}=\theta_{kq}$ and, conversely, $\theta_{hl}
=\theta_{kq}$ implies $\theta_{hk}=\theta_{lq}$.
\endproclaim
\mytheorem{8.3} Any two rotations of a plane $\alpha$ about 
a point $O\in\alpha$ commute: $\theta_{hk}\compos\theta_{lq}
=\theta_{lq}\compos\theta_{hk}$.
\endproclaim
    The theorem~\mythetheorem{8.2} is easily derived from the 
theorem~\mythetheorem{8.1}. It is analogous to the
theorem~\mythetheoremchapter{4.1}{3} in 
Chapter~\uppercase\expandafter{\romannumeral 3}. 
The theorem~\mythetheorem{8.3} is analogous to the
theorem~\mythetheoremchapter{4.3}{3} in
Chapter~\uppercase\expandafter{\romannumeral 3}. Ti means that
the group of pure rotations of a plane $\alpha$ about a point 
$O\in\alpha$ is commutative (Abelian) (see 
definition~\mythedefinitionchapter{4.3}{3} in
Chapter~\uppercase\expandafter{\romannumeral 3}). The following
theorem is derived from the theorem~\mythetheorem{8.1}.
\mytheorem{8.4} The equality $\theta_{hk}=\theta_{lq}$ implies 
the congruence of angles $\angle hk\cong\angle lq$.
\endproclaim
\myexercise{8.1} Using the theorem~\mythetheorem{8.1}, prove the
theorem~\mythetheorem{8.2} and the theorem~\mythetheorem{8.4}.
\endproclaim
\myexercise{8.2} Using the analogy of the rotations $\theta_{hk}$ 
on a plane and the congruent translations $p_{\sssize AB}$ on a
straight line, prove the theorem~\mythetheorem{8.3}.
\endproclaim
\head
\SectionNum{9}{113} Rotation of the space about a straight line.
\endhead
\rightheadtext{\S\,9. Rotation of the space \dots}
\parshape 12 0cm 10.1cm 0cm 10.1cm 0cm 10.1cm 3.5cm 6.6cm 3.5cm 6.6cm
3.5cm 6.6cm 3.5cm 6.6cm 3.5cm 6.6cm 3.5cm 6.6cm 3.5cm 6.6cm 3.5cm 6.6cm
0cm 10.1cm
    Let $\theta_{hk}\!:\alpha\to\alpha$ be the rotation of a plane
$\alpha$ by the angle $\angle hk$ about a point $O$. The plane
$\alpha$ divides the space $\Bbb E$ into two half-spaces $\alpha_{+}$
and $\alpha_{-}$. Let's set $\beta=\alpha$, $\alpha_{+}=\beta_{+}$,
$\alpha_{-}=\beta_{-}$. \vadjust{\vskip 5pt\hbox to 0pt{\kern -5pt
\includegraphics{Oris45.eps}\hss}\vskip -5pt}Applying the
theorem~\mythetheorem{5.5}, we get the mapping $\theta_{hk}\!:\Bbb E\to
\Bbb E$. Let's draw the straight line $c$ passing through the point $O$ 
and being perpendicular to the plane $\alpha$. This line is called the 
{\it axis of the rotation}, while the mapping $\theta_{hk}\!:\Bbb E\to
\Bbb E$ itself is called the {\it rotation by the angle $\angle hk$ 
about the axis $c$}.\par
     The ray $h$ and the axis of the rotation $c$ have the common point
$O$. Let's draw the plane $\gamma$ through the ray $h$ and the line
$c$. Similarly, we draw the plane $\delta$ passing through the ray $k$ 
and the line $c$. Both planes $\gamma$ and $\delta$ are perpendicular
to the plane $\alpha$. They intersect along the line $c$ and define
a dihedral angle with the edge $c$, for which the angle $\angle hk$ 
lying on the plane $\alpha$ is a plane angle.\par
     Each dihedral angle defines a rotation of the space 
about an axis.  Indeed, assume that a dihedral angle with the edge $c$
is given. It is the intersection of two closed half-spaces:
$\overline{\gamma_{+}}\cap\overline{\delta_{-}}$. Denote $c=d$ and
denote by $c_{+}\subset\gamma$ and $d_{+}\subset\delta$ the sides of 
this angle. Then for the mapping $f\!:c\to d$ we take the identical
mapping $\id\!:c\to c$. Applying the theorem~\mythetheorem{5.4}, we
extend it up to a mapping $f\!:\gamma\to\delta$ taking the half-plane
$c_{+}$ to the half-plane $d_{+}$ and taking the half-plane $c_{-}$ to
the half-plane $d_{-}$. Then we apply the theorem~\mythetheorem{5.5}
and extend this mapping up to a mapping $f\!:\Bbb E\to\Bbb E$ taking
$\gamma_{+}$ to $\delta_{+}$ and taking $\gamma_{-}$ to $\delta_{-}$.
\myexercise{9.1} Show that the above construction of the mapping 
$f\!:\Bbb E\to\Bbb E$ on the base of a dihedral angle yields the same
result as the previous construction extending the rotation of a plane
about a point up to a rotation of the space about an axis.
\endproclaim
     The rotations of the space about a fixed axis inherit all of the
properties of the rotations of a plane about a fixed point.
\mytheorem{9.1} Let $\theta_{hk}\!:\Bbb E\to\Bbb E$ be a rotation of
the space about some axis $c$ different from the identical mapping
($h\neq k$). Then the set of stable point of the mapping $\theta_{hk}$
coincides with the rotation axis $c$.
\endproclaim
\mytheorem{9.2} Let $h$ and $k$ be two half-planes with common 
boundary $c$. There are exactly two mappings of congruent translations
$f\!:\Bbb E\to\Bbb E$ which preserves stable the points of the line $c$
and takes the half-plane $h$ to the half-plane $k$. The first of them
$f=\theta_{hk}$ is the rotation about the axis $c$, while the second
one $f=z_\gamma$ is the mirror reflection in the plane $\gamma$ that
contains the bisector of the dihedral angle formed by the half-planes
$h$ and $k$.
\endproclaim
{\bf A remark}. The {\it bisector of a dihedral angle\/} is a half-plane
bounded by the edge of this angle and containing the bisector of any
one of its plane angles. The bisector divides a dihedral angle into two
dihedral angles whose plane angles are congruent.\par
     Like in theorem~\mythetheorem{7.2}, in theorem~\mythetheorem{9.2}
there are special cases where the half-planes $h$ and $k$ lie on one
plane. If the half-planes $h$ and  $k$ do coincide, the rotation 
$\theta_{hk}$ is the identical mapping: $\theta_{hk}=\id$. If $h$ and
$k$ are two complementary half-planes on one plane, then the rotation
$\theta_{hk}$ coincides with the mirror reflection in the line $c$
separating these half-planes: $\theta_{hk}=z_c$.
\mytheorem{9.3} The rotations of the space about a fixed axis and the
mirror reflections in planes passing through this axis possess the 
following properties:
$$
\xalignat 2
&\theta_{kl}\compos\theta_{hk}=\theta_{hl},
&&\theta_{hk}\compos z_\alpha=z_\alpha\compos \theta_{kh},\\
&z_\alpha\compos z_\alpha=\id,&&z_\alpha\compos z_\beta=\theta_{hk},
\text{ \ where \ }h\subset\beta,\,k=z_\alpha(h).
\endxalignat
$$
\endproclaim
\mytheorem{9.4} The equality $\theta_{hk}=\theta_{lq}$ implies
$\theta_{hl}=\theta_{kq}$ and, conversely, the equality $\theta_{hl}
=\theta_{kq}$ implies $\theta_{hk}=\theta_{lq}$.
\endproclaim
\mytheorem{9.5} Any two rotations about the same axis commute: 
$\theta_{hk}\compos\theta_{lq}=\theta_{lq}\compos\theta_{hk}$.
\endproclaim
\myexercise{9.2} Compare the theorems~\mythetheorem{9.1}, 
\mythetheorem{9.2}, \mythetheorem{9.3}, \mythetheorem{9.4},
and \mythetheorem{9.5} with the corresponding theorems for 
the rotations of a plane. Suggest your scheme for proving them.
\endproclaim
\mytheorem{9.6} Let $\theta_{hk}$ and $\theta_{lq}$ be two rotations
whose axes $c_1$ and $c_2$ do not coincide but intersect at some point
$O$. Then the composition $\theta_{hk}\compos\theta_{lq}$ is some
rotation $\theta_{rp}$ about a third axis $c_3$ passing through the
point $O$.
\endproclaim
\demo{Proof} Let's draw the plane $\beta$ passing through the 
intersecting lines $c_1$ and $c_2$. The line $c_1$ divides this 
plane into two half-planes. Let's denote one of these half-planes 
through $h$ and denote $k=\theta_{hk}(h)$. The half-planes $h$ and 
$k$ define a dihedral angle with the edge $c_1$. Let's denote
by $\alpha$ the plane containing its bisector. Then the mirror 
reflection $z_\alpha$ takes the half-plane $h$ to the half-plane
$k$. Therefore we can apply the fourth relationship from the
theorem~\mythetheorem{9.3} in order to expand the rotation 
$\theta_{hk}$ into the composition of two mirror reflections: 
$\theta_{hk}=z_\alpha\compos z_\beta$. We use this expansion 
in the following calculations:
$$
\theta_{hk}\compos\theta_{lq}=z_\alpha\compos z_\beta\compos
\theta_{lq}=z_\alpha\compos (z_\beta\compos\theta_{lq})=
z_\alpha\compos(\theta_{ql}\compos z_\beta).
$$
By construction the plane $\beta$ contains the axis $c_2$ of the 
second rotation, therefore, we used the second relationship
from the theorem~\mythetheorem{9.3} in the form $z_\beta\compos
\theta_{lq}=\theta_{ql}\compos z_\beta$. Now, applying the above
trick to the rotation $\theta_{ql}$, we get the expansion
$\theta_{ql}=z_\gamma\compos z_\beta$. For the initial composition
$\theta_{hk}\compos\theta_{lq}$ it yields
$$
\theta_{hk}\compos\theta_{lq}=z_\alpha\compos(z_\gamma\compos
z_\beta\compos z_\beta)=z_\alpha\compos z_\gamma\compos(z_\beta
\compos z_\beta)=z_\alpha\compos z_\gamma.
$$
The plane $\alpha$ contains the line $c_1$, but it does not contain
$c_2$. Similarly, the plane $\gamma$ contains $c_2$, but it does not
contain the line $c_1$. Therefore, these planes do not coincide, but
they have a common point $O$. Let's denote by $c_3$ the line arising
as their intersection. It is clear that $O\in c_3$. Now due to the 
fourth relationship from the theorem~\mythetheorem{9.3} for the
composition $z_\alpha\compos z_\beta$ we derive $z_\alpha\compos
z_\beta=\theta_{rp}$, where $\theta_{rp}$ is some rotation
about the axis $c_3$. Hence, $\theta_{hk}\compos\theta_{lq}
=\theta_{rp}$. The theorem is proved.
\qed\enddemo
     Let's combine the first relationship from the 
theorem~\mythetheorem{9.3} and the theorem~\mythetheorem{9.6}.
As a result we get the following result.
\mytheorem{9.7} The composition of two rotations of the space
whose axes have a common point $O$ is a rotation about an axis
passing through the point $O$.
\endproclaim
\head
\SectionNum{10}{116} The theorem on the decomposition of rotations.
\endhead
\rightheadtext{\S\,10. Decomposition of rotations.}
     Above we considered several types of the mappings of congruent
translation of the space. \pagebreak They are rotations about axes, 
mirror reflections in planes or in a lines, and inversions with respect 
to various points. Their common property is that they have stable points.
\mydefinition{10.1} A mapping of congruent translation of the space 
$f\!:\Bbb E\to\Bbb E$ having a stable point $O$, i\.\,e\. such that
$f(O)=O$, is called a {\it generalized rotation of the space about 
the point $O$}.
\enddefinition
     Let $f$ and $g$ be two mappings of congruent translation of the
space. Then their composition $f\compos g$ is obviously a mapping of
congruent translation of the space. If $f$ and $g$ are generalized
rotations about a point $O$, then $f\compos g$ is also a generalized
rotation about this point.
\mytheorem{10.1} Any generalized rotation $f\!:\Bbb E\to\Bbb E$ 
of the space about a point $O$ is either a rotation about some axis
passing through the point $O$, i\.\,e\. $f=\theta_{hk}$, or the
composition of such a rotations and a mirror reflection in some plane
containing the point $O$, i\.\,e\. $f=z_\alpha\compos\theta_{hk}$.
\endproclaim
     The term {\it rotation\/} in this theorem is treated so that 
the identical mapping and any reflection in a line are assumed to be
rotations (see the remark to the theorem~\mythetheorem{9.2}). In order 
to prove the theorem~\mythetheorem{10.1} we need two auxiliary  lemmas.
\mylemma{10.1} Let $M$ be some point lying on the plane of a triangle 
$ABC$. If for a point $X$ in the space the relationships $[XA]\cong [MA]$,
$[XB]\cong [MB]$, and $[XC]\cong [MC]$ are fulfilled, then the point
$X$ coincides with $M$.
\endproclaim
\demo{Proof} Assume that $X\neq M$. Then from the relationships
$[XA]\cong [MA]$, $[XB]\cong [MB]$, $[XC]\cong [MC]$ due to the
theorem~\mythetheorem{2.1} we derive that the  points $A$, $B$, and
$C$ lie on the plane of perpendicular bisectors of the segment $[MX]$.
Let's denote this plane through $\alpha$. Then due to the
theorem~\mythetheorem{6.1} we have $X=z_\alpha(M)$, where $z_\alpha$ 
is the mirror reflection in the plane $\alpha$. But, according to the
premise of the theorem, the point $M$ lies on the plane $\alpha$.
Therefore, $z_\alpha(M)=M$. Hence, $X=M$ despite the assumption that 
$X\neq M$. Thus, the lemma is proved.
\qed\enddemo
\mylemma{10.2} Let $M$ be some point not lying on the plane of a 
triangle $ABC$. In the space there is exactly one point $X$ 
different from $M$ and such that the relationships $[XA]\cong [MA]$, 
$[XB]\cong [MB]$, $[XC]\cong [MC]$ are fulfilled. It is the mirror
image of the point $M$ in the plane of the triangle $ABC$.
\endproclaim
\demo{Proof} Let's denote through $\alpha$ the plane of the triangle
$ABC$ and let $X=z_\alpha(M)$. The relationships $[XA]\cong [MA]$,
$[XB]\cong [MB]$, and $[XC]\cong [MC]$ for $X$ follow from the fact 
that the mirror reflection $z_\alpha$ is a congruent translation of 
the space. The existence of the required point $X$ is proved.\par
     Let's prove the uniqueness of the point $X\neq M$ satisfying the
relationships $[XA]\cong [MA]$, $[XB]\cong [MB]$, and $[XC]\cong [MC]$.
From these relationships and from the theorem~\mythetheorem{2.1}
we derive that the points $A$, $B$, and $C$ lie on the plane of
perpendicular bisectors for the segment $[MX]$. Hence, due to the 
theorem~\mythetheorem{6.1} we get $X=z_\alpha(M)$. This formula fixes
the point $X$ uniquely.
\qed\enddemo
\demo{Proof of the theorem~\mythetheorem{10.1}} Assume that the mapping
$f$ being a generalized rotation about the point $O$ is different from
the identical mapping. Let's show that either this mapping or its
composition with the mirror reflection in some plane has a stable point 
$Z$ different from $O$.\par
     Let's study various pairs of points $X$ and $f(X)$, where 
$X\neq O$. If $X=f(X)$, then the required stable point is found. 
If $X\neq f(X)$, we construct the plane of perpendicular bisectors 
for the segment $[Xf(X)]$. From $[OX]\cong [f(O)f(X)]$, since $O$ is
a stable point, we derive $[OX]\cong [Of(X)]$. Then, according to the
theorem~\mythetheorem{2.1} the point $O$ lies on the plane of 
perpendicular bisectors for the segment $[Xf(X)]$. In other words,
all planes of perpendicular bisectors for the segments of the
form $[Xf(X)]$ have the common point $O$. There is a case where all
such planes do coincide. In this case we denote by $\alpha$ the common
plane of perpendicular bisectors for all segments of the form $[Xf(X)]$. 
Then $f$ is the mirror reflection in the plane $\alpha$. Thus we have
$f=z_\alpha\compos\id$, where the identical mapping is treated as a
special case of the rotation about an axis.\par
     Now let's consider the case where there are two points $X$ and $Y$
for which the plane of perpendicular bisectors of the segments $[Xf(X)]$ 
and $[Yf(Y)]$ do not coincide. These planes have the common point $O$, 
hence, they intersect along some straight line $a$. Let's choose some 
point $Z\neq O$ on this line. Then
$$
\xalignat 3
&[ZX]\cong [Zf(X)],&&[ZY]\cong [Zf(Y)],&&[ZO]\cong [Zf(O)].
\endxalignat
$$
The first two relationship follow from the theorem~\mythetheorem{2.1},
while the last one is the trivial consequence of $f(O)=O$. From the fact
that $f$ is a congruent translation of the space we derive
$$
\xalignat 3
&[ZX]\cong [f(Z)f(X)],&&[ZY]\cong [f(Z)f(Y)],&&[ZO]\cong [f(Z)f(O)].
\endxalignat
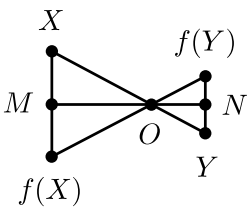
$$\par
\parshape 15 4cm 6cm 4cm 6cm 4cm 6cm 4cm 6cm 4cm 6cm 
4cm 6cm 4cm 6cm 4cm 6cm 4cm 6cm 4cm 6cm 4cm 6cm 4cm
6cm 4cm 6cm 4cm 6cm 0cm 10cm
\noindent Let's denote $\tilde Z=f(Z)$ and, comparing the above
\vadjust{\vskip 5pt\hbox to 0pt{\kern -5pt
\includegraphics{Oris46.eps}\hss}\vskip -5pt}sets 
of relationships, we find
$$
\aligned
[Zf(X)]&\cong [\tilde Zf(X)],\\
[Zf(Y)]&\cong [\tilde Zf(Y)],\\
[Zf(O)]&\cong [\tilde Zf(O)].
\endaligned
\quad
\mytag{10.1}
$$
In order to apply one of the lemmas~\mythelemma{10.1} or \mythelemma{10.2}, 
let's prove that the points $X$, $Y$, and $O$ do not lie on one straight
line. If we assume that these points lie on one straight line, then the
points $f(X)$, $f(Y)$, and $f(O)$ also lie on one straight line (see
theorem~\mythetheorem{5.1}) and we have the situation shown on one of
the figures~10.1. All of the five points $X$, $Y$, $f(X)$, $f(Y)$, and
$O$ lie on one plane. From $[OX]\cong [Of(X)]$ and $[OY]\cong [Of(Y)]$
we conclude that the triangles $XOf(X)$ and $YOf(Y)$ are isosceles. Their
medians $[OM]$ and $[ON]$ lie on one straight line because they are 
bisectors of one angle or bisectors of of two vertical angles. They are
the height of the corresponding triangles at the same time. For this 
reason the line $MN$ is the common perpendicular bisector of the segments
$[Xf(X)]$ and $[Yf(Y)]$ that lie on one plane. Through the point $O$ we
draw the perpendicular to the plane of the triangles $XOf(X)$ and $YOf(Y)$.
Then through this perpendicular and through the line $MN$ one can draw
the plane $\beta$ which (due to the theorem~\mythetheorem{3.1}) is the
common plane of perpendicular bisectors for the segments $[Xf(X)]$ and
$[Yf(Y)]$. But we consider the case where the planes of perpendicular
bisectors for the segments $[Xf(X)]$ and $[Yf(Y)]$ do not coincide. Due
to this contradiction, our initial assumption that the points $X$, $Y$,
and $O$ lie on one straight line is invalid.\par
     Now, having proved that the points $X$, $Y$, and $O$ do not lie
on one line, we return to the relationships \mythetag{10.1}. Due to the
theorem~\mythetheorem{5.1} the points $f(X)$, $f(Y)$, and $f(O)$ also 
do not lie on one straight line. Let's denote by $\beta$ the plane of
the triangle $f(X)f(Y)f(O)$. If $Z\in\beta$, then relying upon the
relationships \mythetag{10.1}, we apply the lemma~\mythelemma{10.1}. 
It yields $Z=\tilde Z=f(Z)$, i\.\,e\. $Z\neq O$ is a required stable
point of the mapping $f$. If $Z\notin\beta$, then the 
lemma~\mythelemma{10.2} is applicable. In this case $Z=\tilde Z$ or
the points $Z$ and $\tilde Z$ are mirror symmetric with respect to
the plane $\beta$. If $Z=\tilde Z$, then $Z$ again is a stable point
for $f$. Otherwise, if $Z\neq\tilde Z$, then $Z$ is a stable point
for the mapping $z_\beta\compos f$, which is also a generalized
rotation about the point $O$.\par
     Let's denote $g=f$ in the first case and denote $g=z_\beta
\compos f$ in the second case. From $g(O)=O$ and $g(Z)=Z$ due to the 
theorem~\mythetheoremchapter{2.1}{3} from
Chapter~\uppercase\expandafter{\romannumeral 3} we conclude that all
points of the line $OZ$ are stable under the action of the mapping
$g$. Let $h$ be some arbitrary half-plane having the line $OZ$ as its
boundary. Let's denote $k=g(h)$ and apply the theorem~\mythetheorem{9.2}.
Due to this theorem $g$ is the rotation about the axis $OZ$ taking the 
half-plane $h$ to the half-plane $k$, or $g$ is the mirror reflection
in the plane $\delta$ containing the bisector of the dihedral angle
formed by the half-planes $h$ and $k$. Hence, for the initial mapping 
$f$ we have the following four possible expansions:
$$
\xalignat 4
&f=\theta_{hk},&&f=z_\beta\compos\theta_{hk},
&&f=z_\delta,&&f=z_\beta\compos z_\delta.
\endxalignat
$$
If the first or the second case takes places, the proof is over. In 
the third case we can write $f=z_\delta\compos\id$, therefore, in this 
case the proof is also completed. The rest is the fourth case. The 
planes $\beta$ and $\delta$ have the common point $O$. Hence, either 
$\beta=\delta$ or these planes intersect along some line $b$. If 
$\delta=\beta$, then $f=z_\beta\compos z_\beta=\id$. Otherwise, if 
$\delta\neq\beta$, we apply the fourth relationship from the 
theorem~\mythetheorem{9.3}. For the mapping $f$ it yields $f=z_\beta
\compos z_\beta=\theta_{lq}$, where $\theta_{lq}$ is the rotation
about the line $b=\beta\cap\delta$. Thus, the theorem~\mythetheorem{10.1}
is proved.
\qed\enddemo
\head
\SectionNum{11}{121} The total rotation group and the group of pure
rotations of the space.
\endhead
\rightheadtext{\S\,11. The total rotation group and \dots}
    The theorem~\mythetheorem{10.1} proved in previous section has a very
important consequence. It means that any mapping of generalized rotation
about a point is bijective. Indeed, a rotation about an axis and a mirror
reflection both are bijective mappings, while the composition of two
bijective mappings is also a bijective mapping. From this fact we derive
that the set of all generalized rotations of the space about a fixed 
point $O$ is a group with respect to the composition. This group is 
called the {\it total rotation group of the space\/} about a fixed 
point $O$.\par
     According to the theorem~\mythetheorem{9.6} the set of rotations 
about various axes passing through a fixed point $O$ is also a group 
with respect to the composition. This group is called the {\it group
of pure rotations of the space\/} about a fixed point $O$.\par
     The theorem~\mythetheorem{10.1} determines the division of 
generalized rotations into {\it even\/} and {\it odd\/} ones. Pure
rotations belong to even rotations, while mirror reflections and their
compositions with pure rotations are odd rotations. The same generalized
rotation cannot be even and odd simultaneously. Indeed, $f=\theta_{hk}$ 
and $f=z_\beta\compos\theta_{lq}$ would imply $z_\beta=\theta_{hk}\compos
\theta_{ql}=\theta_{rp}$, which is impossible.\par
\myexercise{11.1} Show that the composition of two even rotations 
and the composition of two odd rotations both are even rotations, 
while the composition of an even rotation and an odd rotation is 
an odd rotation.
\endproclaim
\head
\SectionNum{12}{122} Orthogonal projection onto a straight line.
\endhead
\rightheadtext{\S\,12. Orthogonal projection \dots}
      Let $a$ be some straight line in the space. Choosing some 
arbitrary point $X$ not lying on the line $a$, we drop the perpendicular
$[XY]$ from the point $X$ onto the line $a$. The foot of this
perpendicular (the point $Y$) is called the {\it orthogonal projection\/}
of the point $X$ onto the line $a$. According to the
theorem~\mythetheoremchapter{6.5}{3} from
Chapter~\uppercase\expandafter{\romannumeral 3}, once a point $X$ is
given, its projection $Y$ is fixed uniquely. Hence, we can define a 
mapping $\pi_a\!:\Bbb E\to a$. For a point $Y\in a$ we set $\pi_a(Y)=Y$.
The mapping $\pi_a$ is called the {\it orthogonal projection\/} onto
the line $a$.\par
     The orthogonal projection is not a congruent translation. Moreover,
two different points $X_1\neq X_2$ can be taken to one point 
$\pi_a(X_1)=\pi_a(X_2)$ under this mapping. Let $\alpha$ be some plane
perpendicular to the line $a$. According to the 
definition~\mythedefinition{1.1}, such a plane intersects the line 
$a$ at some point $A$. Comparing this definition with the construction 
of the mapping $\pi_a$, we see that all points of the plane $\alpha$
are taken to the point $A$ by the projection $\pi_a$.\par
     Let $\alpha$ and $\beta$ be two different planes perpendicular 
to the line $a$ and intersecting the line $a$ at the points $A$ and 
$B$. Such planes have no common points. Indeed, the existence of a point 
$M\in\alpha\cap\beta$ would mean that $\pi_a(M)=A$ and $\pi_a(M)=B$. But
the result of orthogonal projection of the point $M$ onto the line $a$ 
is defined uniquely. Therefore, $\alpha\cap\beta=\varnothing$.
\mytheorem{12.1} Let $\pi_a$ be the orthogonal projection onto a
straight line $a$. If the projections of some two points $A$ and 
$B$ do coincide, then the whole line $AB$ is projected onto one point 
$C=\pi_a(A)=\pi_a(B)$ of the line $a$.
\endproclaim
\demo{Proof} The relationship $C=\pi_a(A)$ means that $C=A$ or
$C$ is the foot of the perpendicular dropped from the point $A$ onto
the line $a$. In both cases the point $A$ lies on the plane $\alpha$,
passing through the point $C$ and being perpendicular to the line $a$ 
(see theorems~\mythetheorem{1.2} and \mythetheorem{1.4}). Similar 
considerations yield $B\in\alpha$. Hence the whole line $AB$ lies
on the plane $\alpha$ which is projected onto the single point 
$C\in a$.
\qed\enddemo
\mytheorem{12.2} Let $\pi_a$ be the orthogonal projection onto a
straight line $a$ and assume that $b$ is some straight line which
is not projected onto a single point of the line $a$. Then for the 
points of the line $b$ the following propositions are valid:
\roster
\item $A\neq B$ implies $\pi_a(A)\neq\pi_a(B)$;
\item $(A\blacktriangleright B\blacktriangleleft C)$ implies
      $(\pi_a(A)\blacktriangleright\pi_a(B)\blacktriangleleft\pi_a(C))$.
\endroster
\endproclaim
\demo{Proof} The first item of this theorem is an immediate consequence
of the previous theorem~\mythetheorem{12.1}.\par
\parshape 16 0cm 10cm 0cm 10cm 
5.2cm 4.8cm 5.2cm 4.8cm 5.2cm 4.8cm 5.2cm 4.8cm 5.2cm 4.8cm 5.2cm 4.8cm
5.2cm 4.8cm 5.2cm 4.8cm 5.2cm 4.8cm 5.2cm 4.8cm 5.2cm 4.8cm 5.2cm 4.8cm
5.2cm 4.8cm 0cm 10cm 
     Let's consider the second item of the theorem. Assume 
that $(A\blacktriangleright B\blacktriangleleft C)$. Let's use the
\vadjust{\vskip 5pt\hbox to 0pt{\kern -5pt
\includegraphics{Oris47.eps}\hss}\vskip -5pt}following 
notations for the projections of the points $A$, $B$, $C$:
$$
\align
&\tilde A=\pi_a(A),\\
&\tilde B=\pi_a(B),\\
&\tilde C=\pi_a(C).
\endalign
$$
Let's draw the plane $\beta$ passing through the point $\tilde B\in a$
and being perpendicular to the line $a$. Such a plane does exist and it
is unique (see theorem~\mythetheorem{1.2}). The points $B$ and $\tilde B$
lie on this plane. From $\tilde A\neq\tilde B$ and from $\tilde C\neq
\tilde B$ we conclude that the four points $A$, $C$, $\tilde A$, and
$\tilde C$ do not lie on the plane $\beta$. The plane $\beta$ divides 
the set of points not lying on this plane into two open half-spaces
$\beta_{+}$ and $\beta_{-}$. The half-spaces $\beta_{+}$ and $\beta_{-}$
arise as the equivalence classes, where two points $X$ and $Y$ are treated
to be equivalent if $X=Y$ or if the segment $[XY]$ does not intersect the
plane $\beta$ (see \S\,6 in 
Chapter~\uppercase\expandafter{\romannumeral 2}).\par
     In our case the points $A$ and $\tilde A$ lie on the plane $\alpha$
passing through the point $\tilde A\in a$ and being perpendicular to the
line $a$. The plane $\alpha$ does not intersect the plane $\beta$. 
Therefore either $A=\tilde A$ or, if $A\neq\tilde A$, the segment
$[A\tilde A]$ lies on the plane $\alpha$ and does not intersect the
plane $\beta$. Hence, $A\sim\tilde A$. In a similar way we get $C\sim
\tilde C$. But $(A\blacktriangleright B\blacktriangleleft C)$ implies
that the segment $[AC]$ intersects the plane $\beta$ at the point $B$. 
Hence the points $A$ and $C$ are not equivalent. Then from $A\sim\tilde A$
and $C\sim\tilde C$ we get that the points $\tilde A$ and $\tilde C$ 
are also not equivalent. Therefore the segment $[\tilde A\tilde C]$ 
intersects the plane $\beta$ at the point $\tilde B$. Hence, we
immediately derive the required relationship 
$(\tilde A\blacktriangleright\tilde B\blacktriangleleft\tilde C)$.
\qed\enddemo
     As a corollary of the theorem~\mythetheorem{12.2} we get that
if the mapping $\pi_a$ does not project a line $b$ onto one point,
then any ray lying on the line $b$ is projected onto a ray ands any 
segment of the line $b$ is projected onto a segment.
\head
\SectionNum{13}{124} Orthogonal projection onto a plane.
\endhead
\rightheadtext{\S\,13. Orthogonal projection onto a plane.}
     Let $\alpha$ be some plane and let $X$ be a point outside
this plane. Let's drop the perpendicular from the point $X$  onto
the plane $\alpha$ and denote by $\pi_\alpha(X)$ the foot of this
perpendicular. For the points $Y$ lying on the plane $\alpha$
we set $\pi_\alpha(Y)=Y$. Due to the theorem~\mythetheorem{3.4}
this construction defines a mapping $\pi_\alpha\!:\Bbb E\to\alpha$
which is called the {\it orthogonal projection\/} onto the plane
$\alpha$.\par
     The orthogonal projection onto a plane $\pi_\alpha$ is not a
congruent translation. It can take some different points to one point.
Let $Y$ be a point on the plane $\alpha$ and let $a$ be the perpendicular
to the plane $\alpha$ passing through the point $Y$ (see
theorem~\mythetheorem{1.3}). Then all points of the line $a$ and only
these points are projected onto the point $Y\in\alpha$.
\mytheorem{13.1} Let $\pi_\alpha$ be the orthogonal projection on some
plane $\alpha$. If for some two points $A$ and $B$ their projections
do coincide, then the whole straight line $AB$ is projected onto one 
point $C=\pi_a(A)=\pi_a(B)$ on the plane $\alpha$. In this case the 
line $AB$ coincides with the perpendicular to the plane $\alpha$ passing
through the point $C$.
\endproclaim
\demo{Proof} The relationship $C=\pi_\alpha(A)$ means that $C=A$ or
$C$ is the foot of the perpendicular dropped from the point $A$ onto
the plane $\alpha$. In any of these two cases the point $A$ lies on the
perpendicular $a$ to the plane $\alpha$ passing through the point $C$
(see theorem~\mythetheorem{1.3}). In a similar way we prove that $B\in a$.
Hence, the line $AB$ coincides with the line $a$ which is projected 
onto one point $C\in\alpha$.\qed\enddemo
\mytheorem{13.2} Let $\pi_\alpha$ be the orthogonal projection onto
a plane $\alpha$ and assume that $b$ is some straight line which is not
projected onto one point. Then the line $b$ is projected onto some line
$a$ lying on the plane $\alpha$ and for the points of $b$ the following
propositions are valid:
\roster
\item $A\neq B$ implies $\pi_a(A)\neq\pi_a(B)$;
\item $(A\blacktriangleright B\blacktriangleleft C)$ implies
      $(\pi_a(A)\blacktriangleright\pi_a(B)\blacktriangleleft
      \pi_a(C))$.
\endroster
\endproclaim
\demo{Proof} \parshape 11 4.2cm 5.8cm 4.2cm 5.8cm 4.2cm 5.8cm 
4.2cm 5.8cm 4.2cm 5.8cm 4.2cm 5.8cm 4.2cm 5.8cm 4.2cm 5.8cm 
4.2cm 5.8cm 4.2cm 5.8cm 
0cm 10cm 
    The line $b$ is not projected onto one point. Therefore, there 
are at least two points $A$ and $B$ on this line whose projections 
are different. Let $\tilde A=\pi_\alpha(A)$ and 
$\tilde B=\pi_\alpha(B)$ and denote by $a$ the line $\tilde A
\tilde B$. Then we draw the perpendiculars $m$ and $n$ to the plane $\alpha$
\vadjust{\vskip 5pt\hbox to 0pt{\kern -5pt
\includegraphics{Oris48.eps}\hss}\vskip -5pt}passing 
through the points $\tilde A$ and $\tilde B$ respectively.
According to the theorem~\mythetheorem{3.2} there is a plane $\beta$
containing both these perpendiculars. This plane is perpendicular to 
the plane $\alpha$. It intersects the plane alpha along the line $a$ 
since $\tilde A$ and $\tilde B$ belong to $\alpha\cap\beta$.\par
     The points $A$ and $B$ lie on the perpendiculars $m$ and $n$, hence,
they belong to the plane $\beta$. This yields $b\subset\beta$. Let $X$ 
be some arbitrary point of the plane $\beta$ not lying on  the line $a$. 
Within the plane $\beta$ we drop the perpendicular $[X\tilde X]$ from the
point $X$ onto the line $a$. From $\alpha\perp\beta$ due to the
definitions~\mythedefinition{3.1} and \mythedefinition{1.1} and due to
the theorem~\mythetheoremchapter{6.5}{3} from 
Chapter~\uppercase\expandafter{\romannumeral 3} we derive that the 
segment $[X\tilde X]$ is the perpendicular dropped from the point 
$X$ onto the plane $\alpha$. Hence, the whole plane $\beta$ is projected
onto the line $a$ and the restriction of the mapping $\pi_\alpha\!:\Bbb E
\to\alpha$ to the plane $\beta$ coincides with the restriction of the
mapping $\pi_a\!:\Bbb E\to a$ to the same plane. This observation 
reduces the theorem~\mythetheorem{13.2} to the theorem~\mythetheorem{12.2}
proved in previous section.
\qed\enddemo
     As an auxiliary result in proving the theorem~\mythetheorem{13.2}
we have shown that the plane $\beta$ perpendicular to the plane $\alpha$
is projected onto the line $a=\alpha\cap\beta$ under the projection
$\pi_\alpha$. This result can be strengthened. 
\mytheorem{13.3} A plane $\beta$ is projected onto a line $a\subset
\alpha$ by the mapping $\pi_\alpha\!:\Bbb E\to\alpha$ if an only if 
$\beta\perp\alpha$.
\endproclaim
    The following theorem is well-known. It is called the theorem
{\it on three perpendiculars}.
\mytheorem{13.4} A straight line $b$ intersecting a plane $\alpha$
at a point $O$ is perpendicular to a line $c$ lying on this plane
and passing through the point $O$ if and only if its projection
$a=\pi_\alpha(b)$ is perpendicular to $c$.
\endproclaim
\demo{Proof} If the line $b$ lies on the plane $\alpha$, then
$b=a$. In this case proposition of the theorem is trivial.\par
     Assume that $b$ intersects the plane $\alpha$ at a point $O$, 
but does not lie on this plane. Let's choose some point 
$B\neq O$ on the line $b$ and denote $A=\pi_\alpha(B)\in a$. The
line $AB$ is the perpendicular to the plane $\alpha$ passing through
the point $A$. We draw the other perpendicular to the plane $\alpha$
through the point $O$ and we denote it by $d$. According to the
theorem~\mythetheorem{3.2}, there is a plane $\beta$ containing both
perpendiculars to the plane $\alpha$. This plane is perpendicular to
$\alpha$. The points $O$ and $B$ lie on the plane $\beta$, hence
$b\subset\beta$. Similarly, from $O\in\beta$ and $A\in\beta$ we 
find that $a\subset\beta$.\par
      Since $d\perp\alpha$, the line $d$ is perpendicular to any line
lying on the plane $\alpha$ and passing through the point $O$. In
particular, $d\perp c$. Now if $b\perp c$, then from $c\perp b$ and
$c\perp d$ due to the theorem~\mythetheorem{1.1} we derive $c\perp\beta$.
Hence, $c\perp a$.\par
     Conversely, if $c\perp a$, then, complementing this condition
with the condition $c\perp d$, we again get the orthogonality $c\perp
\beta$. Hence, $c\perp b$. The proof is over.\qed\enddemo
\mytheorem{13.5} Let $\pi_\alpha$ be the orthogonal projection onto
a plane $\alpha$ and let $\beta$ be some plane which is not projected
onto a line under the projection $\pi_\alpha$. Then for the points
of the plane $\beta$ the following propositions are valid:
\roster
\rosteritemwd=3pt
\item $A\neq B$ implies $\pi_a(A)\neq\pi_a(B)$;
\item if the points $A$, $B$, and $C$ do not lie on one straight
      line, then their projections also do not lie on one straight
      line;
\item if the points $A$, $B$, and $C$ lie on one straight line, 
      then their projections also lie on one straight line and 
      $(A\blacktriangleright B\blacktriangleleft C)$ implies
      $(\pi_\alpha(A)\blacktriangleright\pi_\alpha(B)
      \blacktriangleleft\pi_\alpha(C))$.
\item if $b$ is some straight line on the plane $\beta$ and if the
      points $A$ and $C$ lie on different sides of the line $b$, then
      their projections $\pi_\alpha(A)$ and $\pi_\alpha(C)$ lie on
      different sides of the line $a=\pi_\alpha(b)$ on the plane
      $\alpha$.
\endroster
\endproclaim
\demo{Proof} If we assume that the projections of two distinct points 
$A$ and $B$ do coincide, then the line $AB$ connecting them is 
perpendicular to the plane $\alpha$. The plane $\beta$ comprising the
perpendicular to the plane $\alpha$ is perpendicular to $\alpha$ (see 
definition~\mythedefinition{3.1}). Then it is projected onto a line 
despite to the premise of the theorem. This contradiction proves the 
item \therosteritem{1}.\par
     Let $A$, $B$, and $C$ be three points of the plane $\beta$ not
lying on one straight line. We denote by $\tilde A$, $\tilde B$, and
$\tilde C$ their projections:
$$
\xalignat 3
&\tilde A=\pi_\alpha(A),
&&\tilde B=\pi_\alpha(B),
&&\tilde C=\pi_\alpha(C).
\endxalignat
$$
Die to the item \therosteritem{1}, which is already proved, the 
coincidences $\tilde A=\tilde B$, $\tilde B=\tilde C$, and
$\tilde A=\tilde C$ are impossible. Let's draw two perpendiculars 
to the plane $\alpha$ trough the points $\tilde A$ and $\tilde B$.
According to the theorem~\mythetheorem{3.2}, there is a plane
$\gamma$ containing both these perpendiculars. This plane $\gamma$
is perpendicular to $\alpha$. It contains the line $\tilde A
\tilde B$.\par
     Now if we assume that the points $\tilde A$, $\tilde B$,
$\tilde C$ lie on one straight line, then the point $\tilde C$
belong to the line $\tilde A\tilde B$ which lies on the plane
$\gamma$. From $\gamma\perp\alpha$ we derive that the perpendicular
to the plane $\alpha$ passing through the point $\tilde C$ lies
on the plane $\gamma$. Hence, all of the three points $A$, $B$, and
$C$ lie on the plane $\gamma$. This fact implies $\beta=\gamma$ 
and $\beta\perp\alpha$. But this is not possible since $\beta$ is not
projected onto a line according to the premise of the theorem. The
contradiction obtained proves the second item of the 
theorem~\mythetheorem{13.5}.\par
     Let's proceed to the third item. In this case the points $A$, $B$,
and $C$ lie on one straight line. We denote this line by $b$. Due to
the item \therosteritem{1} of the theorem, which is already proved, the
line $b$ is not projected onto one point. Therefore, the third item of the
theorem follows from the theorem~\mythetheorem{13.2}.\par
     Assume that the points $A$ and $C$ lie on the plane $\beta$ on
different sides of the line $b$. Then the segment $[AC]$ crosses the
line $b$ at some its interior point $B$, i\.\,e\. $(A\blacktriangleright B
\blacktriangleleft C)$. Passing to the projections $\tilde A$, $\tilde B$,
and $\tilde C$, due to the item~\therosteritem{3} we conclude that the
segment $[\tilde A\tilde C]$ lying on the plane $\alpha$ intersects the
line $a=\pi_\alpha(b)$ at some its interior point $\tilde B$. Thus, the
fourth item and the theorem~\mythetheorem{13.5} in whole are proved.
\qed\enddemo
      The theorem~\mythetheorem{13.5} has an important corollary. If the
planes $\alpha$ and $\beta$ are not perpendicular, then the mapping
$\pi_\alpha\!:\Bbb E\to\alpha$ takes each line $b$ lying on the plane
$\beta$ to some line $a$ on the plane $\alpha$ so that the half-planes
bounded by the line $b$ are projected onto the half-planes on
$\alpha$ bounded by the line $a$.
\myexercise{13.1} Prove the \pagebreak theorem~\mythetheorem{13.3}.
\endproclaim
\head
\SectionNum{14}{129} Translation by a vector along a straight line.
\endhead
\rightheadtext{\S\,14. Translation by a vector \dots}
     Let $a$ be some line in a plane $\alpha$. It divides the plane 
$\alpha$ into two half-planes $a_{+}$ and $a_{-}$. Let $\beta=\alpha$ 
and $b=a$. Assume also that $b_{+}=a_{+}$ and $b_{-}=a_{-}$, while for 
the mapping $f\!: a\to a$ we take $p_{\bold c}$, where $\bold c$ is some
slipping vector on the line $a$. The mapping $p_{\bold c}$, which is
called a congruent translation by a vector, was defined above in \S\,3 
of Chapter~\uppercase\expandafter{\romannumeral 3}. Applying the
theorem~\mythetheorem{5.4}, we can extend $p_{\bold c}$ up to the
mapping $p_{a\bold c}\!:\alpha\to\alpha$, which is called the 
{\it translation of the plane $\alpha$ by the vector $\bold c$ along
the line $a$}.\par
     The plane $\alpha$ divides the space into two half-spaces
$\alpha_{+}$ and $\alpha_{-}$. Let's denote $\beta=\alpha$, 
$\beta_{+}=\alpha_{+}$, and $\beta_{-}=\alpha_{-}$. For the mapping
$f\!:\alpha\to\alpha$ we choose $p_{a\bold c}$, then we apply
the theorem~\mythetheorem{5.5}. As a result we get the mapping
$p_{a\bold c}:\Bbb E\to\Bbb E$, which is called the {\it translation 
of the space by the vector $\bold c$ along the line $a$}.\par
      The plane $\alpha$ containing the line $a$ is an auxiliary object
in constructing the  mapping $p_{a\bold c}:\Bbb E\to\Bbb E$. This fact
is explained by the following theorem.
\mytheorem{14.1} The restriction of the mapping $p_{a\bold c}:\Bbb E
\to\Bbb E$ to any plane $\beta$ containing the line $a$ coincides with 
the translation of the plane $\beta$ by the vector $\bold c$ along the 
line $a$.
\endproclaim
\demo{Proof} \parshape 16 0cm 10cm 0cm 10cm 
5.8cm 4.2cm 5.8cm 4.2cm 5.8cm 4.2cm 5.8cm 4.2cm 5.8cm 4.2cm
5.8cm 4.2cm 5.8cm 4.2cm 5.8cm 4.2cm 5.8cm 4.2cm 5.8cm 4.2cm 
5.8cm 4.2cm 5.8cm 4.2cm 0cm 10cm 0cm 10cm
     Let's consider the mapping $p_{a\bold c}:\Bbb E\to\Bbb E$ 
constructed with the use of the auxiliary plane $\alpha$ containing the
line $a$. Let $\beta$ be some other plane containing the line $a$ and let 
\vadjust{\vskip 5pt\hbox to 0pt{\kern -5pt
\includegraphics{Oris49.eps}\hss}\vskip -5pt}$X$ be some 
point of the plane $\beta$ not lying on the plane $\alpha$. 
For the beginning we study the case, where the plane $\beta$ is not 
perpendicular to the plane $\alpha$. Let's denote by $\alpha_{+}$ the
half-space bounded by the plane $\alpha$, and containing the point $X$. 
By $b_{+}$ we denote the half-plane on the plane $\beta$ bounded by the
line $a$ and containing the point $X$. Then $b_{+}=\beta\cap\alpha_{+}$.
Let $Y=p_{a\bold c}(X)$. Following the general scheme of constructing the
extensions for the mappings of congruent translation (see \S\,5), in order
to fix the point $Y$ we drop the perpendicular from the point $X$ onto
the plane $\alpha$. As a result we get the point $\tilde X=\pi_\alpha(X)$
and the segment $[X\tilde X]$. We apply the mapping $p_{a\bold c}\!:\alpha
\to\alpha$ to the point $\tilde X$ and get the point $\tilde Y=p_{a\bold c}
(\tilde X)$. Let's draw the perpendicular to the plane $\alpha$ through
the point $\tilde Y$. On this perpendicular we mark a point $Y$ such that
$Y\in a_{+}$ and $[Y\tilde Y]\cong [X\tilde X]$. These conditions fix the
point $Y=p_{a\bold c}(X)$ uniquely.\par
\vskip 0pt plus 1pt minus 1pt
    Remember that the mapping $p_{a\bold c}\!:\alpha\to\alpha$ itself
is an extension of the mapping $p_{\bold c}\!:a\to a$ from the line
$a$ to the plane $\alpha$. Due to the theorem~\mythetheorem{13.5}
under the orthogonal projection $\pi_\alpha\!:\Bbb E\to\alpha$ the 
half-plane $b_{+}$ is mapped to some half-plane with the boundary $a$. 
We denote it $a_{+}$. From $\tilde X=\pi_\alpha(X)$ we get $\tilde X
\in a_{+}$. In order to fix the point $\tilde Y=p_{a\bold c}(\tilde X)$
we apply the procedure of extending the mapping  $p_{\bold c}\!:a\to a$ 
(see theorem~\mythetheorem{5.4}, its proof, and comments preceding
this theorem in \S\,5). On the plane $\alpha$ we drop the perpendicular 
from the point $\tilde X$ onto the line $a$. We denote by $A$ the foot 
of such perpendicular. Let $B=p_{\bold c}(A)$. Then the vector
$\overrightarrow{AB}$ is a geometric realization for the slipping 
vector $\bold c$ on the line $a$ (see \S\,3 and \S\,4 in
Chapter~\uppercase\expandafter{\romannumeral 3}). Having fixed the 
point $B$, on the plane $\alpha$ we draw the perpendicular to the
line $a$ through this point. Then on this perpendicular we choose 
a point $\tilde Y$ such that $\tilde Y\in a_{+}$ and 
$[B\tilde Y]\cong [A\tilde X]$. These two conditions fix the point
$\tilde Y=p_{a\bold c}(\tilde X)$ uniquely. By construction this point
appears to be the orthogonal projection of the point $Y$ onto the
plane $\alpha$, i\.\,e\. $\tilde Y=\pi_\alpha(Y)$.\par
\vskip 0pt plus 1pt minus 1pt
     Now let's consider the planes of the triangles $AX\tilde X$ and 
$BY\tilde Y$, which are congruent due to $[X\tilde X]\cong [Y\tilde Y]$,
$[A\tilde X]\cong [B\tilde Y]$ and since the angles $\angle A\tilde XX$ 
and $\angle B\tilde YY$ both are right angles. Let's denote the first of 
these two planes by $\gamma$ and the second one by $\delta$. The plane
$\gamma$ contains the line $X\tilde X$, which is perpendicular to the 
plane $\alpha$. Hence, $\gamma\perp\alpha$ (see 
definition~\mythedefinition{3.1}). Due to $\tilde X=\pi_\alpha(X)$ 
the line $A\tilde X$ is the projection of the line $AX$. Therefore we
can apply the theorem on three perpendiculars.(see
theorem~\mythetheorem{13.4} above). Due to this theorem since the lines
$A\tilde X$ and $a$ are perpendicular, we derive $a\perp AX$. Now from
$a\perp A\tilde X$ and $a\perp AX$ we get $a\perp\gamma$. Thus, the
half-planes $b_{+}$ and $a_{+}$ define a dihedral angle with the edge 
$a=b$, while the plane $\gamma$ cut the angle $\angle XA\tilde X$
being a plane angle for this dihedral angle.\par
     Let's repeat the above considerations with respect to the plane
$\delta$. This yields $a\perp\delta$, i\.\,e\. the plane $\delta$ 
cut the other plane angle of the dihedral angle formed by half-planes
$b_{+}$ and $a_{+}$. From the congruence of the triangles $AX\tilde X$ 
and $BY\tilde Y$ we derive the congruence of the angles $\angle YB
\tilde Y$ and $\angle XA\tilde X$. Now, applying the 
theorem~\mythetheorem{4.1} and taking into account the fact that the 
point $Y$ lies on the half-plane $\alpha_{+}$, we conclude that 
the angle $\angle YB\tilde Y$ coincides with the plane angle of
the dihedral angle formed by the half-planes $b_{+}$ and $a_{+}$. 
Hence, $Y\in b_{+}\subset\beta$. Thus, the segments $[AX]$ and $[BY]$ 
lie on the plane $\beta$ on one side of the line $b=a$. They both are
perpendicular to the line $a$ and congruent to each other since the
triangles $AX\tilde X$ and $BY\tilde Y$ are congruent. We write the
obtained results for $X$ and $Y$ as follows:
\roster
\item the points $X$ and $Y$ lie on the plane $\beta$ on one side of
      the line $a$;
\item the segments $[AX]\cong [BY]$ are perpendicular to $a$;
\item the vector $\overrightarrow{AB}$ is a geometric realization of 
      the slipping vector $\bold c$.
\endroster
The conditions listed above mean exactly that the point $Y$ is produced
from the point $X$ by applying the mapping $p_{a\bold c}\!:\beta\to
\beta$ that translate the plane $\beta$ by the vector $\bold c$ along
the line $a$.\par
     The case where $\beta\perp\alpha$ appears to be much simpler than
the case we have already considered. Here the points $X$ and $Y$ are
projected to the points $A$ and $B$ on the line $a$, i\.\,e\. 
$A=\pi_\alpha(X)$ and $B=\pi_\alpha(Y)$. By construction of the
mapping $p_{a\bold c}\!:\Bbb E\to\Bbb E$ this means that the conditions
\therosteritem{2} and \therosteritem{3} are fulfilled. The relationship
$Y\in\beta$ follows from $\beta\perp\alpha$ and $[YB]\perp\alpha$ due
to the definition~\mythedefinition{3.1}. The fact that the points $X$ 
and $Y$ lie on one side of the line $a$ follows from the fact that they 
are on one side of the plane $\alpha$ by construction of the mapping
$p_{a\bold c}\!:\Bbb E\to\Bbb E$. Hence the condition \therosteritem{1} 
in the case $\beta\perp\alpha$ is also fulfilled and the point $Y$ is
the result of applying the translation $p_{a\bold c}\!:\beta\to\beta$ 
to the point $X$. The theorem is proved.
\qed\enddemo
       The mappings of translation along some line $a$ inherit the
properties of the translations on this line considered in \S\,3 and
in \S\,4 of 
Chapter~\uppercase\expandafter{\romannumeral 3}. They are commutative
and the equality $\bold d=\bold c+\bold e$ for the vectors $\bold d$,
$\bold c$, and $\bold e$ implies
$$
p_{a\bold d}=p_{a\bold c}\compos p_{a\bold e}.
$$
The translation by the zero vector appears to be the identical mapping:
$p_{a\bold 0}=\id$. The translations along nonzero vectors constitute
a separate class of the mappings of congruent translation. For a vector 
$\bold c\neq\bold 0$ the mapping $p_{a\bold c}$ cannot be reduced to a
generalized rotation about some point. This fact is derived from the
following theorem.
\mytheorem{14.2} For $\bold c\neq\bold 0$ the mapping 
$p_{a\bold c}\!:\Bbb E\to\Bbb E$ has no stable points at all.
\endproclaim
     Let $\bold c$ be some nonzero slipping vector on the line 
$a$. Let's choose a point $A$ not lying on this line. Then let's
draw the plane $\alpha$ passing through the point $A$ and the line 
$a$. Denote $B=p_{a\bold c}(A)$. According to the
theorem~\mythetheorem{14.1}, the point $B$ belongs to the plane 
$\alpha$. Let's denote by $b$ the line $AB$ and consider the vector
$\overrightarrow{BA}$. We denote by $\bold d$ the slipping vector
on the line $b$ corresponding to the geometric vector 
$\overrightarrow{BA}$. It is obvious that $A=p_{b\bold d}(B)$,
because of which the point $A$ appears to be a stable point for the
composition $f=p_{b\bold d}\compos p_{a\bold c}$. According to
the theorem~\mythetheorem{10.1} the mapping $f$ is either a rotation
of the space about some axis or the composition of such a rotation with
the mirror reflection in some plane. As we shall see, in our case only 
the firs option is possible. The mapping 
$f=p_{b\bold d}\compos p_{a\bold c}$
is a rotation about the axis passing through the point $A$ and being
perpendicular to the plane $\alpha$. This proposition can be
strengthened a little bit.
\mytheorem{14.3} Let $p_{a\bold c}$ and $p_{b\bold d}$ be two mappings 
of translation by two vectors $\bold c$ and $\bold d$ along two
lines $a$ and $b$. If their composition $f=p_{b\bold d}\compos p_{a
\bold c}$ has a stable point $O$, then the lines $a$ and $b$ lie 
on some plane $\alpha$, while the mapping $f$ itself is a rotation
about an axis passing through the point $O$ and being perpendicular to
the plane $\alpha$.
\endproclaim
\mytheorem{14.4} Let $p_{a\bold m}$, $p_{b\bold n}$, and $p_{c\bold k}$ 
be three mappings of translation by the vectors $\bold m$, $\bold n$, 
and $\bold k$ along the lines $a$, $b$, and $c$. If their composition
$f=p_{c\bold k}\compos p_{b\bold n}\compos p_{a\bold m}$ has a stable 
point $O$, then this composition is a rotation about some axis passing
through the point $O$.
\endproclaim
     We do not prove the theorems~\mythetheorem{14.3} and 
\mythetheorem{14.4} here. Their proofs are left to the reader as 
exercises. The matter is that upon formulating the axiom~\mytheaxiom{A20}
the statements of these theorems simplify substantially. The proofs of the
corresponding simplified theorems based on the axiom~\mytheaxiom{A20}
are given in Chapter~\uppercase\expandafter{\romannumeral 6}
\myexercise{14.1} Prove the theorem~\mythetheorem{14.2} on the base of the
theorem~\mythetheorem{14.1}.
\endproclaim
\myexercise{14.2} Prove the theorems~\mythetheorem{14.3} and
\mythetheorem{14.4}
\endproclaim
      The theorems~\mythetheorem{14.3} and \mythetheorem{14.4} differ 
only in the number of translations in the composition. As appears, this
number could be arbitrary. Using the theorems~\mythetheorem{9.6}, 
\mythetheorem{14.2}, and \mythetheorem{14.4}, by induction one can prove
the following fact.
\mytheorem{14.5} Let $f=f_1\compos\ldots\compos f_n$ be the composition 
of $n$ translations by $n$ vectors along $n$ lines. If such a composition
has a stable point $O$, then it is a rotation about some axis passing 
through the point $O$.
\endproclaim
\head
\SectionNum{15}{134} Motions and congruence of complicated 
geometric forms.
\endhead
\rightheadtext{\S\,15. Motions and congruence \dots}
\mytheorem{15.1} Each mapping of congruent translation of the
space $f\!:\Bbb E\to\Bbb E$ is the composition $f=g\compos
p_{a\bold c}$, where $g$ is a generalized rotation about some
point $O$, while $p_{a\bold c}$ is a translation by some 
vector $\bold c$ along some straight line $a$.
\endproclaim
\demo{Proof} If the mapping $f$ has a stable point $O$, then the
mapping $f$ itself is a generalized rotation about this point. In
this case we choose the zero vector $\bold c=\bold 0$ on an arbitrary
line $a$ and assign $g=f$. From $p_{a\bold c}=\id$ we get the required
expansion $f=g\compos p_{a\bold c}$ for the initial mapping $f$.\par
      Assume that $f$ has no stable points. Let's choose some arbitrary
point $A$ and denote $O=f(A)$. Then $O\neq A$. Let's draw the line
$AO$ and denote it by $a$. Then we consider the vector $\bold c=
\overrightarrow{AO}$. The translation $p_{a\bold c}\!:\Bbb E\to\Bbb E$
also maps the point $A$ to the point $O$. This mapping is bijective
and the translation by the opposite vector $\overrightarrow{OA}=
-\bold c$ along the line $a$ is the inverse mapping for it. Let's
denote $g=f\compos p^{-1}_{a\bold c}$. It is easy to verify that
the point $O$ is a stable point for the mapping $g$, i\.\,e\. $g$ is
a generalized rotation about the point $O$. From $g=f\compos 
p^{-1}_{a\bold c}$ one easily derives the required expansion
$f=g\compos p_{a\bold c}$ for $f$.
\qed\enddemo
    Using the theorems~\mythetheorem{10.1} and \mythetheorem{15.1},
one finds that an arbitrary congruent translation of the space $f$
admits an expansion of one of the following two sorts: 
$f=\theta_{hk}\compos p_{a\bold c}$ or $f=z_\beta\compos\theta_{hk}
\compos p_{a\bold c}$. Due to the existence of such expansions the
mapping $f$ is bijective since mirror reflections, rotations about 
axes, and translation along lines all are bijective mappings.\par
     Let $f$ and $g$ be two mappings of congruent translation. 
Since $f$ is bijective, we calculate its inverse mapping $f^{-1}$
and construct the composition $g'=f\compos g\compos f^{-1}$. Such
a procedure is called the {\it conjugation} of the mapping $g$ by
means of the mapping $f$. The mapping $g'$ obtained as a result of
conjugation is called a conjugate mapping for $g$.
\mytheorem{15.2} Let $f$ and $g$ are two mappings of congruent
translation and let $g'=f\compos g\compos f^{-1}$ be obtained 
through the conjugation of the mapping $g$ by means of $f$. In
this case
\roster
\item if $g=z_\beta$ is a mirror reflection in a plane $\beta$,
      then $g'$ is a mirror reflection in the plane $f(\beta)$;
\item if $g=\theta_{hk}$ is a rotation about an axis $a$ by an
      angle $\angle hk$, then $g'$ is the rotation about the
      axis $f(a)$ by an angle congruent to $\angle hk$;
\item if $g=p_{a\bold c}$ is a translation by some vector $\bold c$ 
      along some line $a$, then $g'$ is the translation by the 
      vector $f(\bold c)$ along the line $f(a)$.
\endroster
\endproclaim
\myexercise{15.1} Prove the theorem~\mythetheorem{15.2}. For this
purpose use the fact that if $g$ maps a point $A$ to a point $B$, 
then its conjugate mapping $g'=f\compos g\compos f^{-1}$ maps 
the point $f(A)$ to the point $f(B)$.
\endproclaim
\mydefinition{15.1} A congruent translation of the space $f$, which 
expands into the composition of a rotation about some axis and a
translation along some line ($f=\theta_{hk}\compos p_{a\bold c}$),
is said to be {\it even}. If there is a mirror reflection in a plane in
the expansion of $f$ (i\.\,e\. $f=z_\beta\compos\theta_{hk}\compos
p_{a\bold c}$), then $f$ is said to be {\it odd}.
\enddefinition
     A mapping of congruent translation cannot be even and odd at the 
same time. For the mappings which are generalized rotations about a fixed 
point this fact was proved in \S\,11. In general case it should be proved 
separately. Let's prove it by contradiction. Assume that $f$ is a mapping
of congruent translation possessing some expansions of two sorts:
$$
\xalignat 2
&\hskip -2em
\quad f=g_{\sssize A}\compos p_{a\bold c},
&&f=\tilde g_{\sssize B}\compos p_{b\bold d}.
\quad
\mytag{15.1}
\endxalignat
$$
Here $g_{\sssize A}$ is an even rotation about some fixed point $A$,
while $\tilde g_{\sssize B}$ is an odd rotation about some fixed 
point $B$. From the relationships \mythetag{15.1} we derive $g_{\sssize A}
\compos p_{a\bold c}\compos p_{b\bold e}=\tilde g_{\sssize B}$, where
$\bold e$ is the vector opposite to the vector $\bold d$.\par
     For the beginning let's study the case $A=B$. Here we have 
$g_{\sssize A}\compos p_{a\bold c}\compos p_{b\bold e}
=\tilde g_{\sssize A}$. If this relationship is fulfilled, the
composition of translations $p_{a\bold c}\compos p_{b\bold e}$
has the stable point $A$. Applying the theorem~\mythetheorem{14.3},
we get $p_{a\bold c}\compos p_{b\bold e}=\theta_{hk}$, where
$\theta_{hk}$ is the rotation about some axis passing through the
point $A$. The composition of $g_{\sssize A}$ with such a rotation
does not change the parity, therefore the equality $g_{\sssize A}
\compos\theta_{hk}=\tilde g_{\sssize A}$ contradicts the fact 
that $\tilde g_{\sssize A}$ is odd. This contradiction proves 
that the relationships \mythetag{15.1} cannot be fulfilled simultaneously
in the case $A=B$.\par
     Now assume that $A\neq B$. In this case we consider the translation
along the line $AB$ by the vector $\bold s=\overrightarrow{AB}$. Let's
denote it $p_{r\bold s}$. Then $B=p_{r\bold s}(A)$. Using this equality, 
we set $g_{\sssize B}=p_{r\bold s}\compos g_{\sssize A}\compos
p^{-1}_{r\bold s}$. The mapping $g_{\sssize B}$ produced from 
$g_{\sssize A}$ by conjugation is a generalized rotation about a fixed
point $B$, having the same parity as the mapping $g_{\sssize A}$, i\.\,e\.
it is even. This fact follows from the theorem~\mythetheorem{15.2}. From
$g_{\sssize B}=p_{r\bold s}\compos g_{\sssize A}\compos p^{-1}_{r\bold s}$
for the initial mapping $g_{\sssize A}$ we derive
$$
g_{\sssize A}=p^{-1}_{r\bold s}\compos g_{\sssize B}\compos p_{r\bold s}
=g_{\sssize B}\compos(g^{-1}_{\sssize B}\compos p^{-1}_{r\bold s}
\compos g_{\sssize B})\compos p_{r\bold s}.
$$
The mapping $g^{-1}_{\sssize B}\compos p^{-1}_{r\bold s}\compos
g_{\sssize B}$ is produced from $p^{-1}_{r\bold s}$ through 
conjugation by means of $g^{-1}_{\sssize B}$. According to the
theorem~\mythetheorem{15.2} it is a translation by a vector $\bold v$ 
along some line $u$. Then for $g_{\sssize A}$ we have 
$g_{\sssize A}=g_{\sssize B}\compos p_{u\bold v}\compos p_{r\bold s}$. 
Let's substitute this formula into the relationship
$g_{\sssize A}\compos p_{a\bold c}\compos p_{b\bold e}=\tilde
g_{\sssize B}$, which follows from \mythetag{15.1}. Then
$$
\hskip -2em
p_{u\bold v}\compos p_{r\bold s}\compos p_{a\bold c}\compos
p_{b\bold e}=g^{-1}_{\sssize B}\compos\tilde g_{\sssize B}.
\mytag{15.2}
$$
Due to the relationship \mythetag{15.2} the composition of four
translations in the left hand side of this relationship has the
stable point $B$. Applying the theorem~\mythetheorem{14.5} with
$n=4$, we get that this composition is a rotation about some axis
passing through the point $B$. Hence, $\tilde g_{\sssize B}
=g_{\sssize B}\compos\theta_{hk}$, which contradicts the initial 
assumption that $\tilde g_{\sssize B}$ is odd. Thus, in the case
$A\neq B$ the relationships \mythetag{15.1} cannot be fulfilled
simultaneously either.\par
\mydefinition{15.2} A congruent translation of the space $f\!:
\Bbb E\to\Bbb E$ is called a {\it motion}, if it is even.
\enddefinition
     The congruence of segments and angles are basic concepts, they
enter the statements of the axioms. The congruence of triangles is
a derived concept. Passing from triangles to more complicated geometric
forms, we could formulate the definitions of congruence for each 
particular form. However, the concept of congruent translations enable us
to do it at once.
\mydefinition{15.3} Two geometric forms $\Phi_1$ and $\Phi_2$ are 
called {\it congruent}, if there is a congruent translation $f\!:\Bbb E
\to\Bbb E$ performing one-to-one correspondence for the points of these
forms.
\enddefinition
     Since the congruent translations are divided into even and odd ones,
we can sharpen our concept of congruence for arbitrary geometric forms.
\mydefinition{15.4} Two geometric forms $\Phi_1$ and $\Phi_2$ are called
{\it strictly congruent} if there is a motion performing one-to-one 
correspondence for the points of these forms.
\enddefinition
\mydefinition{15.5} Two geometric forms $\Phi_1$ and $\Phi_2$ are 
called {\it mirror congruent} or being {\it mirror images of each other}, 
if there is an odd congruent translation $f\!:\Bbb E\to\Bbb E$ 
performing one-to-one correspondence for the points of these forms.
\enddefinition
\myexercise{15.2} Show that for triangles the 
definitions~\mythedefinition{15.3}, \mythedefinition{15.4}, and
\mythedefinition{15.5} are equivalent to the definition
\mythedefinitionchapter{5.1}{3} from
Chapter~\uppercase\expandafter{\romannumeral 3}.
\endproclaim
\newpage
\setfirstpage
\topmatter
\title\chapter{5}
Axioms of continuity.
\endtitle
\leftheadtext{CHAPTER \uppercase\expandafter{\romannumeral 5}.
AXIOMS OF CONTINUITY.}
\endtopmatter
\document
\chapternum=5
\head
\SectionNum{1}{138} Comparison of straight line segments.
\endhead
\rightheadtext{\S\,1. Comparison of straight line segments.}
\mydefinition{1.1} Let $[AB]$ and $[CD]$ be two straight line segments.
We say that the segment $[AB]$ is {\it smaller\/} than the segment
$[CD]$, and write it as $[AB]<[CD]$, if in the interior of the
segment $[CD]$ there is a point $E$ such that $[AB]\cong [CE]$.
\enddefinition
\mytheorem{1.1} The binary relation of comparison for straight line 
segments introduced in the definition~\mythedefinition{1.1} possesses
the following five properties:
\roster
\item the condition $[AB]<[CD]$ excludes $[AB]\cong [CD]$;
\item the condition $[AB]<[CD]$ excludes $[CD]<[AB]$;
\item $[AB]<[CD]$ and $[CD]<[EF]$ imply $[AB]<[EF]$;
\item if $[\tilde A\tilde B]\cong [AB]$ and $[\tilde C\tilde D]\cong
      [CD]$, then $[AB]<[CD]$ implies  $[\tilde A\tilde B]<[\tilde C
      \tilde D]$;
\item for arbitrary two non-congruent segments $[AB]$ 
      and $[CD]$ one of the two conditions $[AB]<[CD]$ or $[CD]<[AB]$ 
      is always fulfilled.
\endroster
\endproclaim
\demo{Proof} The first item of the theorem is a direct consequence 
of the axiom~\mytheaxiom{A13}. Indeed, if the conditions $[AB]<[CD]$ 
and $[AB]\cong [CD]$ are fulfilled simultaneously, then it would mean
that on the ray $[CD\rangle$ there are two points $D$ and $E$ such 
that $[AB]\cong [CD]$ and $[AB]\cong [CE]$, which contradicts the 
axiom~\mytheaxiom{A13}.\par
     Let's prove the second item of the theorem by contradiction. 
Assume that both conditions $[AB]<[CD]$ and $[CD]<[AB]$ are fulfilled 
simultaneously. Then in the interior of the segment $[CD]$ there is
a point $E$ and in the interior of the segment $[AB]$ there is a point 
$F$ such that the following relationships are fulfilled:
$$
\xalignat 2
&[AB]\cong [CE],&&[CD]\cong [AF].
\endxalignat
$$
Using the relationship $[CD]\cong [AF]$, we define the mapping of 
congruent translation $f$ such that it maps the line $CD$ to the
line $AB$ and such that $f(C)=A$ and $f(D)=F$ (see 
theorems~\mythetheoremchapter{2.1}{3}, \mythetheoremchapter{2.2}{3},
and \S\,3 of Chapter~\uppercase\expandafter{\romannumeral 3}). Let's
denote $\tilde E=f(E)$. Then from $(C\blacktriangleright E
\blacktriangleleft D)$ we get $(A\blacktriangleright\tilde E
\blacktriangleleft F)$. The segments $[A\tilde E]$ and $[CE]$ are
congruent. As a result on the ray $[AB\rangle$ we get two points
$\tilde E\neq B$ such that $[A\tilde E]\cong [CE]$ and $[AB]\cong
[CE]$, which contradicts the axiom~\mytheaxiom{A13}. The contradiction
obtained shows that the conditions $[AB]<[CD]$ and $[CD]<[AB]$ cannot
be fulfilled simultaneously.\par
    Now let's consider the third item of the theorem. The relationships 
$[AB]<[CD]$ and $[CD]<[EF]$ mean that in the interiors of the segments 
$[CD]$ and $[EF]$ there are two points $M$ and $N$ such that the following
conditions are fulfilled:
$$
\xalignat 2
&[AB]\cong [CM],&&[CD]\cong [EN].
\endxalignat
$$
Using the second of these two conditions, we construct a congruent 
translation $f$ mapping the line $CD$ to the line $EF$ and such
that $f(C)=E$ and $f(D)=N$. By means of the mapping $f$ we define
the point $K=f(M)$ in the interior of the segment $[CN]$. For this
point $[EK]\cong [CM]$. From $[EK]\cong [CM]$ and from 
$[AB]\cong [CM]$ we derive $[AB]\cong [EK]$. From
$(E\blacktriangleright K\blacktriangleleft N)$ and $(E\blacktriangleright
N\blacktriangleleft F)$ we get $(E\blacktriangleright K
\blacktriangleleft F)$ (see lemma~\mythelemmachapter{3.2}{2} in
Chapter~\uppercase\expandafter{\romannumeral 2}). In other words, the
point $K$ lies in the interior of the segment $[EF]$ and $[AB]\cong [EF]$. 
The required relationship $[AB]<[EF]$ is proved.\par
     Let's prove the fourth item of the theorem. From $[AB]<[CD]$ we
derive the existence of a point $E$ in the interior of the segment 
$[CD]$ such that $[AB]\cong [CE]$. Using $[CD]\cong [\tilde C\tilde D]$,
we construct a congruent translation $f$ of the line $CD$ to the line
$\tilde C\tilde D$ such that $f(C)=\tilde C$ and $f(D)=\tilde D$. 
Applying $f$ to the point $E$, we define the point $\tilde E=f(E)$ 
lying in the interior of the segment $[\tilde C\tilde D]$ such that 
$[CE]\cong [\tilde C\tilde E]$. Combining this relationship with
$[\tilde A\tilde B]\cong [AB]$ and with $[AB]\cong [CE]$, we get 
$[\tilde A\tilde B]\cong [\tilde C\tilde E]$. The relationship 
$[\tilde A\tilde B]<[\tilde C\tilde D]$ is proved.\par
      In order to prove the fifth item of the theorem, applying the
axiom~\mytheaxiom{A13}, we mark a point $E$ on the ray $[CD\rangle$ 
such that the segment $[AB]$ is congruent to the segment $[CE]$. The
coincidence $D=E$ is impossible since it would mean the congruence
of the segments $[AB]$ and $[CD]$. Therefore the point $E$ lies
in the interior of the segment $[CD]$ or outside this segment. In the
first case $[AB]<[CD]$. In the second case $[CD]<[CE]$. Complementing 
this relationship with $[AB]\cong [CE]$ and applying the fourth item 
of the theorem, which is already proved, we get $[CD]<[AB]$. The theorem 
is proved.
\qed\enddemo
     The properties \therosteritem{1}--\therosteritem{3} of the comparison
for segments are very similar to the corresponding properties of the order
relation (see \S\,3 in Chapter~\uppercase\expandafter{\romannumeral 1}). 
The only difference is that instead of {\tencyr\char '074}$[AB]<[CD]$ 
excludes $[AB]=[CD]${\tencyr\char '076} the theorem~\mythetheorem{1.1} 
here says {\tencyr\char '074}$[AB]<[CD]$ excludes $[AB]\cong [CD]${\tencyr
\char '076}. This means the comparison of segments is an order relation 
not in the set of segments, but in the factorset consisting of classes
of congruent segments. The next two properties \therosteritem{4} and 
\therosteritem{5} mean that this factorset is linearly ordered.\par
     The comparison relation $[AB]<[CD]$ is sometimes written as
$[CD]>[AB]$. In this case we say that the segment $[CD]$ is 
{\it bigger\/} than the segment $[AB]$. The writings $[AB]\leqslant [CD]$
and $[CD]\geqslant [AB]$ mean that one or the two conditions $[AB]<[CD]$ 
or $[AB]\cong [CD]$ is fulfilled.
\mytheorem{1.2} If two points $C$ and $D$ lie in the interior of the
segment $[AB]$, then $[CD]<[AB]$.
\endproclaim
\demo{Proof} Let's introduce the order relation on the line $AB$ by setting
$A\prec B$ (see \S\,4 in Chapter~\uppercase\expandafter{\romannumeral 2}).
If $C$ and $D$ are in the interior of the segment $[AB]$, then one of the
following two conditions is fulfilled:
$$
\xalignat 2
&A\prec C\prec D\prec B,
&&A\prec D\prec C\prec B.
\endxalignat
$$
Let's consider the case where the first condition is fulfilled. The 
case where the second condition is fulfilled is reduced to this case 
by exchanging the notations of the points $C$ and $D$. From the 
relationship $A\prec C\prec D\prec B$ we derive $(A\blacktriangleright C
\blacktriangleleft D)$ and $(A\blacktriangleright D\blacktriangleleft B)$. 
From the first of these two relationships we get $[CD]<[AD]$, while from 
the second one we derive $[AD]<[AB]$. Now, applying the third item of the
theorem~\mythetheorem{1.1}, we obtain the required relationship
$[CD]<[AB]$. The theorem is proved.
\qed\enddemo
\head
\SectionNum{2}{141} Comparison of angles.
\endhead
\rightheadtext{\S\,2. Comparison of angles.}
\mydefinition{2.1} Let $\angle hk$ and $\angle lq$ be two arbitrary
angles. We say that the $\angle hk$ is {\it smaller\/} than the angle
$\angle lq$ and write it as $\angle hk<\angle lq$ if there is a ray 
$m$ coming out from the vertex of the angle $\angle lq$ and lying 
inside it so that $\angle hk\cong\angle lm$.
\enddefinition
\mytheorem{2.1} The binary relation of comparison for angles possesses
the following five properties:
\roster
\item the condition $\angle hk<\angle lq$ excludes 
      $\angle hk\cong\angle lq$;
\item the condition $\angle hk<\angle lq$ excludes 
      $\angle lq<\angle hk$;
\item $\angle hk<\angle lq$ and $\angle lq<\angle mn$ imply
      $\angle hk<\angle mn$;
\item if $\angle \tilde h\tilde k\cong\angle hk$ and if 
      $\angle \tilde l\tilde q\cong\angle lq$, then the relationship
      $\angle hk<\angle lq$ implies $\angle \tilde h\tilde k
      <\angle \tilde l\tilde q$;
\item for any two non-congruent angles $\angle hk$ and 
      $\angle lq$ one of the two conditions $\angle hk<\angle lq$ 
      or $\angle lq<\angle hk$ is fulfilled .
\endroster
\endproclaim
\mytheorem{2.2} If the rays $l$ and $q$ coming out from the vertex
of the $\angle hk$ lie inside this angle, then $\angle lq<\angle hk$.
\endproclaim
\myexercise{2.1} By analogy to theorems~\mythetheorem{1.1} and
\mythetheorem{1.2} prove the theorem~\mythetheorem{2.1} and
the theorem~\mythetheorem{2.2}.
\endproclaim
     In the case of angles we have two reference angles, they are 
the right angle and the straight angle. Each angle can be enclosed 
into some straight angle, therefore, any angle is smaller than any
straight angle. Comparing angles with the right angle we divide
those different from a right angle into two sets. An angle smaller 
than a right angle is called an {\it acute angle}. An angle bigger
than a right angle is called an {\it obtuse angle}. According to
the theorem~\mythetheorem{2.1} any acute angle is smaller than any
obtuse angle.
\mydefinition{2.2} An angle adjacent to an internal angle of a triangle
is called an {\it external angle} of this triangle.
\enddefinition
\mytheorem{2.3} In a triangle any internal angle is smaller than any
external angle not adjacent with it.
\endproclaim
\demo{Proof} \parshape 12 4cm 6cm 4cm 6cm 4cm 6cm 4cm 6cm 4cm 6cm 
4cm 6cm 4cm 6cm 4cm 6cm 4cm 6cm 4cm 6.1cm 4cm 6.1cm 0cm 10cm
     Let $ABC$ be a triangle . We extend its side $[BC]$ by 
drawing the ray $[CE\rangle$. The angle $\angle ACE$ adjacent to 
the angle $\angle ACB$ is an external angle of the triangle 
$ABC$. \vadjust{\vskip 5pt\hbox to 0pt{\kern -5pt
\includegraphics{Oris50.eps}\hss}\vskip -5pt}We compare 
it with the internal angle $\angle CAB$ not adjacent to $\angle ACE$. 
The line $AC$ divides the plane of the triangle $ABC$ into two
half-planes, the points $E$ and $B$ lying on different sides of this 
line. Applying the axiom~\mytheaxiom{A16}, in the half plane containing
the point $B$ we draw the ray $[AD\rangle$ so that the angle $\angle CAD$ 
is congruent to the angle $\angle ACE$. If we assume that
$\angle CAD\cong\angle CAB$ or $\angle CAD<\angle CAB$, then the ray
$[AD\rangle$ coincides with the ray $[AB\rangle$ or it lies inside the
angle $\angle CAB$. In both cases the ray $[AD\rangle$ should intersect
the segment $[BC]$ (see lemma~\mythelemmachapter{6.2}{2} in
Chapter~\uppercase\expandafter{\romannumeral 2}). But in our case the
angles $\angle CAD$ and $\angle ACE$ are inner crosswise lying angles
at the intersections of the line $AC$ with two lines $AD$ and $BC$. 
They are congruent $\angle CAD\cong\angle ACE$, therefore, due to
the theorem~\mythetheoremchapter{8.1}{3} from 
Chapter~\uppercase\expandafter{\romannumeral 3} the lines $AD$ and
$BC$ cannot intersect. The contradiction obtained proves the required 
relationship $\angle CAB<\angle ACE$.
\qed\enddemo
\mytheorem{2.4} In a triangle a bigger side is opposite to a bigger 
angle and, conversely, a bigger angle is opposite to a bigger side.
\endproclaim
\demo{Proof}\parshape 9 4cm 6cm 4cm 6cm 4cm 6cm 4cm 6cm 4cm 6cm 
4cm 6cm 4cm 6cm 4cm 6cm 0cm 10cm
      Let $ABC$ be a triangle. For the beginning we prove that  
the relationship $[AC]>[BC]$ implies the relationship 
$\angle ABC>\angle BAC$. \vadjust{\vskip 5pt\hbox to 0pt{\kern 0pt
\includegraphics{Oris51.eps}\hss}\vskip -5pt}From 
$[AC]>[BC]$ we conclude that in the interior of the segment $[AC]$ 
there is a point $D$ such that the segment $[DC]$ is congruent to the
segment $[BC]$. Let's draw the segment $[BD]$ and consider the 
triangle $BCD$. Due to $[DC]\cong [BC]$ it is isosceles. Hence,
$\angle BDC\cong\angle DBC$.\par
     By construction the ray $[BD\rangle$ is inside the angle $ABC$. 
Therefore $\angle ABC>\angle DBC$. On the other hand, the angle
$\angle BDC$ is an external angle of the triangle $ABD$. Applying the
theorem~\mythetheorem{2.3}, we get $\angle BDC>\angle BAC$. As a result
we have three relationships
$$
\xalignat 3
&\angle ABC>\angle DBC,&&\angle BDC\cong\angle DBC,
&&\angle BDC>\angle BAC.
\endxalignat
$$
Applying the theorem~\mythetheorem{2.1}, from these relationships we 
derive $\angle ABC>\angle BAC$. The first proposition of the theorem 
is proved.\par
     We prove the converse proposition by contradiction. Assume 
that the relationship $\angle ABC>\angle BAC$ is fulfilled, but the 
relationship $[AC]>[BC]$ is not fulfilled. Then $[AC]\cong [BC]$ 
or $[AC]<[BC]$. If $[AC]\cong [BC]$, then the triangle $ABC$ is
isosceles and $\angle ABC\cong\angle BAC$, which contradicts 
$\angle ABC>\angle BAC$. If $[AC]<[BC]$, then due to the first
proposition of the theorem, which is already proved, we have
$\angle ABC<\angle BAC$, which contradicts $\angle ABC>\angle BAC$. 
Hence, $\angle ABC>\angle BAC$ implies $[AC]>[BC]$. The theorem
is proved.\qed\enddemo
     Assume that a point $B$ lies \pagebreak in the interior of a 
segment $[AC]$ between the points $A$ and $C$. The segment $[AC]$ 
in this case is composed of two segments $[AB]$ and $[BC]$ (see 
theorems~\mythetheoremchapter{3.1}{2} and \mythetheoremchapter{3.2}{2}
in Chapter~\uppercase\expandafter{\romannumeral 2}). It is sometimes 
called the {\it sum\/} of the segments $[AB]$ and $[BC]$. The segment
$[BC]$ is called the {\it difference\/} of the segments $[AC]$ and
$[AB]$.\par
     Having two arbitrary segments $[MN]$ and $[PQ]$, we can draw
two segments congruent to them on one line so that they lie on
different sides of some point $B$ on this line:
$$
\xalignat 2
&[MN]\cong [AB],&&[PQ]\cong [BC].
\endxalignat
$$
Then the segment $[AC]$ is the sum of the segments $[MN]$ and $[PQ]$. 
Such a segment is not unique, however, all segments representing the 
sum of the segments $[MN]$ and $[PQ]$ are congruent to each other.\par
     Let $[MN]<[PQ]$. Let's draw a segment $[AC]$ congruent to
$[PQ]$ on some line. Then $[MN]<[AC]$, therefore, in the interior
of the segment $[AC]$ there is a point $B$ such that $[AB]\cong [MN]$. 
For such a point $B$ the segment $[BC]$ is the difference of the
segments $[PQ]$ and $[MN]$.
\mytheorem{2.5} In an arbitrary triangle the sum of any two sides is 
bigger than the third side.
\endproclaim
\demo{Proof} \parshape 11 4.5cm 5.5cm 4.5cm 5.5cm 4.5cm 5.5cm
4.5cm 5.5cm 4.5cm 5.5cm 4.5cm 5.5cm 4.5cm 5.5cm 4.5cm 5.5cm
4.5cm 5.5cm 4.5cm 5.5cm 0cm 10cm
    Let $ABC$ be some triangle. The point $B$ divides the line
$AB$ into two rays. \vadjust{\vskip 5pt\hbox to 0pt{\kern 0pt
\includegraphics{Oris52.eps}\hss}\vskip -5pt}Let's 
choose the ray opposite to the ray $[BA\rangle$ and mark the 
segment $[B\tilde C]$ congruent to the segment $[BC]$ on this 
ray. The ray $[CB\rangle$ crosses the segment $[A\tilde C]$, 
therefore it lies inside the angle $\angle AC\tilde C$ (see
lemma~\mythelemmachapter{6.2}{2} in
Chapter~\uppercase\expandafter{\romannumeral 2}). This yields
the formula $\angle AC\tilde C>\angle BC\tilde C$. On the other
hand, in the isosceles triangle $CB\tilde C$ the angles 
$\angle BC\tilde C$ and $\angle B\tilde CC$ are congruent. 
Hence, $\angle AC\tilde C>\angle A\tilde CC$. Applying the
theorem~\mythetheorem{2.4} to the triangle $AC\tilde C$, we 
get $[A\tilde C]>[AC]$. By construction the segment $[A\tilde C]$ 
is the sum of the segments $[AB]$ and $[BC]$. It is bigger that
the segment $[AC]$. The theorem is proved.\qed\enddemo
\vskip 0pt plus 1pt minus 1pt
\mytheorem{2.6} In any triangle at least two angles are acute.
\endproclaim
\vskip 0pt plus 1pt minus 1pt
\demo{Proof} \parshape 11 4.5cm 5.5cm 4.5cm 5.5cm 4.5cm 5.5cm 
4.5cm 5.5cm 4.5cm 5.5cm 4.5cm 5.5cm 4.5cm 5.5cm 4.5cm 5.5cm 
4.5cm 5.5cm 4.5cm 5.5cm 0cm 10cm
     Let $ABC$ be a triangle, the angle $\angle ABC$ of which is
obtuse. \vadjust{\vskip 5pt\hbox to 0pt{\kern 0pt
\includegraphics{Oris53.eps}\hss}\vskip -5pt}Let's 
draw the perpendicular to the line $AB$ through the point $B$. 
Since an obtuse angle is bigger than a right angle, 
this perpendicular crosses the side $AC$ at some interior point
$\tilde C$. As a result we get the triangle $AB\tilde C$ in which
the obtuse angle is replaced by the right angle, while the angle
at the vertex $A$ remains unchanged. If this angle was initially
obtuse, repeating the procedure, we can replace id by a right 
angle either. Thus, we conclude: if a triangle has two obtuse 
angle or if it has an obtuse angle and a right angle, it can 
be transformed to a triangle with two right angles. However, 
we know that a triangle with two right angles is impossible 
(see theorem~\mythetheoremchapter{8.2}{3} in
Chapter~\uppercase\expandafter{\romannumeral 3}). 
\qed\enddemo
\vskip 0pt plus 1pt minus 1pt
      Due to the theorem~\mythetheorem{2.6} all triangles are divided
into three sorts: {\it acute-angular triangles\/}, {\it right-angular
triangles}, and {\it obtuse-angular triangles}. In an acute-angular 
triangle all angles are acute; a right-angular triangle has one
right angle and two acute angles; an obtuse-angular triangle has
one obtuse angle, while two other angles are acute.\par
\vskip 0pt plus 1pt minus 1pt
     In a right-angular triangle the side opposite to the right angle
is called the {\it hypotenuse}. Other two sides are called {\it legs}.
Due to the theorems~\mythetheorem{2.4} and \mythetheorem{2.6} in a
right-angular triangle the hypotenuse is bigger \pagebreak
than any one of two legs.
\head
\SectionNum{3}{146} Axioms of real numbers.
\endhead
\rightheadtext{\S\,3. Axioms of real numbers.}
      The comparison of segments and angles introduced in 
definitions~\mythedefinition{1.1} and \mythedefinition{2.1} yields
very rough concept of what are the sizes of objects being compared.
More precise knowledge of a size require the concept of length for
segments and a quantitative measure for angles. In order to introduce
these concepts we need to use some facts from the theory of real
numbers.\par
       Real numbers constitute the set $\Bbb R$, where two
algebraic operations are defined --- the {\it addition\/} and
the {\it multiplication}. Basic properties of real numbers are
formulated in seventeen axioms R1--R17. The whole theory of real
numbers is deduced from these axioms.
\myaxiom{R1} The addition of real numbers is commutative, i\.\,e\.
$a+b=b+a$ for all $a$ and $b$ in $\Bbb R$.
\endproclaim
\myaxiom{R2} The addition of real numbers is associative, i\.\,e\.
$(a+b)+c=a+(b+c)$ for all $a$, $b$, and $c$ in $\Bbb R$.
\endproclaim
\myaxiom{R3} There is a number $0$ in $\Bbb R$, which is called zero,
such that $a+0=a$ for all $a\in\Bbb R$.
\endproclaim
      One can prove that zero is unique. Indeed, if we assume that 
there is another number $\tilde 0$ with the same property, then 
$\tilde 0+0=\tilde 0$ and $0+\tilde 0=0$. From the axiom~\mytheaxiom{R1}
we derive that they do coincide: $\tilde 0=\tilde 0+0=0+\tilde 0=0$.
\myaxiom{R4} For any number $a\in\Bbb R$ there is an opposite number 
$a'\in\Bbb R$ such that $a+a'=0$.
\endproclaim
     For any number $a\in\Bbb R$ its opposite number $a'$ is unique.
If we assume that there is another opposite number $\tilde a'$ for $a$,
then from the axioms~\mytheaxiom{R1}, \mytheaxiom{R2}, and \mytheaxiom{R3}
we derive $a'=\tilde a'$:
$$
a'=a'+0=a'+(a+\tilde a')=(a'+a)+\tilde a'=0+\tilde a'=\tilde a'.
$$
\myexercise{3.1} Let $a'$ be the number \pagebreak 
opposite to the number $a\in\Bbb R$. Prove that the number opposite 
to $a'$ coincides with the number $a$, i\.\,e\. $(a')'=a$. 
\endproclaim
\myexercise{3.2} Prove that the number opposite to zero coincides with
itself, i\.\,e\.  $0'=0$.
\endproclaim
     The concept of an opposite number is the base for introducing the
{\it subtraction}: $b-a=b+a'$. Moreover, the opposite number is usually
denoted $a'=-a$ for the sake of uniformity of notations.
\myexercise{3.3} Prove the following relationships:
\vadjust{\vskip -6pt}
$$
\xalignat 2
&a-b=-(b-a),&&(a-b)+c=a+(c-b),\\
&(a+b)-c=a+(b-c),&&(a-b)-c=a-(b+c).
\endxalignat
$$
\endproclaim
\myaxiom{R5} The multiplication of real numbers is commutative,
i\.\,e\. $a\cdot b=b\cdot a$ for all $a$ and $b$ in $\Bbb R$.
\endproclaim
\myaxiom{R6} The multiplication of real numbers is associative, 
i\.\,e\. $(a\cdot b)\cdot c=a\cdot (b\cdot c)$ for all $a$, $b$, 
and $c$ in $\Bbb R$.
\endproclaim
\myaxiom{R7} There is a number $1$, which is called one, such that
$a\cdot 1=a$ for all $a\in\Bbb R$.
\endproclaim
\myaxiom{R8} For any nonzero number $a\neq 0$ there is an inverse
number $a^*\in\Bbb R$ such that $a\cdot a^*=1$.
\endproclaim
\myexercise{3.4} Prove that the number one is unique. Also prove 
that for any number $a\neq 0$ its inverse number $a^*$ is unique.
\endproclaim
     The concept of an inverse number is the base for introducing the
{\it division}: $b:a=b\cdot a^*$. There are several ways for denoting 
the operation of division for real numbers:
$$
b:a=b/a=\frac{b}{a}=b\cdot a^{-1}.
$$
The last form writing the quotient $b/a$ is due to the notation
$a^*=a^{-1}$ for an inverse number.
\myaxiom{R9} The multiplication and the addition of real numbers 
are related trough the distributivity law: $(a+b)\cdot c=a\cdot c+b\cdot c$.
\endproclaim
\myexercise{3.5} Prove the distributivity law relating the multiplication
and the subtraction: $(a-b)\cdot c=a\cdot c-b\cdot c$. Also derive the
following rules for operating with fractions:
$$
\xalignat 2
&\frac{b}{a}+\frac{d}{c}=\frac{b\cdot c+a\cdot d}{a\cdot c},
&&\frac{b}{a}\cdot\frac{d}{c}=\frac{b\cdot d}{a\cdot c}.
\endxalignat
$$
\endproclaim
    The set of real numbers $\Bbb R$ is equipped with a binary relation 
of order with respect to which it is a linearly ordered set (see \S\,3 
in Chapter~\uppercase\expandafter{\romannumeral 1}).
\myaxiom{R10} For any two numbers $a$ and $b$ in $\Bbb R$ at least 
one of the three condition $a<b$, $a=b$, or $b<a$ is fulfilled.
\endproclaim
\myaxiom{R11} The condition $a<b$ excludes $a=b$.
\endproclaim
\myaxiom{R12} The condition $a<b$ excludes $b<a$.
\endproclaim
\myaxiom{R13} The conditions $a<b$ and $b<c$ imply $a<c$.
\endproclaim
      The following two axioms bind the relation of order with algebraic 
operations of addition and multiplication for real numbers.
\myaxiom{R14} The condition $a<b$ implies $a+c<b+c$ for any number
$c\in\Bbb R$.
\endproclaim
\myaxiom{R15} The conditions $a<b$ and $0<c$ imply $a\cdot c<b\cdot c$.
\endproclaim
\myexercise{3.6} Prove that the relation of order in the set of real 
numbers possesses the following properties:
\roster
\item $a>b$ implies $-a<-b$;
\item $a<b$ and $c<0$ imply $a\cdot c>b\cdot c$;
\item $1>0$;
\item $a>b>0$ implies $b^{-1}>a^{-1}>0$.
\endroster
\endproclaim
     Positive integers are obtained by successive adding the unity:
$2=1+1$, $3=1+1+1=2+1$ etc. They constitute the set $\Bbb N=\{1,\,2,
\,3,\,\ldots\}$. Complementing $\Bbb N$ with the number zero and with
the numbers opposite to positive integers, we get the set of all
integers $\Bbb Z=\{\ldots,\,-2,\,-1,\,0,\,1,\,2,\,\ldots\}$. Fractions
$n/m$, where $n\in\Bbb Z$ and $m\in\Bbb N$, constitute the set of rational 
numbers $\Bbb Q$. We consider also the set fractions of the special form
$$
r=\frac{n}{2^m}\text{, \ where \ }n\in\Bbb Z\text{ \ and \ }
m\in\Bbb N.
$$
Such fractions constitute the set of {\it binary-rational numbers}.
THe set of binary-rational numbers is closed with respect to the
addition, subtraction, and multiplication, but it is not closed with
respect to the division.
\myaxiom{R16} For any real number $\xi\in\Bbb R$ there is a positive
integer $n\in\Bbb N$ such that $n>\xi$.
\endproclaim
     The axiom~\mytheaxiom{R16} is known as the Archimedes axiom. It 
means that each real number has an upper estimate in the set of 
positive  integers. This estimate can be strengthened.
\mytheorem{3.1} For any real number $\xi$ there is a positive 
integer $n\in\Bbb N$ such that $-n<\xi<n$.
\endproclaim
\demo{Proof}  If $\xi=0$, we can take $n=1$. If $\xi>0$, then the
number $n$ given by the Archimedes axiom~\mytheaxiom{R16} provides
the estimate $-n<\xi<n$. If $\xi<0$, we apply the Archimedes axiom
to the number $-\xi$. The resulting number $n$ in this case 
provides the estimate $-n<\xi<n$.
\qed\enddemo
\mytheorem{3.2} For any real number $\xi>0$ there is a positive
integer $m$ such that $2^{-m}<\xi$.
\endproclaim
\demo{Proof} Let's consider the number $\xi^{-1}$ and apply the
Archi\-medes axiom~\mytheaxiom{R16} to it. As a result we get a
positive integer $m$ such that $\xi^{-1}<m$. Let's use the inequality
$m<2^m$ which is fulfilled for all positive integers $m\in\Bbb N$. 
It is easily proved by induction on $m$. Now from $\xi^{-1}<2^m$ and
$\xi>0$ we derive $\xi>2^{-m}$. The theorem is proved.
\qed\enddemo
\myaxiom{R17} Let $\{a_n\}_{n\sssize\in\Bbb N}$ and
$\{b_n\}_{n\sssize\in\Bbb N}$ be two monotonic sequences of real 
numbers such that
$$
a_1\leqslant\ldots\leqslant a_n\leqslant a_{n+1}\leqslant\ldots
\leqslant b_{n+1}\leqslant b_n\leqslant\ldots\leqslant b_1.
$$
Then there is a real number $\xi\in\Bbb R$ separating these
sequences, i\.\,e\. such that $a_n\leqslant\xi\leqslant b_n$ 
for all $n\in\Bbb N$.
\endproclaim
The axiom~\mytheaxiom{R17} is known as Cantor's axiom. It 
expresses the property of completeness of real numbers. 
Cantor's axiom~\mytheaxiom{R17} lies in the base of proving many
well-known facts of mathematical analysis (see details in \mycite{6}).
\head
\SectionNum{4}{150} Binary rational approximations of real numbers.
\endhead
\rightheadtext{\S\,4. Binary rational approximations \dots}
     Let $p$, $q$, and $\xi$ be three real numbers. We say that
the real numbers $p$ and $q$ {\it approximate\/} the number real 
$\xi$ if 
$$
p\leqslant\xi\leqslant q.
$$
The number $p$ is called the {\it lower estimate}, while $q$ is
the {\it upper estimate}. The difference $q-p$ is called the 
{\it accuracy} of the approximation.\par
     Let's study the problem of approximating real numbers with
binary rational numbers. Let's denote by $I(m,k)$ the following
semi-open intervals within the set of real numbers:
$$
I(m,k)=\left[\frac{k}{2^m},\frac{k+1}{2^m}\right)
=\left\{x\in\Bbb R:\quad\frac{k}{2^m}\leqslant x
<\frac{k+1}{2^m}\right\}.
$$
The integer $k\in\Bbb Z$ is called the {\it number} of such interval,
while the number $m\in\Bbb N$ is its {\it order}.
\mytheorem{4.1} For any fixed $m\in\Bbb N$ each real number $\xi$ 
is trapped exactly in one interval $I(m,k)$.
\endproclaim
\demo{Proof} The intervals $I(m,k_1)$ and $I(m,k_2)$ with different
numbers $k_1\neq k_2$ do not intersect each other. Therefore, once 
the number $\xi$ is in some interval, it cannot be in any other.\par
     Applying the theorem~\mythetheorem{3.1} to the number $\xi$, we
find that there is an integer $n$, such that $\xi$ is in some
interval $(-n,n)$ enclosed between the numbers $-n$ and $n$. But for this
interval we have
$$
(-n,n)\subset\bigcup^{2^mn}_{k=-2^mn}I(m,k).
$$
Therefore, the number $\xi$ is certainly contained in one of the intervals
$I(m,k)$ whose union contains the interval $(-n,n)$.
\qed\enddemo
     The theorem~\mythetheorem{4.1} proved just above means that for 
any real number $\xi$ and for any positive integer $m$ there is a unique 
integer $k_m$ defining two binary-rational numbers $a_m$ and $b_m$ such 
that
$$
\hskip -2em
\qquad \frac{k_m}{2^m}=a_m\leqslant\xi<b_m=\frac{k_m+1}{2^m}.
\mytag{4.1}
$$
The numbers $a_m$ and $b_m$ are called the {\it binary-rational
approximations\/} of the order $m$ for the number $\xi$.\par
\mylemma{4.1} The sequences $\{a_m\}_{m\sssize\in\Bbb N}$ and
$\{b_m\}_{m\sssize\in\Bbb N}$ of binary-rational approximations
of a given real number $\xi$ are monotonic: $a_{m+1}\geqslant a_m$ 
and $b_{m+1}\leqslant b_m$. The sequence $\{a_m\}_{m\sssize\in\Bbb N}$
stabilizes if and only if the number $\xi$ itself is a binary-rational
number.
\endproclaim
\demo{Proof} When advancing the number $m$ by one the interval $I(m,k)$
divides into two mutually non-intersecting intervals of the next order:
$I(m,k)=I(m+1,2k)\cup I(m+1,2k+1)$. Therefore, when passing from $m$ to
$m+1$ we obtain $k_{m+1}=2k_m$ or $k_{m+1}=2k_m+1$. For $a_{m+1}$ and
$b_{m+1}$ this yields
$$
\xalignat 2
&a_{m+1}=a_{m},&&b_{m+1}=b_m-\frac{1}{2^{m+1}},\\
\displaybreak
&a_{m+1}=a_m+\frac{1}{2^{m+1}},&&b_{m+1}=b_m.
\endxalignat
$$
In any of these two cases the monotony conditions $a_{m+1}\geqslant a_m$
and $b_{m+1}\leqslant b_m$ are fulfilled.\par
     Let's consider the case where the number $\xi$ is a binary-rational.
Then $\xi$ can be written as a fraction
$$
\hskip -2em
\xi=\frac{q}{2^n}=\frac{q\cdot 2^{m-n}}{2^m}.
\mytag{4.2}
$$
From \mythetag{4.2} for $k_m$, where $m\geqslant n$, we get 
$k_m=q\cdot 2^{m-n}$. For the numbers $a_m$ and $b_m$, where
$m\geqslant n$, this yields
$$
\xalignat 2
&a_m=\xi,&&b_m=\xi+\frac{1}{2^m}.
\endxalignat
$$
This means that the sequence $\{a_m\}_{m\sssize\in\Bbb N}$ stabilizes
at the value $a_m=\xi$ for $m\geqslant n$.\par
     Assume, conversely, that the sequence $\{a_m\}_{m\sssize\in\Bbb N}$
stabilizes for $m\geqslant n$. Let $a=a_n$ be the value at which this
sequence stabilizes. Then from \mythetag{4.1} we derive
$$
\hskip -2em
a\leqslant\xi<a+\frac{1}{2^m}.
\mytag{4.3}
$$
The inequalities \mythetag{4.3} are easily transformed to
$$
\frac{1}{2^m}>\xi-a\geqslant 0.
$$
If $\xi\neq a$, then $\xi-a\neq 0$ and we get a contradiction with
the theorem~\mythetheorem{3.2}. Hence, the number $\xi$ coincides
with the binary-rational number $a$. The lemma is proved.
\qed\enddemo\nopagebreak
\myexercise{4.1} Show that the second sequence $\{b_m\}_{m\sssize\in
\Bbb N}$ cannot stabilize for any choice of the number $\xi$.
\endproclaim
     Let $\{a_m\}_{m\sssize\in\Bbb N}$ and $\{b_m\}_{m\sssize\in
\Bbb N}$ be the sequences of binary-rational approximations of 
the number $\xi$. Due to the lemma~\mythelemma{4.1} we can apply
Cantor's axiom~\mytheaxiom{R17}. Let $\tilde\xi$ is a number whose
existence is declared in this axiom. Then we have two inequalities
$$
\xalignat 2
&a_m\leqslant\xi<b_m,
&&a_m\leqslant\tilde\xi\leqslant b_m.
\endxalignat
$$
From these inequalities  we derive the estimate for the modulus of the
difference of two numbers $\xi$ and $\tilde\xi$:
$$
|\xi-\tilde\xi|\leqslant b_m-a_m=\frac{1}{2^m}.
$$
If we assume that $\xi\neq\tilde\xi$, then we immediately get a
contradiction to the theorem~\mythetheorem{3.2}. Therefore, $\xi
=\tilde\xi$. Thus, we have an important conclusion: any real number 
$\xi$ is uniquely fixed by two sequences of its binary-rational
approximations \mythetag{4.1}.\par
\mytheorem{4.2} Let $\xi<\eta$ be two real numbers and let $a_m$, $b_m$,
$\tilde a_m$, and $\tilde b_m$ be the sequences of their binary-rational
approximations. Then there is an integer number $n$ such that 
$b_m<\tilde a_m$ for all integers $m>n$.
\endproclaim
\demo{Proof} Let's consider the relationships \mythetag{4.1}. From these
relationships for the numbers $a_m$ and $b_m$ we derive the estimates
$$
\hskip -2em
\aligned
b_m=a_m+\frac{1}{2^m}\leqslant\xi+\frac{1}{2^m},\\
\vspace{1ex}
a_m=b_m-\frac{1}{2^m}>\xi-\frac{1}{2^m}.
\endaligned
\mytag{4.4}
$$
The estimates analogous to \mythetag{4.4} can be written for the numbers
$\tilde a_m$ and $\tilde b_m$ either. Let's write such estimates for 
$b_m$ and $\tilde a_m$:
$$
\xalignat 2
&\tilde a_m>\eta-\frac{1}{2^m},
&&b_m\leqslant\xi+\frac{1}{2^m}.
\endxalignat
$$
Subtracting the second inequality from the first one, we get
$$
\hskip -2em
\tilde a_m-b_m>(\eta-\xi)-\frac{1}{2^{m-1}}.
\mytag{4.5}
$$
Relying on the theorem~\mythetheorem{3.2}, we choose a number $n$ such
that 
$$
\hskip -2em
\eta-\xi>\frac{1}{2^{n-1}}.
\mytag{4.6}
$$
From \mythetag{4.5} and from \mythetag{4.6} for the difference
$\tilde a_m-b_m$ we get the lower estimate 
$$
\tilde a_m-b_m>\frac{1}{2^{n-1}}-\frac{1}{2^{m-1}}.
$$
For $m>n$ it yields the required inequality $b_m<\tilde a_m$. The 
theorem is proved.
\qed\enddemo
\head
\SectionNum{5}{154}  The Archimedes axiom and Cantor's axiom in geometry.
\endhead
\rightheadtext{\S\,5. The Archimedes axiom and Cantor's axiom \dots}
     Assume that a segment $[OE]$ on some line $a$ is given. Let's 
define the positive direction on this line by means of the vector
$\bold c=\overrightarrow{OE}$. Let's denote $A(0)=O$, $A(1)=E
=p_{\bold c}(A(0))$ and construct the sequence of points 
$A(2)=p_{\bold c}(A(1))$, $A(3)=p_{\bold c}(A(2))$ etc by applying
the congruent translation by the vector $\bold c$ repeatedly. Let's
enumerate with the negative integers the points obtained by several 
subsequent translations of the point $O$ by the opposite vector 
$-\bold c$: $A(-1)=p_{-\bold c}(A(0))$, $A(-2)=p_{-\bold c}(A(-1))$
etc. This sequence of points enumerated with integer numbers 
appears to be monotonic (see definition in \S\,4 of
Chapter~\uppercase\expandafter{\romannumeral 2}):
$$
\ldots\prec A(-2)\prec A(-1)\prec A(0)\prec A(1)\prec A(2)\prec
\ldots\,.
$$
The segments connecting the neighboring points in this sequence
are congruent to each other:
$$
[A(-n)A(-n+1)]\cong\,\ldots\,\cong [OE]\cong\,\ldots\,\cong [A(n-1)A(n)]
$$
Such a monotonic sequence of points is called {\it equidistant}. The 
segment $[A(0)A(2)]$ is obtained by adding two segments --- the
segment $[A(0)A(1)]$ and the segment $[A(1)A(2)]$, each being congruent 
to $[OE]$. We say that it is obtained by {\it duplicating\/} the segment
$[OE]$. The segment $[A(0)A(3)]$ is obtained by {\it triplicating\/} the
segment $[OE]$, the next segment $[A(0)A(4)]$ is obtained by 
{\it quadrupling\/} etc. Let's visualize these words as follows:
$$
[A(0)A(2)]\cong 2\cdot [OE],\,\ldots,\,
[A(0)A(n)]\cong n\cdot [OE].
$$
\myaxiom{A18} For any two segments $[AB]$ and $[OE]$ there is 
a positive integer $n\in\Bbb N$ such that $[AB]<[CD]\cong n
\cdot [OE]$, where the segment $[CD]$ is obtained as the sum 
of $n$ replicas of the segment $[OE]$.
\endproclaim
     The axiom~\mytheaxiom{A18} is known as the Archimedes axiom. It is
a geometric version of the axiom~\mytheaxiom{R16} considered in \S\,3.
\par
     Let's consider the segment $[OE]$ on the line $a$ again. According
to the theorem~\mythetheoremchapter{7.2}{3} from
Chapter~\uppercase\expandafter{\romannumeral 3}, in the interior of 
this segment one can find a point $E_1$ being its center. Then we 
have $[OE]\cong 2\cdot [OE_1]$. Let's write this fact as follows:
$$
[OE_1]\cong\frac{1}{2}\cdot [OE].
$$
Applying the operation of bisection once more, we find the center 
of the segment $[OE_1]$. We denote it $E_2$. Repeating this procedure 
many times, we get a series of segments each of which is twice as small 
as the previous one:
$$
[OE_1]\cong\frac{1}{2}\cdot [OE],\,\ldots,\,[OE_m]\cong
\frac{1}{2^m}\cdot [OE].
$$
The following theorem is a geometric version of the 
theorem~\mythetheorem{3.2}.
\mytheorem{5.1} For any two segments $[AB]$ and $[OE]$ there is a positive
integer $m$ such that 
$$
[AB]>[CD]\cong\frac{1}{2^m}\cdot [OE],
$$
where the segment $[CD]$ is obtained as a result of $m$-tuple bisection
of the segment $[OE]$.
\endproclaim
\demo{Proof} If we assume that the proposition of the theorem is not valid,
then for any positive integer $m$ and for $n=2^m$ we have
$$
[AB]<[CD]\cong\frac{1}{2^m}\cdot [OE]
\text{ \ or \ }[OE]>2^m\cdot [AB],
$$
which contradicts the Archimedes axiom~\mytheaxiom{A18}. This contradiction
proves the theorem.\qed\enddemo
\myaxiom{A19} Let $\{[A_nB_n]\}_{\sssize n\in\Bbb N}$ be a sequence
of segments on some straight line such that
$$
[A_{n+1}B_{n+1}]\subset [A_nB_n]\text{ \ for all \ }
n\in\Bbb N.
$$
Then the intersection of these segments is not empty and there is a point
$X$ belonging to all of them.
\endproclaim
The axiom~\mytheaxiom{A19} is called Cantor's axiom. It is a geometric
version of Cantor's axiom~\mytheaxiom{R17} from \S\,3.
\head
\SectionNum{6}{156} The real axis.
\endhead
\rightheadtext{\S\,6. The real axis.}
     Let's consider again the line $a$ and the segment $[OE]$ on it. 
In \S\,5 we have constructed the sequence of points 
$\{A(n)\}_{n\sssize\in\Bbb Z}$ on this line enumerated by integer 
numbers. The sequence $\{A(n)\}_{n\sssize\in\Bbb Z}$ is equidistant 
and monotonic, it is such that the inequality $m<n$ implies
$A(m)\prec A(n)$. It defines a {\it coordinate network\/} (or a 
{\it gauge}) on the line $a$. We bisect each of the segments connecting 
the neighboring points $A(n)$ and $A(n+1)$, then we enumerate the centers
of these segments by means of half-integers. As a result we get the doubly
dense gauge of points:
$$
\ldots\prec A(-1)\prec A(-1/2)\prec A(0)\prec A(1/2)\prec A(1)
\prec\ldots\,.
$$
The procedure of doubling the density of points on the line $a$ can 
be performed repeatedly again and again. On the second step we use 
the numbers of the forms $k/4$, then of the form $k/8$, $k/16$ etc.
Repeating this procedure infinitely many times, we get the set of
points enumerated with binary-rational numbers. Let's denote this
set through $\Cal A$, while the set of binary-rational numbers through
$\Bbb D$. Then, putting into correspondence the point $A(r)$ to a number
$r\in\Bbb D$, we get  a mapping 
$$
\hskip -2em
A\!:\Bbb D\to\Cal A.
\mytag{6.1}
$$
The mapping \mythetag{6.1} is bijective. For this mapping 
$$
\hskip -2em
p<q\text{ \ implies \ }A(p)\prec A(q).
\mytag{6.2}
$$
The property \mythetag{6.2} follows from the way in which the mapping $A$ 
is constructed. It is easy understand it if we write the numbers $p$ and
$q$ in the form brought to a common denominator:
$$
\xalignat 2
&p=\frac{k_1}{2^m},
&&q=\frac{k_2}{2^m}.
\endxalignat
$$
\mytheorem{6.1} The mapping \mythetag{6.1} can be extended up to
a bijective mapping $A\!:\Bbb R\to a$ from the set of real numbers
$\Bbb R$ to the line $a=OE$ preserving the property \mythetag{6.2}.
\endproclaim
\demo{Proof} First of all we prove the existence of a mapping 
that extends \mythetag{6.1}. We do it by constructing this mapping 
directly. Let $\xi\in\Bbb R$ be some real number. According to 
the results of \S\,5 for each real number $\xi$ there are two
sequences of its binary-rational approximations $a_m$ and $b_m$, 
which are fixed by this number uniquely. Applying the mapping
\mythetag{6.1} to the numbers $a_m$ and $b_m$, we get two
sequences of points $A_m=A(a_m)$ and $B_m=A(b_m)$. Due to the
lemma~\mythelemma{4.1} and due to the property \mythetag{6.2} 
of the mapping $A$ for these sequences we have 
$$
[A_{m+1}B_{m+1}]\subset [A_mB_m].
$$
Hence we can apply Cantor's axiom~\mytheaxiom{A19}. Due to this
axiom there is a point $X$ belonging to all segments $[A_mB_m]$.
\par
     The point $X$ belonging to all segments $[A_mB_m]$ simultaneously
is unique. Indeed, from $b_m-a_m=2^{-m}$ and by construction of the
points $A_m$ and $B_m$ we have 
$$
[A_mB_m]\cong\frac{1}{2^m}\cdot [OE].
$$
The existence of a second point $\tilde X$ on the line $a$ belonging 
to all segments $[A_mB_m]$ would mean
$$
[X\tilde X]<[A_mB_m]\cong\frac{1}{2^m}\cdot [OE],
$$
which contradicts the theorem~\mythetheorem{5.1}. The existence and 
uniqueness of the point $X$ determined by a real number $\xi$ through
its binary-rational approximations defines the required mapping 
$A\!:\Bbb R\to a$ if we set $A(\xi)=X$. If the number $\xi$ 
binary-rational, then the sequence $a_m$ stabilizes: $a_m=\xi$ for
$m>m_0$. Therefore, $A_m=X$ for $m>m_0$, and we get that the restriction
of $A\!:\Bbb R\to a$ to the set of binary-rational numbers coincides
with \mythetag{6.1}.\par
     Let $\xi$ and $\eta$ be two real numbers and let $\xi<\eta$. 
Let's denote through $a_m$, $b_m$, $\tilde a_m$, and $\tilde b_m$ their
binary-rational approximations. From the theorem~\mythetheorem{4.2}
we derive the existence of a positive integer number $m$ such that
$$
a_m<b_m<\tilde a_m<\tilde b_m.
$$
Let $A_m=A(a_m)$, $B_m=A(b_m)$, $\tilde A_m=A(\tilde a_m)$, and
$\tilde B_m=A(\tilde b_m)$. Let's use the property \mythetag{6.2}
of the mapping \mythetag{6.1} and get
$$
\hskip -2em
A_m\prec B_m\prec\tilde A_m\prec\tilde B_m.
\mytag{6.3}
$$
But the point $X=A(\xi)$ by construction belongs to the segment 
$[A_mB_m]$, while the point $Y=A(\eta)$ belongs to $[\tilde A_m
\tilde B_m]$. Therefore, we can sharpen the relationships 
\mythetag{6.3}:
$$
A_m\preccurlyeq A(\xi)\preccurlyeq B_m\prec\tilde A_m
\preccurlyeq A(\eta)\preccurlyeq\tilde B_m.
$$
From these relationships we already can derive the required 
formula $A(\xi)\prec A(\eta)$. Thus, for the mapping $A\!:\Bbb R
\to a$ the following condition is fulfilled:
$$
\hskip -2em
\xi<\eta\text{ \ implies \ }A(\xi)\prec A(\eta).
\mytag{6.4}
$$
As an immediate consequence of the condition \mythetag{6.4} we get
that the mapping $A\!:\Bbb R\to a$ is injective. In order to prove
its bijectivity now it is sufficient to prove its surjectivity.
We formulate and prove this fact as a separate theorem.
\qed\enddemo
\mytheorem{6.2} For any point $X$ on the line $a=OE$ there is a real
number $\xi$ such that $X=A(\xi)$.
\endproclaim
\demo{Proof} Let's consider the intervals $I(m,k)$ with the use of 
which in \S\,4 we have defined binary-rational approximations of
real numbers. Their boundaries are defined by the numbers
$$
\xalignat 2
&a_{mk}=\frac{k}{2^m},&&b_{mk}=\frac{k+1}{2^m}.
\endxalignat
$$
Let's consider the analogous semi-open intervals on the line $a$:
$$
\hskip -2em
\varTheta(m,k)=[A(a_{mk})A(b_{mk})).
\mytag{6.5}
$$
For each fixed $m$ the intervals with different numbers $k_1\neq k_2$ 
do not intersect, while their union covers the whole line $a$. This
fact follows from the Archimedes axiom~\mytheaxiom{A18} (compare
it with the theorem~\mythetheorem{4.1} for real numbers). Hence,
for each $m$ there exists exactly one interval $\varTheta(m,k)$ that
contains the point $X$. Let's denote by $k_m$ the number of this
interval:
$$
\xalignat 2
&a_m=\frac{k_m}{2^m},&&b_m=\frac{k_m+1}{2^m}.
\endxalignat
$$
The intervals \mythetag{6.5} are so that each $\varTheta(m,k)$ is the 
union of two intervals of the next order:
$$
\varTheta(m,k)=\varTheta(m+1,2k)\cup\varTheta(m+1,2k+1).
$$
Hence, $k_{m+1}=2k_m$ or $k_{m+1}=2k_m+1$. For the numbers $a_m$ and 
$b_m$ this fact means $a_m\leqslant a_{m+1}<b_{m+1}\leqslant b_m$. 
Therefore, Cantor's axiom~\mytheaxiom{R17} is applicable to the 
sequences $a_m$ and $b_m$. It yields the existence of a real number
$\xi$ satisfying the inequalities $a_m\leqslant\xi\leqslant b_m$. 
Such a number $\xi$ is unique since the existence of another number
$\tilde\xi$ satisfying the inequalities $a_m\leqslant\tilde\xi
\leqslant b_m$ would lead to the relationship that contradicts
the theorem~\mythetheorem{3.2}:
$$
|\xi-\tilde\xi|\leqslant b_m-a_m=\frac{1}{2^m}.
$$
The existence and uniqueness of a number $\xi$ given by a point $X$ 
on the line $a=OE$ means that we have constructed a mapping from
the line $a$ to the set of real numbers:
$$
\hskip -2em
\xi\!: a\to\Bbb R.
\mytag{6.6}
$$\par
     Let's prove that $A(\xi)=X$. For this purpose we sharpen the
inequalities $a_m\leqslant\xi\leqslant b_m$. The coincidence $\xi=b_m$ 
in these inequalities is excluded since due to $\xi\leqslant b_{m+1}
\leqslant b_m$ it would mean that the sequence $b_m$ stabilizes, i\.\,
e\. $b_m=b=\xi$ for $m>m_0$. For the point $X\in a$ this stabilization
would yield
$$
\hskip -2em
X\in [A(b-2^{-m})A(b))\text{ \ for all \ } m>m_0,
\mytag{6.7}
$$
which is impossible since the intersection of all semi-open intervals 
$[A(b-2^{-m})A(b))$ in \mythetag{6.7} is empty. Thus, for the number
$\xi$ and the binary-rational numbers $a_m$ and $b_m$ the inequalities
$a_m\leqslant\xi<b_m$ are fulfilled. Hence, $\xi\in I(m,k_m)$ and the
numbers $a_m$ and $b_m$ coincides exactly with the binary-rational 
approximations of the number $\xi$, which are used in constructing the
points $A(\xi)$. Now from $X\in [A(a_m)A(b_m)]$ for all $m\in\Bbb N$
we derive the required relationship $X=A(\xi)$.
\qed\enddemo
\myexercise{6.1} Show that the intersection of all semi-open intervals
$[A(b-2^{-m})A(b))$ from \mythetag{6.7} is an empty set. For this purpose
consider the segments $[A(b-2^{-m})A(b)]$ and, relying upon Cantor's
axiom, prove that the intersection of these segments consists of the 
single point $A(b)$.
\endproclaim
    In proving the theorem~\mythetheorem{6.2}, we not only have 
finished the proof of bijectivity of the mapping $A\!:\Bbb R\to a$ 
in the theorem~\mythetheorem{6.1}, but have constructed its inverse
mapping \mythetag{6.6}. In constructing both mappings $A$ and $\xi$
we used substantially the point $O$ and the vector 
$\overrightarrow{OE}$ on the line $a$.
\mydefinition{6.1} If a point $O$ and some vector 
$\bold e=\overrightarrow{OE}$ on a straight line $a$ are given,
the define a {\it Cartesian coordinate system\/} on this line. 
The point $O$ is called the {\it origin}, while the vector $\bold e$ 
is called the {\it basis vector}.
\enddefinition
\mydefinition{6.2} Let $(O,\bold e)$ be a Cartesian coordinate 
system on some line $a$. Then the number $\xi=\xi(X)$ is called 
the {\it coordinate\/} of a point $X\in a$, while the vector
$\overrightarrow{OX}$ connecting the origin with the point $X$ 
is called the {\it radius-vector\/} of this point.
\enddefinition
      If a point $X\in a$ is given, its coordinate is fixed uniquely
and, conversely, if the coordinate of a point $X$ is given, the point
$X\in a$ is fixed uniquely. This fact follows from the bijectivity 
of the mappings $A$ and $\xi=A^{-1}$.\par
    Using the mappings $A$ and $\xi$, we can identify real numbers 
with the points of some straight line. This is a visual image for 
understanding the set $\Bbb R$. For this reason the set of real 
numbers $\Bbb R$ itself is sometimes called the {\it real axis}.
\head
\SectionNum{7}{161} Measuring straight line segments.
\endhead
\rightheadtext{\S\,7. Measuring straight line segments.}
     Assume that a straight line segment $[OE]$ is given. Then we 
introduce the Cartesian coordinate system with the origin $O$ and
with the basis vector $\bold e=\overrightarrow{OE}$ on the line
$a=OE$. We use the line as a ruler and the segment $[OE]$ as a
gauge unit on this ruler. Let $[PQ]$ is some arbitrary segment. 
On the ray $[OE\rangle$ we choose a point $M$ such that $[OM]\cong
[PQ]$. Then the point is $M$ associated with the number $\xi(M)$
being its coordinate. From $O\prec M$ due to the relationship
\mythetag{6.4} we find that $\xi=\xi(M)$ is a positive number.
The positive number $\xi$ is called the {\it length\/} of the
segment $[PQ]$ measured relative to the gauge unit $[OE]$ on the
line $a$. This fact is written as follows:
$$
[PQ]\cong \xi\cdot [OE]\text{ \ or \ }|PQ|=\xi\cdot |OE|.
$$
If the reference segment $[OE]$ is fixed, the number $\xi$ can be 
taken for the length of the segment $[PQ]$. In this case we write 
$$
|PQ|=\xi.
$$
If the segments $[AB]$ and $[CD]$ are congruent, then the points 
$M$ and $\tilde M$ on the ray $[OE\rangle$ of the reference line 
$a$ determined by the segments $[AB]$ and $[CD]$ do coincide. For 
their lengths this yields $|AB|=|CD|$. Conversely, the equality 
of lengths $|AB|=|CD|$ imply $M=\tilde M$, $[OM]\cong [AB]$, and
$[OM]\cong [CD]$, which yields $[AB]\cong [CD]$. Let's formulate 
this result as a theorem.
\mytheorem{7.1} Two segments $[AB]$ and $[CD]$ are congruent if 
and only if their lengths measured with respect to the same reference
segment $[OE]$ are equal.
\endproclaim
\mytheorem{7.2} The relationship $[AB]<[CD]$ for two segments $[AB]$ 
and $[CD]$ is equivalent to the inequality $|AB|<|CD|$ for their
lengths measured with respect to some fixed reference segment $[OE]$.
\endproclaim
\demo{Proof} The procedure of measuring associates the segments $[AB]$ 
and $[CD]$ with two points $M$ and $\tilde M$ on the ray $[OB\rangle$
such that $[OM]\cong [AB]$ and $[O\tilde M]\cong [CD]$. The relationship
$[AB]<[CD]$ is equivalent to $[OM]<[O\tilde M]$ and to  
$$
\hskip -2em
O\prec M\prec\tilde M.
\mytag{7.1}
$$
Hence, due to the property \mythetag{6.4} of the mapping $A$ we get
$\xi(M)<\xi(\tilde M)$ \pagebreak which means $|AB|<|CD|$.\par
      Conversely, from $|AB|<|CD|$ we derive $\xi(M)<\xi(\tilde M)$,
which leads to the relationship \mythetag{7.1}. From \mythetag{7.1} 
and from $[OM]\cong [AB]$ with $[O\tilde M]\cong [CD]$ we get 
$[AB]<[CD]$.
\qed\enddemo
\mytheorem{7.3} Assume that a point $B$ lies in the interior of a 
segment $[AC]$. Then for the lengths of the segments $[AB]$, $[BC]$,
and $[AC]$ measured with respect to some fixed reference segment
$[OE]$ we have the equality $|AC|=|AB|+|BC|$.
\endproclaim
\demo{Proof} Let's consider the procedure of measuring the lengths 
of the segments $[AB]$ and $[AC]$ with respect to some reference 
segment $[OE]$. On the ray $[OE\rangle$ we mark the points $M$ and 
$\tilde M$ such that $[OM]\cong [AB]$ and $[O\tilde M]\cong [AC]$. 
Then the point $M$ lies in the interior of the segment $[O\tilde M]$, 
which follows from $[AB]<[AC]$. Moreover, $[M\tilde M]\cong [BC]$, 
which is due to the axiom~\mytheaxiom{A15}.\par
     Let $\xi=\xi(M)$ and $\tilde\xi=\xi(\tilde M)$. Then $\xi<\tilde\xi$
and for the lengths of the segments $[AB]$ and $[AC]$ we have
$$
\xalignat 2
&\hskip -2em
|AB|=\xi,&&|AC|=\tilde\xi.
\mytag{7.2}
\endxalignat
$$
The numbers $\xi$ and $\tilde\xi$ are associated with the sequences
$a_m$, $b_m$, $\tilde a_m$,and $\tilde b_m$ of their binary-rational
approximations. Applying the theorem~\mythetheorem{4.2}, we get 
$b_m<\tilde a_m$ for $m>m_0$. Hence, we write
$$
a_m\leqslant\xi<b_m<a_m\leqslant\tilde\xi<b_m.
$$
Now let's take into account that $A(\xi)=M$ and $A(\tilde\xi)=\tilde M$, 
then let's apply the property \mythetag{6.4} of the mapping $A$. This
yields 
$$
A(a_m)\preccurlyeq M\prec A(b_m)\prec A(\tilde a_m)\preccurlyeq
\tilde M\prec A(\tilde b_m).
$$
In other words, the segment $[M\tilde M]$ comprises the segment $[A(b_m)
A(\tilde a_m)]$ and is enclosed into the segment $[A(a_m)A(\tilde b_m)]$.
Let's apply the theorem~\mythetheorem{7.2} in this situation. It yields
$$
\hskip -2em
|A(b_m)A(\tilde a_m)|<|M\tilde M|<|A(a_m)A(\tilde b_m)|.
\mytag{7.3}
$$
Binary-rational numbers $a_m$, $b_m$, $\tilde a_m$, and $\tilde b_m$ are
determined by two integer numbers $k_m$ and $\tilde k_m$ (see formulas
\mythetag{4.1}). From 
$$
\xalignat 2
&b_m=\frac{k_m+1}{2^m},&&\tilde a_m=\frac{\tilde k_m}{2^m}
\endxalignat
$$
we conclude that the segment $[A(b_m)A(\tilde a_m)]$ is composed 
of congruent segments each of which can be obtained by mens of
the $m$-fold bisection of the segment $[OE]$. The number of such 
segments is equal to $\tilde k_m-k_m-1$. Therefore, we have
$$
|A(b_m)A(\tilde a_m)|=\frac{\tilde k_m-k_m-1}{2^m}=\tilde a_m-b_m.
$$
Similarly for the segment $[A(a_m)A(\tilde b_m)]$ we get
$$
|A(a_m)A(\tilde b_m)|=\frac{\tilde k_m-k_m+1}{2^m}=\tilde b_m-a_m.
$$
Now from \mythetag{7.3} we obtain the following estimates for the length
of the segment $[M\tilde M]$:
$$
\tilde a_m-b_m<|\tilde MM|<\tilde b_m-a_m.
$$
Let's remember that $[M\tilde M]\cong [BC]$ and use the inequalities
$a_m\leqslant\xi<b_m$ and $\tilde a_m\leqslant\tilde\xi<\tilde b_m$.
This yields
$$
\aligned
&\tilde\xi-\xi-|BC|<\tilde b_m-a_m-(\tilde a_m-b_m)=2^{m+1},\\
&\tilde\xi-\xi-|BC|>\tilde a_m-b_m-(\tilde b_m-a_m)=-2^{-m+1}.
\endaligned
\mytag{7.4}
$$
If we take into account \mythetag{7.2}, the inequalities \mythetag{7.4}
can be written as the following estimates:
$$
-2^{-m+1}<|AC|-|AB|-|BC|<2^{-m+1}.
$$
Since $m>m_0$ is an arbitrary integer greater than $m_0$, these estimates
yield $|AC|=|AB|+|BC|$. The theorem is proved.
\qed\enddemo
\mydefinition{7.1} The function associating each segment $[PQ]$ with 
some positive number $\xi([PQ])$ is called a {\it length function\/}, 
if the following conditions are fulfilled:
\roster
\item $\xi([OE])=1$ for some reference segment $[OE]$;
\item $[AB]\cong [CD]$ implies $\xi([AB])=\xi([CD])$;
\item if a point $B$ lies in the interior of the segment $[AC]$,
      then $\xi([AC])=\xi([AB])+\xi([BC])$.
\endroster
\enddefinition
\mytheorem{7.4} A length function satisfying the above conditions
\therosteritem{1}--\therosteritem{3} in the 
definition~\mythedefinition{7.1} does exist. It is unique if 
some reference segment $[OE]$ is fixed.
\endproclaim
\demo{Proof} Indeed, one of the length functions, satisfying the 
conditions \therosteritem{1}--\therosteritem{3} from the
definition~\mythedefinition{7.1} was constructed above. Let's
denote it $\xi$. Let $\eta$ be some other such function satisfying
the same conditions \therosteritem{1}-\therosteritem{3}. Let's show 
that these functions do coincide: $\xi([PQ])=\eta([PQ])$. The property
\therosteritem{2} means that it is sufficient to consider the segments
of the form $[OM]$, where $M$ is some point on the ray $[OE\rangle$
comprising the segment $[OE]$.\par
     If $M=E$, then the equality $\xi([OE])=\eta([OE])=1$ follows 
from the condition \therosteritem{1}. Assume that the segment $[AB]$ 
is obtained by $m$-fold bisection of the segment $[OE]$ and by 
subsequent composing $k$ replicas of the resulting segment. Then
$$
[AB]\cong\frac{k}{2^m}\cdot [OE].
$$
Using the condition  \therosteritem{3}, now one can derive the 
relationships
$$
\xalignat 2
&\xi([AB])=\frac{k}{2^m},&&\eta([AB])=\frac{k}{2^m}.
\endxalignat
$$
Thus, we have proved $\xi([AB])=\eta([AB])$ for a segment $[AB]$ 
being binary-rational multiple of the reference segment $[OE]$.\par
     The condition \therosteritem{3} in the
definition~\mythedefinition{7.1} has another important consequence.
It leads to $\eta([AB])<\eta([AC])$ for any two segments $[AB]$ and
$[AC]$ such that $[AB]<[AC]$. This property of the function $\eta$
is derived from the formula 
$$
\eta([AC])=\eta([AB])+\eta([BC])
$$
since $\eta([BC])$ in this formula is positive. Let's use this
property in order to prove $\xi([OM])=\eta([OM])$ for a point
$M$ on the ray $[OE\rangle$. According to the
theorem~\mythetheorem{6.2}, we  have $M=A(\alpha)$, where $\alpha$ 
is some positive real number. Let's denote by $a_m$ and $b_m$ the
binary-rational approximations of this number. Then 
$a_m\leqslant\alpha<b_m$ and we have 
$$
A(a_m)\preccurlyeq M\prec A(b_m).
$$
From the inclusions $[OA(a_m)]\subset [OM]\subset [OA(b_m)]$ we derive
$$
\eta([OA(a_m)])\leqslant\eta([OM])\leqslant\eta([OA(b_m)]).
$$
But $\eta([OA(a_m)])=a_m$ and $\eta([OA(b_m)])=b_m$ since we have already
proved that $\eta$ and $\xi$ do coincide for segments being binary-rational
multiples of the segment $[OE]$. Hence, we get
$$
\alpha-\frac{1}{2^m}<a_m\leqslant\eta([OM])\leqslant b_m
\leqslant\alpha+\frac{1}{2^m}.
$$
These inequalities yield $\eta([OM])=\alpha=\xi([OM])$ since $m\in\Bbb N$
is an arbitrary positive integer. The theorem is proved.
\qed\enddemo
\mytheorem{7.5} If the length of a segment $[AB]$ measured with 
respect to a reference segment $[OE]$ is equal to $\xi$ and if
the length of the the reference segment $[OE]$ measured with 
respect to another reference segment $[\tilde O\tilde E]$ is 
equal to $\eta$, then the length of the segment $[AB]$ measured 
with respect to the second reference segment $[\tilde O\tilde E]$ 
is equal to the product $\tilde\xi=\xi\cdot\eta$.
\endproclaim
\myexercise{7.1} Derive the theorem~\mythetheorem{7.5} from the
theorem~\mythetheorem{7.4}.
\endproclaim
\head
\SectionNum{8}{167} Similarity mappings for straight lines. 
Multiplication of vectors by a number.
\endhead
\rightheadtext{\S\,8. Similarity mappings \dots}
     Assume that on a straight line $a$ a Cartesian coordinate system 
with the origin $O$ and with the basis vector $\bold e=\overrightarrow{OE}$
is given. Let's consider another straight line $b$ on which another 
Cartesian coordinate system with the origin $Q$ and with the basis vector
$\bold h=\overrightarrow{QH}$ is given. We define a mapping $f\!:a\to b$
as follows. For a point $X\in a$ we take its coordinate $\xi=\xi(X)$ and
then we associate to $X$ the point $Y\in b$ with exactly the same
coordinate $\xi$ . In the other words, $f$ is the composition of two
mappings
$$
\xalignat 2
&\xi\!:a\to\Bbb R,&&A\!:\Bbb R\to b.
\endxalignat
$$
Such a mapping $f=A\compos\xi$ is called a {\it similarity mapping}.
The ratio of the lengths of reference vectors 
$$
k=|QH|/|OE|
$$
is called the {\it similarity factor\/} for this mapping. Any similarity
mapping $f\!:a\to b$ is bijective. The inverse mapping $f^{-1}\!:b\to a$
is also a similarity mapping, its similarity factor is $k^{-1}$. The
composition of two similarity mappings $f\!:a\to b$ and $g\!:b\to c$ is
a similarity mapping $g\compos f\!:a\to c$. Its similarity factor is
equal to the product of similarity factors of the mappings $f$ and $g$.
\mylemma{8.1} Assume that on a line $a$ a Cartesian coordinate system
with the origin $O$ and the basis vector $\bold e=\overrightarrow{OE}$
is given. Then the following propositions are valid:
\roster
\rosteritemwd=0pt
\item the vector $\overrightarrow{AB}$ is codirected with the 
      vector $\overrightarrow{OE}$ if and only if for the points
      $B$ and $A$ the difference of their coordinates is positive:
      $\xi(B)-\xi(A)>0$;
\item the length of the segment $[AB]$ measured with respect to the
      reference segment $[OE]$ is given by the formula
$$
|AB|=|\xi(B)-\xi(A)|.
$$
\endroster
\endproclaim
\myexercise{8.1} Consider all possible dispositions of the points 
$A$ and $B$ relative to the point $O$ on the line $a$ and, using
the theorem~\mythetheorem{7.3}, prove the lemma~\mythelemma{8.1}.
\endproclaim
\mytheorem{8.1} Assume that $a$ and $b$ are two straight lines
with the vectors $\overrightarrow{OE}$ and $\overrightarrow{QH}$ 
on them. Let's define the positive directions on the lines $a$ 
and $b$ by means of the vectors $\overrightarrow{OE}$ and
$\overrightarrow{QH}$ and then consider the similarity mapping 
$f\!:a\to b$ defined by them. This mapping
\roster
\rosteritemwd=0pt
\item preserves the relation of precedence for points, i\.\,e\. 
      $X\prec Y$ implies $f(X)\prec f(Y)$;
\item multiplies the lengths of segments by $k$, where $k$ 
      is the similarity factor: $|f(X)f(Y)|=k\cdot |XY|$.
\endroster
\endproclaim
\mytheorem{8.2} Let $f\!:a\to b$ be a similarity mapping. Then for
some arbitrary points $A$, $B$, $C$, and $D$ on the line $a$ and for
their images $A'=f(A)$, $B'=f(B)$, $C'=f(C)$, and $D'=f(D)$ on the
line $b$ the equality $\overrightarrow{AB}=\overrightarrow{CD}$ 
implies the equality $\overrightarrow{A'B'}=\overrightarrow{C'D'}$.
\endproclaim
\myexercise{8.2} Using the theorem~\mythetheorem{7.5} and the
lemma~\mythelemma{8.1}, pro\-ve the theorem~\mythetheorem{8.1} and 
derive the theorem~\mythetheorem{8.2} from it.
\endproclaim
    A special sort of similarity mappings arise if the lines $a$ and
$b$ do coincide: $a=b$. Assume that on a straight line $a$ two vectors
$\overrightarrow{OE}$ and $\overrightarrow{OH}$ are given. Then we
have two Cartesian coordinate systems with the common origin $O$ on 
this line. They define a similarity mapping $f\!:a\to a$. Such a mapping
is called a {\it homothety\/} with the center $O$. The mapping $f$
takes a point $X$ with the coordinate $\xi(X)$ in the first coordinate
system to the point $Y=f(X)$ with the coordinate $\tilde\xi(Y)=
\xi(X)$ in the second coordinate system. Using the 
theorem~\mythetheorem{8.1}, it is easy to calculate the coordinate 
of the point $Y$ in the first coordinate system:
$$
\xi(Y)=\cases
\hphantom{-}(|OH|/|OE|)\cdot\xi(X)&\text{for \ }\overrightarrow{OH}
\upuparrows\overrightarrow{OE},\\
-(|OH|/|OE|)\cdot\xi(X)&\text{for \ }\overrightarrow{OH}
\updownarrows\overrightarrow{OE}.
\endcases
$$
For the homothety mapping $f$ we introduce the numeric parameter, 
which is called the {\it homothety factor}:
$$
k=\cases
\hphantom{-}|OH|/|OE|&\text{if \ \ }\overrightarrow{OH}
\upuparrows\overrightarrow{OE},\\
-|OH|/|OE|&\text{if \ \ }\overrightarrow{OH}
\updownarrows\overrightarrow{OE}.
\endcases
$$
Then the homothety mapping $f\!:a\to a$ can be defined using only one
Cartesian coordinate system: a point $X$ with the coordinate $\xi$ is
taken to the point $Y=f(X)$ with the coordinate $k\cdot\xi$.\par
\mydefinition{8.1} Assume that $\overrightarrow{AB}$ is a vector on
some straight line $a$. A vector $\overrightarrow{CD}$ with the 
length $|CD|=|k|\cdot |AB|$ on the line $a$ is called the {\it product
of the vector $\overrightarrow{AB}$ by a number $k\neq 0$} if it is 
codirected to $\overrightarrow{AB}$ for $k>0$ and is oppositely
directed to $\overrightarrow{AB}$ for $k<0$. We express this fact by
writing $\overrightarrow{CD}=k\cdot\overrightarrow{AB}$.
\enddefinition
\mytheorem{8.3} A point $Y$ on a straight line $a$ is the image of a 
point $X$ on the same line under the homothety mapping with the center 
at the point $O$ and with the homothety factor $k$ if and only if its
radius-vector $\overrightarrow{OY}$ is obtained from the radius vector
$\overrightarrow{OX}$ through multiplying it by the number $k$.
\endproclaim
\myexercise{8.3} Prove the theorem~\mythetheorem{8.3} by choosing
some Cartesian coordinate system in the line $a$.
\endproclaim
     The definition~\mythedefinition{8.1} does not fix a definite
position of the vector $\overrightarrow{CD}=k\cdot\overrightarrow{AB}$ 
on a line, the vector $\overrightarrow{CD}$ is determined up to the
replacement of it by any other vector equal to it in the sense of the
definition~\mythedefinitionchapter{4.1}{3} from 
Chapter~\uppercase\expandafter{\romannumeral 3}. If
$\overrightarrow{GF}=\overrightarrow{AB}$, then 
$k\cdot\overrightarrow{GF}=k\cdot\overrightarrow{AB}$. This means that
in the definition~\mythedefinition{8.1} the product 
$\overrightarrow{CD}=k\cdot\overrightarrow{AB}$ is defined only 
as a slipping vector $\bold c=\overrightarrow{CD}$ obtained through 
multiplying another slipping vector $\bold b=\overrightarrow{AB}$ 
by the number $k$. If $k=0$ or if $\bold b=\bold 0$, the product 
$k\cdot\bold b$ is assumed to be equal to zero vector:
$$
\align
&0\cdot\bold b=\bold 0\text{ \ for any vector  \ }\bold b,\\
\displaybreak
&k\cdot\bold 0=\bold 0\text{ \ for any number \ }k\in\Bbb R.
\endalign
$$
\mytheorem{8.4} The operations of addition and multiplication by a 
number for slipping vectors possess the following properties:
\roster
\rosteritemwd=5pt
\item commutativity of addition: $\bold a+\bold b=\bold b+\bold a$;
\item associativity of addition: $(\bold a+\bold b)+\bold c=
      \bold a+(\bold b+\bold c)$;
\item there is a zero vector $\bold 0$ such that $\bold a+\bold 0
      =\bold a$ for an arbitrary vector $\bold a$;
\item for any vector $\bold a$ there is an opposite vector $\bold a'$
      such that $\bold a+\bold a'=\bold 0$;
\item distributivity of multiplication by a number with respect 
      to addition of vectors: 
      $k\cdot (\bold a+\bold b)=k\cdot\bold a+k\cdot\bold b$;
\item distributivity of multiplication by a number with respect to 
      addition of numbers: $(k+q)\cdot\bold a=k\cdot\bold a+q\cdot
      \bold a$;
\item associativity of multiplication: 
      $(k\cdot q)\cdot\bold a=k\cdot(q\cdot\bold a)$;
\item the property of the numeric unity: $1\cdot\bold a=\bold a$.
\endroster
\endproclaim
\myexercise{8.4} Prove those propositions of the 
theorem~\mythetheorem{8.4} which are not yet proved.
\endproclaim
\head
\SectionNum{9}{170} Measuring angles.
\endhead
\rightheadtext{\S\,9. Measuring angles.}
     The numeric measure of angles is constructed approximately in the 
same way as the length of segments. The only difference is that here  
the propositions like the Archimedes axiom~\mytheaxiom{A18} and Cantor's
axiom~\mytheaxiom{A19} should be proved.
\mytheorem{9.1} Let $\{\angle h_nk_n\}_{\sssize n\in\Bbb N}$ be s sequence
of angles lying on one plane and having a common vertex $O$. Assume that 
the following relationships are fulfilled:
$$
\angle h_{n+1}k_{n+1}\subset\angle h_nk_n\text{ \ for all \ }
n\in\Bbb N.
$$
Then there is a ray $l$ coming out from the point $O$ and belonging to the
intersection of all these angles.
\endproclaim
\demo{Proof}\parshape 18 0cm 10cm 0cm 10cm 0cm 10cm 0cm 10cm 5cm
5cm 5cm 5cm 5cm 5cm 5cm 5cm 5cm 5cm 5cm 5cm 5cm 5cm 5cm 5cm 5cm 5cm 5cm 5cm
5cm 5cm 5cm 5cm 5cm 5cm 0cm 10cm
Let's mark a point $A_1$ on the ray $h_1$ and a point $B_1$ on the ray
$k_1$. Then we connect them with the segment $[A_1B_1]$. According to 
the lemma~\mythelemmachapter{6.2}{2} from
Chapter~\uppercase\expandafter{\romannumeral 2}, the rays $h_n$ and 
$k_n$ intersect the segment $[A_1B_1]$. We denote the intersection 
points through $A_n$ and $B_n$ respectively. Such points form a
sequence of segments such that
$$
[A_{n+1}B_{n+1}]\subset [A_nB_n].
$$
In this situation Cantor's axiom~\mytheaxiom{A19} is applicable. 
\vadjust{\vskip 5pt\hbox to 0pt{\kern 0pt
\includegraphics{Oris54.eps}\hss}\vskip -5pt}It yields the
existence of a point $X$ belonging to all segments $[A_nB_n]$. We draw 
the ray $[OX\rangle$ trough this point $X$ and denote it $l$. According 
to the lemma~\mythelemmachapter{6.2}{2} from
Chapter~\uppercase\expandafter{\romannumeral 2}, this ray lies
in the intersection of all angles $\angle h_nk_n$. The theorem is 
proved.\qed\enddemo
\parshape 14 4cm 6cm 4cm 6cm 4cm 6cm 4cm 6cm
4cm 6cm 4cm 6cm 4cm 6cm 4cm 6cm 4cm 6cm 4cm 6cm
4cm 6cm 4cm 6cm 4cm 6cm 0cm 10cm
     Let's consider an angle $\angle h_0h_n$ with the vertex at some point
$O$. Assume that inside this angle $\angle h_0h_n$ the rays 
$h_1,\,\ldots,\,h_{n-1}$ are drawn so that they form the angles
$$
\angle h_0h_1,\,\ldots,\,\angle h_{n-1}h_n
$$
congruent to each other. \vadjust{\vskip 5pt\hbox to 0pt{\kern 0pt
\includegraphics{Oris55.eps}\hss}\vskip -5pt}In this case we 
say that the angle $\angle h_0h_n$ is $n$ times as bigger than the angle
$\angle h_0h_1$:
$$
\angle h_0h_n\cong n\cdot\angle h_0h_1.
$$
Conversely, for the angle $\angle h_0h_1$ we say that it is obtained from
$\angle h_0h_n$ by dividing into $n$ congruent parts. We write this fact
as
$$
\angle h_0h_1\cong\frac{1}{n}\cdot\angle h_0h_n.
$$
\mylemma{9.1} Assume that in a triangle $ABC$ the bisector $AD$ is drawn. 
Then $\angle ABC>\angle ACB$ implies $[CD]>[BD]$.
\endproclaim
\demo{Proof}\parshape 16 0cm 10cm 0cm 10cm 3cm 7cm 3cm 7cm
3cm 7cm 3cm 7cm 3cm 7cm 3cm 7cm 3cm 7cm 3cm 7cm 3cm 7cm
3cm 7cm 3cm 7cm 3cm 7cm 0cm 10cm 0cm 5.8cm 
     From the relationship $\angle ABC>\angle ACB$, applying the 
theorem~\mythetheorem{2.4}, we get $[AC]>[AB]$. On the ray $[AC\rangle$
we mark a point $E$ so that $[AE]\cong [AB]$. 
\vadjust{\vskip 5pt\hbox to 0pt{\kern -10pt
\includegraphics{Oris56.eps}\hss}\vskip -5pt}Due to 
$[AC]>[AB]$ the point $E$ is an interior point of the segment $[AC]$. 
From the congruence $\angle EAD\cong\angle BAD$ and from
$[AE]\cong [AB]$ we derive that the triangles $EAD$ and $BAD$ are
congruent. Hence, $[DE]\cong [DB]$. The angle $\angle CED$ is adjacent 
to the angle $\angle AED$ which is congruent to $\angle ABC$. Therefore 
the angle $\angle CED$ is congruent to the external angle of the triangle
$ABC$ at the vertex $B$. According to the theorem~\mythetheorem{2.3} the
internal angle of this triangle at the vertex $C$ is smaller than its
external angle at the vertex $B$. 
\vadjust{\vskip 5pt\hbox to 0pt{\kern 175pt
\includegraphics{Oris57.eps}\hss}\vskip -5pt}Hence, we get
$\angle CED>\angle ECD$. Now, applying the theorem~\mythetheorem{2.4} 
to the triangle $CED$, we get $[CD]>[ED]$ which is equivalent to the
relationship $[CD]>[BD]$ due to $[ED]\cong [BD]$. The lemma is proved.
\qed\enddemo
\mytheorem{9.2}\parshape 1 0cm 5.8cm For any two acute angles $\angle hk$ 
and $\angle lq$ there is a positive integer $n$ such that
$$
\angle hk<\angle rp \cong n\cdot\angle lq,
$$
where the angle $\angle rp$ is $n$ times as bigger than the angle $\angle lq$.
\endproclaim
\demo{Proof}\parshape 5 0cm 5.8cm 0cm 5.8cm 0cm 5.8cm 0cm 5.8cm 0cm 10cm 
If the angle $\angle lq$ is bigger th$\angle hk$, then by choosing
$n=1$ we provide the required relationship $\angle hk<n\cdot\angle lq$. 
If $\angle lq\cong\angle hk$, it is sufficient take $n=2$. Therefore,
we consider the case, where the angle $\angle lq$ is smaller than the
angle $\angle hk$. In this case we denote by $A$ the vertex of the 
angle $\angle hk$, mark some point $C$ on the ray $k$ and drop the
perpendicular from the point $C$ onto the line containing the ray
$h$. The foot of this perpendicular lies on the ray $h$ since the
angle $\angle hk$ is acute. We denote it through $B$ (see 
Fig\.~9.4 above).\par
     The further proof of the theorem is by contradiction. We denote
$h=h_0$ and in the half-plane containing the ray $k$ we draw the series 
of rays $h_1,\,h_2,\,\ldots,\,h_n$ so that all angles of the form 
$\angle h_sh_{s+1}$ are congruent to the angle $\angle lq$. If we 
assume that the relationship $\angle hk<n\cdot\angle lq$ is not 
fulfilled for all positive integers $n\in\Bbb N$, we would be able to 
draw infinitely many rays $h_1,\,h_2,\,\ldots,\,h_n,\ldots$, all of 
them lying inside the angle $\angle hk$. Let's denote by $B_1,\,B_2,
\,\ldots,\,B_n,\,\ldots$ the intersection point of these rays and the 
segment $[BC]$ (such points do exist due to the 
lemma~\mythelemmachapter{6.2}{2} from
Chapter~\uppercase\expandafter{\romannumeral 2}). Let's consider the
triangle $AB_{s-1}B_{s+1}$. For $s>1$ the angle $\angle AB_{s-1}B_{s+1}$ 
in this triangle is adjacent to the acute angle $\angle AB_{s-1}B$ in the
rectangular triangle $ABB_{s-1}$. Therefore, it is an obtuse angle. For
$s=1$ the angle $\angle AB_{s-1}B_{s+1}$ is a right angle. In each of these
two cases we have 
$$
\angle AB_{s-1}B_{s+1}>\angle AB_{s+1}B_{s-1}.
$$
The ray $[AB_s\rangle$ is the bisector of the angle 
$\angle B_{s-1}AB_{s+1}$. Therefore we can apply the 
lemma~\mythelemma{9.1}, which yields
$$
[B_{s-1}B_s]<[B_sB_{s+1}].
$$
This means that the lengths of the segments $[B_sB_{s+1}]$ form a
monotonic increasing sequence of numbers. For the segment $[BB_n]$ 
this yields the relationship $[BB_n]>n\cdot [BB_1]$. Now the
assumption that the relationship $\angle hk<n\cdot\angle lq$ is 
not fulfilled for all positive integers $n$ leads to
$n\cdot [BB_1]<[BC]$ for all $n\in\Bbb N$. But this contradicts
the Archimedes axiom~\mytheaxiom{A18}. The contradiction obtained
completes the proof of the theorem.
\qed\enddemo
     The construction of a gauge for measuring angles does not differ
from that for segments on a straight line. The natural restriction
that all angles are enclosed into a straight angle, defines the natural
choice of the reference angle. In order to have an acute reference 
angle we take some straight angle, bisect it twice, and then assign the
value of $\pi/4$ radians to the resulting angle. Here $\pi=3.14\ldots$
is the well-known irrational number arising as the area of a unit circle.
\mytheorem{9.3} Each angle $\angle hk$ is associated with some real
number $\xi(\angle hk)=\widehat{hk}$ from the interval  $0<\xi\leqslant\pi$ 
so that the following conditions are fulfilled:
\roster
\item a straight angle is associated with the number $\pi$;
\item if $\angle hk\cong\angle lq$, then $\xi(\angle hk)=\xi(\angle lq)$;
\item if a ray $l$ lies inside an angle $\angle hk$ and divides it into
      two angles, then $\xi(\angle hk)=\xi(\angle hl)+\xi(\angle lk)$.
\endroster
\endproclaim
\myexercise{9.1} Using the analogy to segments, completes the details 
required for proving the theorem~\mythetheorem{9.3} and prove it.
\endproclaim
\newpage
\setfirstpage
\topmatter
\title\chapter{6}
Axiom of parallels.
\endtitle
\leftheadtext{CHAPTER~\uppercase\expandafter{\romannumeral 6}. 
THE AXIOM OF PARALLELS.}
\endtopmatter
\document
\chapternum=6
\head
\SectionNum{1}{175} The axiom of parallels 
and the classical Euclidean geometry.
\endhead
\rightheadtext{\S\,1. The axiom of parallels \dots}
     The fifth group of Euclid's axioms consists of one
axiom~\mytheaxiom{A20}. It is formulated as follows.
\myaxiom{A20} For any point $O$ not lying on a straight line $a$
there is exactly one straight line passing through the point $O$
and being parallel to the line $a$.
\endproclaim
\noindent The axiom~\mytheaxiom{A20} played an important role 
in the history of science. Multiple attempts to prove it by deriving
from the other axioms lasted more than 2000 years. However, they did 
not succeed. To the contrary, giving up the idea to prove it, people
had discovered new non-Euclidean geometries. For the firs time this
was done by Lobachevsky, Bolyai, and Gauss.\vadjust{\vskip 5pt
\hbox to 0pt{\kern -5pt
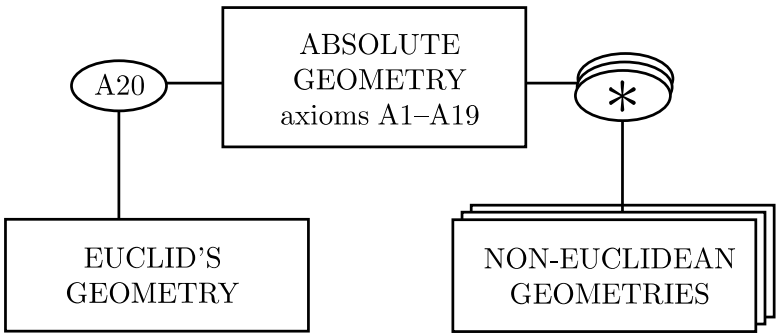
\includegraphics{Oris58.eps}\hss}\vskip -5pt}\par
\vskip 130pt\vphantom{.}
     The axioms~A1--A19 and their consequences considered in 
Chapters~\uppercase\expandafter{\romannumeral 1}--\uppercase
\expandafter{\romannumeral 5} constitute the so called 
{\tencyr\char '074}absolute geometry{\tencyr\char '076}. They 
are valid in classical Euclidean geometry and they remain valid 
in its non-Euclidean variations, where the axiom~\mytheaxiom{A20} 
is replaced by some propositions not equivalent to it. To the 
contrary, in Chapter~\uppercase\expandafter{\romannumeral 6} we 
consider those results which are specific to Euclidean geometry 
only. Non-Euclidean geometries are beyond the scope of this book.
\mytheorem{1.1} Assume that $a$ and $b$ are two straight lines
lying on one plane and intersecting with a third straight line
$c$ at the points $A$ and $B$. The lines $a$ and $b$ are parallel
if and only if the inner crosswise lying angles at the points $A$ 
and $B$ are congruent.
\endproclaim
\demo{Proof}\parshape 10 4cm 6cm 4cm 6cm 4cm 6cm 4cm 6cm 4cm 6cm
4cm 6cm 4cm 6cm 4cm 6cm 4cm 6cm 0cm 10cm
     The direct proposition of the theorem is already proved (see 
theorem~\mythetheoremchapter{8.1}{3} in
Chapter~\uppercase\expandafter{\romannumeral 3}).
\vadjust{\vskip 5pt\hbox to 0pt{\kern 0pt
\includegraphics{Oris59.eps}\hss}\vskip -5pt}Let's prove 
the converse proposition. Assume that $a\parallel b$. We mark some 
point $D\neq B$ on the line $b$. The line $c$ divides the plane of 
the lines $a$ and $b$ into two half-planes. On that half-plane which 
contains the point $D$ we draw the ray $[B\tilde D\rangle$ so that 
$\angle AB\tilde D\cong\angle BAC$. Applying  the
theorem~\mythetheoremchapter{8.1}{3} from
Chapter~\uppercase\expandafter{\romannumeral 3} to  the line 
$\tilde b=B\tilde D$, we get $\tilde b\parallel a$. If the line 
$\tilde b$ would be different from $b$, we would have two lines
passing through the point $B$ and being parallel to the line $a$, 
which contradicts the axiom~\mytheaxiom{A20}. Hence, $\tilde b=b$ and 
$\angle ABD\cong\angle BAC$. The theorem~\mythetheorem{1.1} is proved.
\qed\enddemo
\mytheorem{1.2} Let $a\neq b$ be two parallel lines lying on a plane
$\alpha$. If a line $c\neq b$ lying on the plane $\alpha$ intersect
the line $b$ at some point $B$, then it intersects the line $a$ at
some other point $A$.
\endproclaim
\demo{Proof} If we assume that $c$ does not intersect the line $a$, 
then, according to the definition~\mythedefinitionchapter{8.1}{3}
from Chapter~\uppercase\expandafter{\romannumeral 3}, these lines 
are parallel: $c\parallel a$. As a result we get two lines $b$ and
$c$, passing through the point $B$ and parallel to the line $a$, 
which contradicts the axiom~\mytheaxiom{A20}. Hence, the line $c$ 
intersects the line $a$ at some point $A\neq B$.
\qed\enddemo
\mytheorem{1.3} Assume that $a$ and $b$ are two parallel lines lying
on a p[lane $\alpha$. Then a perpendicular to the line $a$ drawn on 
the plane $\alpha$ is a perpendicular to the line $b$ either.
\endproclaim
    The theorem~\mythetheorem{1.3} is a simple consequence of the
theorems~\mythetheorem{1.1} and \mythetheorem{1.2}. It does not
require a separate proof.
\mytheorem{1.4} In Euclidean geometry the relation of parallelism 
of straight lines is reflexive, symmetric, and transitive, because 
of which it is an equivalence relation.
\endproclaim
\demo{Proof} \parshape 3 0cm 10cm 0cm 10cm 5cm 5cm
     Reflexivity and symmetry of parallelism of straight lines follows
immediately from its definition (see 
definition~\mythedefinitionchapter{8.1}{3} in
Chapter~\uppercase\expandafter{\romannumeral 3}). 
\vadjust{\vskip 5pt\hbox to 0pt{\kern -5pt
\includegraphics{Oris60.eps}\hss}\vskip -5pt}Let's prove its 
transitivity. Assume that $a\parallel b$ and $b\parallel c$. If $a=b$,
$b=c$, or $c=a$, the relationship $a\parallel c$ is obviously fulfilled.
Therefore, we consider the case, where the lines $a$, $b$, and $c$ are
mutually distinct. From $a\parallel b$ we conclude that there is a plane
containing the lines $a$ and $b$. We denote this plane by $\alpha$.
Similarly, we denote by $\beta$ the plane containing the lines $b$ and
$c$.\par
\parshape 2 5cm 5cm 0cm 10cm
     Let's choose a point $B$ on the line $b$ and draw the plane $\gamma$
perpendicular to the line $b$ through this point. At the intersection
of the planes $\alpha$ and $\gamma$ we get the perpendicular to the line
$b$. According to the theorem~\mythetheorem{1.2}, this perpendicular 
crosses the line $a$ at some point $A$. According to the 
theorem~\mythetheorem{1.3} the line $AB$ is perpendicular to the line
$a$. By construction the line $b$ is a perpendicular to the plane
$\gamma$, while the plane $\alpha$ contains this line $b$. Hence, 
$\alpha\perp\gamma$ (see definition~\mythedefinitionchapter{3.1}{4} in
Chapter~\uppercase\expandafter{\romannumeral 4}). Due to the
theorem~\mythetheoremchapter{3.1}{4} from 
Chapter~\uppercase\expandafter{\romannumeral 4} the plane contains
the perpendicular to the plane $\gamma$ passing through the point $A$.
This perpendicular coincides with the line $a$ due to the
theorem~\mythetheoremchapter{6.3}{3} from
Chapter~\uppercase\expandafter{\romannumeral 3} since $a\perp AB$.
\par
     Thus, $a\perp\gamma$. A similar situation arises on the plane
$\beta$. The line $BC$ obtained as the intersection of the planes 
$\gamma$ and $\beta$ appears to be perpendicular to the lines $b$ and
$c$. Then from $b\perp\gamma$ we derive $\beta\perp\gamma$ and 
$c\perp\gamma$. But any two perpendiculars to the same plane are
parallel (see theorem~\mythetheoremchapter{3.3}{4} in
Chapter~\uppercase\expandafter{\romannumeral 4}). Therefore from
$a\perp\gamma$ and $c\perp\gamma$ we derive the required relationship
$a\parallel c$.\par
     Note that there is a special case, where the planes $\alpha$ and 
$\beta$ do coincide. In this case the line $AB$ coincides with the line
$BC$. However, the above considerations remain valid for this case too.
The theorem is proved.
\qed\enddemo
\head
\SectionNum{2}{178} Parallelism of a straight line and a plane.
\endhead
\rightheadtext{\S\,2. Parallelism of a line and a plane.}
\mydefinition{2.1} A straight line $a$ is said to be {\it parallel\/} 
to a plane $\alpha$ if it lies on the plane $\alpha$ or if it does
not intersect this plane. The parallelism of $a$ and $\alpha$ is
written as $a\parallel\alpha$.
\enddefinition
\mytheorem{2.1} A straight line $a$ is parallel to a plane $\alpha$ 
if and only if it is parallel to some straight line $b$ lying on the
plane $\alpha$.
\endproclaim
\demo{Proof} The case where the line $a$ lies on the plane $\alpha$
is trivial. In tis case we can choose $b=a$, upon which both propositions 
of the theorem (direct and inverse ones) appear to be obviously valid.
\par
     Let's consider the case where the line $a$ does not lies on the 
plane $\alpha$. Assume that $a\parallel\alpha$. Let's choose some
arbitrary point $A$ on the plane $\alpha$. Then let's draw a plane
$\beta$ through the line $a$ and the point $A$. At the intersection 
of the planes $\beta$ and $\alpha$ we get a straight line $b$. The
line $b$ does not intersect the line $a$ since it lies on the plane
$\alpha$ that has no common points with the line $a$. On the other
hand, the line $b$ lies on the same plane $\beta$ as the line $a$.
Hence, we have $b\parallel a$.\par
     Conversely, assume that the line $a$ is parallel to some line
$b$ lying on the plane $\alpha$. Let's draw a plane $\beta$ through
these two parallel lines. Then $\alpha\cap\beta=b$. If the line
$a$ would intersect the plane $\alpha$, the intersection point would
lie on the line $b$. But the lines $a$ and $b$ have no common points
since they are parallel and do not coincide. Hence, $a\parallel\alpha$.
\qed\enddemo
      Note that the definition~\mythedefinition{2.1} and the
theorem~\mythetheorem{2.1} can be formulated in absolute geometry
either. The proof of the theorem~\mythetheorem{2.1} does not use the
axiom~\mytheaxiom{A20}. However, all other theorems below in this 
section are valid only in Euclidean geometry.
\mytheorem{2.2} For two straight lines $a$ and $b$ and for a plane 
$\alpha$ the conditions $a\parallel b$ and $b\parallel\alpha$ imply
$a\parallel\alpha$.
\endproclaim
\myexercise{2.1} Derive the theorem~\mythetheorem{2.2} as a consequence
of the theorems~\mythetheorem{1.4} and \mythetheorem{2.1}.
\endproclaim
\mytheorem{2.3} Let $a$ be some line parallel to a plane $\alpha$ and 
let $O$ be some point on this plane. If a line $b$ passes through the 
point $O$ and if $b\parallel a$, then the line $b$ lies on the plane
$\alpha$.
\endproclaim
\demo{Proof} Let's begin with the case where the line $a$ lies on the
plane $\alpha$. If $O\in a$, then the lines $a$ and $b$ have the common
point $O$. In this case $b\parallel a$ implies $b=a$ (see 
definition~\mythedefinitionchapter{8.1}{3} in 
Chapter~\uppercase\expandafter{\romannumeral 3}). Hence, $b\subset\alpha$.
\par
    If $a\subset\alpha$, but $O\notin a$, the lines $a$ and $b$ 
do not coincide. Two parallel lines $a$ and $b$, according to the 
definition~\mythedefinitionchapter{8.1}{3} in 
Chapter~\uppercase\expandafter{\romannumeral 3}, lie on one plane.
Let's denote this plane by $\beta$. The planes $\alpha$ and 
$\beta$ both contain the line $a$ and the point $O$ not lying 
on this line. Therefore they coincide: $\beta=\alpha$. Hence,
$b\subset\alpha$.\par
    And finally, let's consider the case where the line $a$ does not
lie on the plane $\alpha$. Since $a\parallel\alpha$, it does not 
intersect the plane $\alpha$. The line $b$ has the common point
$O$ with the plane $\alpha$. Hence, $a\neq b$. Since $b\parallel a$,
there is a plane $\beta$ passing through the lines $a$ and $b$ 
(see definition~\mythedefinitionchapter{8.1}{3} in 
Chapter~\uppercase\expandafter{\romannumeral 3}). In this case the 
plane $\beta$ is different from the plane $\alpha$. These two planes
$\alpha$ and $\beta$ have the common point $O$. Therefore, they
intersect along some line $\tilde b$ that contains the point $O$ and
has no common points with the line $a$. Hence, $\tilde b\parallel a$
since non-intersecting straight lines $a$ and $\tilde b$ lie on one
plane $\beta$. For the last step in our proof we apply the 
axiom~\mytheaxiom{A20} which says that there is a unique line $b$
passing through the point $O$ and being parallel to $a$. Hence,
$b=\tilde b$ and $b\subset\alpha$ since $\tilde b\subset\alpha$ by
construction.\qed\enddemo
\mytheorem{2.4} Assume that $a$ and $b$ are arbitrary two non-parallel
straight lines. Then there is a unique plane $\beta$ passing through the
line $b$ and being parallel to the line $a$.
\endproclaim
\demo{Proof} For the beginning we prove the existence of a required 
plane. From $a\nparallel b$ we derive $a\neq b$. On the line $b$ we
choose a point $B$ not lying on the line $a$. Then we draw the line
$\tilde a$ parallel to $a$ through this point $B$. The lines
$\tilde a$ and $b$ have one common point $B$, but they do not coincide
(since $\tilde a=b$ would mean $a\parallel b$). There is a plane
containing both of two such lines. Let's denote this plane by $\beta$.
Since $\tilde a\parallel a$ and $\tilde a\subset\beta$, applying the
theorem~\mythetheorem{2.1}, we derive that $a\parallel\beta$.\par
     Now let's prove the uniqueness of the plane $\beta$ constructed
just above. Assume that $\tilde\beta$ is some other plane containing
the line $b$ and being parallel to $a$. Then $B\in\tilde\beta$. The 
line $\tilde a$ by construction is parallel to the line $a$ and it
passes through the point $B$. Due to the theorem~\mythetheorem{2.3}
we have $\tilde a\subset\tilde\beta$. Hence, the plane $\tilde\beta$
contains two intersecting at the point $B$ but not coinciding straight
lines $\tilde a$ and $b$. So does the plane $\beta$. Hence, we have
$\tilde\beta=\beta$. The theorem is proved.\qed\enddemo
\mytheorem{2.5} Let $a$ be some straight line. If two distinct planes
$\alpha$ and $\beta$ are parallel to the line $a$ and if they intersect
along a line $b$, then $b\parallel a$.
\endproclaim
\demo{Proof} If $b$ is the line at the intersection of the planes 
$\alpha$ and $\beta$. As stated in the theorem, there are two distinct
planes $\alpha$ and $\beta$ passing through the line $b$ and being 
parallel to the line $a$. If the line would be not parallel to $a$, 
then, according to the theorem~\mythetheorem{2.4}, there would be only
one such plane. These considerations show that $b\parallel a$.
\qed\enddemo
\head
\SectionNum{3}{181} Parallelism of two planes.
\endhead
\rightheadtext{\S\,3. Parallelism of two planes.}
\mydefinition{3.1} Two planes $\alpha$ and $\beta$ are called 
{\it parallel}, if they coincide $\alpha=\beta$ or if they have
no common points.
\enddefinition
    For denoting the binary relation of parallelism of planes we use
the same sign as for the parallelism of straight lines. We write
$\alpha\parallel\beta$.
\mytheorem{3.1} Assume that a plane $\alpha$ is parallel to a plane
$\beta$. If a plane $\gamma$ intersects both planes $\alpha$ and
$\beta$ along two lines $a$ and $b$ respectively, then $a\parallel b$.
\endproclaim
\demo{Proof} If $\alpha=\beta$, the parallelism of the lines 
$a\parallel b$ follows from their coincidence $a=b$.\par
     Now assume that $\alpha\neq\beta$. Then from $\alpha\parallel\beta$
we get that the planes $\alpha$ and $\beta$ do not intersect. Hence, the
lines $a=\alpha\cap\gamma$ and $b=\beta\cap\gamma$ also do not intersect. 
They lie on one plane $\gamma$, therefore they are parallel. The theorem
is proved.\qed\enddemo
\mytheorem{3.2} Let $\alpha\neq\beta$ be two planes intersecting a line 
$c$ at the points $A$ and $B$. Then $c\perp\alpha$ and $c\perp\beta$
imply $\alpha\parallel\beta$.
\endproclaim
\demo{Proof} If $A=B$ we would have two planes $\alpha\neq\beta$
passing through the point $A\in c$ and being perpendicular to the
line $c$. However, this is impossible due to the 
theorem~\mythetheoremchapter{1.2}{4} from 
Chapter~\uppercase\expandafter{\romannumeral 4}. Hence, $A\neq B$.\par
     If we assume that the planes $\alpha$ and $\beta$ intersect and if
we denote by $C$ some point from their intersection $\alpha\cap\beta$, 
then would have two perpendiculars $[CB]$ and $[CA]$ dropped from the
point $C\notin c$ onto the line $c$. However, this contradicts the
theorem~\mythetheoremchapter{6.5}{3} from
Chapter~\uppercase\expandafter{\romannumeral 3}. Hence, the planes 
$\alpha$ and $\beta$ are parallel. The theorem is proved.
\qed\enddemo
      The definition of parallelism for two planes can be formulated
in absolute geometry either. The theorems~\mythetheorem{3.1} and
\mythetheorem{3.2} do not use the axiom~\mytheaxiom{A20}, they are
valid in absolute geometry. But the other theorems below in this 
section can be proved only in Euclidean geometry.
\mytheorem{3.3} If two intersecting straight lines $a$ and $b$ 
on a plane $\alpha$ are parallel to intersecting lines $\tilde a$ 
and $\tilde b$ on another plane $\tilde\alpha$, then the planes
$\alpha$ and $\tilde\alpha$ are parallel.
\endproclaim
\demo{Proof} If $\alpha=\tilde\alpha$, the parallelism of $\alpha$ 
and $\tilde\alpha$ follows from their coincidence (see 
definition~\mythedefinition{3.1}). Therefore, it is sufficient to 
consider the case where the planes $\alpha$ and $\tilde\alpha$ do not
coincide.\par
     Since the line $a$ lies on the plane $\alpha$, we have
$a\parallel\alpha$ (see definition~\mythedefinition{2.1}). The 
parallelism $a\parallel\tilde\alpha$ follows from the parallelism
of the lines $a$ and $\tilde a$ and from the fact that the line
$\tilde a$ lies on the plane $\tilde\alpha$ (see
theorem~\mythetheorem{2.1}). Hence, $\alpha$ and $\tilde\alpha$ are
two non-coinciding planes parallel to the line $a$. If we assume that
they intersect along some line $c$, then from the 
theorem~\mythetheorem{2.5} we derive $a\parallel c$. In a similar way 
we derive $b\parallel c$, hence, applying the theorem~\mythetheorem{1.4},
we get $a\parallel b$. However, two intersecting, but not coinciding
lines cannot be parallel (see definition~\mythedefinitionchapter{8.1}{3}
in Chapter~\uppercase\expandafter{\romannumeral 3}). The contradiction
obtained shows that the planes $\alpha$ and $\tilde\alpha$ do not
coincide. Hence, they are parallel. The theorem is proved.
\qed\enddemo
\mytheorem{3.4} For any point $O$ not lying on a plane $\alpha$
there is exactly one plane passing through this point and being 
parallel to the plane $\alpha$.
\endproclaim
\demo{Proof} Let's choose some point $B$ on the plane $\alpha$ and
draw two different lines $a$ and$b$ through this point on the plane
$\alpha$. Then through the point $O$ we draw two lines $\tilde a$ 
and $\tilde b$ parallel to $a$ and $b$ respectively. The existence
and uniqueness of such lines $\tilde a$ and $\tilde b$ follow from
the axiom~\mytheaxiom{A20}. We know, that there is a plane $\beta$
passing through the pair of non-coinciding lines $\tilde a$ and 
$\tilde b$ crossing at the point $O$. This plane $\beta$ is parallel
to the plane $\alpha$ due to the theorem~\mythetheorem{3.3}.\par
     Let's prove the uniqueness of the plane $\beta$ constructed
just above. Let's consider some plane $\tilde\beta$ passing through
the point $O$ and being parallel to the plane $\alpha$. Such a plane
has no common points with $\alpha$. Therefore, for two lines $a$ and
$b$ lying on the plane $\alpha$ we have $a\parallel\tilde\beta$ and 
$b\parallel\tilde\beta$. Therefore, we can apply the 
theorem~\mythetheorem{2.3} to the plane $\tilde\beta$, to the point
$O$ and to the lines $a$ and $\tilde a$. According to this theorem,
$\tilde a\parallel a$ implies $\tilde a\subset\tilde\beta$. Similarly,
$\tilde b\parallel b$ implies $\tilde b\subset\tilde\beta$. Hence,
the plane $\tilde\beta$ passes through the lines $\tilde a$ and 
$\tilde b$ which define the plane $\beta$. This yields $\tilde\beta
=\beta$. Thus, the theorem is proved.
\qed\enddemo
\mytheorem{3.5} Let $\alpha\neq\beta$ be two parallel planes. If a
plane $\gamma\neq\beta$ intersects the plane $\beta$ along a line
$b$, then it it intersects the plane $\alpha$ along some line $a$.
\endproclaim
\demo{Proof} If we assume that $\gamma$ does not intersect 
$\alpha$, then $\gamma$ is parallel to $\alpha$ according to 
the definition~\mythedefinition{3.1}. Assume that $B$ is some 
point on the line $b$ produced as the intersection of the planes
$\beta$ and $\gamma$. Then we have two planes $\beta$ and $\gamma$
passing through the point $B$ and being parallel to the plane
$\alpha$, which contradicts the theorem~\mythetheorem{3.4}. 
The contradiction obtained shows that the plane $\gamma$ intersects
the plane $\alpha$ along some line $a$.
\qed\enddemo
\mytheorem{3.6} Let $\alpha\neq\beta$ be two parallel planes. If the
line $c$ intersects the plane $\beta$, but does not lies on it,
then the line $c$ intersects the plane $\alpha$ too.
\endproclaim
\demo{Proof} Let's denote by $B$ the intersection point for the plane 
$\beta$ and the line $c$. Then we choose some point $C\neq B$ on the
plane $\beta$ and draw a plane $\gamma$ through the line $c$ and the
point $C$. At the intersection of the planes $\beta$ and $\gamma$
we get the line $b=BC$. According to the theorem~\mythetheorem{3.5},
the plane $\gamma$ intersects the plane $\alpha$ along some
line $a$. From $\alpha\parallel\beta$ by applying the 
theorem~\mythetheorem{3.1} we get $a\parallel b$. The lines $a$, $b$,
and $c$ lie on one plane $\gamma$, therefore we can apply the 
theorem~\mythetheorem{1.2} to them. According to this theorem, the line
$c$ intersecting the line $b$ at the point $B$ intersects the line
$a$ at some point $A$. But $a\subset\alpha$, hence, the point $A$
is the intersection point for the line $c$ and the plane $\alpha$.
\qed\enddemo
\mytheorem{3.7} Assume that $a\neq b$ are two parallel straight lines.
If a plane $\gamma$ intersects the line $b$, but does not contain
this line, then it intersects the line $a$ too.
\endproclaim
\myexercise{3.1} Derive the theorem~\mythetheorem{3.7} as a direct
consequence of the theorem~\mythetheorem{2.2}.
\endproclaim
\mytheorem{3.8} Let $\alpha\neq\beta$ be two parallel planes. Then 
any perpendicular to the plane $\alpha$ is a perpendicular to the 
plane $\beta$.
\endproclaim
\demo{Proof} Assume that the line $c$ is a perpendicular to the plane
$\alpha$ at a point $A\in\alpha$. According to the 
theorem~\mythetheorem{3.7}, the line $c$ crosses the plane $\beta$ at
some point $B$.\par
     Let's prove that $c\perp\beta$. For this purpose we consider
some line $b$ lying on the plane $\beta$ and passing through 
the point $B$. There is a plane $\gamma$ passing through the lines
$b$ and $c$. It intersects the plane $\beta$ along the line $b$. 
Let's denote by $a$ the line produced at the intersection of the 
planes $\gamma$ and $\alpha$. The line $a$ passes through the point 
$A$. From the theorem~\mythetheorem{3.1} we derive $a\parallel b$. 
The lines $a$, $b$, and $c$ lie on one plane $\gamma$, therefore, we
can apply the theorem~\mythetheorem{1.3} to them. Due to this theorem
$a\perp c$ implies $b\perp c$. Thus, the line $c$ appears to be
perpendicular to an arbitrary line $b$ lying on the plane $\beta$ and
passing through the point $B$. This fact yields the required result
$c\perp\beta$.\qed\enddemo
\mytheorem{3.9} Let $\alpha\neq\beta$ be two parallel planes. If a plane
$\gamma$ is perpendicular to the plane $\alpha$, then it is perpendicular
to the plane $\beta$ either.
\endproclaim
\myexercise{3.2} Derive the theorem~\mythetheorem{3.9} as a consequence
of the theorems~\mythetheorem{3.5} and \mythetheorem{3.8}.
\endproclaim
\head
\SectionNum{4}{184} The sum of angles of a triangle.
\endhead
\rightheadtext{\S\,4. The sum of angles of a triangle.}
\mytheorem{4.1} The sum of angles in an arbitrary triangle is equal 
to a straight angle.
\endproclaim
\demo{Proof} \parshape 13 0cm 10cm0cm 10cm 3.5cm 6.5cm 3.5cm 6.5cm 
3.5cm 6.5cm 3.5cm 6.5cm 3.5cm 6.5cm 3.5cm 6.5cm 3.5cm 6.5cm 3.5cm 6.5cm
3.5cm 6.5cm 3.5cm 6.5cm 0cm 10cm
     Let's choose some arbitrary triangle $ABC$. Let's draw a line
parallel to the side $AC$ through the vertex $B$ in this triangle. 
\vadjust{\vskip 5pt\hbox to 0pt{\kern 0pt
\includegraphics{Oris61.eps}\hss}\vskip -5pt}On this line 
we mark two point $D$ and $E$ on different sides of the point $B$. The
points $A$ and $C$ and the triangle $ABC$ in whole lie on one side with
respect to the lune $DE$ since the segment $[AC]$ lies on the line 
parallel to the line $DE$ and, therefore, cannot intersect this line. 
Let's apply the theorem~\mythetheorem{1.1} to the line $AB$ intersecting
two parallel lines $AC$ and $DE$. It yields the congruence of the angles
$\angle DBA\cong\angle CAB$. Similarly, applying the
theorem~\mythetheorem{1.1} to the line $CB$, we get $\angle EBC
\cong\angle BCA$. The angles $\angle DBA$, $\angle ABC$, and $\angle EBC$
have the common vertex $B$ and compose the straight angle $\angle DBE$.
Now, taking into account the above congruence relationships for angles, 
we can write 
$$
\angle CAB+\angle ABC+\angle BCA\cong\angle DBE.
$$
This is the very relationship which means that the sum of internal angles
of the triangle $ABC$ coincide with a straight angle.\qed\enddemo
\head
\SectionNum{5}{185} Midsegment of a triangle.
\endhead
\rightheadtext{\S\,5. Midsegment of a triangle.}
     Let $ABC$ be a triangle. Let's mark the centers of the sides $[AB]$ 
and $[BC]$ in this triangle. Let's denote them by $M$ and $N$ respectively.
The segment $[MN]$ is called the {\it midsegment\/} of the triangle $ABC$.
\mytheorem{5.1} The midsegment $[MN]$ connecting the centers of the sides
$[AB]$ and $[BC]$ in a triangle $ABC$ is parallel to the side $[AC]$ of
this triangle and $[AC]\cong 2\cdot [MN]$.
\endproclaim
\demo{Proof}\parshape 6 0cm 10cm 0cm 10cm 0cm 10cm 0cm 10cm
0cm 10cm 4.2cm 5.8cm 
     Let $M$ be the center of the side $[AB]$ in a triangle $ABC$. 
We draw the line parallel to the side $[AC]$ through this point $M$.
Applying Pasch's axiom~\mytheaxiom{A12}, it is easy to show that 
such line crosses the side $[BC]$ at some interior point $N$. Then 
we draw the line parallel to the side $[BC]$ through the point $M$. 
\vadjust{\vskip 5pt\hbox to 0pt{\kern -5pt
\includegraphics{Oris62.eps}\hss}\vskip -5pt}At the 
intersection of this line with the side $[AC]$ we get a point $K$
lying in the interior of the segment $[AC]$. Let's connect the 
points $N$ and $K$ with the segment $[NK]$. On the lines 
$MN$ and $MK$ we mark two points $D$ and $E$ for the sake of 
convenience.\par
\parshape 5 4.2cm 5.8cm 4.2cm 5.8cm 4.2cm 5.8cm 4.2cm 5.8cm 0cm 10cm
     The angles $\angle MBN$ and $\angle BME$ are inner crosswise 
lying angles at the intersections of the line $AB$ with two parallel
lines $BC$ and $MK$. The angles $\angle BME$ and $\angle AMK$ are
vertical angles. Hence, applying the theorem~\mythetheorem{1.1}, 
we get $\angle MBN\cong\angle AMK$. Now let's consider the angles 
$\angle MAK$ and $\angle DMA$. They are inner crosswise lying 
angles at the intersections of the line $AB$ with two parallel lines 
$AC$ and $MN$. If we take into account that $\angle DMA$ and 
$\angle BMN$ are vertical angles, we get $\angle MAK\cong\angle BMN$.
From the following three relationships
$$
\xalignat 3
&[AM]\cong [MB],
&&\angle AMK\cong\angle MBN,
&&\angle MAK\cong\angle BMN
\endxalignat
$$
we derive the congruence of the triangles $AMK$ and $MBN$ (see
theorem~\mythetheoremchapter{5.2}{3} in
Chapter~\uppercase\expandafter{\romannumeral 3}). Due to this 
congruence we get
$$
\xalignat 2
&\hskip -2em
[AK]\cong [MN],
&&[MK]\cong [BN].
\mytag{5.1}
\endxalignat
$$
Let's complement the congruences \mythetag{5.1} with one more 
relationship $\angle NMK\cong\angle MNB$. This relationship is 
derived if we consider the inner crosswise lying angles at the 
intersections of the line $MN$ with two parallel lines $BC$ and
$MK$. Now, applying the theorem~\mythetheoremchapter{5.1}{3} 
from Chapter~\uppercase\expandafter{\romannumeral 3}, we find
that the triangles $MBN$ and $NKM$ are congruent.\par
     At the intersections of the line $NK$ with two parallel lines
$BC$ and $MK$ we get the inner crosswise lying angles $\angle MKN$ 
and $\angle KNC$. Similarly, at the intersections of the line
$NK$ wih two parallel lines $AC$ and $MN$ we get the inner crosswise 
lying angles $\angle CKN$ and $\angle KNM$. Now, from the 
relationships
$$
\xalignat 3
&[NK]\cong [KN],
&&\angle MKN\cong\angle KNC,
&&\angle CKN\cong\angle KNM
\endxalignat
$$
we derive the congruence of the triangles $NKM$ and $KNC$. Thus,
we see that the four triangles $AMK$, $MBN$, $NKM$, and $KNC$ on
Fig\.~5.1 are congruent to each other. Hence, we have
$$
\xalignat 2
&[BN]\cong [NC],
&&[AK]\cong [KC]\cong [MN].
\endxalignat
$$
The first of these relationships means that the segment $[MN]$ lying
on the line parallel to the line $AC$ is the midsegment of the initial
triangle $ABC$. The second relationship is equivalent to $[AC]\cong 
2\cdot [MN]$. The theorem is proved.
\qed\enddemo
     In a triangle $ABC$ there are three midsegments $[MN]$, $[NK]$, 
and $[KM]$ parallel to the sides $[AC]$, $[BA]$, and $[CB]$ of 
this triangle respectively. They divide the triangle $ABC$ into 
four triangles $BMN$, $NKC$, $MAK$, and $KNM$, whose sides are 
twice as smaller than the corresponding sides of the triangle 
$ABC$. The angles of these triangles are congruent to the 
corresponding angles of the triangle $ABC$.
\head
\SectionNum{6}{187} Midsegment of a trapezium.
\endhead
\rightheadtext{\S\,6. Midsegment of a trapezium.}
\parshape 12 0cm 10cm 0cm 10cm 0cm 10cm 0cm 10cm 0cm 10cm 0cm 10cm 
0cm 10cm 0cm 10cm 0cm 10cm 0cm 10cm 0cm 10cm 4.2cm 5.8cm 
     Assume that on two parallel lines $a\neq b$ two segments 
$[AB]$ and $[CD]$ are marked. Let's connect the points $A$, $B$, $C$,
and $D$ with four segments $[DA]$, $[DB]$, $[CA]$, and $[CB]$. The
points $A$ and $B$ lies on one side of the line $CD$ since the segment
does not intersect the line $CD$ (see \S\,5 in 
Chapter~\uppercase\expandafter{\romannumeral 2}). The rays $[DA\rangle$ 
and $[DB\rangle$ cannot coincide since in this case the segment $[AB]$
would lie on the line $DA$ intersecting the line $CD$. Hence, we
conclude: one of the two angles $\angle CDA$ and $\angle CDB$ lies
inside the other. For the sake of certainty assume that 
$\angle CDB<\angle CDA$ as it is shown on Fig\.~6.1. In this case
we can apply the lemma~\mythelemmachapter{6.2}{2} from
Chapter~\uppercase\expandafter{\romannumeral 2}. From this lemma
we derive that the ray $[DB\rangle$ intersects the segment $[AC]$
\vadjust{\vskip 5pt\hbox to 0pt{\kern -5pt
\includegraphics{Oris63.eps}\hss}\vskip -5pt}at some its interior 
point $O$.\par
\parshape 7 4.2cm 5.8cm 4.2cm 5.8cm 4.2cm 5.8cm 4.2cm 5.8cm 
4.2cm 5.8cm 4.2cm 5.8cm 0cm 10cm
     The point $O$ lies on the ray $[DB\rangle$, therefore, the following
two mutual dispositions of the points $D$, $B$, and $O$ are possible:
$$
\aligned
&(D\blacktriangleright O\blacktriangleleft B),\\
&(D\blacktriangleright B\blacktriangleleft O).
\endaligned
\mytag{6.1}
$$
The coincidence $O=B$ is impossible since in this case the segment
$[AB]$ would lie on the line $AC$, which is not parallel to the
line $CD$. If we assume that thew point $B$ lies between the points 
$D$ and $O$, then the line $AB$ intersects the side $[DO]$ in the
triangle $DOC$, but it does not intersect the side $[OC]$. Applying
Pasch's axiom~\mytheaxiom{A12} in this situation, we find that the
line $AB$ intersects the cide $[CD]$ in the triangle $DOC$, which
contradicts the parallelism of the lines $AB$ and $CD$. This 
contradiction means that only the first disposition of points $D$, 
$B$, and $O$ in \mythetag{6.1} can be actually implemented, i\.\,e\.
$O$ is an interior point for both segments $[AC]$ and $[BD]$.\par
     Now it is easy to show that the segments $[DA]$ and $[CB]$ do
not intersect. For this purpose it is sufficient to consider the 
triangle $DAO$ and the line $CB$. The line $CB$ does not intersect
the sides $[DO]$ and $[OA]$ in the triangle $DAO$. Hence, according
to Pasch's axiom~\mytheaxiom{A12} it cannot intersect the third side 
$[DA]$ of this triangle. These considerations prove the following
theorem.
\mytheorem{6.1} For any two segments $[AB]$ and $[CD]$ lying on
two parallel lines $a\neq b$ exactly one of two segments $[DA]$ or
$[DB]$ intersects exactly one of two segments $[CA]$ or $[CB]$ at
some interior point $O$ for both intersecting segments.
\endproclaim
\parshape 15 0cm 10cm 0cm 10cm 0cm 10cm 0cm 10cm 4.2cm 5.8cm
4.2cm 5.8cm 4.2cm 5.8cm 4.2cm 5.8cm 4.2cm 5.8cm 4.2cm 5.8cm 
4.2cm 5.8cm 4.2cm 5.8cm 4.2cm 5.8cm 4.2cm 5.8cm 0cm 10cm
     According to the theorem~\mythetheorem{6.1}, the segments $[AB]$ 
and $[CD]$, lying on two parallel lines can be complemented up to a 
closed polygonal line with no self-intersections. 
\vadjust{\vskip 5pt\hbox to 0pt{\kern -5pt
\includegraphics{Oris64.eps}\hss}\vskip -5pt}For the situation
shown on Fig\.~6.2 this line is $ABCD$. It bounds a part of the plane
consisting of two triangles $ABC$ and $ADC$, which  intersect along their
common side $[AC]$. In general case any closed polygonal line without
self-intersections bounds some set of points on a plane which can be
represented as the union of finite number of triangles none two of which
have common interior points. This fact is known as Jordan's theorem. Its
proof can be found in the book \mycite{7}.\par
     A set of points on a plane bounded by a closed polygonal line 
is called a {\it polygon}. The segments of such a line are called 
{\it sides\/} of a polygon. By the number of sides polygons are divided
into triangles, quadrangles, pentagons, hexagons etc.
\mydefinition{6.1} A quadrangle two sides of which lie on two parallel 
lines is called a {\it trapezium}. Parallel sides of a trapezium are 
called {\it bases}, other two sides are called {\it lateral sides}.
\enddefinition
\mydefinition{6.2} The segment connecting centers of lateral sides of a 
trapezium is called the {\it midsegment\/} of this trapezium.
\enddefinition
    Let $ABCD$ be a trapezium. The segments $[AC]$ and $[BD]$ are called
{\it diagonals\/} of this trapezium. According to the 
theorem~\mythetheorem{6.1} they intersect at some $O$ lying in the interior
of each of these two segments (see Fig\.~6.1).
\mytheorem{6.2} The midsegment $[MK]$ of a trapezium $ABCD$ is parallel to
its bases $[AB]$ and $[CD]$ and such that the relationship 
$2\cdot [MK]\cong [AB]+[CD]$ is fulfilled.
\endproclaim
\demo{Proof} Let's consider the pair of triangles $ABC$ and $CDA$, 
which compose the trapezium $ABCD$ (see Fig\.~6.2). Let's draw their
midsegments $[NK]$ and $[NM]$. Then we apply the theorem~\mythetheorem{5.1}
to them. This yields the relationships
$$
\xalignat 2
&NK\parallel AB,
&&NM\parallel CD.
\endxalignat
$$
Now from $AB\parallel CD$ we derive that both lines $NM$ and $NK$ 
passing through the point $N$ are parallel to the line $CD$. Due 
to the axiom~\mytheaxiom{A20} these two lines should coincide. 
This means that the points $M$, $N$, and $K$ lie on one straight 
line $MK$ parallel to the bases of our trapezium.\par
     Note that the segment $[MK]$ is composed of two segments $[MN]$ 
and $[NK]$. For these segments from the theorem~\mythetheorem{5.1}
we derive $2\cdot [MN]\cong [CD]$ and $2\cdot [NK]\cong [AB]$. This
yields the required relationship $2\cdot [MK]\cong [AB]+[CD]$ for
the segment $[MK]$. The theorem is proved.
\qed\enddemo
\head
\SectionNum{7}{190} Parallelogram.
\endhead
\rightheadtext{\S\,7. Parallelogram.}
\mydefinition{7.1} A trapezium lateral sides of which are parallel
is called a {\it parallelogram}.
\enddefinition
\mytheorem{7.1} A trapezium is a parallelogram if and only if its bases
are congruent.
\endproclaim
\demo{Proof} Let's consider a trapezium $ABCD$. Its diagonals
$[AC]$ and $[BD]$ intersect at some point $O$ lying in the interior 
of both of these two segments (see Fig\.~6.1). From the parallelism of 
bases $AB\parallel CD$ in our trapezium we derive the congruence of the
following inner crosswise lying angles:
$$
\xalignat 2
&\angle BAC\cong\angle ACD,&&\angle ABD\cong\angle  BDC.
\endxalignat
$$
If we complement these relationships with the congruence of the bases 
$[AB]\cong [CD]$, then due to the theorem~\mythetheoremchapter{5.2}{3}
from Chapter~\uppercase\expandafter{\romannumeral 3} we find that the
triangles $AOB$ and $COD$ are congruent. Hence, we have the following
relationships:
$$
\xalignat 2
&\hskip -2em
[AO]\cong [CO],
&&[BO]\cong [DO].
\mytag{7.1}
\endxalignat
$$
The angles $\angle AOD$ and $\angle BOC$ are congruent since they 
are vertical angles. Applying the theorem~\mythetheoremchapter{5.1}{3} 
from Chapter~\uppercase\expandafter{\romannumeral 3} and taking into
account the relationships \mythetag{7.1}, we find that the triangles
$AOD$ and $COB$ are congruent. Hence, for the inner crosswise lying 
angles $\angle ADB$ and $\angle DBC$ we get $\angle ADB\cong\angle DBC$,
which yields $AD\parallel BC$. Hence, $ABCD$ is a parallelogram.\par
     Conversely, assume that the quadrangle $ABCD$ is a parallelogram.
Then from $AB\parallel CD$ and $AD\parallel BC$ we derive 
$$
\xalignat 2
&\angle BAC\cong\angle ACD,&&\angle BCA\cong\angle CAD.
\endxalignat
$$
We complement these relationships with the trivial relationship 
$[AC]\cong [CA]$ and we apply the theorem~\mythetheoremchapter{5.2}{3}
from Chapter~\uppercase\expandafter{\romannumeral 3}. As a result we 
derive the congruence of the triangles $ABC$ and $CDA$. Hence, we have
$[AB]\cong [CD]$ and $[AD]\cong [BC]$. The theorem is proved.
\qed\enddemo
    In proving the theorem~\mythetheorem{7.1} we have proved the 
following two important additional facts:
\roster
\item opposite sides in any parallelogram are congruent, i\.\,e\.
      if $ABCD$ is a parallelogram, then $[AB]\cong [CD]$ and
      $[AD]\cong [BC]$;
\item diagonals of any parallelogram intersect each other and the 
      intersection point divides them into halves, i\.\,e\. if 
      $ABCD$ is a parallelogram and if $[AC]\cap[BD]=O$, then
      $[AO]\cong [OC]$ and $[BO]\cong [OD]$.
\endroster
\mytheorem{7.2} A quadrangle $ABCD$ is a parallelogram if and only 
if the opposite sides of this quadrangle are congruent, i\.\,e\. if
$[AB]\cong [CD]$ and if $[AD]\cong [BC]$.
\endproclaim
\mytheorem{7.3} A quadrangle $ABCD$ is a parallelogram if and only 
if its diagonals intersect each other at some interior point $O$ 
that divides them into halves.
\endproclaim
\myexercise{7.1} Considering various dispositions of the points $A$ 
and $C$ relative to the line $BD$, prove the theorems~\mythetheorem{7.2}
and \mythetheorem{7.3}.
\endproclaim
\head
\SectionNum{8}{192} Codirected and equal vectors in the space.
\endhead
\rightheadtext{\S\,8. Codirected and equal vectors \dots}
    The concept of codirectedness was introduced above in
Chapter~\uppercase\expandafter{\romannumeral 2}. However, it was
applicable only to vectors lying on one straight line. The concept
of equality was also applicable only to vectors lying on one
straight line (see definition~\mythedefinitionchapter{4.2}{2} 
in Chapter~\uppercase\expandafter{\romannumeral 2}
and definition~\mythedefinitionchapter{4.1}{3} in
Chapter~\uppercase\expandafter{\romannumeral 3}). Here we extend
these concepts for the case of vectors not lying on one line.
\mydefinition{8.1} Two vectors $\overrightarrow{AB}$ and
$\overrightarrow{CD}$ not lying on one straight line are called
{\it codirected\/} if
\roster
\rosteritemwd=5pt
\item they lie on parallel straight lines;
\item the segment $[BD]$ connecting their ending points does not
      intersect the segment $[AC]$ connecting their initial points.
\endroster
\enddefinition
\mydefinition{8.2} Two vectors $\overrightarrow{AB}$ and
$\overrightarrow{CD}$ are called {\it equal\/} if they are
codirected and if the segment $[AB]$ is congruent to the
segment $[CD]$.
\enddefinition
     The definitions~\mythedefinition{8.1} and 
\mythedefinition{8.2} can be formulated in absolute geometry either.
However, only in Euclidean geometry the axiom~\mytheaxiom{A20} providing
the transitivity of the parallelism relation (see
theorem~\mythetheorem{1.4}) makes these definitions reasonable.\par
     It is easy to see that the relation of codirectedness of vectors
introduced by the definition~\mythedefinition{8.1} and by the
definition~\mythedefinitionchapter{4.2}{2} in 
Chapter~\uppercase\expandafter{\romannumeral 2} is reflexive and symmetric.
In order to prove its transitivity one should study several special cases.
\mylemma{8.1} Let $A_0\prec A_1\prec A_2\prec A_3$ be a monotonic sequence
of points on a straight line $a$. If a vector $\overrightarrow{MN}$ lying
on another line $b\neq a$ is codirected to the vector
$\overrightarrow{A_1A_3}$, then it is codirected to each of the vectors
$\overrightarrow{A_0A_3}$ and $\overrightarrow{A_2A_3}$.
\endproclaim
\demo{Proof} \parshape 16 0cm 10cm 0cm 10cm 0cm 10cm 
0cm 10cm 0cm 10cm 4.5cm 5.5cm 4.5cm 5.5cm 4.5cm 5.5cm 
4.5cm 5.5cm 4.5cm 5.5cm 4.5cm 5.5cm 4.5cm 5.5cm 4.5cm 5.5cm 
4.5cm 5.5cm 4.5cm 5.5cm 0cm 10cm 
     The codirectedness of the vectors $\overrightarrow{MN}$ and
$\overrightarrow{A_1A_3}$ means that the quadrangle $MNA_3A_1$ is
a trapezium. Lateral sides of a trapezium does not intersect each 
other, therefore, the line $NA_3$ has no common points 
with the side $[MA_1]$ \vadjust{\vskip 5pt\hbox to -5pt{\kern -5pt
\includegraphics{Oris65.eps}\hss}\vskip -5pt}in the 
triangle $MA_1A_2$. This line has no common points with the side 
$[A_1A_2]$ in the triangle $MA_1A_2$ either. Applying Pasch's
axiom~\mytheaxiom{A12}, we conclude that the line $NA_3$ cannot 
intersect the side $[MA_2]$ in this triangle. In other words, the 
segments $[MA_2]$ and $[NA_3]$ do not intersect each other. Hence, 
the vectors $\overrightarrow{MN}$ and $\overrightarrow{A_2A_3}$ 
are codirected. The relationship
$\overrightarrow{MN}\upuparrows\overrightarrow{A_0A_3}$ is proved 
in a similar way by considering the triangle $MA_0A_1$. The lemma 
is proved.
\qed\enddemo
\mylemma{8.2} Let $A_1\prec A_2\prec A_3\prec A_4$ be a monotonic 
sequence of points on a line $a$. If a vector $\overrightarrow{MN}$
lying on another line $b\neq a$ is codirected to the vector
$\overrightarrow{A_1A_3}$, then it is codirected to each of the
vectors $\overrightarrow{A_1A_2}$ and $\overrightarrow{A_1A_4}$.
\endproclaim
\myexercise{8.1} Prove the lemma~\mythelemma{8.2} using considerations
similar to those used in proving the lemma~\mythelemma{8.1}.
\endproclaim
\mylemma{8.3} Assume that a vector $\overrightarrow{MN}$ lies on
a line $b$, while the vectors $\overrightarrow{AB}$ and
$\overrightarrow{CD}$ lie on another line $a\neq b$. Then the
relationships $\overrightarrow{MN}\upuparrows\overrightarrow{AB}$ 
and $\overrightarrow{AB}\upuparrows\overrightarrow{CD}$ imply 
$\overrightarrow{MN}\upuparrows\overrightarrow{CD}$.
\endproclaim
\demo{Proof} Let's enumerate the points $A$, $B$, $C$, $D$ on the 
line $a$ so that they form a monotonic sequence $A_1\prec A_2\prec 
A_3\prec A_4$. Let $A=A_i$, $B=A_j$, $C=A_p$, $D=A_q$ so that $i<j$. 
Then $\overrightarrow{AB}\upuparrows\overrightarrow{CD}$ implies 
$p<q$. Let's denote $r=\min\{i,p\}$ and $s=\max\{j,q\}$. Then we
have the following implications:
$$
\align
&\overrightarrow{MN}\upuparrows\overrightarrow{A_iA_j}\quad\Rightarrow
\quad\overrightarrow{MN}\upuparrows\overrightarrow{A_rA_j}\quad
\Rightarrow\quad\overrightarrow{MN}\upuparrows\overrightarrow{A_rA_s},\\
&\overrightarrow{MN}\upuparrows\overrightarrow{A_rA_s}\quad\Rightarrow
\quad\overrightarrow{MN}\upuparrows\overrightarrow{A_rA_q}\quad
\Rightarrow\quad\overrightarrow{MN}\upuparrows\overrightarrow{A_pA_q}.
\endalign
$$
This sequence of implications is obtained as a result of applying the
lemmas~\mythelemma{8.1} and \mythelemma{8.2} and taking into account 
the inequalities $r\leqslant i<j\leqslant s$ and $r\leqslant p<q
\leqslant s$. In the end of this sequence of implications we get the
required relationship $\overrightarrow{MN}\upuparrows\overrightarrow{CD}$.
Thus, the lemma is proved.\qed\enddemo
\mylemma{8.4} Assume that a vector $\overrightarrow{MN}$ lies on 
a line $b$, while $\overrightarrow{AB}$ and $\overrightarrow{CD}$
are two vectors lying on another line $a\neq b$. Then
$\overrightarrow{AB}\upuparrows\overrightarrow{MN}$ and 
$\overrightarrow{MN}\upuparrows\overrightarrow{CD}$ imply 
$\overrightarrow{AB}\upuparrows\overrightarrow{CD}$.
\endproclaim
\demo{Proof} Let's prove the lemma by contradiction. Assume that
the conditions $\overrightarrow{AB}\upuparrows \overrightarrow{MN}$ and
$\overrightarrow{MN}\upuparrows\overrightarrow{CD}$ are fulfilled, while
the condition $\overrightarrow{AB}\upuparrows\overrightarrow{CD}$ is
not fulfilled. Then the vector $\overrightarrow{CD}$ is codirected to
the vector $\overrightarrow{BA}$, which is an opposite vector for 
the vector $\overrightarrow{AB}$. From $\overrightarrow{MN}\upuparrows
\overrightarrow{CD}$ and $\overrightarrow{CD}\upuparrows
\overrightarrow{BA}$, applying the lemma~\mythelemma{8.3}, we get
$\overrightarrow{MN}\upuparrows\overrightarrow{BA}$. But two 
relationships $\overrightarrow{AB}\upuparrows \overrightarrow{MN}$
and $\overrightarrow{MN}\upuparrows \overrightarrow{BA}$ cannot be
fulfilled simultaneously because of the theorem~\mythetheorem{6.1}. 
The contradiction obtained proves the lemma.
\qed\enddemo
\mytheorem{8.1} The relation of codirectedness of vectors in the space 
is transitive, i\.\,e\. the relationships $\overrightarrow{AB}\upuparrows
\overrightarrow{CD}$ and $\overrightarrow{CD}\upuparrows\overrightarrow{EF}$
imply the relationship $\overrightarrow{AB}\upuparrows\overrightarrow{EF}$.
\endproclaim
\demo{Proof} The case where all of the vectors $\overrightarrow{AB}$,
$\overrightarrow{CD}$, and $\overrightarrow{EF}$ lie on one straight
line was considered in \S\,2 of
Chapter~\uppercase\expandafter{\romannumeral 2}. The case where some
two of these three vectors lie on one straight line is described by
the lemmas~\mythelemma{8.3} and \mythelemma{8.4}. The rest is the case 
of general position where these vectors lie on three distinct straight
lines $a$, $b$, and $c$. In this case the relationships 
$\overrightarrow{AB}\upuparrows\overrightarrow{CD}$ and
$\overrightarrow{CD}\upuparrows\overrightarrow{EF}$ imply the
parallelism of the corresponding lines $a\parallel b$ and
$b\parallel c$. Hence, due to the theorem~\mythetheorem{1.4}
we have $a\parallel c$.\par
     Let's consider the vector $\overrightarrow{AB}$ lying on the line 
$a$. Let's draw the plane $\alpha$ perpendicular to the line $a$ through
the point $A$. According to the theorem~\mythetheorem{3.7}, this plane
intersects the lines $b$ and $c$. We denote the intersection points
through $C'$ and $E'$. respectively. In a similar way, we draw the plane
$\beta\perp a$ through the point $B$. At the intersections of this plane
$\beta$ with the lines $b$ and $c$ we find the point $D'$ and $F'$
respectively. Due to the theorem~\mythetheorem{3.2} the planes $\alpha$
and $\beta$ are parallel. Therefore, the segment $[AC']$ does not
intersect the segment $[BD']$, the segment $[C'E']$ does not intersect
the segment $[D'F']$, and the segment $[AE']$ does not intersect the
segment $[BF']$. As a result we get the following three relationships:
$$
\xalignat 3
&\overrightarrow{AB}\upuparrows\overrightarrow{C'D'},
&&\overrightarrow{C'D'}\upuparrows\overrightarrow{E'F'},
&&\overrightarrow{AB}\upuparrows\overrightarrow{E'F'}.
\qquad
\mytag{8.1}
\endxalignat
$$
Let's combine the relationship $\overrightarrow{AB}\upuparrows
\overrightarrow{CD}$ with the first relationship \mythetag{8.1} 
and take into account that the vectors $\overrightarrow{CD}$ and
$\overrightarrow{C'D'}$ lie on one straight line $b$. Applying
the lemma~\mythelemma{8.4}, we get $\overrightarrow{C'D'}\upuparrows
\overrightarrow{CD}$. Let's combine this relationship with the
second relationship in \mythetag{8.1} and apply the 
lemma~\mythelemma{8.3}. As a result we get $\overrightarrow{CD}
\upuparrows\overrightarrow{E'F'}$. Being combined with
$\overrightarrow{CD}\upuparrows\overrightarrow{EF}$, upon applying
the lemma~\mythelemma{8.4}, this relationship yields
$\overrightarrow{EF}\upuparrows\overrightarrow{E'F'}$. The last
step is to combine $\overrightarrow{EF}\upuparrows
\overrightarrow{E'F'}$ with the third relationship
\mythetag{8.1} and apply the lemma~\mythelemma{8.3}. As a result
we get the required relationship
$\overrightarrow{AB}\upuparrows\overrightarrow{EF}$. The theorem
is proved.\qed\enddemo
    Thus, we have proved that the relation of codirectedness for
vectors in the space is reflexive, symmetric, and transitive, 
i\.\,e\. it is an equivalence relation. The relation of congruence 
for segments possesses the same properties. Now from the 
definition~\mythedefinition{8.2} we derive the following theorem.
\mytheorem{8.2} In Euclidean geometry the relation of equality for 
vectors in the space is reflexive, symmetric, and transitive, because 
of which it is an equivalence relation.
\endproclaim
     The classes of mutually equal vectors in Euclidean geometry are
called {\it free vectors}. Geometric vectors composing a class are
called {\it geometric realizations\/} of a free vector.
\mytheorem{8.3} For any free vector $\bold a$ and for any point $A$
there is a geometric realization $\overrightarrow{AB}$ of $\bold a$ 
with the initial point $A$.
\endproclaim
     The theorem~\mythetheorem{8.3} approved the above terminology.
A free vector is called {\tencyr\char '074}free{\tencyr\char '076}
since it can be realized at any point of the space without any 
limitations.
\myexercise{8.2} Let $\overrightarrow{CD}$ be some geometric 
realization of a free vector $\bold a$. Using this vector 
$\overrightarrow{CD}$, prove the theorem~\mythetheorem{8.3} through
constructing a required geometric realization $\overrightarrow{AB}$ 
of the vector $a$ with the initial point $A$
\endproclaim
\mydefinition{8.3} Free vectors $\bold a$ and $\bold b$ are called
{\it collinear\/} and are written $\bold a\parallel\bold b$ if some 
of their geometric realizations $\overrightarrow{AB}$ and
$\overrightarrow{CD}$ lie on parallel straight lines. If
$\overrightarrow{AB}\upuparrows\overrightarrow{AB}$, then the vectors
$\bold a$ and $\bold b$ are said to be {\it codirected}. In the case
where the vectors $\bold a$ and $\bold a$ are collinear, but not
codirected, they are said to be {\it oppositely directed}.
\enddefinition
\myexercise{8.3} Prove the correctness of the 
definition~\mythedefinition{8.3} by showing that the properties 
of collinearity and codirectedness of free vectors $\bold a$ and
$\bold b$ do not depend on any particular choice of their geometric
realizations $\overrightarrow{AB}$ and $\overrightarrow{CD}$.
\endproclaim
\mytheorem{8.4} For any four points $A$, $B$, $C$, and $D$ the 
equality $\overrightarrow{AB}=\overrightarrow{CD}$ implies 
$\overrightarrow{AC}=\overrightarrow{BD}$ and, conversely, the
equality $\overrightarrow{AC}=\overrightarrow{BD}$ implies 
$\overrightarrow{AB}=\overrightarrow{CD}$.
\endproclaim
\demo{Proof} In the case where the points $A$, $B$, $C$, and $D$ 
lie on one straight line the theorem~\mythetheorem{8.4} reduces to
the theorem~\mythetheoremchapter{4.1}{3} from
Chapter~\uppercase\expandafter{\romannumeral 3}. Therefore, we consider
the case where the points $A$, $B$, $C$, and $D$ do not lie on one 
straight line. The condition $\overrightarrow{AB}=\overrightarrow{CD}$ 
yields the relationships
$$
\xalignat 2
&AB\parallel CD,&&[AB]\cong [CD]
\endxalignat
$$
and the condition $[AC]\cap [BD]=\varnothing$. In this case we can 
apply the theorem~\mythetheorem{7.1} which means that the quadrangle
$ACDB$ is a parallelogram. Hence, with the use of the 
theorem~\mythetheorem{7.2} we get
$\overrightarrow{AC}=\overrightarrow{BD}$. Conversely, \pagebreak 
the relationship $\overrightarrow{AB}=\overrightarrow{CD}$ is derived 
from the relationship $\overrightarrow{AC}=\overrightarrow{BD}$ in a 
quite similar way.
\qed\enddemo
\head
\SectionNum{9}{197} Vectors and parallel translations.
\endhead
\rightheadtext{\S\,9. Vectors and parallel translations.}
     Translations by some vectors along straight lines were defined
in \S\,14 of Chapter~\uppercase\expandafter{\romannumeral 4}, 
In Euclidean geometry due to the axiom~\mytheaxiom{A20} one can
specify the properties of these mappings making their description
substantially more detailed.
\mytheorem{9.1} If $\bold c$ is a slipping vector on a straight line
$a$ and if $\overrightarrow{AB}$ is some its geometric realization
lying on the line $a$, then the relationship $p_{a\bold c}(C)=D$ is
fulfilled if and only if the vectors $\overrightarrow{AB}$ and
$\overrightarrow{CD}$ are equal in the sense of the
definition~\mythedefinition{8.2}.
\endproclaim
\demo{Proof} If the point $C$ lies on the line $a$, then the point $D$
also lies on the line $a$. In this case the theorem~\mythetheorem{9.1}
is reduced to the theorem~\mythetheoremchapter{4.2}{3} from
Chapter~\uppercase\expandafter{\romannumeral 3}.\par
\parshape 13 0cm 10cm 0cm 10cm 0cm 10cm 0cm 10cm 4.2cm 5.8cm
4.2cm 5.8cm 4.2cm 5.8cm 4.2cm 5.8cm 4.2cm 5.8cm 4.2cm 5.8cm 4.2cm 5.8cm
4.2cm 5.8cm 0cm 10cm
     Let's consider the case where the point $C$ does not lie on the
line $a$. The points $A$, $B$, $C$, and $D$ lie on one plane (see
theorem~\mythetheoremchapter{14.1}{4} in
Chapter~\uppercase\expandafter{\romannumeral 4}). Let's denote this 
plane by $\alpha$. Taking into account that $B=p_{a\bold c}(A)$, we
mark one more point $F=p_{a\bold c}(B)$ on the line $a$.
\vadjust{\vskip 5pt\hbox to -5pt{\kern 0pt
\includegraphics{Oris66.eps}\hss}\vskip -5pt}Then we 
draw the lines $AC$ and $BD$ on the plane $\alpha$ and mark a point
$E$ on the line $BD$ as shown on Fog\.~9.1. Since $p_{a\bold c}$ is
a mapping of congruent translation, now from the relationships
$$
\hskip -2em
\aligned
p_{a\bold c}(A)&=B,\\
p_{a\bold c}(C)&=D,\\
p_{a\bold c}(B)&=F
\endaligned
\mytag{9.1}
$$
we derive the congruence of angles $\angle CAB\cong\angle DBF$ 
and the congruence of segments $[AC]\cong[BD]$. Then, using the 
congruence of the vertical angles $\angle DBF$ and $\angle ABE$, 
we get $\angle CAB\cong\angle ABE$. But $\angle CAB$ and $\angle DBF$ 
are the inner crosswise lying angles arising at the intersections of
the lines $AC$ and $BD$ with the line $a$. Hence, from $\angle CAB
\cong\angle ABE$ we derive the parallelism of lines $AC\parallel BD$. 
By construction of the mapping $p_{a\bold c}$ the points $C$ and
$D$ lie on one side of the line $a$ (see \S\,14 in 
Chapter~\uppercase\expandafter{\romannumeral 4}). Therefore the
segments $[AB]$ and $[CD]$ do not intersect and the quadrangle 
$ACDB$ is the trapezium with the bases $[AC]$ and $[BD]$. Moreover,
let's take into account the relationship $[AC]\cong [BD]$ derived
from \mythetag{9.1} and apply the theorem~\mythetheorem{7.1}. 
According to this theorem the quadrangle $ACDB$ is a parallelogram. 
Hence, we immediately get $\overrightarrow{AB}=\overrightarrow{CD}$.
\par
     Now, conversely, assume that $\overrightarrow{AB}=
\overrightarrow{CD}$. Then the quadrangle $ACDB$ is a parallelogram, 
because of which $[AC]\cong [BD]$ and $\angle CAB\cong\angle DBF$.
Remember that the mapping $p_{a\bold c}$ is constructed as an
extension of the mapping $p_{\bold c}\!:a\to a$ from the line $a$
to the plane $\alpha$, then from this plane to the whole space (see
\S\,14 in Chapter~\uppercase\expandafter{\romannumeral 4}). It maps
the half-planes $a_{+}$ to $a_{+}$ and $a_{-}$ to $a_{-}$. Then from
the relationships
$$
\xalignat 2
&[AC]\cong [BD],&&\angle CAB\cong\angle DBF,\\
&p_{a\bold c}(A)=B,&&p_{a\bold c}(B)=F
\endxalignat
$$
and from the fact that the point $C$ and $D$ lies on one side with
respect to the line $a$, we get $p_{a\bold c}(C)=D$. Both propositions
of the theorem are proved.
\qed\enddemo
\mytheorem{9.2} Let $\overrightarrow{AB}$ be a geometric 
realization of a slipping vector $\bold c$ on a line $a$ and 
let $\overrightarrow{CD}$ be a geometric realization of
another slipping vector $\bold d$ on a line $b$. Then the
relationship $p_{a\bold c}=p_{b\bold d}$ is fulfilled if and
only if $\overrightarrow{AB}=\overrightarrow{CD}$.
\endproclaim
     The theorem~\mythetheorem{9.2} is easily derived from the 
theorem~\mythetheorem{9.1} if we take into account the symmetry
and transitivity of the relation of equality of vectors. This
theorem shows that in Euclidean geometry the $a$ plays an auxiliary
role in constructing the mapping $p_{a\bold c}$. This line cal be
replaced by any other line $b$ parallel to $a$ if we replace the
slipping vector $\bold c$ \pagebreak by the free vector $\bold c$.
Therefore, in Euclidean geometry the mapping $p_{a\bold c}$ is denoted 
as $p_{\bold c}$. Here it is called the {\it parallel translation\/}
by the vector $c$.
\mytheorem{9.3} For arbitrary two points $A$ and $B$ in the space 
there is exactly one parallel translation $p_{\bold c}$ taking the 
point $A$ to the point $B$.
\endproclaim
\demo{Proof} In order to prove the existence of the required parallel
translation $p_{\bold c}$ it is sufficient to take the free vector 
$\bold c$ whose geometric realization is $\overrightarrow{AB}$. Now,
if we assume that $p_{\bold d}(A)=B$ for some free vector $\bold d$
with geometric realization $\overrightarrow{CD}$, then from the
theorem~\mythetheorem{9.1} we derive 
$\overrightarrow{AB}=\overrightarrow{CD}$. Hence, $\bold c=\bold d$, 
which proves the uniqueness of the required parallel translation
$p_{\bold c}$.\qed\enddemo
Due to the theorem~\mythetheorem{9.3} we can use the notation
$p_{\sssize AB}$ in order to designate the parallel translation
$p_{\bold c}$ taking the point $A$ to the point $B$. In this
case $\overrightarrow{AB}$ appears to be a geometric realization 
for the free vector $\bold c$.
\head
\SectionNum{10}{199} The group of parallel translations.
\endhead
\rightheadtext{\S\,10. The group of parallel translations.}
\mytheorem{10.1} The mapping $f\!:\Bbb E\to\Bbb E$ is a parallel 
translation if and only if for any two points $X$ and $Y$ the
equality $\overrightarrow{Xf(X)}=\overrightarrow{Yf(Y)}$ is valid.
\endproclaim
\demo{Proof} If $f=p_{\bold c}$ is a parallel translation by a
vector $\bold c$ taking a point $A$ to another point $B$, then
the equality $\overrightarrow{Xf(X)}=\overrightarrow{Yf(Y)}$ is
derived from the equalities 
$$
\xalignat 2
&\overrightarrow{Xf(X)}=\overrightarrow{AB},
&&\overrightarrow{Yf(Y)}=\overrightarrow{AB},
\endxalignat
$$
which follow from the theorem~\mythetheorem{9.1}.\par
     Conversely, assume that $f$ is a mapping such that for any
two points $X$ and $Y$ the equality $\overrightarrow{Xf(X)}=
\overrightarrow{Yf(Y)}$ is valid. Let's fix some point $A$ and 
denote $B=f(A)$. Let $\bold c=\overrightarrow{AB}$. \pagebreak
Then for any point $X\in\Bbb E$ we have the equalities
$$
\xalignat 2
&\hskip -2em
\overrightarrow{Xf(X)}=\overrightarrow{AB},
&&\overrightarrow{Xp_{\bold c}(X)}=\overrightarrow{AB}.
\mytag{10.1}
\endxalignat
$$
The first of the equalities \mythetag{10.1} follows from 
$\overrightarrow{Xf(X)}=\overrightarrow{Yf(Y)}$ by substituting
$Y=A$, the second one is the consequence of the 
theorem~\mythetheorem{9.1}. From \mythetag{10.1} we derive
$\overrightarrow{Xf(X)}=\overrightarrow{Xp_{\bold c}(X)}$, which
yields $f(X)=p_{\bold c}(X)$ for all $X\in\Bbb E$. Hence, the 
mappings $f$ and $p_{\bold c}$ do coincide. The theorem is proved.
\qed\enddemo
\mytheorem{10.2} The composition of two parallel translations is
a parallel translation.
\endproclaim
\demo{Proof}\parshape 13 0cm 10cm 0cm 10cm 4.2cm 5.8cm 4.2cm 5.8cm 
4.2cm 5.8cm 4.2cm 5.8cm 4.2cm 5.8cm 4.2cm 5.8cm 4.2cm 5.8cm 
4.2cm 5.8cm 4.2cm 5.8cm 4.2cm 5.8cm 0cm 10cm
     Let $\bold a$ and $\bold b$ be two free vectors defining two
parallel translations $p_{\bold a}$ and $p_{\bold b}$. Let $X$ and $X'$
be two arbitrary points in the space. 
\vadjust{\vskip 5pt\hbox to -5pt{\kern 0pt
\includegraphics{Oris67.eps}\hss}\vskip -5pt}We denote
$$
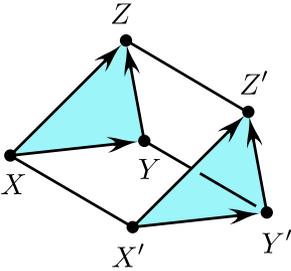
\xalignat 2
&Y=p_{\bold a}(X),&&Z=p_{\bold b}(Y),\\
&Y'=p_{\bold a}(X'),&&Z'=p_{\bold b}(Y').
\endxalignat
$$
If $f=p_{\bold b}\compos p_{\bold a}$, then $Z=f(X)$ and $Z'=f(X')$.
Let's apply the theorem~\mythetheorem{10.1} to the mapping $p_{\bold a}$.
This yields the equality $\overrightarrow{XY}=\overrightarrow{X'Y'}$. 
Similarly, applying the theorem~\mythetheorem{10.1} to the
mapping $p_{\bold b}$, we get $\overrightarrow{YZ}=\overrightarrow{Y'Z'}$.
Now let's use the theorem~\mythetheorem{8.4}. Due to this theorem 
$\overrightarrow{XY}=\overrightarrow{X'Y'}$ implies
$\overrightarrow{XX'}=\overrightarrow{YY'}$ and 
$\overrightarrow{YZ}=\overrightarrow{Y'Z'}$ implies
$\overrightarrow{YY'}=\overrightarrow{ZZ'}$. Using the transitivity
of the relation of equality for vectors, we get
$\overrightarrow{XX'}=\overrightarrow{ZZ'}$. Applying the
theorem~\mythetheorem{8.4} once more, we derive the relationship 
$\overrightarrow{XZ}=\overrightarrow{X'Z'}$. Let's write this relationship
as 
$$
\hskip -2em
\overrightarrow{Xf(X)}=\overrightarrow{X'f(X')}.
\mytag{10.2}
$$
Now, since $X$ and $X'$ are two arbitrary points, \pagebreak
applying the theorem!\mythetheorem{10.1} to the relationship
\mythetag{10.2}, we prove that the mapping $f=p_{\bold b}\compos 
p_{\bold a}$ is a parallel translation.\qed\enddemo
\mytheorem{10.3} Any two mappings of parallel translation are
commutative: $p_{\bold a}\compos p_{\bold b}=p_{\bold b}\compos
p_{\bold a}$.
\endproclaim
\demo{Proof} Let's choose some point $A$ in the space and denote
$B=p_{\bold b}(A)$, $D=p_{\bold a}(B)$, $C=p_{\bold a}(A)$. Applying
the theorem~\mythetheorem{10.1} to the mapping $p_{\bold a}$, we get
$\overrightarrow{AC}=\overrightarrow{BD}$. Hence, we can apply the
theorem~\mythetheorem{8.4}. From this theorem we derive
$\overrightarrow{AB}=\overrightarrow{CD}$. Then $p_{\bold b}(C)=D$,
which follows from the theorem~\mythetheorem{9.1}. As a result for the
mappings $p_{\bold a}\compos p_{\bold b}$ and
$p_{\bold b}\compos p_{\bold a}$ we get
$$
\xalignat 2
&\hskip -2em
p_{\bold a}\compos p_{\bold b}(A)=D,
&&p_{\bold b}\compos p_{\bold a}(A)=D.
\mytag{10.3}
\endxalignat
$$
According to the theorem~\mythetheorem{10.3} the composition
$p_{\bold a}\compos p_{\bold b}$ and the composition
$p_{\bold b}\compos p_{\bold a}$ both are parallel translations
and, as we see in \mythetag{10.3}, they both take the point 
$A$ to the point $D$. Due to the theorem~\mythetheorem{9.3} they
should coincide: $p_{\bold a}\compos p_{\bold b}=p_{\bold b}\compos
p_{\bold a}$.\qed\enddemo
\parshape 7 0cm 10cm 0cm 10cm 0cm 10cm 0cm 10cm 0cm 10cm 0cm 10cm
4.5cm 5.5cm  
     The theorem~\mythetheorem{10.2} shows that in Euclidean geometry
the set of parallel translations is closed with respect to the 
composition. It is a group with respect to this operation (see
definition~\mythedefinitionchapter{4.2}{3} in
Chapter~\uppercase\expandafter{\romannumeral 3}). The unity of this
group is the identical mapping $\id=p_{\bold 0}$ interpreted as the
parallel translation by the zero vector. According to the
theorem~\mythetheorem{10.3} the group of parallel translations is
an Abelian group.\par
\parshape 11 4.5cm 5.5cm 4.5cm 5.5cm 4.5cm 5.5cm 4.5cm 5.5cm 4.5cm 5.5cm 
4.5cm 5.5cm 4.5cm 5.5cm 4.5cm 5.5cm 4.5cm 5.5cm 4.5cm 5.5cm 0cm 10cm
     The theorems~\mythetheorem{9.1}, \mythetheorem{9.2}, and 
\mythetheorem{9.3} establishing one-to-one correspond\-ence between
free vectors and parallel tran\-slations
\vadjust{\vskip 5pt\hbox to -5pt{\kern 0pt
\includegraphics{Oris68.eps}\hss}\vskip -5pt}enable 
us to define the operation of addition for vectors by means of the 
formula
$$
p_{\bold a}\compos p_{\bold b}=p_{\bold a+\bold b}.
$$
Since $p_{\bold a}$ commute with $p_{\bold b}$, the addition of vectors
is a commutative operation, i\.\,e\. $\bold a+\bold b=\bold b+\bold a$.
Passing from the free vectors $\bold b$ and $\bold a$ to their geometric
realizations $\overrightarrow{AB}$ and $\overrightarrow{BD}$, we obtain
$$
\hskip -2em
\overrightarrow{AB}+\overrightarrow{BD}=\overrightarrow{AD},
\mytag{10.4}
$$
where $\overrightarrow{AD}$ is a geometric realization of the vector
$\bold c=\bold a+\bold b$. The relationship \mythetag{10.4} is called
the {\it triangle rule\/} for adding vectors. The triangle $ABC$ on
Fig\.~10.2 can be complemented up to a parallelogram, which yields
$$
\hskip -2em
\overrightarrow{AB}+\overrightarrow{AC}=\overrightarrow{AD}.
\mytag{10.5}
$$
The relationship \mythetag{10.5} is known as the {\it parallelogram
rule\/} for adding vectors.\par
     In Chapter~\uppercase\expandafter{\romannumeral 4} we have
formulated three theorems describing the properties of congruent
translations by vectors along straight lines in absolute geometry.
These are the theorems~\mythetheoremchapter{14.3}{4}, 
\mythetheoremchapter{14.4}{4}, and \mythetheoremchapter{14.5}{4}. 
In Euclidean geometry they are formulated as follows.
\mytheorem{10.4} Let $p_{\bold c}$ and $p_{\bold d}$ be two parallel
translations. If their composition $f=p_{\bold c}\compos p_{\bold d}$ 
has a stable point $O$, then it is the identical mapping: $f=\id$.
\endproclaim
\mytheorem{10.5} If the composition of three or more parallel
translation has a stable point, then it is the identical mapping.
\endproclaim
     The rotations mentioned in the 
theorems~\mythetheoremchapter{14.3}{4}, \mythetheoremchapter{14.4}{4}, 
and  \mythetheoremchapter{14.5}{4} in Euclidean geometry turn to the
trivial rotation by the zero angle, which coincide with the identical
mapping. The proof of the theorems~\mythetheorem{10.4} and
\mythetheorem{10.5} is obvious since the composition of any number of
parallel translations is a parallel translation. A parallel translation
with a stable point is the translation by the zero vector. It coincides
with the identical mapping.
\head
\SectionNum{11}{202} Homothety and similarity.
\endhead
\rightheadtext{\S\,11. Homothety and similarity.}
     The homothety and similarity mappings on straight lines were
introduced in \S\,8 of Chapter~\uppercase\expandafter{\romannumeral 5}.
The Homothety mapping can be defined in the whole space either. Let's
choose some point $O$, which is called the {\it center\/} of homothety,
and some real number, which is called the {\it homothety factor}. The
homothety mapping itself  $h_{k\sssize O}\!:\Bbb E\to\Bbb E$ is defined 
as follows:
\roster
\rosteritemwd=0pt
\item for the point $O$ we set $h_{k\sssize O}(O)=O$;
\item if $X\neq O$, we draw the line $OX$, take the vector 
      $\overrightarrow{OX}$ on this line, multiply it by the number
      $k$, lay the vector $\overrightarrow{OY}=k\cdot
      \overrightarrow{OX}$ on the line $OX$, and then assign
      $Y=h_{k\sssize O}(X)$.
\endroster
The above manipulations defining a homothety mapping can be performed 
in Euclidean and in absolute geometries. However, in absolute geometry 
we could not prove any properties of a homothety mapping that make it
worth considering.\par
     Let $h_{k\sssize O}$ be the homothety mapping with the center
at the point $O$ and with the factor $k$. For $k=1$ this mapping
coincides with the identical mapping $h_{k\sssize O}=\id$, while
for $k=-1$ it turns to be the inversion $h_{k\sssize O}=i_{\sssize O}$.
Moreover, if $k=p\cdot q$, then $h_{k\sssize O}=h_{p\sssize O}\compos
h_{q\sssize O}$. In particular, if $k=-q$, then $h_{k\sssize O}=
h_{q\sssize O}\compos i_{\sssize O}$. Therefore, we often can consider
homotheties with positive factors $k>0$ only.
\mytheorem{11.1} Let $f=h_{k\sssize O}$ be the homothety with the 
factor $k$ and with the center at a point $O$. Then for any two
points $X$ and $Y$ the relationship $|f(X)f(Y)|=|k|\cdot|XY|$ is
fulfilled and the line connecting the points $f(X)$ and $f(Y)$ is
parallel to the line connecting the points $X$ and $Y$.
\endproclaim
\demo{Proof} In the case where the points $X$, $Y$, and $O$ lie on
one straight line the proposition of the theorem~\mythetheorem{11.1}
follows from the theorem~\mythetheoremchapter{8.1}{5} in
Chapter~\uppercase\expandafter{\romannumeral 5}.\par
     Let's consider the case where the points $X$ and $Y$ do not lie
on one straight line. Let's draw the lines $OX$ and $OY$. To each real
number $k$ we associate two points $X(k)$ and $Y(k)$. We define these
points in the following way:
$$
\xalignat 2
&X(k)=h_{k\sssize O}(X),
&&Y(k)=h_{k\sssize O}(Y).
\endxalignat
$$
Then $X=X(1)$ and $Y=Y(1)$. For the beginning we consider the integer
values $k=1,\,2,\,4,\,8,\ldots$ and the rational values $k=1/2,\,1/4,
\,1/8,\,\ldots$ being integer exponentials of the number two. 
\vadjust{\vskip 5pt\hbox to -5pt{\kern 10pt
\includegraphics{Oris69.eps}\hss}\vskip 125pt}Note that the 
segment $[X(1/2)Y(1/2)]$
is the midsegment of the triangle $XOY$, the segment $[XY]$ is the 
midsegment
in the triangle $X(2)OY(2)$, the segment $[X(2)Y(2)]$ is the midsegment
in the triangle $X(4)OY(4)$ and so on. Hence, the segments of the form
$[X(2^q)Y(2^q)]$ are parallel, their lengths are given by the formula
$|X(2^q)Y(2^q)|=2^q\cdot|XY|$. In the other words, the proposition of
the theorem is valid for the homotheties with the factors $k=2^q$,
where $q\in\Bbb Z$.\par
     In the second step we prove the proposition of the theorem for all
binary-rational values $k>0$. Each such $k$ is represented as $k=2^q
\cdot n$, where$n$ is some odd positive integer. Let's convert $n$ into
the binary format, i\.\,e\. into the system with the base $2$:
$$
n=a_m\cdot 2^m+a_{m-1}\cdot 2^{m-1}+\ldots+a_1\cdot 2+a_0.
$$
Here $a_0=1$, wile the numbers $a_1,\,\ldots,\,a_m$ take the values $0$ 
or $1$. Let's denote my $\mu(n)$ the number of unities in the binary 
representation of the number $n$ and perform the induction on this
number. If $\mu(n)=1$, then $n=1$ since $n$ is odd. In this case
$k=2^q$. For such a case the proposition of the theorem was already 
proved in the first step.\par
      Assume that the proposition of the theorem is proved for all 
$k=2^q\cdot n$, such that $\mu(n)<s$. Let's consider some odd positive
integer $n$ such that $\mu(n)=s$. Then we have
$$
\hskip -2em
n=2^p\cdot\tilde n+1=\frac{2^{p+1}\cdot\tilde n+2}{2},
\mytag{11.1}
$$
where $p>0$, while $\tilde n$ is odd and $\mu(\tilde n)=s-1$. Let's
denote $k_1=2^{q+1}\cdot 1$ and $k_2=2^{q+p+1}\cdot\tilde n$. For $k$ 
from \mythetag{11.1} we derive
$$
\hskip -2em
k=\frac{k_1+k_2}{2}.
\mytag{11.2}
$$
For the numbers $k_1$ and $k_2$ the proposition of the theorem
is fulfilled by the inductive hypothesis. Due to \mythetag{11.2} 
we conclude that the segment $[X(k)Y(k)]$ is the midsegment of the 
trapezium $X(k_1)X(k_2)Y(k_2)Y(k_1)$. Hence, the line $X(k)Y(k)$
is parallel to the line $XY$, while the length of the segment
$[X(k)Y(k)]$ is calculated in the following way:
$$
\gathered
|X(k)Y(k)|=\frac{|X(k_1)Y(k_1)|+|X(k_2)Y(k_2)|}{2}=\\
\vspace{1.5ex}
\qquad=\frac{k_1\cdot |XY|+k_2\cdot |XY|}{2}=\frac{(k_1+k_2)
\cdot |XY|}{2}.
\endgathered
\quad
\mytag{11.3}
$$
From \mythetag{11.2} and \mythetag{11.3} we derive $|X(k)Y(k)|=k\cdot
|XY|$. Thus, the inductive step from $\mu(n)<s$ to $\mu(n)=s$ is performed.
It means that we have proved the proposition of the theorem for all 
positive binary-rational values of $k$.\par
    Now let's consider some positive real value of $k$. Let $a_m$ and
$b_m$ be the binary-rational approximations of the number $k$:
$$
\hskip -2em
\frac{p_m}{2^m}=a_m\leqslant k<b_m=\frac{p_m+1}{2^m}.
\mytag{11.4}
$$
Then for $|X(k)Y(k)|$ from the triangle inequality we derive
$$
|X(k)Y(k)|<|X(k)X(b_m)|+|X(b_m)Y(b_m)|+|Y(b_m)Y(k)|
$$
(see theorem~\mythetheoremchapter{2.5}{5} in
Chapter~\uppercase\expandafter{\romannumeral 5}). But the lengths of
the segments in the right hand side of the above inequality are known:
$$
\align
&\hskip -2em
|X(k)X(b_m)|=(b_m-k)\cdot|OX|,\\
&\hskip -2em
|X(b_m)Y(b_m)|=b_m\cdot |XY|,
\mytag{11.5}\\
&\hskip -2em
|Y(b_m)Y(k)|=(b_m-k)\cdot |OY|.
\endalign
$$
Due to \mythetag{11.4} we have $b_m-k<2^{-m}$. Combining this inequality
with \mythetag{11.5}, we can transform the above estimate for the length
of the segment $[X(k)Y(k)]$ to the following form:
$$
|X(k)Y(k)|< k\cdot |XY|+2^{-m}\cdot(|OX|+|XY|+|YO|).
$$
In a similar way, using the triangle inequality, we derive
$$
|X(a_m)Y(a_m)|\leqslant |X(a_m)X(k)|+|X(k)Y(k)|+|Y(k)Y(a_m)|,
$$
which can be then transformed to
$$
|X(k)Y(k)|\geqslant k\cdot |XY|-2^{-m}\cdot(|OX|+|XY|+|YO|).
$$\par
\parshape 2 0cm 10cm 5.2cm 4.8cm \noindent Now let's take into
account that $m$ is an arbitrary positive integer. For $|X(k)Y(k)|$. 
This yields the required relationship $|X(k)Y(k)|=k\cdot |XY|$.\par
\parshape 8 5.2cm 4.8cm 5.2cm 4.8cm 5.2cm 4.8cm 5.2cm 4.8cm 5.2cm 4.8cm
5.2cm 4.8cm 5.2cm 4.8cm 0cm 10cm
     The rest is to prove that the lines $X(k)Y(k)$ and $XY$ are parallel.
\vadjust{\vskip 5pt\hbox to -5pt{\kern 5pt
\includegraphics{Oris70.eps}\hss}\vskip -5pt}We do it by
contradiction. Assume that the lines $X(k)Y(k)$ and $XY$ are not
parallel. Let's draw the line parallel to the line $XY$ through the
point $Y(k)$ and denote by $\tilde X$ the intersection point of the
lines $a$ and $OX$. The point $\tilde X$ is associated with some
number $q$ such that $\tilde X=X(q)$. From the inequalities 
\mythetag{11.4} we find that the points $X(k)$ and $Y(k)$ lie on
the lateral sides of the trapezium $X(a_m)X(b_m)Y(b_m)Y(a_m)$. 
The line $a$ parallel to the bases of this trapezium and crossing 
the lateral side $[Y(b_m)Y(a_m)]$ should necessarily cross the other
lateral side $[X(a_m)X(b_m)]$. This fact follows from the 
Pasch's axiom~\mytheaxiom{A12} applied to the triangles $Y(a_m)X(b_m)
Y(b_m)$ and $X(a_m)Y(a_m)X(b_m)$. Hence, the point $X(q)$, like the 
point $X(k)$, lies within the segment $[X(a_m)X(b_m)]$. This yields 
the following estimate:
$$
|X(k)X(q)|\leqslant |X(a_m)X(b_m)|=2^{-m}\cdot |OX|.
$$
If our assumption is valid, then $X(k)\neq X(q)$ and we come to 
a contradiction with the theorem~\mythetheoremchapter{5.1}{5} 
from Chapter~\uppercase\expandafter{\romannumeral 5}.
This contradiction proves the coincidence $X(q)=X(k)$, which means
that the line $X(k)Y(k)$ is parallel to the line $XY$.\par
    Thus, we have proved the theorem~\mythetheorem{11.1} for any
real number $k>0$. It is easy to extend it to the case $k<0$ since
for $k=-q$ the homothety $h_{k\ssize O}$ with the negative factor 
$k<0$ is the composition of the homothety $h_{q\ssize O}$ and the
inversion $i_{\sssize O}$.
\qed\enddemo
\mytheorem{11.2} Let $f=h_{k\sssize O}$ be the homothety with the
factor $k$ and with the center at some point $O$. In this case
\roster
\item if some three points $X$, $Y$, and $Z$ lie on one straight line, 
      then their images $f(X)$, $f(Y)$, and $f(Z)$ lie on one
      straight line and $(X\blacktriangleright Y\blacktriangleleft Z)$
      implies $(f(X)\blacktriangleright f(Y)\blacktriangleleft f(Z))$;
\item if the points $X$, $Y$, and $Z$ do not lie on one straight line
      then the points $f(X)$, $f(Y)$, and $f(Z)$ also do not lie
      on one straight line.
\endroster
\endproclaim
\demo{Proof} Assume that the points $X$, $Y$, and $Z$ lie on one 
straight line so that $(X\blacktriangleright Y\blacktriangleleft Z)$. 
Then we have the equality $|XY|+|YZ|=|XZ|$, which follows from the 
theorem~\mythetheoremchapter{7.3}{5} in
Chapter~\uppercase\expandafter{\romannumeral 5}. Applying the 
theorem~\mythetheorem{11.1}, we derive the equality 
$|f(X)f(Y)|+|f(Y)f(Z)|=|f(X)f(Z)|$. If we assume that the points 
$f(X)$, $f(Y)$, and $f(Z)$ do not lie on one straight line, then we
get the triangle inequality $|f(X)f(Y)|+|f(Y)f(Z)|>|f(X)f(Z)|$,
which is not compatible with $|f(X)f(Y)|+|f(Y)f(Z)|=|f(X)f(Z)|$.
\par
      Let's prove that the relationship $(X\blacktriangleright
Y\blacktriangleleft Z)$ implies the analogous relationship 
$(f(X)\blacktriangleright f(Y)\blacktriangleleft f(Z))$. If it is 
not so, the following two dispositions of points are possible:
$$
\xalignat 2
&(f(Y)\blacktriangleright f(X)\blacktriangleleft f(Z)),
&&(f(X)\blacktriangleright f(Z)\blacktriangleleft f(Y)).
\endxalignat
$$
In the first case we have $|f(X)f(Z)|<|f(Y)f(Z)|$, while in the second
we have $|f(X)f(Z)|<|f(X)f(Y)|$. It is easy to see that none of these
inequalities can be fulfilled simultaneously with the equality
$|f(X)f(Y)|+|f(Y)f(Z)|=|f(X)f(Z)|$, which follows from 
$(X\blacktriangleright Y\blacktriangleleft Z)$. Therefore, we 
actually have the required disposition of points where
$(f(X)\blacktriangleright f(Y)\blacktriangleleft f(Z))$.\par
      Now, conversely, assume that the points $X$, $Y$, and $Z$ do not
lie on one straight line and assume that their images $f(X)$, $f(Y)$,
and $f(Z)$ lie on one straight line so that $(f(X)\blacktriangleright
f(Y)\blacktriangleleft f(Z))$. Then $|f(X)f(Y)|+|f(Y)f(Z)|=|f(X)f(Z)|$.
Using the equalities $|f(X)f(Y)|=|k|\cdot|XY|$, $|f(Y)f(Z)|=|k|\cdot|YZ|$,
$|f(X)f(Z)|=|k|\cdot|XZ|$, which follow from the
theorem~\mythetheorem{11.1}, for the initial points $X$, $Y$, and $Z$ 
we get $|XY|+|YZ|=|XZ|$. However, the points $X$, $Y$, and $Z$ do not 
lie on one straight line and satisfy the triangle inequality 
$|XY|+|YZ|>|XZ|$. This inequality is not compatible with $|XY|+|YZ|=|XZ|$.
The contradiction obtained proves that the points $f(X)$, $f(Y)$, and
$f(Z)$ do not lie on one line either.\qed\enddemo
\mytheorem{11.3} Let $f=h_{k\sssize O}$ be the homothety with the factor 
$k$ and with the center at a point $O$. If three points $X$, $Y$, and $Z$
do no lie on one straight line, then $\angle XYZ\cong\angle f(X) f(Y)f(Z)$.
\endproclaim
\demo{Proof} Assume that $k>0$. Let's draw some arbitrary ray coming out 
from the center of the homothety $h_{k\sssize O}$ and a point $A$ on it
so that the relationship $[OA]\cong [XY]$ is fulfilled. Then we draw another 
ray $q$ coming out from the point $O$ so that $\angle hq\cong
\angle XYZ$. Let's mark a point $B$ on the ray $k$ so that $[OB]\cong [YZ]$.
The three points $O$, $A$, and $B$ form a triangle $AOB$ congruent
to the triangle $XYZ$ (see theorem~\mythetheoremchapter{5.1}{3} in
Chapter~\uppercase\expandafter{\romannumeral 3}). This fact yields
the following relationships:
$$
\xalignat 3
&|XY|=|OA|,&&|YZ|=|OB|,&&|XZ|=|AB|.
\qquad
\mytag{11.6}
\endxalignat
$$
Let's apply the homothety mapping $f=h_{k\sssize O}$ to the points
$O$, $A$, $B$, $X$, $Y$, and $Z$. The point $O$ is a stable point,
i\.\,e\. $f(O)=O$, while the points $A$ and $B$ are taken to the
points $f(A)$ and $f(B)$ lying on the same rays $h$ and $k$ as the
initial points $A$ and $B$. Hence, we have the following relationship: 
$$
\hskip -2em
\angle f(A)Of(B)=\angle hk=\angle AOB\cong\angle XYZ.
\mytag{11.7}
$$
The images of the points $X$, $Y$, and $Z$ form the triangle 
$f(X)f(Y)f(Z)$ whose sides are congruent to the sides of the triangle
$f(A)Of(B)$. This fact follows from
$$
\aligned
&|f(X)f(Y)|=|f(A)O|,\\
&|f(Y)f(Z)|=|f(B)O|,\\
&|f(X)f(Z)|=|f(A)f(B)|,
\endaligned
$$
which, in turn, are derived by applying the theorem~\mythetheorem{11.1}
to the relationships \mythetag{11.6}. Now, applying the 
theorem~\mythetheoremchapter{5.5}{3} from
Chapter~\uppercase\expandafter{\romannumeral 3}, we get the congruence
of the triangles $f(X)f(Y)f(Z)$ and $f(A)Of(B)$. Then 
$\angle f(X)f(Y)f(Z)\cong\angle f(A)Of(B)$. Combining this fact with
\mythetag{11.7}, we get the required relationship $\angle f(X)f(Y)f(Z)
\cong\angle XYZ$.\par
     The case $k<0$ is reduced to the case $k>0$ since for $k=-q$ we have
the the relationship $h_{k\sssize O}=h_{q\sssize O}\compos i_{\sssize O}$,
where $i_{\sssize O}$ is the inversion mapping. The inversion 
$i_{\sssize O}$ is a congruent translation, it maps each angle to a
congruent angle.\qed\enddemo
     The theorems~\mythetheorem{11.1}, \mythetheorem{11.2}, and
\mythetheorem{11.3} show that under a homothety straight lines are 
mapped to straight lines (being parallel to the initial ones), 
segments are mapped to segments, and rays are mapped to rays. 
The quantitative measures of angles are preserved. Therefore,
the homothety preserves the orthogonality of lines. Applying the
theorem~\mythetheoremchapter{1.5}{4} from
Chapter~\uppercase\expandafter{\romannumeral 4}, we can prove that
that under a homothety a plane is mapped to a plane, a half-plane
is mapped to a half-plane, and a half-space --- to a half-space. 
These properties of a homothety are identical to those of a 
congruent translation.\par
     Note that the mapping $f=h_{k\sssize O}$ is bijective, the 
inverse mapping $f^{-1}$ for it is the homothety $h_{q\sssize O}$ 
with the factor $q=1/k$.
\mydefinition{11.1} A mapping $f\!:\Bbb E\to\Bbb E$ is called a
{\it similarity mapping\/} if it admits an expansion
$f=h_{k\sssize O}\compos\varphi$, where $\varphi$ is a mapping
of congruent translation, while $h_{k\sssize O}$ is a homothety
with the coefficient $k\neq 0$. The number $|k|$ is called 
the {\it similarity factor\/} of such a mapping.
\enddefinition
\myexercise{11.1} Show that the composition of two similarity
mappings is a similarity mapping.
\endproclaim
\myexercise{11.2} Show that the mapping inverse to a similarity
mapping is a similarity mapping as well.
\endproclaim
\mydefinition{11.2} Two geometric forms $\Phi_1$ and $\Phi_2$ are 
called {\it similar\/} if there is a similarity mapping 
$f\!:\Bbb E\to\Bbb E$ establishing a one-to-one correspondence betwen
the points of these two forms.
\enddefinition
\mytheorem{11.4} If for two triangles $ABC$ and $\tilde A\tilde B
\tilde C$ the conditions $|AB|:|\tilde A\tilde B|=|AC|:|\tilde A
\tilde C|$ and $\angle BAC\cong\angle\tilde B\tilde A\tilde C$ are
fulfilled, then the triangle $ABC$ is similar to the triangle
$\tilde A\tilde B\tilde C$.
\endproclaim
\mytheorem{11.5} If some two angles of a triangle $ABC$ are congruent 
to the corresponding angles of another triangle $\tilde A\tilde B
\tilde C$, then the triangle $ABC$ is congruent to the triangle 
$\tilde A \tilde B\tilde C$.
\endproclaim
\mytheorem{11.6} If for triangles $ABC$ and $\tilde A\tilde B\tilde C$
the conditions $|AB|:|\tilde A\tilde B|=|AC|:|\tilde A\tilde C|=|BC|
:|\tilde B\tilde C|$ are fulfilled, then the triangle $ABC$ is similar
to the triangle $\tilde A\tilde B\tilde C$.
\endproclaim
\myexercise{11.3} The theorems~\mythetheorem{11.4}, \mythetheorem{11.5},
and \mythetheorem{11.6} are known as the similarity criteria for triangles. 
Prove these theorems relying upon the 
theorems~\mythetheoremchapter{5.1}{3}, \mythetheoremchapter{5.2}{3}, 
and \mythetheoremchapter{5.5}{3} from 
Chapter~\uppercase\expandafter{\romannumeral 3}.
\endproclaim
\head
\SectionNum{12}{211} Multiplication of vectors by a number.
\endhead
\rightheadtext{\S\,12. Multiplication of vectors by a number.}
     The multiplication of slipping vectors by a number is given 
by the definition~\mythedefinitionchapter{8.1}{5} in
Chapter~\uppercase\expandafter{\romannumeral 5}. In Euclidean geometry
this operation can be extended to the set of free vectors.
\mydefinition{12.1} Assume that a geometric vector 
$\overrightarrow{AB}$ is given. The vector $\overrightarrow{CD}$ 
with the length $|CD|=|k|\cdot |AB|$ is called the {\it product of 
the vector $\overrightarrow{AB}$ by the number $k\neq 0$} and written as
$\overrightarrow{CD}=k\cdot\overrightarrow{AB}$ if it is codirected
to $\overrightarrow{AB}$ for $k>0$ and if it is oppositely
directed to $\overrightarrow{AB}$ for $k<0$.
\enddefinition
     The result of multiplying a vector $\overrightarrow{AB}$ by
a number $k$ is not unique since the definition~\mythedefinition{12.1}
does not fix the position of the vector $\overrightarrow{CD}$ in the
space. However, there is the following theorem.
\mytheorem{12.1} If two vectors $\overrightarrow{CD}$ and 
$\overrightarrow{C'D'}$ are obtained through multiplying some vector
$\overrightarrow{AB}$ by a number $k$, then they are equal in the
sense of the definition~\mythedefinition{8.2}.
\endproclaim
\demo{Proof} Indeed, the lengths of the segment $[CD]$ and $[C'D']$ 
are equal to each other since $|CD|=|k|\cdot |AB|$ and $|C'D'|=|k|
\cdot |AB|$. Hence, we have the congruence of segments $[CD]\cong
[C'D']$.\par
     For $k>0$ we have 
$\overrightarrow{AB}\upuparrows\overrightarrow{CD}$ and 
$\overrightarrow{AB}\upuparrows\overrightarrow{C'D'}$. These 
two relationships yield $\overrightarrow{CD}\upuparrows
\overrightarrow{C'D'}$. For $k<0$ from $\overrightarrow{AB}
\updownarrows\overrightarrow{CD}$ and $\overrightarrow{AB}
\updownarrows\overrightarrow{C'D'}$ we get the relationships 
$\overrightarrow{BA}\upuparrows\overrightarrow{CD}$ and  
$\overrightarrow{BA}\upuparrows\overrightarrow{C'D'}$, 
which also yield $\overrightarrow{CD}\upuparrows
\overrightarrow{C'D'}$. The conditions
$\overrightarrow{CD}\upuparrows\overrightarrow{C'D'}$ and
$[CD]\cong [C'D']$ are the very conditions that mean
$\overrightarrow{CD}=\overrightarrow{C'D'}$ in the sense of the
definition~\mythedefinition{8.2}.
\qed\enddemo
\mytheorem{12.2} The equality $\overrightarrow{AB}=\overrightarrow{A'B'}$
implies the equality $\overrightarrow{CD}=\overrightarrow{C'D'}$ for 
vectors $\overrightarrow{CD}$ and $\overrightarrow{C'D'}$ obtained through
multiplying $\overrightarrow{AB}$ and $\overrightarrow{A'B'}$ by a number
$k$.
\endproclaim
\demo{Proof} The equality $\overrightarrow{AB}=\overrightarrow{A'B'}$
means that the vectors $\overrightarrow{AB}$ and $\overrightarrow{A'B'}$
are codirected, while their lengths are equal:
$$
\xalignat 2
&\hskip -2em
\overrightarrow{AB}\upuparrows\overrightarrow{A'B'},
&&|AB|=|A'B'|.
\mytag{12.1}
\endxalignat
$$
From the relationships $\overrightarrow{CD}=k\cdot\overrightarrow{AB}$
and $\overrightarrow{C'D'}=k\cdot\overrightarrow{A'B'}$ for $k>0$ we
derive the relationships
$$
\xalignat 2
&\hskip -2em
\overrightarrow{CD}\upuparrows\overrightarrow{AB},
&&|CD|=|k|\cdot |AB|,\\
\vspace{-1.2ex}
&&&\mytag{12.2}\\
\vspace{-1.2ex}
&\hskip -2em
\overrightarrow{C'D'}\upuparrows\overrightarrow{A'B'},
&&|C'D'|=|k|\cdot |A'B'|
\endxalignat
$$
for the vectors $\overrightarrow{CD}$ and $\overrightarrow{C'D'}$.
Then from \mythetag{12.1} and \mythetag{12.2} we derive that the
vectors $\overrightarrow{CD}$ and $\overrightarrow{C'D'}$ are
codirected, while their lengths are equal. This fact yields
$\overrightarrow{CD}=\overrightarrow{C'D'}$. The case $k<0$ differs
from the case $k>0$ only in that the vectors $\overrightarrow{AB}$ 
and $\overrightarrow{A'B'}$ in \mythetag{12.2} are replaced by the
opposite vectors $\overrightarrow{BA}$ and $\overrightarrow{B'A'}$. 
The theorem is proved.\qed\enddemo
     The theorems~\mythetheorem{12.1} and \mythetheorem{12.2} show
that the multiplication by a number is a correctly defined unambiguous 
operation in the set of free vectors. To a vector $\bold b$ and to a 
number $k\in\Bbb R$ it associates some definite vector $\bold c=k\cdot
\bold b$. For $k=0$ of for $\bold b=\bold 0$ the product $k\cdot\bold b$
is taken to be equal to the zero vector:
$$
\aligned
&0\cdot\bold b=\bold 0\text{ \ for any vector  \ }\bold b,\\
&k\cdot\bold 0=\bold 0\text{ \ for any number \ }k\in\Bbb R.
\endaligned
$$
\mytheorem{12.3} The operation of addition and the operation of
multiplication by a number in the set of free vectors possesses
the following properties:
\roster
\rosteritemwd=5pt
\item commutativity of addition: $\bold a+\bold b=\bold b+\bold a$;
\item associativity of addition: $(\bold a+\bold b)+\bold c=
      \bold a+(\bold b+\bold c)$;
\item there is a zero vector $\bold 0$ such that $\bold a+\bold 0
      =\bold a$ for an arbitrary vector $\bold a$;
\item for any vector $\bold a$ there is an opposite vector $\bold a'$
      such that $\bold a+\bold a'=\bold 0$;
\item distributivity of multiplication by a number with respect 
      to the addition of vectors: 
      $k\cdot (\bold a+\bold b)=k\cdot\bold a+k\cdot\bold b$;
\item distributivity of multiplication by a number with respect to 
      the addition of numbers: $(k+q)\cdot\bold a=k\cdot\bold a+q\cdot
      \bold a$;
\item associativity of multiplication: 
$(k\cdot q)\cdot\bold a=k\cdot(q\cdot\bold a)$;
\item the property of the numeric unity: $1\cdot\bold a=\bold a$.
\endroster
\endproclaim
\myexercise{12.1} Prove those propositions of the 
theorem~\mythetheorem{12.3} which are not already proved.
\endproclaim
\newpage
\setfirstpage
\topmatter
\title
REFERENCES
\endtitle
\endtopmatter
\document
\rightheadtext{References.}
\leftheadtext{References.}

\par
\newpage
\setfirstpage
\topmatter
\title
Contacts
\endtitle
\endtopmatter
\document
\line{\vtop{\hsize 5cm
{\bf Address: }
\medskip\noindent
Ruslan A. Sharipov,\newline
Math. Department,\newline
Bashkir State University,\newline
32 Frunze street,\newline
Ufa 450074, Russia
\medskip
{\bf Phone:}\medskip
\noindent
+7-(347)-273-67-18 (Office)\newline
+7-(917)-476-93-48 (Cell)
}\hss
\vtop{\hsize 4.3cm
{\bf Home address:}\medskip\noindent
Ruslan A. Sharipov,\newline
5 Rabochaya street,\newline
Ufa 450003, Russia
\vskip 1cm
{\bf E-mails:}\medskip
\noindent
r-sharipov\@mail.ru\newline
R\hskip 0.5pt\_\hskip 1.5pt Sharipov\@ic.bashedu.ru
}
}
\bigskip
{\bf URL's:}\medskip
\noindent
\myhref{http://www.geocities.com/r-sharipov/}
{http:/\negskp/www.geocities.com/r-sharipov}\newline
\myhref{http://www.freetextbooks.boom.ru/}
{http:/\negskp/www.freetextbooks.boom.ru}\newline
\myhref{http://sovlit2.narod.ru/}
{http:/\negskp/sovlit2.narod.ru}\newline
\par
\newpage
\setfirstpage
\topmatter
\title
Appendix
\endtitle
\endtopmatter
\document
\rightheadtext{List of publications.}
\leftheadtext{List of publications.}

\hphantom{1111}

\Refs\nofrills{}
\ref\myrefno{1}\by Backelman~I.~Ya.\book Higher geometry
\publ {\tencyr\char '074}Prosveshchenie{\tencyr\char '076}
publishers\publaddr Mos\-cow\yr 1967
\endref
\ref\myrefno{2}\by Efimov~N.~V.\book Higher geometry
\publ {\tencyr\char '074}Nauka{\tencyr\char '076}
publishers\publaddr Moscow\yr 1971
\endref
\ref\myrefno{3}\by Pogorelov~A.~V.\book Geometry
\publ {\tencyr\char '074}Prosveshchenie{\tencyr\char '076}
publishers\publaddr Moscow\yr 1993
\endref
\ref\myrefno{4}\by Hilton~P.~J. \& Wylie~S.\book Homology theory:
an introduction to algebraic topology\publ Cambridge University Press
\publaddr Cambridge\yr 1960
\endref
\ref\myrefno{5}\by Kostrikin~A.~I.\book \  Introduction \ to 
\ algebra \publ \ \ {\tencyr\char '074}Nauka{\tencyr\char '076} 
\ publishers\publaddr \ Mos\-cow\yr 1977
\endref
\ref\myrefno{6}\by Kudryavtsev~L.~D.\book Course of mathematical 
analysis, \rm Vol\.~\uppercase\expandafter{\romannumeral 1},
\uppercase\expandafter{\romannumeral 2}
\publ {\tencyr\char '074}Visshaya Shkola{\tencyr\char '076}
publishers\publaddr Moscow\yr 1985
\endref
\ref\myrefno{7}\by Alexandrov~A.~D.\book \ Convex \ polyhedra
\publ \ {\tencyr\char '074}Gostehizdat{\tencyr\char '076}
\ publishers\publaddr \ Mos\-cow\yr 1950
\endref
\endRefs

\Refs\nofrills{List
of publications by the author\\ for the period 1986--2006.}
\refstyle C
{\bf Part 1. Soliton theory.}\medskip
\ref\myrefno{1}\by Sharipov R. A.\paper Finite-gap analogs of $N$-multiplet 
solutions of the KdV equation\jour Uspehi Mat. Nauk\vol 41\issue 5\yr 1986 
\pages 203--204
\endref
\ref\myrefno{2}\by Sharipov R. A.\paper Soliton multiplets of the
Korteweg-de Vries equation\jour Dokladi AN SSSR\vol 292\yr 1987
\issue 6\pages 1356--1359
\endref
\ref\myrefno{3}\by Sharipov R. A.\paper Multiplet solutions of 
the Kadomtsev-Petviashvili equation on a finite-gap background
\jour Uspehi Mat. Nauk\vol 42\yr 1987\issue 5\pages 221--222
\endref
\ref\myrefno{4}\by Bikbaev R. F. \& Sharipov R. A.\paper Magnetization 
waves in Landau-Lifshits model\jour Physics Letters A\vol 134\yr 1988
\issue 2\pages 105-108\moreref see
\myhref{http://arxiv.org/abs/solv-int/9905008}{solv-int/9905008}
\endref 
\ref\myrefno{5}\by Bikbaev R. F. \& Sharipov R. A.\paper Assymptotics 
as $t\to\infty$ for a solution of the Cauchy problem for the Korteweg-de 
Vries equation in the class of potentials with finite-gap behaviour as
$x\to\pm\infty$\jour Theor\. and Math\. Phys\.\vol 78\yr 1989\issue 3
\pages 345--356 
\endref 
\ref\myrefno{6}\by Sharipov R. A.\paper On integration of the Bogoyavlensky
chains\jour Mat\. zametki\vol 47\yr 1990\issue 1\pages 157--160
\endref
\ref\myrefno{7}\by Cherdantsev I. Yu. \& Sharipov R. A.\paper Finite-gap
solutions of the Bul\-lough-Dodd-Jiber-Shabat equation\jour Theor\. and 
Math\. Phys\.\vol 82\yr 1990\issue 1\pages 155--160
\endref
\ref\myrefno{8}\by Cherdantsev I. Yu. \&  Sharipov R. A.\paper Solitons 
on a finite-gap background in Bullough-Dodd-Jiber-Shabat model\jour 
International\. Journ\. of Modern Physics A\vol 5\yr 1990\issue 5
\pages 3021--3027\moreref see
\myhref{http://arxiv.org/abs/math-ph/0112045}{math-ph/0112045}
\endref
\ref\myrefno{9}\by Sharipov R. A. \& Yamilov R. I.\paper Backlund
transformations and the construction of the integrable boundary value 
problem for the equation\linebreak $u_{xt}=e^u-e^{-2u}$\inbook 
{\tencyr\char '074}Some problems of mathematical physics and asymptotics 
of its solutions{\tencyr\char '076}\publ Institute of mathematics BNC UrO
AN SSSR\publaddr Ufa\yr 1991\pages 66--77\moreref see
\myhref{http://arxiv.org/abs/solv-int/9412001}{solv-int/9412001}
\endref
\ref\myrefno{10}\by Sharipov R. A.\paper Minimal tori in 
five-dimensional sphere in $\Bbb C^3$\jour Theor\. and Math\. 
Phys\.\vol 87\yr 1991\issue 1\pages 48--56\moreref see
\myhref{http://arxiv.org/abs/math.DG/0204253}{math.DG/0204253}
\endref
\ref\myrefno{11}\by Safin S. S. \& Sharipov R. A.\paper Backlund 
autotransformation for the equation $u_{xt}=e^u-e^{-2u}$\jour Theor\. 
and Math\. Phys\.\vol 95\yr 1993\issue 1\pages 146--159 
\endref
\ref\myrefno{12}\by Boldin A. Yu. \& Safin S. S. \& Sharipov R. A.
\paper On an old paper of Tzitzeika and the inverse scattering 
method\jour Journal of Mathematical Physics\vol 34\yr 1993\issue 12
\pages 5801--5809
\endref
\ref\myrefno{13}\by Pavlov M. V. \& Svinolupov S. I. \& Sharipov R. A. 
\paper Invariant criterion of integrability for a system of equations 
of hydrodynamical type\inbook {\tencyr\char '074}Integrability in dynamical 
systems{\tencyr\char '076}\publ Inst. of Math. UrO RAN\publaddr Ufa
\yr 1994\pages 27--48\moreref\jour Funk\. Anal\. i Pril\.\vol 30\yr 1996
\issue 1\pages 18--29\moreref see 
\myhref{http://arxiv.org/abs/solv-int/9407003}{solv-int/9407003}
\endref
\ref\myrefno{14}\by Ferapontov E. V. \& Sharipov R. A.\paper On
conservation laws of the first order for a system of equations of
hydrodynamical type\jour Theor\. and Math\. Phys\.\vol 108\yr 1996
\issue 1\pages 109--128
\endref
\medskip{\bf Part 2. Geometry of the normal shift.}\medskip
\ref\myrefno{1}\by Boldin A. Yu. \& Sharipov R. A.\paper Dynamical 
systems accepting the normal shift\jour Theor\. and Math\. Phys\.
\vol 97\yr 1993\issue 3\pages 386--395\moreref see
\myhref{http://arxiv.org/abs/chao-dyn/9403003}{chao-}
\myhref{http://arxiv.org/abs/chao-dyn/9403003}{dyn/9403003}
\endref
\ref\myrefno{2}\by Boldin A. Yu. \& Sharipov R. A.\paper Dynamical 
systems accepting the normal shift\jour Dokladi RAN\vol 334\yr 1994
\issue 2\pages 165--167
\endref
\ref\myrefno{3}\by Boldin A. Yu. \& Sharipov R. A.\paper Multidimensional 
dynamical systems accepting the normal shift\jour Theor\. and Math\. Phys\.
\vol 100\yr 1994\issue 2\pages 264--269\moreref see
\myhref{http://arxiv.org/abs/patt-sol/9404001}{patt-sol/9404001}
\endref
\ref\myrefno{4}\by Sharipov R. A.\paper Problem of metrizability for 
the dynamical systems accepting the normal shift\jour Theor\. and Math\.
Phys\.\vol 101\yr 1994\issue 1\pages 85--93\moreref see
\myhref{http://arxiv.org/abs/solv-int/9404003}{solv-int/9404003}
\endref
\ref\myrefno{5}\by Sharipov R. A.\paper Dynamical systems accepting the 
normal shift\jour Uspehi Mat\. Nauk\vol 49\yr 1994\issue 4\page 105
\moreref see \myhref{http://arxiv.org/abs/solv-int/9404002}
{solv-int/9404002}
\endref
\ref\myrefno{6}\by Boldin A. Yu. \& Dmitrieva V. V. \& Safin S. S. 
\& Sharipov R. A.\paper Dynamical systems accepting the normal shift on an
arbitrary Riemannian manifold\inbook {\tencyr\char '074}Dynamical systems
accepting the normal shift{\tencyr\char '076}\publ Bashkir State
University\publaddr Ufa\yr 1994\pages 4--19\moreref see also\nofrills
\jour Theor\. and Math\. Phys\.\vol 103\yr 1995\issue 2\pages 256--266
\nofrills\moreref and \myhref{http://arxiv.org/abs/hep-th/9405021}
{hep-th/9405021}
\endref
\ref\myrefno{7}\by Boldin A. Yu. \& Bronnikov A. A. \& Dmitrieva V. V. 
\& Sharipov R. A.\paper Complete normality conditions for the dynamical
systems on Riemannian manifolds\inbook {\tencyr\char '074}Dynamical 
systems accepting the normal shift{\tencyr\char '076}\publ Bashkir State
University\yr 1994\pages 20--30\moreref see also\nofrills\jour Theor\. 
and Math\. Phys\.\vol 103\yr 1995\issue 2\pages 267--275\nofrills
\moreref and \myhref{http://arxiv.org/abs/astro-ph/9405049}
{astro-ph/9405049}
\endref
\ref\myrefno{8}\by Sharipov R. A.\paper Higher dynamical systems accepting 
the normal shift\inbook {\tencyr\char '074}Dynamical systems accepting the
normal shift{\tencyr\char '076}\publ Bashkir State University\yr 1994
\pages 41--65
\endref
\ref\myrefno{9}\by Bronnikov A. A. \& Sharipov R. A.\paper Axially
symmetric dynamical systems accepting the normal shift in $\Bbb R^n$
\inbook {\tencyr\char '074}Integrability in dynamical systems{\tencyr
\char '076}\publ Inst\. of Math\. UrO RAN\publaddr Ufa\yr 1994
\pages 62--69
\endref
\ref\myrefno{10}\by Sharipov R. A.\paper Metrizability by means of
a conformally equivalent metric for the dynamical systems\inbook 
{\tencyr\char '074}Integrability in dynamical systems{\tencyr\char
'076}\publ Inst\. of Math\. UrO RAN\publaddr Ufa\yr 1994\pages 80--90
\moreref see also\nofrills\jour Theor\. and Math\. Phys\.\vol 103
\yr 1995\issue 2\pages 276--282
\endref
\ref\myrefno{11}\by Boldin A. Yu. \& Sharipov R. A.\paper On the 
solution of the normality equations for the dimension $n\geqslant 3$
\jour Algebra i Analiz\vol 10\yr 1998\issue  4\pages 31--61\moreref
see also \myhref{http://arxiv.org/abs/solve-int/9610006}
{solve-int/9610006}
\endref
\ref\myrefno{12}\by Sharipov R. A.\book Dynamical systems admitting 
the normal shift, \rm Thesis for the degree of Doctor of Sciences in 
Russia\publ \myhref{http://arxiv.org/abs/math.DG/0002202}
{math.DG/0002202}\publaddr Electronic archive \myEarXivlink\yr 2000
\pages 1--219
\endref
\ref\myrefno{13}\by Sharipov R. A.\paper Newtonian normal shift in
multidimensional Riemannian geometry\jour Mat\. Sbornik\vol 192
\yr 2001\issue 6\pages 105--144\moreref see also
\myhref{http://arxiv.org/abs/math.DG/0006125}{math.DG}
\myhref{http://arxiv.org/abs/math.DG/0006125}{/0006125}
\endref
\ref\myrefno{14}\by Sharipov R. A.\paper Newtonian dynamical systems
admitting the normal blow-up of points\jour Zap\. semin\. POMI
\vol 280\yr 2001\pages 278--298\moreref see also 
\myhref{http://arxiv.org/abs/math.DG/0008081}{math.DG/0008081}
\endref
\ref\myrefno{15}\by Sharipov R. A.\paper On the solutions of the weak 
normality equations in multidimensional case\jour
\myhref{http://arxiv.org/abs/math.DG/0012110}{math.DG/0012110}
in Electronic archive \myhref{http://arxiv.org}{http:/\negskp/}
\myhref{http://arxiv.org}{arxiv.org}\yr 2000\pages 1--16
\endref
\ref\myrefno{16}\by Sharipov R. A.\paper First problem of globalization 
in the theory of dynamical systems admitting the normal shift of 
hypersurfaces\jour International Journal of Mathematics and Mathematical
Sciences\vol 30\yr 2002\issue 9\pages 541--557\moreref see also
\myhref{http://arxiv.org/abs/math.DG/0101150}{math.DG/0101150}
\endref
\ref\myrefno{17}\by Sharipov R. A.\paper Second problem of globalization 
in the theory of dyna\-mical systems admitting the normal shift of 
hypersurfaces\jour \myhref{http://arxiv.org/abs/math.DG/0102141}
{math.DG} \myhref{http://arxiv.org/abs/math.DG/0102141}{/0102141}
in Electronic archive \myEarXivlink\yr 2001\pages 1--21 
\endref
\ref\myrefno{18}\by Sharipov R. A.\paper A note on Newtonian, Lagrangian,
and Hamiltonian dynamical systems in Riemannian manifolds\jour
\myhref{http://arxiv.org/abs/math.DG/0107212}{math.DG/0107212} in
Electronic archive \myEarXivlink\yr 2001\pages 1--21
\endref
\ref\myrefno{19}\by Sharipov R. A.\paper Dynamical systems admitting 
the normal shift and wave equations\jour Theor\. and Math\. Phys\.
\vol 131\yr 2002\issue 2\pages 244--260\moreref see also
\myhref{http://arxiv.org/abs/math.DG/0108158}{math.DG/0108158}
\endref
\ref\myrefno{20}\by Sharipov R. A.\paper Normal shift in general 
Lagrangian dynamics\jour \myhref{http://arxiv.org/abs/math.DG/0112089}
{math.DG} \myhref{http://arxiv.org/abs/math.DG/0112089}{/0112089} 
in Electronic archive \myEarXivlink\yr 2001\pages 1--27
\endref
\ref\myrefno{21}\by Sharipov R. A\paper Comparative analysis for a pair 
of dynamical systems one of which is Lagrangian\jour 
\myhref{http://arxiv.org/abs/math.DG/0204161}{math.DG/0204161}
in Electronic archive \myhref{http://arxiv.org}{http:/\negskp/}
\myhref{http://arxiv.org}{arxiv.org}\yr 2002\pages 1--40
\endref
\ref\myrefno{22}\by Sharipov R. A.\paper On the concept of a normal 
shift in non-metric geometry\jour 
\myhref{http://arxiv.org/abs/math.DG/0208029}{math.DG/0208029}
in Electronic archive \myEarXivlink\yr 2002\pages 1--47
\endref
\ref\myrefno{23}\by Sharipov R. A.\paper $V$-representation for 
the normality equations in geometry of generalized Legendre 
transformation\jour
\myhref{http://arxiv.org/abs/math.DG/0210216}{math.DG/0210216}
in Electronic archive \myEarXivlink\yr 2002\pages 1--32
\endref
\ref\myrefno{24}\by Sharipov R. A.\paper On a subset of the normality
equations describing a generalized Legendre transformation\jour
\myhref{http://arxiv.org/abs/math.DG/0212059}{math.DG/0212059}
in Electronic ar\-chive \yr 2002\pages 1--19
\endref
\medskip{\bf Part 3. Several complex variables.}\medskip
\ref\myrefno{1}\by Sharipov R. A. \& Sukhov A. B. On $CR$-mappings 
between algebraic Cauchy-Riemann manifolds and the separate algebraicity 
for holomorphic functions\jour Trans\. of American Math\. Society
\vol 348\yr 1996\issue 2\pages 767--780\moreref see also\nofrills
\jour Dokladi RAN\vol 350\yr 1996\issue 4\pages 453--454
\endref
\ref\myrefno{2}\by Sharipov R. A. \& Tsyganov E. N. On the separate 
algebraicity along families of algebraic curves\book Preprint of Baskir 
State University\publaddr Ufa\yr 1996\pages 1-7\moreref see also\nofrills 
\jour Mat\. Zametki\vol 68\yr 2000\issue 2\pages 294--302
\endref
\medskip{\bf Part 4. Symmetries and invariants.}\medskip
\ref\myrefno{1}\by Dmitrieva V. V. \& Sharipov R. A.\paper On the
point transformations for the second order differential equations
\jour \myhref{http://arxiv.org/abs/solv-int/9703003}{solv-int/9703003}
in Electronic archive \myEarXivlink\yr 1997\pages 1--14 
\endref
\ref\myrefno{2}\by Sharipov R. A.\paper On the point transformations 
for the equation $y''=P+3\,Q\,y'+3\,R\,{y'}^2+S\,{y'}^3$\jour
\myhref{http://arxiv.org/abs/solv-int/9706003}{solv-int/9706003}
in Electronic archive \myhref{http://arxiv.org}
{http:/\negskp/}\linebreak\myhref{http://arxiv.org}{arxiv.org}
\yr 1997\pages 1--35\moreref see also\nofrills\jour Vestnik BashGU
\vol 1(I)\yr 1998\pages 5--8
\endref
\ref\myrefno{3}\by Mikhailov O. N. \& Sharipov R. A.\paper On the 
point expansion for a certain class of differential equations of 
the second order\jour Diff\. Uravneniya\vol 36\yr 2000\issue 10
\pages 1331--1335\moreref see also 
\myhref{http://arxiv.org/abs/solv-int/9712001}{solv-int/9712001}
\endref
\ref\myrefno{4}\by Sharipov R. A.\paper Effective procedure of 
point-classification for the equation $y''=P+3\,Q\,y'+3\,R\,{y'}^2
+S\,{y'}^3$\jour \myhref{http://arxiv.org/abs/math.DG/9802027}
{math.DG/9802027} in Electronic archive \myEarXivlink\yr 1998
\pages 1--35
\endref
\ref\myrefno{5}\by Dmitrieva V. V. \& Gladkov A. V. \& Sharipov R. A.
\paper On some equations that can be brought to the equations of 
diffusion type\jour Theor\. and Math\. Phys\.\vol 123\yr 2000
\issue 1\pages 26--37\moreref see also 
\myhref{http://arxiv.org/abs/math.AP/9904080}{math.AP/9904080}
\endref
\ref\myrefno{6}\by Dmitrieva V. V. \& Neufeld E. G. \& Sharipov R. A. 
\& Tsaregorod\-tsev~A.~ A.\paper On a point symmetry analysis for 
generalized diffusion type equations\jour \ 
\myhref{http://arxiv.org/abs/math.AP/9907130}{math.AP/9907130} \
in \ Electronic \ archive \ \myEarXivlink \yr 1999\pages 1--52
\endref
\medskip{\bf Part 5. General algebra.}\medskip
\ref\myrefno{1}\by Sharipov R. A\paper Orthogonal matrices with 
rational components in composing tests for High School students
\jour \myhref{http://arxiv.org/abs/math.GM/0006230}{math.GM/0006230}
in Electronic archive \myEarXivlink\yr 2000\pages 1--10
\endref
\ref\myrefno{2}\by Sharipov R. A.\paper On the rational extension 
of Heisenberg algebra\jour \myhref{http://arxiv.org/abs/math.RA/0009194}
{math.} \myhref{http://arxiv.org/abs/math.RA/0009194}{RA/0009194} 
in Electronic archive \myEarXivlink\yr 2000\pages 1--12
\endref
\ref\myrefno{3}\by Sharipov R. A\paper An algorithm for generating
orthogonal matrices with rational elements\jour 
\myhref{http://arxiv.org/abs/cs.MS/0201007}{cs.MS/0201007} in
Electronic archive \myEarXivlink\yr 2002\pages 1--7
\endref
\medskip{\bf Part 6. Condensed matter physics.}\medskip
\ref\myrefno{1}\by Lyuksyutov S. F. \& Sharipov R. A.\paper Note 
on kinematics, dynamics, and thermodynamics of plastic glassy media
\jour\myhref{http://arxiv.org/abs/cond-mat/0304190}{cond-mat/0304190}
in Electronic archive \myEarXivlink\yr 2003\pages 1--19
\endref
\ref\myrefno{2}\by Lyuksyutov S. F. \& Sharipov R. A. \& Sigalov G. 
\& Paramonov P. B.\paper Exact analytical solution for electrostatic 
field produced by biased atomic force microscope tip dwelling above 
dielectric-conductor bilayer\jour 
\myhref{http://arxiv.org/abs/cond-mat/0408247}{cond-}
\myhref{http://arxiv.org/abs/cond-mat/0408247}{mat/0408247}
in Electronic archive \myEarXivlink\yr 2004\pages 1--6
\endref
\ref\myrefno{3}\by Lyuksyutov S. F. \& Sharipov R. A.\paper Separation 
of plastic deformations in polymers based on elements of general nonlinear
theory\jour \myhref{http://arxiv.org/abs/cond-mat/0408433}{cond-mat}
\myhref{http://arxiv.org/abs/cond-mat/0408433}{/0408433}
in Electronic archive \myEarXivlink\yr 2004\pages 1--4
\endref
\ref\myrefno{4}\by Comer J. \& Sharipov R. A.\paper A note on the
kinematics of dislocations in crystals\jour
\myhref{http://arxiv.org/abs/math-ph/0410006}{math-ph/0410006}
in Electronic archive \myEarXivlink\yr 2004\pages 1--15
\endref
\ref\myrefno{5}\by Sharipov R. A.\paper Gauge or not gauge?\nofrills
\jour \myhref{http://arxiv.org/abs/cond-mat/0410552}{cond-mat/0410552}
in Electronic archive \myEarXivlink\yr 2004\pages 1--12
\endref
\ref\myrefno{6}\by Sharipov R. A.\paper Burgers space versus real space 
in the nonlinear theory of dislocations\jour
\myhref{http://arxiv.org/abs/cond-mat/0411148}{cond-mat/0411148}
in Electronic archive \myEarXivlink\yr 2004\pages 1--10
\endref
\ref\myrefno{7}\by Comer J. \& Sharipov R. A.\paper On the geometry 
of a dislocated medium\jour
\myhref{http://arxiv.org/abs/math-ph/0502007}{math-ph/0502007}
in Electronic archive \myEarXivlink\yr 2005\pages 1--17
\endref
\ref\myrefno{8}\by Sharipov R. A.\paper A note on the dynamics and
thermodynamics of dislocated crystals\jour
\myhref{http://arxiv.org/abs/cond-mat/0504180}{cond-mat/0504180} 
in Electronic archive \myEarXivlink\yr 2005\pages 1--18
\endref
\ref\myrefno{9}\by Lyuksyutov S. F. \& Paramonov P. B. \& Sharipov R. A. 
\& Sigalov G.\paper Induced nanoscale deformations in polymers using 
atomic force microscopy\jour Phys\. Rev\. B \vol 70\yr 2004
\issue 174110
\endref
\medskip{\bf Part 7. Tensor analysis.}\medskip
\ref\myrefno{1}\by Sharipov R. A.\paper Tensor functions of tensors 
and the concept of extended tensor fields\jour 
\myhref{http://arxiv.org/abs/math/0503332}{math/0503332}
in Electronic archive \myEarXivlink\yr 2005\pages 1--43
\endref
\ref\myrefno{2}\by Sharipov R. A.\paper Spinor functions of spinors 
and the concept of extended spinor fields\jour
\myhref{http://arxiv.org/abs/math.DG/0511350}{math.DG/0511350}
in Electronic archive \myEarXivlink\yr 2005\pages 1--56
\endref
\ref\myrefno{3}\by Sharipov R. A.\paper Commutation relationships and
curvature spin-tensors for extended spinor connections\jour
\myhref{http://arxiv.org/abs/math.DG/0512396}{math.DG/0512396}
in Electronic archive \myEarXivlink\yr 2005\pages 1-22
\endref
\medskip{\bf Part 8. Particles and fields.}\medskip
\ref\myrefno{1}\by Sharipov R. A.\paper A note on Dirac spinors in 
a non-flat space-time of ge\-neral relativity\jour
\myhref{http://arxiv.org/abs/math.DG/0601262}{math.DG/0601262}
in Electronic archive \myEarXivlink\yr 2006\pages 1--22
\endref
\ref\myrefno{2}\by Sharipov R. A.\paper A note on metric connections 
for chiral and Dirac spi\-nors\jour 
\myhref{http://arxiv.org/abs/math.DG/0602359}{math.DG/0602359}
in Electronic archive \myEarXivlink\yr 2006\pages 1--40
\endref
\ref\myrefno{3}\by Sharipov R. A.\paper On the Dirac equation in a
gravitation field and the secondary quantization\jour
\myhref{http://arxiv.org/abs/math.DG/0603367}{math.DG/0603367}
in Electronic archive \myhref{http://arxiv.org}{http:/\negskp/}
\myhref{http://arxiv.org}{arxiv.org}\yr 2006\pages 1--10
\endref
\ref\myrefno{4}\by Sharipov R. A.\paper The electro-weak and color 
bundles for the Standard Model in a gravitation field\jour
\myhref{http://arxiv.org/abs/math.DG/0603611}{math.DG/0603611}
in Electronic archive\linebreak\myEarXivlink\yr 2006\pages 1--8
\endref
\ref\myrefno{5}\by Sharipov R. A.\paper A note on connections of 
the Standard Model in a gravitation field\jour
\myhref{http://arxiv.org/abs/math.DG/0604145}{math.DG/0604145}
in Electronic archive \myEarXivlink\yr 2006\pages 1--11
\endref
\ref\myrefno{6}\by Sharipov R. A.\paper A note on the Standard Model 
in a gravitation field\jour
\myhref{http://arxiv.org/abs/math.DG/0605709}{math.DG/0605709}
in Electronic archive \myEarXivlink\yr 2006\pages 1--36
\endref
\medskip{\bf Part 9. Textbooks.}\medskip
\ref\myrefno{1}\by Sharipov R. A.\book Theory of representations of 
finite groups\publ Bash-NII-Stroy\publaddr Ufa\yr 1995\moreref
see also e-print
\myhref{http://arxiv.org/abs/math.HO/0612104}{math.HO/0612104}
\endref
\ref\myrefno{2}\by Sharipov R. A\book Course of linear algebra and 
multidimensional geometry\publ Bashkir State University\publaddr
Ufa\yr 1996\moreref see also e-print
\myhref{http://arxiv.org/abs/math.HO/0405323}{math.HO/0405323}
\endref
\ref\myrefno{3}\by Sharipov R. A.\book Course of differential geometry 
\publ Bashkir State University\publaddr Ufa\yr 1996\moreref see also
e-print
\myhref{http://arxiv.org/abs/math.HO/0412421}{math.HO/0412421}
\endref
\ref\myrefno{4}\by Sharipov R. A.\book Classical electrodynamics and 
theory of relativity\publ Bash\-kir State University\publaddr Ufa\yr 1996 
\moreref see also e-print
\myhref{http://arxiv.org/abs/physics/0311011}{physics/0311011}
\endref
\ref\myrefno{5}\by Sharipov R. A.\book Foundations of geometry for 
university students and high-school students\publ Bashkir State 
University\yr 1998\moreref see also e-print 
\myhref{http://arxiv.org/abs/math.HO/0702029}{math.HO/0702029}
\endref
\ref\myrefno{6}\by Sharipov R. A.\book Quick introduction to tensor 
analysis\publ free on-line textbook\yr 2004\moreref see also e-print
\myhref{http://arxiv.org/abs/math.HO/0403252}{math.HO/0403252}
\endref
\endRefs
\enddocument
\end